\documentclass[journal]{new-aiaa}
\usepackage[utf8]{inputenc}
\usepackage{textcomp}

\usepackage{graphicx}
\usepackage{amsmath}
\usepackage[version=4]{mhchem}
\usepackage{siunitx}
\usepackage{longtable,tabularx}
 \usepackage{palatino}
 \usepackage{graphics}
 \usepackage{amsmath,environ,bm}
 \usepackage{geometry}
 \usepackage{multicol}
 \usepackage{float}
  \usepackage[table,xcdraw]{xcolor}
 \usepackage{multirow}
 \usepackage[superscript]{cite}
 \usepackage{varioref}%  smart page, figure, table, and equation referencing
 \usepackage{wrapfig}%   wrap figures/tables in text (i.e., Di Vinci style)
 \usepackage{threeparttable}% tables with footnotes
 \usepackage{dcolumn}%   decimal-aligned tabular math columns
  \newcolumntype{d}{D{.}{.}{-1}}
 \usepackage{nomencl}%   nomenclature generation via makeindex
  \makeglossary
 \usepackage{subfigure}% subcaptions for subfigures
 \usepackage{subfigmat}% matrices of similar subfigures, aka small mulitples
 \usepackage{fancyvrb}%  extended verbatim environments
  \fvset{fontsize=\footnotesize,xleftmargin=2em}
\usepackage{mathtools}
\usepackage[ruled,vlined]{algorithm2e}
\usepackage[normalem]{ulem}

\title{High-resolution solutions of nonlinear, advection dominated problems using a generalized finite element method}

\author{Troy Shilt \footnote{PhD candidate, shilt.4@osu.edu, Department of Mechanical and Aerospace Engineering. The Ohio State University, Columbus, OH, 43210, USA}, Patrick J. O'Hara \footnote{Research Engineer, Air Force Research Laboratory, Wright-Patterson AFB, OH, 45433, USA}, and Jack J. McNamara \footnote{Professor, mcnamara.190@osu.edu, Department of Mechanical and Aerospace Engineering. The Ohio State University, Columbus, OH, 43210, USA}}

\begin{document}
\maketitle

\begin{abstract}
Traditional finite element approaches are well-known to introduce spurious oscillations when applied to advection-dominated problems. We explore alleviation of this issue from the perspective of a generalized finite element formulation, which enables stabilization through an enrichment process. The presented work uses solution-tailored enrichments for the numerical solution of the one-dimensional, unsteady Burgers' equation. Mainly, generalizable exponential and hyperbolic tangent enrichments effectively capture local, steep boundary layer/shock features. Results show natural alleviation of oscillations and return smooth numerical solutions over coarse grids. Additionally, significantly improved error levels are observed compared to Lagrangian finite element methods.
\end{abstract}


\section*{Nomenclature}
\begin{tabbing}
\hspace*{0.7in}\=second column tab \kill
$\textbf{A}$\>= matrix containing asymmetric terms\\ 
$E_{\alpha j}$\>= the $j$th enrichment function defined over node $\alpha$\\
$\textbf{f}_{\Gamma_N}$\>= RHS vector corresponding to Neumann boundary condition terms\\
$g_{\Gamma_D}$\>= Dirichlet boundary conditions\\
$g_{\Gamma_N}$\>= Neumann boundary conditions\\
$h$\>= element size\\
$H^1$\>= first order Hilbert space\\
$\widetilde{H}^1$\>= subset of $H^1$ satisfying prescribed Dirichlet boundary conditions $g_{\Gamma_D}$\\
$\dot{H}^1$\>= subset of $H^1$ which vanish on the Dirichlet boundary $\Gamma_D$\\
$\textbf{K}$\>= stiffness matrix\\
$m_\alpha$\>= dimension of the space $\chi_\alpha$\\
$\textbf{M}$\>= mass matrix\\
$N(h)$\>= total number of nodes in domain $\Omega$, determined by $h$\\
$\mbox{Pe}$\>= Péclet number\\
$\mathbb{R}^n$\>= $n$-dimensional real space\\
$S^{GFEM}$\>= GFEM trial space\\
$t_b$\>= Breaking time / time when a shock first forms in the inviscid Burgers' equation\\
$u$\>= test/solution field\\
$u_h$\>= FEM/GFEM approximation\\
$u_{IC}$\>= initial condition\\
$\widetilde{V}$\>= finite-dimensional subspace of $\widetilde{H}^1$\\
$\dot{V}$\>= finite-dimensional subspace of $\dot{H}^1$\\
$w$\>= trial/weighting functions\\
$w_h$\>= FEM/GFEM trial/weighting functions\\
$\textbf{x}_\alpha$\>= point over which patch $\omega_\alpha$ is defined\\
$\alpha$\>= node in the computational domain\\
$\Gamma$\>= computational domain boundary\\
$\Gamma_D$\>= domain boundary where Dirichlet boundary conditions are prescribed\\
$\Gamma_N$\>= domain boundary where Neumann boundary conditions are prescribed\\
$\nu$\>= kinematic viscosity / diffusion coefficient\\
$\phi_{\alpha j}$\>= the $j$th GFEM shape function corresponding to node $\alpha$\\
$\varphi_\alpha$\>= finite element shape function over node $\alpha$\\
$\chi_\alpha$\>= local approximate space / space of enrichment functions\\
$\omega_\alpha$\>= patch defined over node $\alpha$\\
$\Omega$\>= computational domain\\
$\emptyset$\>= null set\\
\end{tabbing}

\section{Introduction}
There is an increased interest in applying finite element methods (FEM) to fluid dynamic problems due to a desire to obtain computationally efficient solutions of multiscale flows. Generalized/extended finite element methods (G/XFEM) are a promising approach towards this goal due to a high degree of flexibility for incorporating solution-tailored features into the finite element approximation space while maintaining local solution conformity \citep{Oden1998, Duarte2000, Melenk1996, Strouboulis2000,Deshmukh2020}. It has been noted in Gracie et al. \citep{Gracie2009} that the G/XFEM are equivalent approaches, and for the rest of the paper, the authors will adopt the use of the term GFEM to refer to both methods. Previous work by the authors has explored the application of GFEM to different aspects of solving fluid dynamics problems. For example, in \citep{Shilt2020} it is shown that properties of the GFEM naturally mitigate the effect of locking in Stokes flow, a regime where viscous forces are dominant. The authors also previously explored mitigation of spurious oscillations using GFEM for the linear advection-diffusion equation in \citep{Shilt2021}. Here insights are provided on the natural capability of the enrichment process for stabilizing advection-dominated problems. Additionally, unlike traditional stabilized methods (streamline upwind/Petrov-Galerkin method \citep{Brooks1982}, Galerkin/least-squares \citep{Hughes1989, Franca1989}, residual-free bubble methods \citep{Baiocchi1993, Franca1996}), no restrictions are placed on the enrichment selection process, thus allowing the choice of solution-tailored enrichments that enable stable, high-accuracy solutions.

This paper extends our previous work on the advection-diffusion equation and explores solution-tailored enrichments applied to the one-dimensional, unsteady Burgers' equation. The viscous Burgers' equation is identical to the advection-diffusion equation, except the advection coefficient is replaced by the solution variable, $u$, thus resulting in a nonlinear term. This equation was first introduced by Bateman in \citep{Bateman1915} as a relatively simple equation to explore discontinuous solutions as the kinematic viscosity tends towards zero. It was not until many years later that Burgers explored this equation in \citep{Burgers1948} as a nonlinear equation with similar phenomena to turbulence. Nowadays, the Burgers' equation is known to have physical relevance for problems which include: viscous flows, shock theory, gas dynamics, cosmology, traffic flow, and quantum computing \citep{Bonkile2018}. The Burgers' equation has many features similar to the Navier-Stokes equations and is used to clarify the interaction between transient, dissipative, and nonlinear advective terms. Specifically, the Burgers' equation contains an inertial and dissipation range similar to turbulence in the Navier-Stokes equations \citep{Burgers1974,Bayona2017}. As such, numerical simulation of the Burgers' equation presents a challenge when inertial effects dominate the solution, analogous to challenges associated with numerically solving Navier-Stokes equations with high Reynolds numbers. These highly advective problems often demand ultra-fine discretizations to resolve the system's multiscale behavior accurately; otherwise, spurious oscillations arise in the numerical solution. 

The remaining outline of this paper is as follows: first, we discuss the governing equations for the viscous Burgers’ equation, formulation of the GFEM nonlinear system of equations, and linearization using Newton-Raphson. Next, the inviscid Burgers' equation is presented, followed by a discussion on the formation of shocks in the domain and numerical stability. Finally, we present numerical examples for the GFEM solution to the unsteady one-dimensional Burgers' equation, along with a general discussion of the results.

\section{Viscous Burgers' equation}

\subsection{Preliminaries}
Let $\Omega$ be an open set contained in $\mathbb{R}^n$, $n\geq 1$, with a piecewise smooth boundary $\Gamma$. Vector and tensor fields defined on $\Omega$ are in boldface notation with lowercase and uppercase variables, respectively (e.g., vector $\textbf{y}$ and tensor $\textbf{A}$). For prescribing boundary conditions, it is necessary to define $\Gamma = \Gamma_D \cup \Gamma_N$ such that $\Gamma_D \cap \Gamma_N = \emptyset$, where $\Gamma_D$ denotes part of the boundary for prescribed Dirichlet boundary conditions, and $\Gamma_N$ denotes part of the boundary for prescribed Neumann boundary conditions.

\subsection{Governing equation}
The one-dimensional viscous Burgers' equation is the following: find $u$ such that

\begin{equation}
\label{viscous_burger_strong_formulation}
\begin{aligned}
\frac{\partial u}{\partial t}  + u \frac{\partial u}{\partial x} - \nu \frac{\partial^2 u}{\partial x^2} = 0& \quad \mbox{on} \quad \Omega\\
u(x,0) = u_{IC}(x)& \quad \mbox{on} \quad \Omega\\
u(x,t) = g_{\Gamma_D}(x,t)& \quad\mbox{on} \quad \Gamma_D\\
\frac{\partial u(x,t)}{\partial x} = g_{\Gamma_N}(x,t)& \quad \mbox{on} \quad \Gamma_N\\
\end{aligned}
\end{equation}

\noindent where when referring to fluids, $\nu$ is the kinematic viscosity, and $u(x,t)$ is the fluid velocity. The weak formulation of the boundary value problem, Eq.\ref{viscous_burger_strong_formulation}, is obtained by multiplying by weighting functions $w$ and integrating over the domain $\Omega$. The formulation is as follows: find $u \in \widetilde{H}^1$ such that for all $w \in \dot{H}^1$:

\begin{equation}
\label{viscous_burger_weak_formulation}
\begin{aligned}
\int_\Omega \Bigg( w \frac{\partial u}{\partial t} + w u \frac{\partial u}{\partial x} + \nu \frac{\partial w}{\partial x} \frac{\partial u}{\partial x}\Bigg) d \Omega = \nu\int_{\Gamma_N} w \frac{\partial u}{\partial x}  \, d \Gamma_N 
\end{aligned}
\end{equation}
\noindent where: 
\begin{align}
\label{spaces_of_u_and_w}
&\widetilde{H}^1 = \{u \in H^1 |\, u = g_{\Gamma_D}(x,t) \,\, \mbox{on} \,\, \Gamma_D\} \\
&\dot{H}^1 = \{w \in H^1 |\, w = 0 \,\, \mbox{on} \,\, \Gamma_D\}
\end{align}

Note that the above formulation does not detail enforcement of Dirichlet boundary conditions ($g_{\Gamma_D}$). This will be the focus of a later section. The Galerkin formulation is obtained by assuming finite-dimensional approximations of the test and trial functions. Let $\widetilde{V}$ be a finite-dimensional subspace of the space $\widetilde{H}^1$, such that $u_h \in \widetilde{V}$ is a finite-dimensional approximate solution to the weak form of the boundary value problem, Eq.\ref{viscous_burger_weak_formulation}, and similarly define $\dot{V}$ to be a finite-dimensional subspace of the space $\dot{H}^1$. The Galerkin formulation is as follows: find $u_h \in \widetilde{V}$ such that for all $w_h \in \dot{V}$:

\begin{equation}
\label{viscous_burger_galerkin_formulation}
\begin{aligned}
\int_\Omega \Bigg( w_h \frac{\partial u_h}{\partial t} + w_h u_h \frac{\partial u_h}{\partial x} + \nu \frac{\partial w_h}{\partial x} \frac{\partial u_h}{\partial x}\Bigg) d \Omega = \nu\int_{\Gamma_N} w_h \frac{\partial u_h}{\partial x}  \, d \Gamma_N 
\end{aligned}
\end{equation}

\subsection{GFEM approximation space}

Constructing the GFEM approximation space consists of three components: a) patches, b) a partition of unity, and c) local approximation spaces.

\begin{figure}[ht!]
\centering{\includegraphics[width=5in]{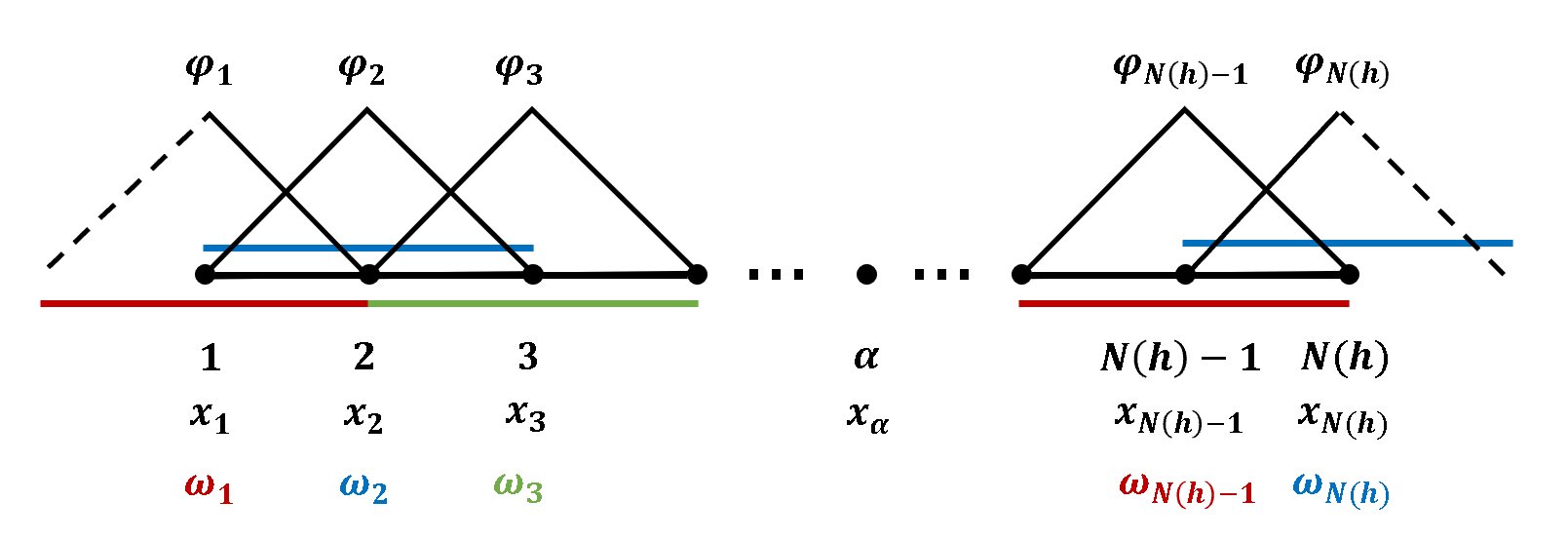}}
\caption{Patches ($\omega_\alpha$) and the partition of unity composed of linear, Lagrangian shape functions ($\varphi_\alpha$) over a one-dimensional, uniformly discretized computational domain ($\Omega$)}
\label{fig:GFEM_patches_and_POU}
\end{figure}

\begin{enumerate}[label=\alph*)]

\item \textit{Patches}: build an open covering defined such that for a parameter $h > 0$:

\begin{equation}
\label{patch_space}
\Omega \subset \bigcup^{N(h)}_{\alpha = 1} \omega_\alpha
\end{equation}

where $\omega_\alpha$ are patches defined over $\textbf{x}_\alpha$, $\alpha = 1, ...\,, N(h)$. Any $\textbf{x} \in \Omega$ belongs to at most $M \leq N(h)$ elements of the set $\{\omega_\alpha\}_{\alpha = 1}^{N(h)}$. In (G)FEM, $\omega_\alpha$ is given by the union of finite elements sharing node $\alpha$ of the finite element mesh covering $\Omega$. Additionally, $N(h)$ is defined to be the total number of nodes in the domain $\Omega$. Fig. \ref{fig:GFEM_patches_and_POU} provides a visual representation of patches typically used in (G)FEM for a one-dimensional domain.

\item \textit{Partition of unity}: let $\{\varphi_\alpha\}_{\alpha = 1}^{N(h)}$ be piecewise $C^0$ functions defined on $\Omega$ satisfying:

\begin{equation}
\label{partition_of_unity}
\sum_{\alpha = 1}^{N(h)} \varphi(\textbf{x}) = 1, \quad \forall \, \textbf{x} \in \Omega
\end{equation}

Then the set $\{\varphi_\alpha\}_{\alpha = 1}^{N(h)}$ forms a partition of unity with respect to the open cover set $\{\omega_\alpha\}_{\alpha = 1}^{N(h)}$. In GFEM, the set $\{\varphi_\alpha\}_{\alpha = 1}^{N(h)}$ is typically chosen as linear, Lagrangian shape functions (see Fig. \ref{fig:GFEM_patches_and_POU}).

\item \textit{Local and approximation spaces}: For each patch $\omega_\alpha$ we associate an $m_\alpha$-dimensional space $\chi_\alpha(\omega_\alpha)$ of functions, denoted the local approximate space, such that:

\begin{equation}
\label{local_approximation_space}
\chi_\alpha = \mbox{span}\{E_{\alpha j},\, 1 \leq j \leq m_\alpha,\, E_{\alpha j} \in H^1\}
\end{equation}

where the functions $E_{\alpha j} \in \chi_\alpha$ are known as \textit{enrichment functions}. It is assumed each $\chi_\alpha$ contain a constant function. This inclusion allows for the set $\{\varphi_\alpha\}_{\alpha = 1}^{N(h)}$ to be contained in the trial space. 

\end{enumerate}

Finally, the GFEM approximation space is given by:

\begin{equation}
\label{GFEM_trial_space}
S^{GFEM}\big(\Omega\big) = \mbox{span}\{\phi_{\alpha j} = \varphi_\alpha E_{\alpha j} \,\mbox{(no sum over } \alpha\mbox{)}, \, 1 \leq \alpha \leq N(h), \, 1 \leq j \leq m_\alpha\}
\end{equation}

\noindent where $\phi_{\alpha j}$ are called the \textit{GFEM shape functions}.

\subsection{Solution of the GFEM system of equations}
Any trial function $u_h \in S^{GFEM} \subset \widetilde{V}$ may be written in vector notation as:

\begin{equation}
\label{GFEM_vector_notation}
u_h = \bm{ \phi }^T \bm{(}x\bm{)} \textbf{c}\bm{(}t\bm{)} 
\end{equation}

\noindent where $\bm{\phi }$ is the vector of GFEM shape functions ($\phi_{\alpha j}$) and $\bm{c}$ is the vector of corresponding weighting coefficients. GFEM test functions $w_h$ are defined identically. Substituting Eq. \ref{GFEM_vector_notation} into the Galerkin formulation Eq. \ref{viscous_burger_galerkin_formulation} results in the following system of equations:

\begin{equation}
\label{GFEM_system_of_equations}
\textbf{M} \dot{\textbf{c}}(t) = - (\textbf{A}\bm{(}t\bm{)} +\textbf{K})\textbf{c}(t) + \textbf{f}_{\Gamma_N} (t)
\end{equation}

\noindent where

\begin{align}
\label{viscous_burgers_equation_GFEM_matrix_terms}
&\textbf{M} = \int_{\Omega} \bm{ \phi } \bm{ \phi }^T \, d\Omega \\
&\textbf{A}\bm{(}t\bm{)} = \int_{\Omega} \bm{ \phi } \bm{\phi}^T \textbf{c}(t) \frac{\partial \bm{ \phi }^T}{\partial x} \, d\Omega \\
&\textbf{K} = \nu \int_{\Omega} \frac{\partial \bm{\phi}}{\partial x} \frac{\partial \bm{\phi}^T}{\partial x} \, d\Omega \\
&\textbf{f}_{\Gamma_N}(t) = \int_{\Gamma_N} \bm{\phi} g_{\Gamma_N} (x,t) \, d\Gamma_N
\end{align}

\subsubsection{Time discretization using Crank-Nicolson method}

The Crank-Nicolson scheme is used for temporal discretization of Eq. \ref{GFEM_system_of_equations}, such that:

\begin{equation}
\label{GFEM_discretization}
2 \textbf{M} \Bigg( \frac{\textbf{c}^{n+1} - \textbf{c}^n}{\Delta t} \Bigg) = \Big[\textbf{f}_{\Gamma_N}(t^{n+1}) - (\textbf{A}\bm{(}\textbf{c}^{n+1}\bm{)} + \textbf{K})\textbf{c}^{n+1}\Big] + \Big[\textbf{f}_{\Gamma_N}(t^{n}) - (\textbf{A}\bm{(}\textbf{c}^{n}\bm{)} +\textbf{K})\textbf{c}^{n}\Big]
\end{equation}

\noindent which may be rearranged as:

\begin{equation}
\label{GFEM_crank_nicolson}
\textbf{G}(\textbf{c}^{n+1}) \textbf{c}^{n+1} = \textbf{g}
\end{equation}

\noindent where $\textbf{G}(\textbf{c}^{n+1}) = \frac{2}{\Delta t} \textbf{M} + \textbf{A}\bm{(}\textbf{c}^{n+1}\bm{)} + \textbf{K}$ and $\textbf{g} = \textbf{f}_{\Gamma_N}(t^{n+1}) + \Big[\textbf{f}_{\Gamma_N}(t^{n}) - (\textbf{A}\bm{(}\textbf{c}^{n}\bm{)} +\textbf{K})\textbf{c}^{n}\Big]$.

\subsubsection{Enforcing Dirichlet boundary conditions}
In the Lagrangian finite element method, the Kronecker delta property of the shape functions allows direct enforcement of Dirichlet boundary conditions by setting the coefficients equal to the desired solution values. However, additional degrees of freedom per node introduced through the GFEM enrichment process makes this approach nontrivial. A straightforward manner of enforcing desired boundary conditions in GFEM is to add a penalty term on both sides of the system of equations Eq. \ref{GFEM_crank_nicolson} such that:

\begin{equation}
\label{GFEM_crank_nicolson_BCs}
\begin{aligned}
[\textbf{M}_{\Gamma_D} + \textbf{G}(\textbf{c}^{n+1})] \textbf{c}^{n+1} = \textbf{g}
 + \textbf{f}_{\Gamma_D}(t)\\
\end{aligned}
\end{equation}

\noindent where $\beta$ is the penalty parameter that is typically very large in relation to the other matrix components, but not so large to cause ill conditioning of system matrices, and:

\begin{align}
\label{BC terms}
&\textbf{M}_{\Gamma_D} = \int_{\Gamma_D} \bm{ \phi } \bm{ \phi }^T \, d\Gamma_D \\
&\textbf{f}_{\Gamma_D}(t) = \int_{\Gamma_D} \bm{\phi} g_{\Gamma_D} (x,t) \, d\Gamma_D
\end{align}

\subsubsection{Iteration of the nonlinear system with Newton-Raphson method}

Equation \ref{GFEM_crank_nicolson_BCs} is a nonlinear system of equations with knowns $\textbf{c}^{n}$ and unknowns $\textbf{c}^{n+1}$. Using the Newton-Raphson method, the solution may be iterated to solve for the solutions at the $n+1$ time step. To do so, assume $\textbf{c}^{n+1}$ = $\textbf{c}^{n} + \bm{\epsilon}$, where solutions at the previous time steps $n$ are an approximation of the $n+1$ solution, and $\bm{\epsilon}$ is a small correction. First, applying this above decomposition to the nonlinear matrix term ($\textbf{A}(t^{n+1})$) product with solution coefficients at the $n+1$ time step ($\textbf{c}^{n+1}$) simplify to the following after neglecting the underlined $\mathcal{O}(\bm{\epsilon}^2)$ terms:

\begin{equation}
\label{GFEM_newton_step1}
\begin{aligned}
\textbf{A}(\textbf{c}^{n+1}) \textbf{c}^{n+1} &\approx \textbf{A}(\textbf{c}^{n}) \textbf{c}^{n} + \Bigg(\textbf{A}(\textbf{c}^{n}) + \int_\Omega \bm{\phi} \bm{\phi}^T \Bigg(\frac{\partial \bm{\phi}}{\partial x}\textbf{c}^n \Bigg) \, d\Omega\Bigg) \bm{\epsilon} + \sout{\mathcal{O}(\epsilon^2)} \\
&\approx \textbf{A}(\textbf{c}^{n}) \textbf{c}^{n} + \Big(\textbf{A}(\textbf{c}^{n}) + \widetilde{\textbf{A}}(\textbf{c}^{n}) \Big) \bm{\epsilon} + \sout{\mathcal{O}(\epsilon^2)}
\end{aligned}
\end{equation}

Substituion of $\textbf{c}^{n+1}$ = $\textbf{c}^{n} + \bm{\epsilon}$ and Eq. \ref{GFEM_newton_step1} into Eq. \ref{GFEM_crank_nicolson_BCs} results in the following system of equations to solve for the corrections $\epsilon$:

\begin{equation}
\label{GFEM_newton_method_indicial}
\begin{aligned}
\widetilde{\textbf{G}} \bm{\epsilon} = \widetilde{\textbf{g}}
\end{aligned}
\end{equation}

\noindent where $\widetilde{\textbf{G}} = \textbf{G}(\textbf{c}^n) + \widetilde{\textbf{A}}(\textbf{c}^{n}) + \textbf{M}_{\Gamma_D}$ and $\widetilde{\textbf{g}} = \textbf{g}
 + \textbf{f}_{\Gamma_D} - (\textbf{G}(\textbf{c}^n) + \textbf{M}_{\Gamma_D}) \textbf{c}^n$. Finally, iterate over Eq. \ref{GFEM_newton_method_indicial} until some residual is converged.
 
\subsubsection{Linear solver}
Note that due to the ill-conditioning of the GFEM formulation (see \citep{Gupta2013}) the iterative algorithm presented in \citep{Duarte2000} and displayed in Algo. \ref{algorithm:linear_solver} is used to solve the potentially indefinite system of equations. For subsequent numerical examples $\epsilon_1 = \epsilon_2 = 10^{-10}$ is used. 

\begin{algorithm}[H]
\label{algorithm:linear_solver}
\textbf{INPUT}: $\tilde{\textbf{A}}, \tilde{\textbf{b}}$, perturbation $\epsilon_1 << 1$, and criterion $\epsilon_2 << 1$\\
\textbf{OUTPUT}: $\textbf{c} = \textbf{c}_i$\\

\SetAlgoLined
%\KwResult{Write here the result }
Initialization:\\
\quad Precondition $\tilde{\textbf{A}} \tilde{\textbf{c}} = \tilde{\textbf{b}}$ to equivalent system $\textbf{A} \textbf{c} = \textbf{b}$ by defining:\\
\quad \quad \quad $T_{ij} = \frac{\delta_{ij}}{\sqrt{\tilde{\textbf{A}}_{ij}}}$\\
\quad \quad \quad $\textbf{A} = \textbf{T} \tilde{\textbf{A}} \textbf{T}$\\
\quad \quad \quad $\textbf{c} = \textbf{T}^{-1} \tilde{\textbf{c}}$\\
\quad \quad \quad $\textbf{b} = \textbf{T} \tilde{\textbf{b}}$\\

\quad Perturbed matrix: $\textbf{A}_{\epsilon} = \textbf{A}+ \epsilon_1 \textbf{I}$;\\
\quad Approximate system of equations solution vector: $\textbf{c}_0 = \textbf{A}_\epsilon^{-1} \textbf{b}$;\\
\quad  Residual error of approximate system of equations: $\textbf{r}_0 = \textbf{b} - \textbf{A}\textbf{c}_0$;\\
\quad Residual error of solution vector: $\textbf{e}_0 = \textbf{c} - \textbf{c}_0 \approx \textbf{A}_\epsilon^{-1} \textbf{r}_0$;\\
 \While{$\Bigg|\frac{\textbf{e}_i \textbf{A} \textbf{e}_i}{\textbf{c}_i \textbf{A} \textbf{c}_i}\Bigg| > \epsilon_2$}{ 
 	$\textbf{r}_i = \textbf{r}_{i-1} - \sum_{i = 0}^{i - 1} \textbf{A} \textbf{e}_i$;\\
 	$\textbf{e}_i = \textbf{A}_\epsilon^{-1} \textbf{r}_i$;\\
 	$\textbf{c}_i = \textbf{c}_0 + \sum_{i = 0}^{i - 1} \textbf{e}_i$;\\
 }
  	\Return{$\tilde{\textbf{c}} = \textbf{T} \textbf{c}_i$}
  	
 \caption{Solution to the system of equations $\tilde{\textbf{A}} \tilde{\textbf{c}} = \tilde{\textbf{b}}$}
\end{algorithm}
\subsubsection{Initial condition}
Solution of Eq. \ref{GFEM_newton_method_indicial} requires an initial solution vector, $\textbf{c}^0$, which approximates the initial value problem $u_h(x,0) = \bm{ \phi }^T \bm{(}x\bm{)} \textbf{c}^0 \approx u_{IC}(x)$. The initial solution vector is obtained by solving the Galerkin formulation of this IVP, such that:

\begin{equation}
\label{eq:ICs_weak}
\textbf{M} \textbf{c}^0= \int_{\Omega} \bm{ \phi } u_{IC}(x) \, d\Omega.
\end{equation}

\section{Inviscid Burgers' Equation}
\subsection{Governing equation}
The inviscid Burgers' equation represents a limiting case where the kinematic viscosity tends toward zero ($\nu \rightarrow 0$). The resulting inviscid Burgers' equation is the following: find $u$ such that:

\begin{equation}
\label{inviscid_burger_strong_formulation}
\begin{aligned}
\frac{\partial u}{\partial t}  + u \frac{\partial u}{\partial x} = 0& \quad \mbox{on} \quad \Omega\\
u(x,0) = u_{IC}(x)& \quad \mbox{on} \quad \Omega\\
u(x,t) = g_{\Gamma_D}(t)& \quad\mbox{on} \quad \Gamma_D\\
\frac{\partial u(x,t)}{\partial x} = g_{\Gamma_N}(t)& \quad \mbox{on} \quad \Gamma_N\\
\end{aligned}
\end{equation}

Using the method of characteristics an implicit solution to Eq. \ref{inviscid_burger_strong_formulation} can be constructed. Readers are directed to \citep{LeVeque1992} for additional details on this procedure. The resulting implicit solution is given by $u(x,t) = u_{IC}(x - ut) = u_{IC}(\xi)$, with a characteristic trajectory $x = ut + \xi$, where $\xi$ is an arbitrary point on the $x$-axis of the $x$-$t$ plane. Note the formulation of the GFEM system of equations for this inviscid case is identical to that of Eq. \ref{GFEM_system_of_equations} with $\textbf{K}$ removed.

\subsection{Formation of shocks}
In the inviscid case, a discontinuity (``shock") will form in the domain if $u_{IC}'(x) < 0$. This work's notion of a shock is assumed to be inclusive of any solution which contains a steep gradient. Therefore, the distinction between the inviscid and viscous Burgers' equation for shock formation is the discontinuity that appears in the domain. Additionally, the time when the discontinuity first occurs is denoted the breaking time and is given by:

\begin{equation}
\label{inviscid_burger_breaking_time}
\begin{aligned}
t_b = \frac{-1}{\min{u_{IC}'(x)}}\\
\end{aligned}
\end{equation}

\begin{figure}[ht!]
\begin{center}
\begin{subfigmatrix}{6}
\subfigure[$b = 1.25$]{\includegraphics[width=3in]{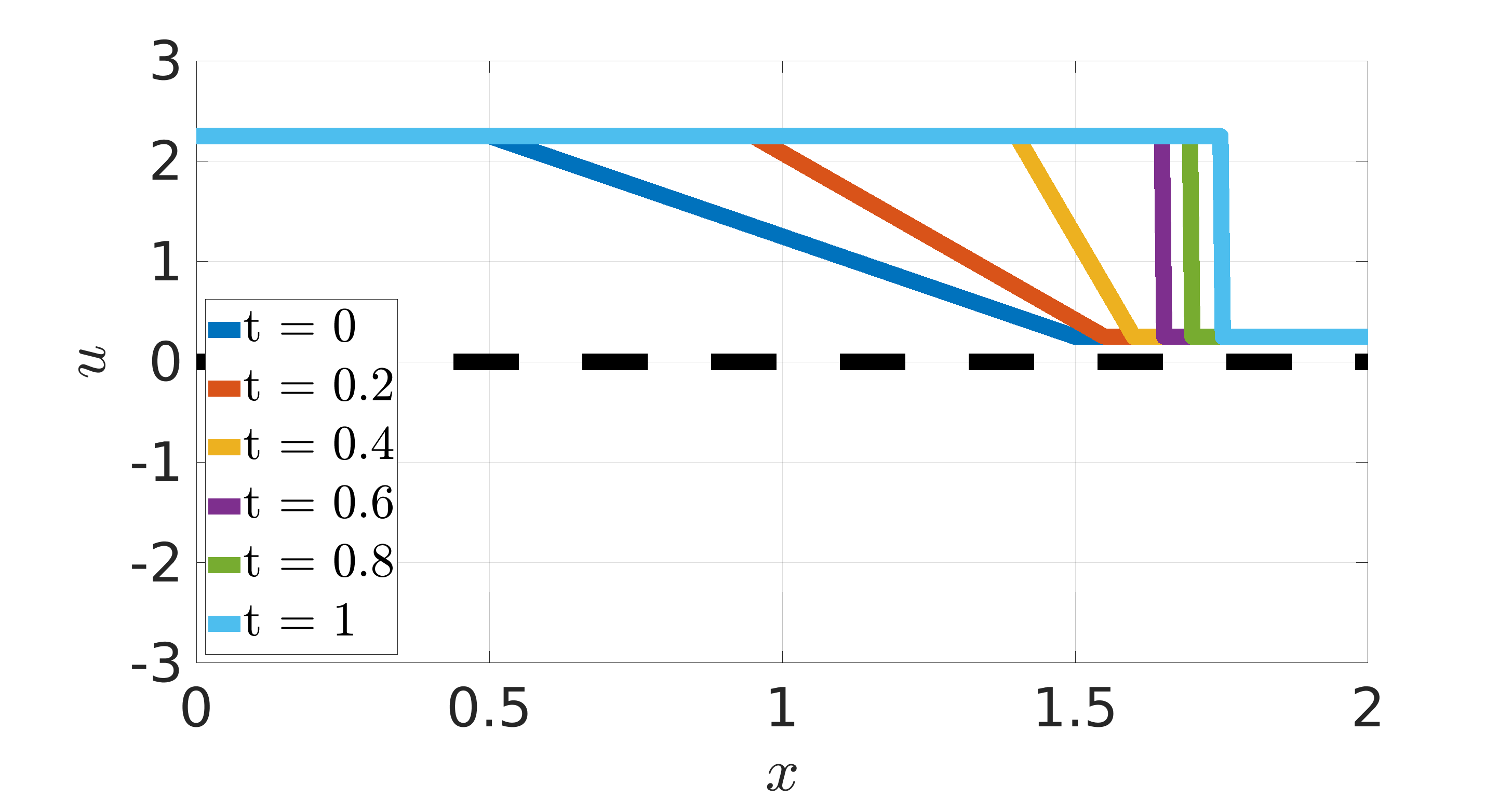}}
\subfigure[$b = 1$]{\includegraphics[width=3in]{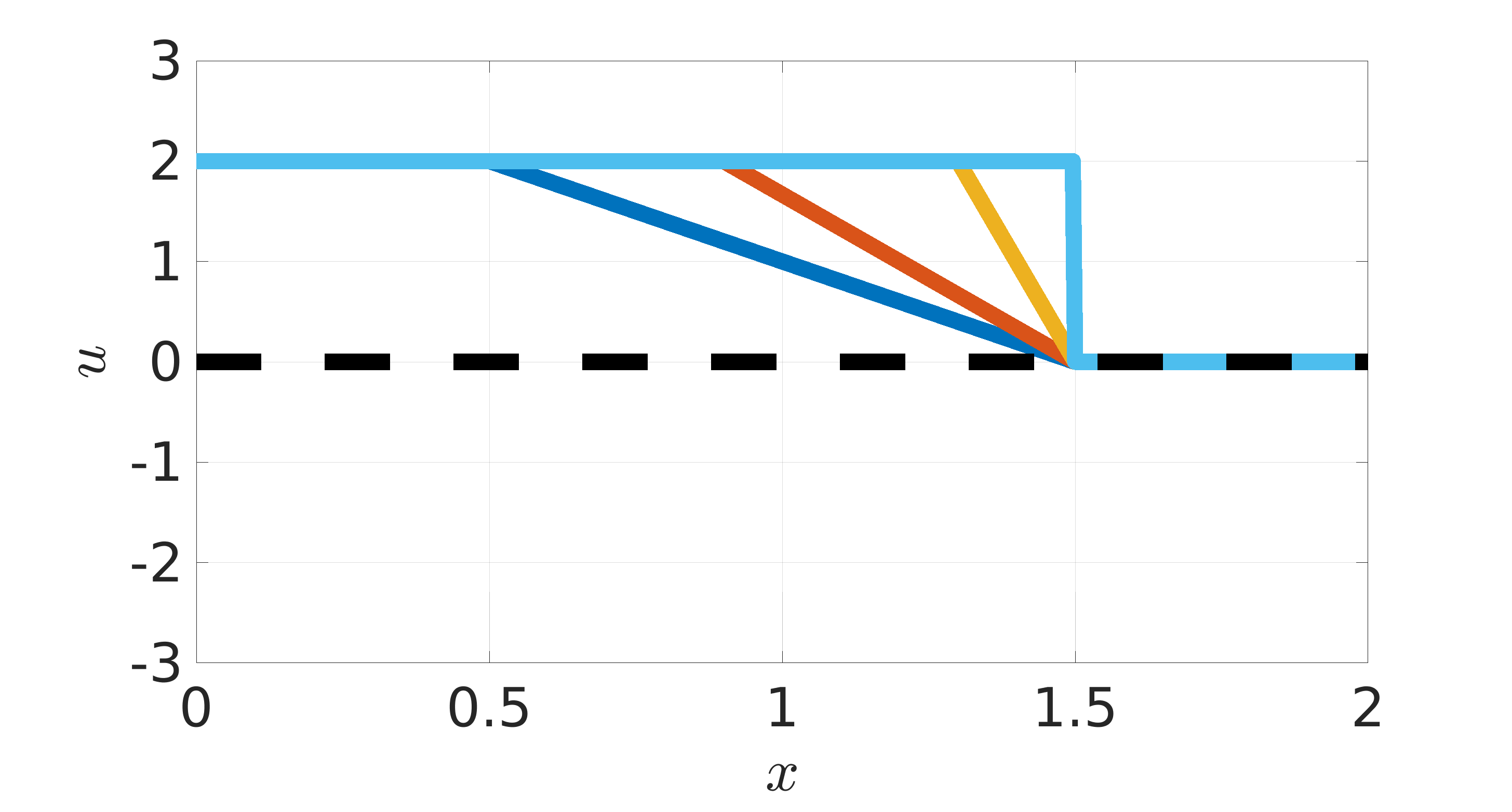}}
\subfigure[$b = 0.5$]{\includegraphics[width=3in]{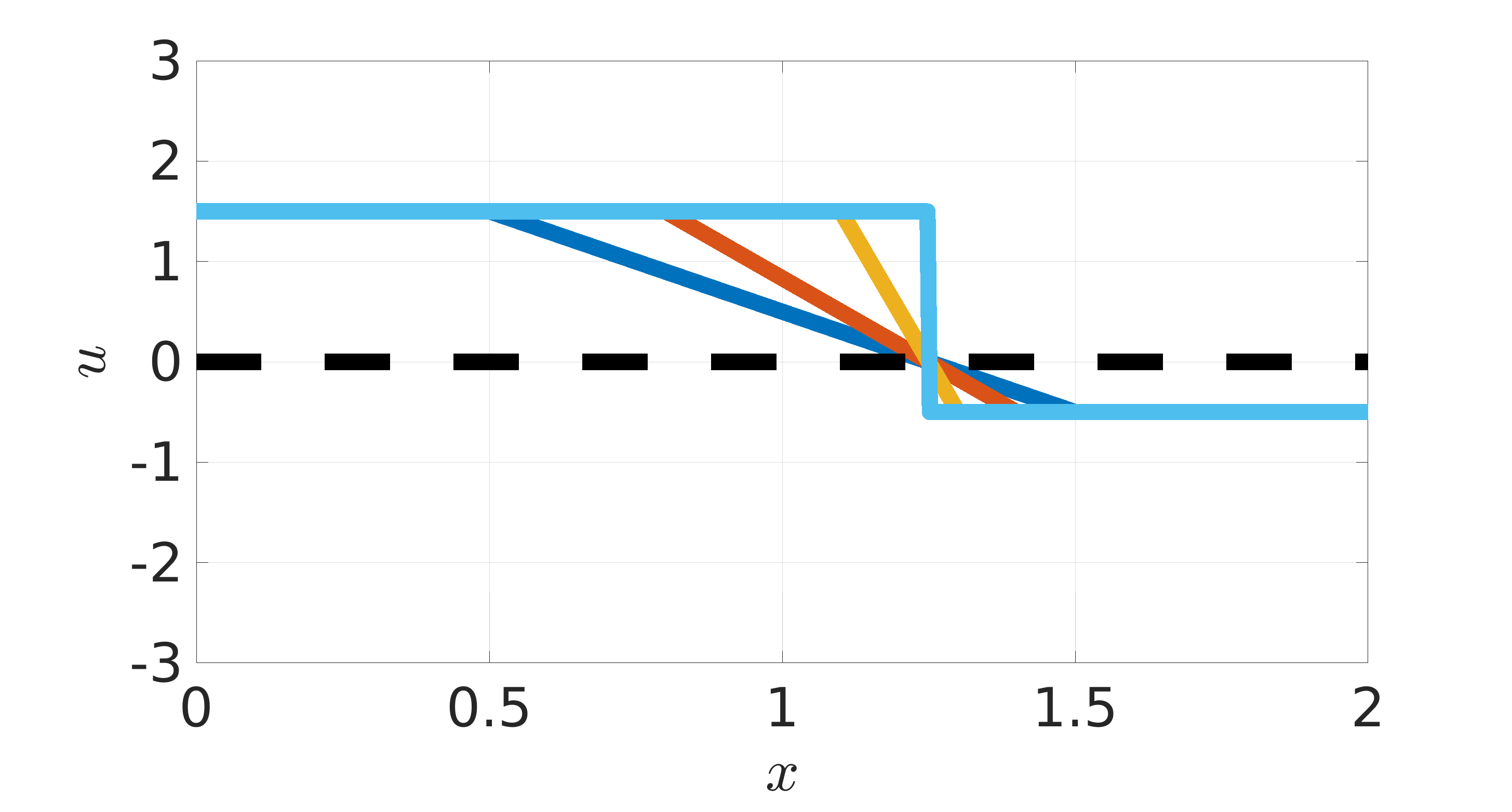}}
\subfigure[$b = 0$]{\includegraphics[width=3in]{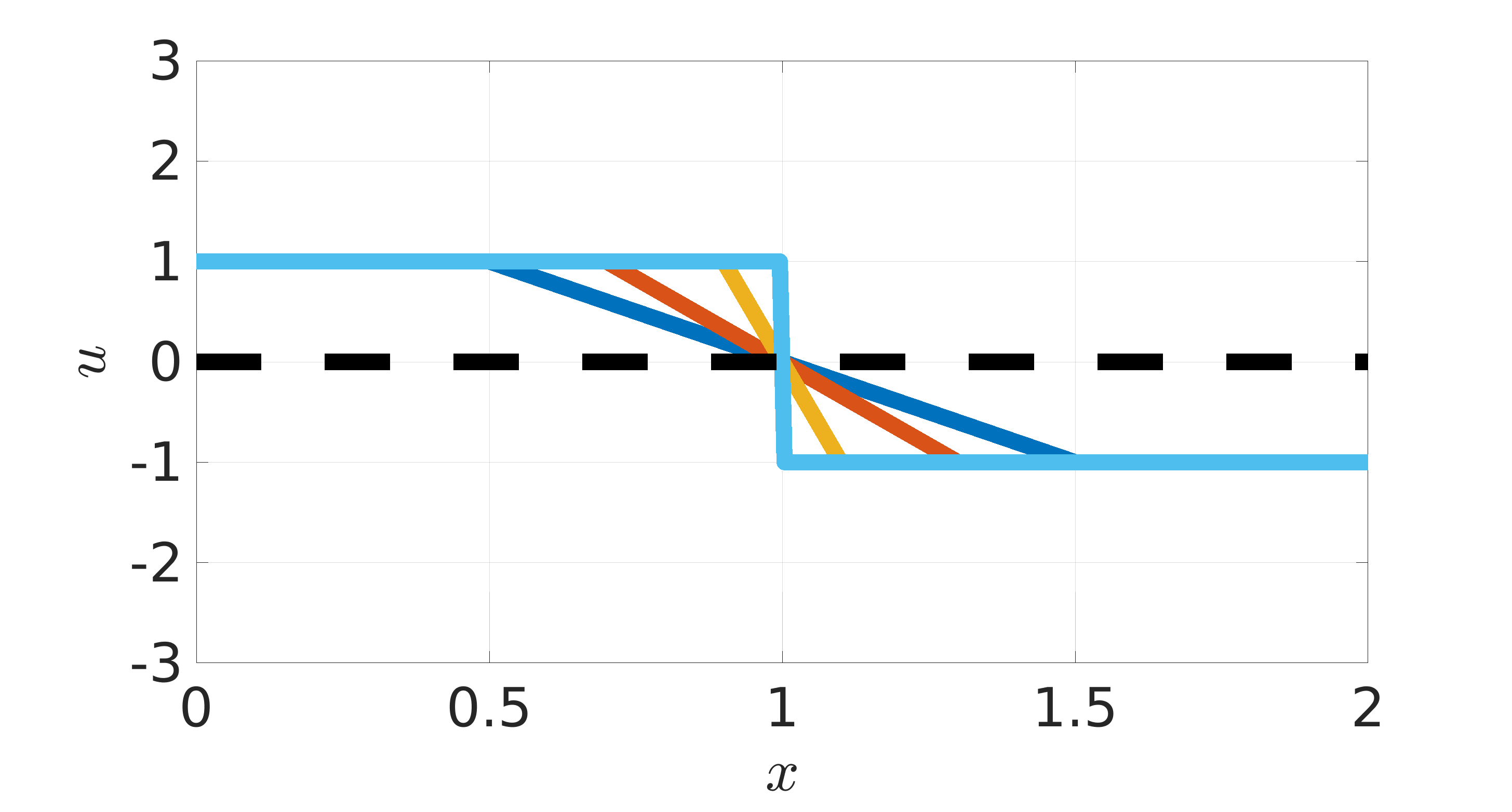}}
\subfigure[$b = -1$]{\includegraphics[width=3in]{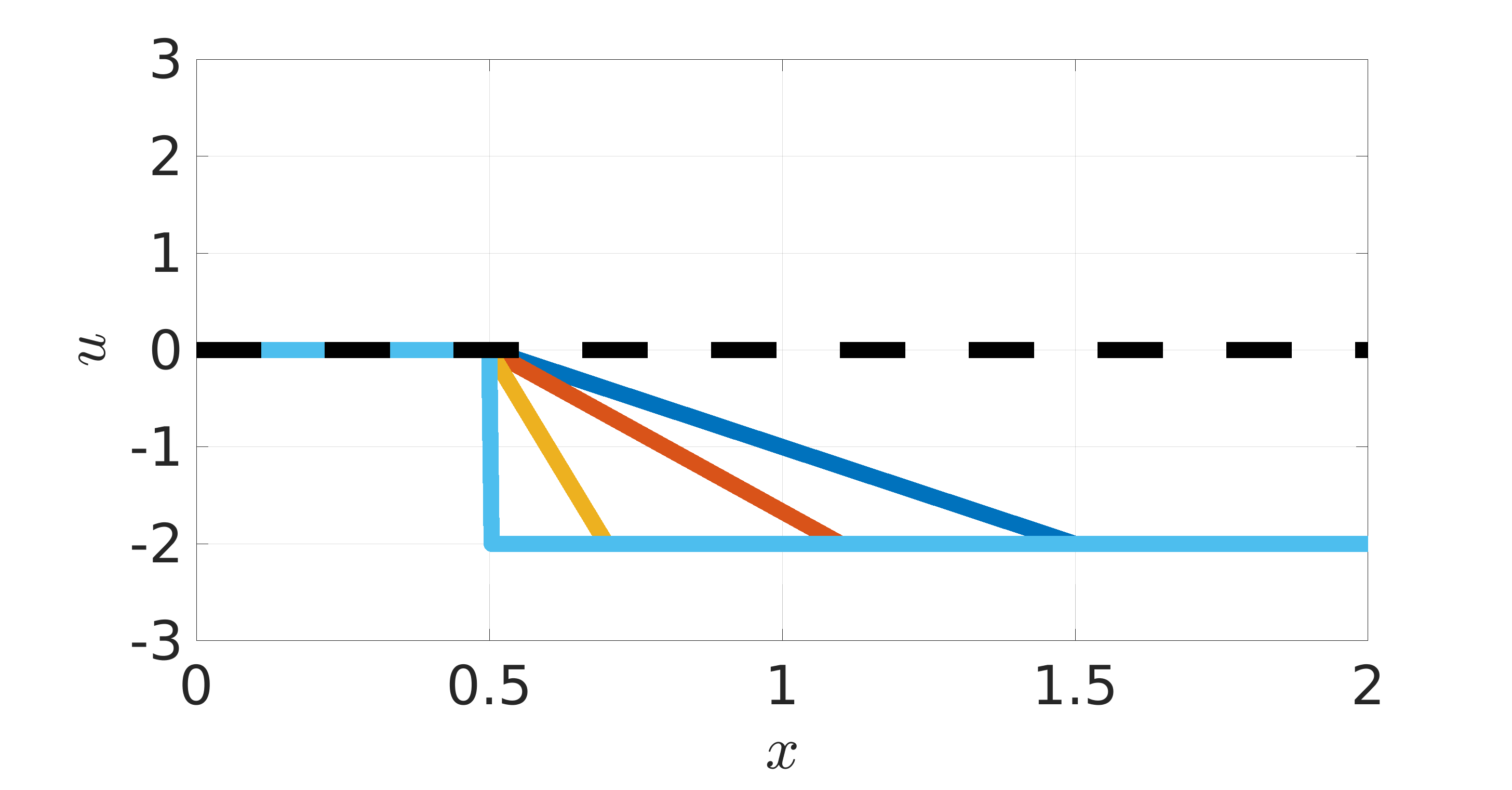}}
\subfigure[$b = -1.25$]{\includegraphics[width=3in]{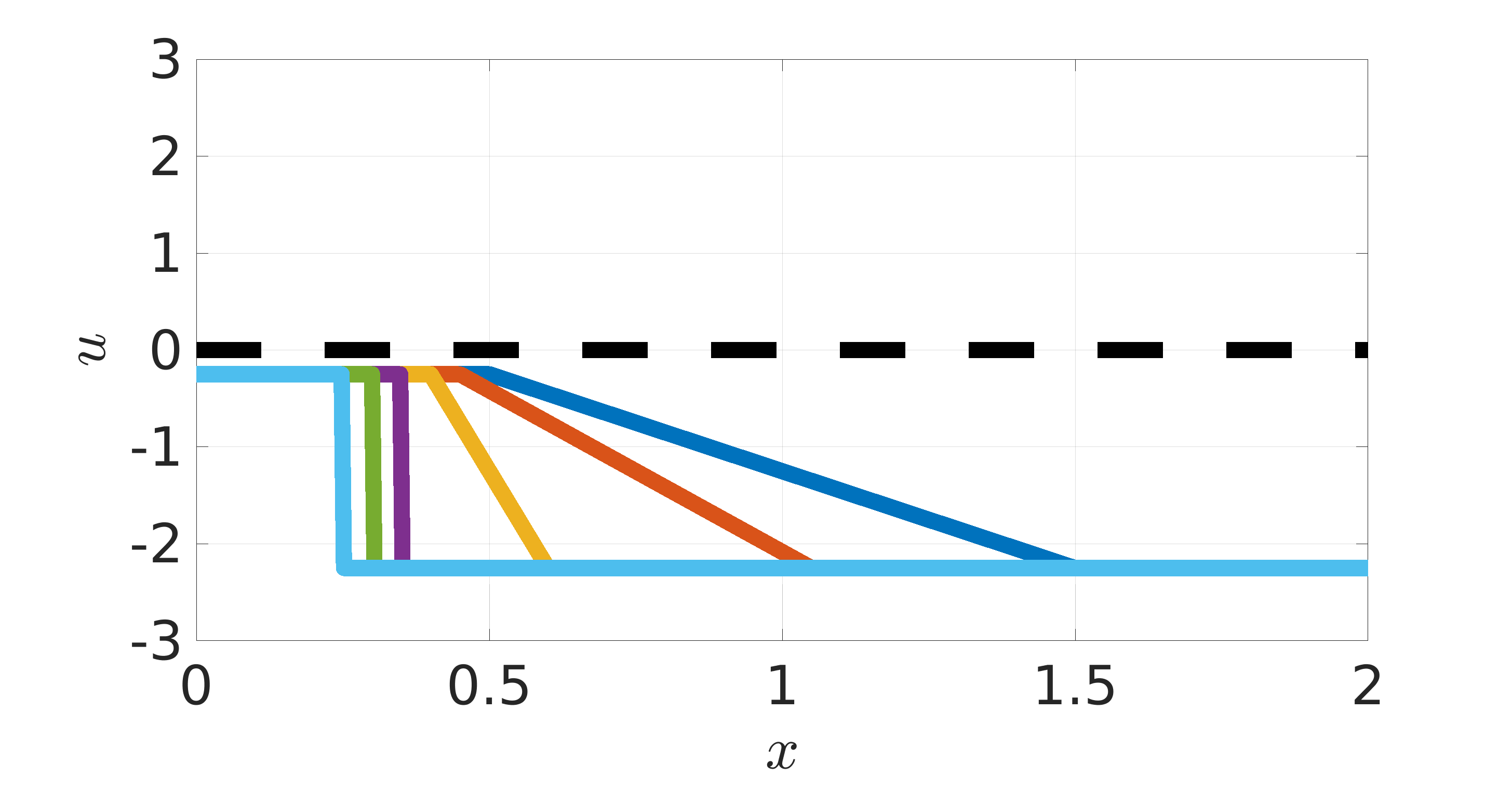}}
\end{subfigmatrix}
\caption{Riemann solution for various initial conditions}
\label{fig:Riemann_problem}
\end{center}
\end{figure}

If $u_{IC}(x)$ crosses the $x-$axis at $x_b$, such that $u_{IC}(x_b) = 0$ and $u_{IC}'(x_b) < 0$, the shock that forms will be stationary at $x_b$. In the case where $u_{IC}'(x) < 0$ and $u_{IC}(x)$ does \textit{not} cross the x-axis, the shock formed will be moving. If $u_{IC}(x) > 0$ the shock will travel in the positive $x$-direction with time, otherwise for $u_{IC}(x) < 0$ the shock will travel in the negative $x$-direction with time. This may be demonstrated by considering a series of Riemann problems represented by the following initial conditions:

\begin{equation}
\label{Riemann_problem_initial_conditions}
u_{IC}(x) = 
\left\{
\begin{array}{ll}
      b + 1 & x \leq \frac{1}{2} \\
      b + 2(1-x) & \frac{1}{2} < x < \frac{3}{2} \\
      b - 1 & \frac{3}{2} \leq x \\
\end{array}
\right.
\end{equation}

\noindent where $b$ is an arbitrary value to translate the initial condition. The resulting solution to Eq. \ref{inviscid_burger_strong_formulation}, with initial conditions Eq. \ref{Riemann_problem_initial_conditions} for values of $b = \{-1.25,\, -1,\, 0,\, 0.5,$ $\, 1,\, 1.25\}$, is shown in Fig \ref{fig:Riemann_problem}. Observe that the initial problem is strictly negative for $b = -1.25$, and the forming shock moves to the left with time. If $b = 1.25$, the initial problem is strictly positive and moves to the right. The remaining values of $b$ cross the $x$-axis at some point and form a stationary wave. Addressing moving shocks is outside of the scope of this paper, as the authors are concerned with addressing stability concerns in the GFEM. Thus, this work considers only stationary shocks for subsequent numerical examples. However, the authors note that previous work on GFEM to solve time-dependent problems indicates a promise to handle transient, local behavior. The authors direct the readers to \citep{Ohara2011} for work on GFEM solution to highly localized sharp, transient thermal gradients as an example of such application.

\subsection{A note on stability}
Recall the viscous Burgers' equation (Eq. \ref{viscous_burger_weak_formulation}) is form identical to the advection-diffusion equation, where the advection coefficient is replaced by the solution variable, $u$. For advection-diffusion equations, the Péclet number is considered for stability of the linear FEM. Specifically, for linear FEM solution over uniform grid size, $h$, Pe = $\frac{a h}{2 \nu} > 1$ results in spurious oscillations, where $a$ and $\nu$ are the rate of advection and rate of diffusion, respectively. Rearranged, the required element size to eliminate spurious oscillations in the numerical solution is determined by $h \leq \frac{2 \nu}{a}$. Using this, a conservative estimate for stability of the Burgers' equation is obtained by replacing $a$ with the absolute maximum value of $u$ at $t = 0$. Specifically, for $\max{|u(x,0)|} = \max{|u_{IC}(x)|}$, then $h \leq \frac{2 \nu}{\max{|u_{IC}(x)|}}$. Note in the limit $\nu \rightarrow 0$, the required grid size for a stable solution in linear FEM is unachievable.

\section{Numerical results}

\noindent This section presents GFEM solutions of the one-dimensional Burgers' equation. For the following examples, please consider:

\begin{enumerate}
\item All enrichments are shift by their nodal values to retain the physical meaning of the standard FEM DOFs at each node, as well as reduce potential linear dependencies between the FEM and GFEM basis.
\item Special consideration is necessary to integrate the non-polynomial enrichment functions accurately. The computational cost of integrating the enrichments is trivial in the following examples since the elemental matrices are not time-dependent. Such, the following work uses a conservative number of Gaussian quadrature points for each grid refinement. For example, we use ten-point Gaussian quadrature on the most refined meshes considered (approx. $\frac{1}{80}$ element size); while we use sixty-point Gaussian quadrature for the coarsest meshes considered (approx. $\frac{1}{10}$ element size). For problems where the elemental matrices are time-dependent, evaluation of the elemental matrices at each time step using Gaussian quadrature may increase costs considerably. More efficient integration strategies may be beneficial for these problems, such as the Gauss-Laguerre quadrature for exponential functions. 
\item \textit{A priori} error estimates are well-known for polynomial approximation spaces: for $\Omega \subset \mathbb{R}^n$ with Lipschitz boundary, a $p$-degree polynomial solution converges in the $L_2$ and $H_1$ norm versus total degrees of freedom at a theoretical convergence rate of $\frac{p+1}{n}$ and $\frac{p}{n}$, respectively. However, for approximation spaces containing solution-tailored enrichments, theoretical convergence rates are not formally developed. Insights into convergence rates for solution-tailored enrichments are provided by considering convergence plots. Unless specified, convergence rates in the $L_2$ and $H_1$ norm versus total degrees of freedom use the finest two grids studied. For a sufficiently smooth solution using polynomial + non-polynomial enrichments, convergence rates are similar to those of the polynomial approximation spaces. An exception is when the numerical solution is of the same order of numerical precision as the reference solution. Same order numerical precision is often the case when using solution-tailored enrichments.
\end{enumerate}

\subsection{Numerical example 1: Boundary layer solution as the kinematic viscosity tends toward zero}
\subsubsection{Problem statement and reference solutions}
Consider the viscous Burgers' equation (Eq. \ref{viscous_burger_strong_formulation}) defined over a unit domain ($\Omega = [0, 1]$) and subject to homogeneous Dirichlet boundary conditions everywhere ($\Gamma = \Gamma_D = 0; \Gamma_N = \emptyset$). The problem formulation is as follows: For $t \in [0, 1]$, find $u$ such that:

\begin{equation}
\label{example1_viscous_burger_strong_formulation}
\begin{aligned}
\frac{\partial u}{\partial t}  + u \frac{\partial u}{\partial x} - \nu \frac{\partial^2 u}{\partial x^2} = 0& \quad \mbox{on} \quad \Omega\\
u(x,0) = \sin{\pi x}& \quad \mbox{on} \quad \Omega\\
u(0,t) = u(1,t) = 0& \quad \mbox{on} \quad \Gamma\\
\end{aligned}
\end{equation}

An analytical Fourier solution to Eq. \ref{example1_viscous_burger_strong_formulation} is obtainable through use of the Hopf-Cole transformation, as detailed in \citep{Kutluay1999}. The resulting analytical Fourier solution is:

\begin{equation}
\label{example1_analytical_fourier_solution}
u(x,t) = 2\pi\nu \frac{\sum_{n = 1}^{\infty} a_n e^{-n^2\pi^2\nu t} n \sin{n\pi x} }{a_0 + \sum_{n = 1}^{\infty} a_n e^{-n^2\pi^2\nu t} \cos{n\pi x}}
\end{equation}

\noindent with Fourier coefficients, $a_n$:

\begin{equation}
\label{example1_fourier_coefficients}
\begin{aligned}
&a_0 = \int_0^1 e^{- \frac{1}{2 \pi \nu} (1-\cos{\pi x})} \, dx\\
&a_n = 2 \int_0^1 e^{- \frac{1}{2 \pi \nu} (1-\cos{\pi x})} \cos{n \pi x} \, dx
\end{aligned}
\end{equation}

\begin{figure}[ht!]
\begin{center}
\begin{subfigmatrix}{4}
\subfigure[$\nu = \frac{1}{10}$]{\includegraphics[width=1.6in]{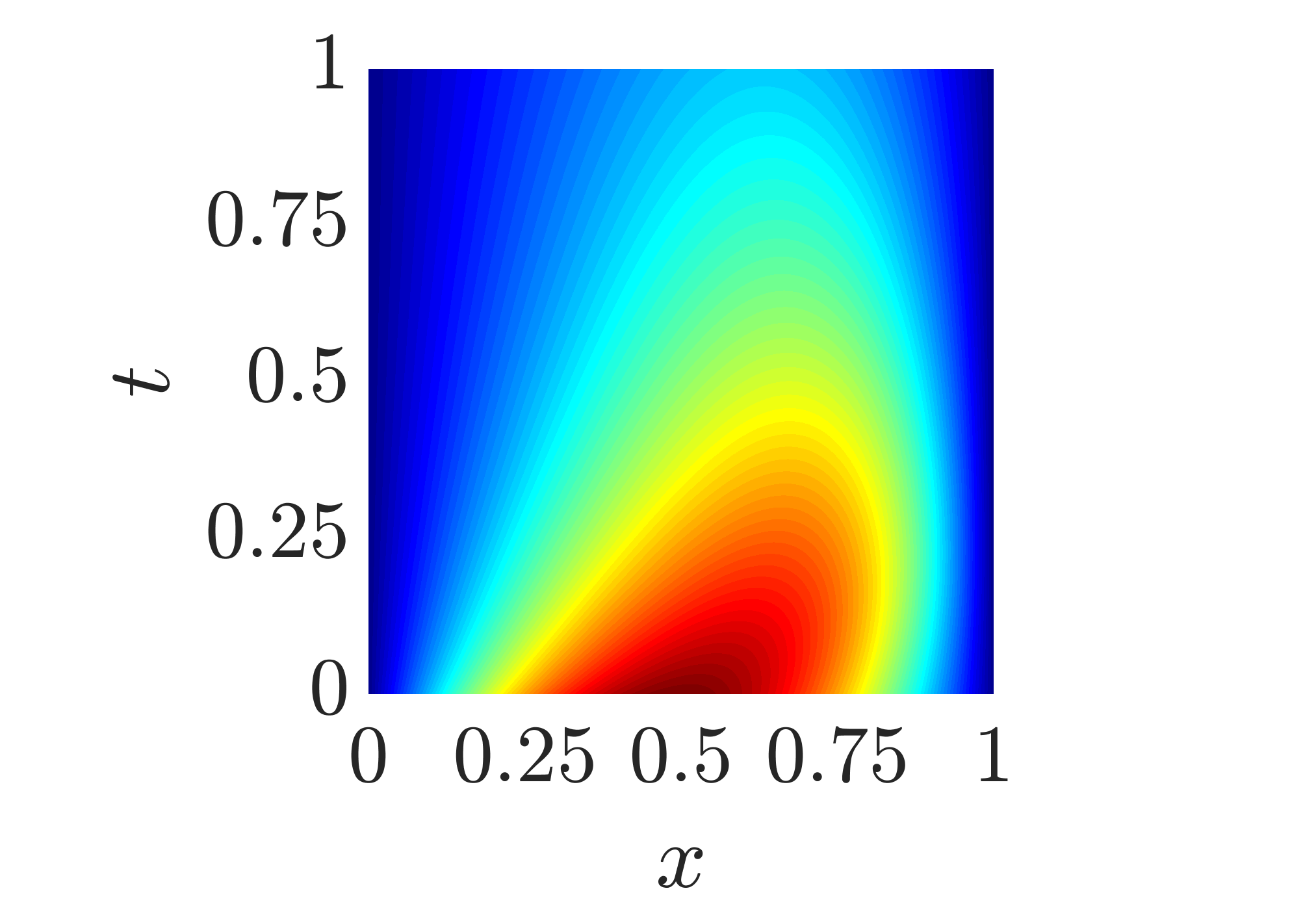}}
\subfigure[$\nu = \frac{1}{50}$]{\includegraphics[width=1.6in]{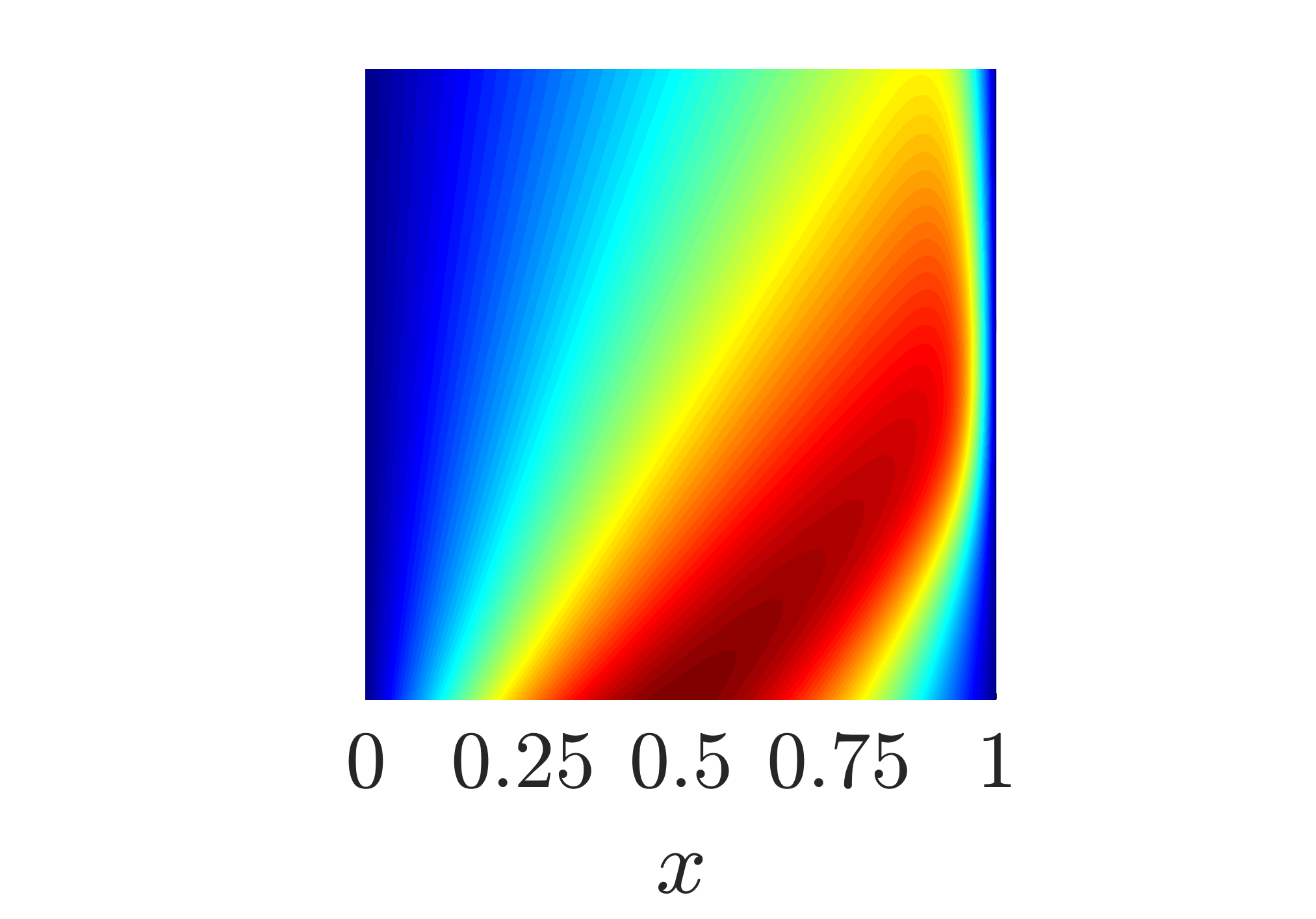}}
\subfigure[$\nu = \frac{1}{100}$]{\includegraphics[width=1.6in]{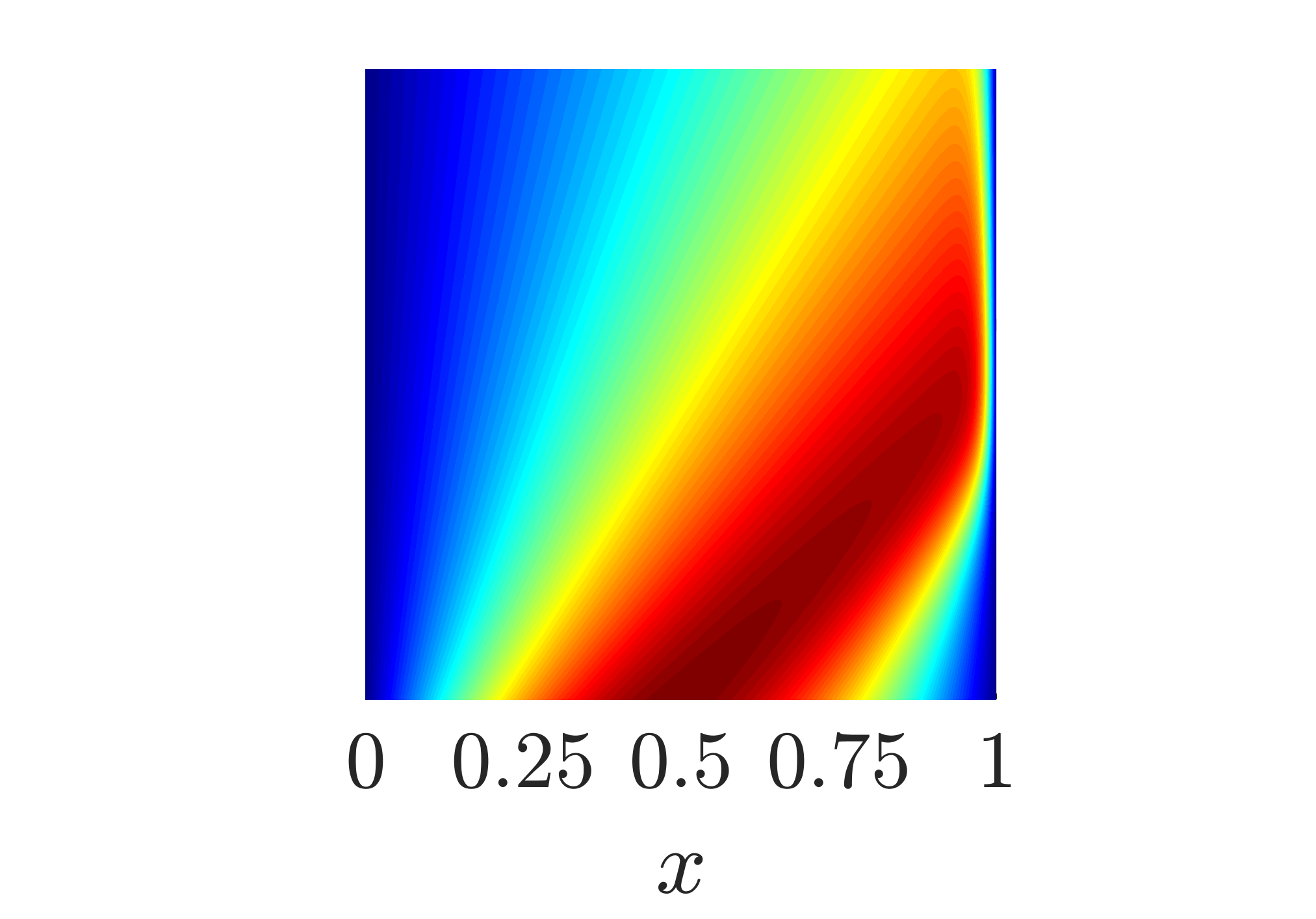}}
\subfigure[$\nu = 0$]{\includegraphics[width=1.6in]{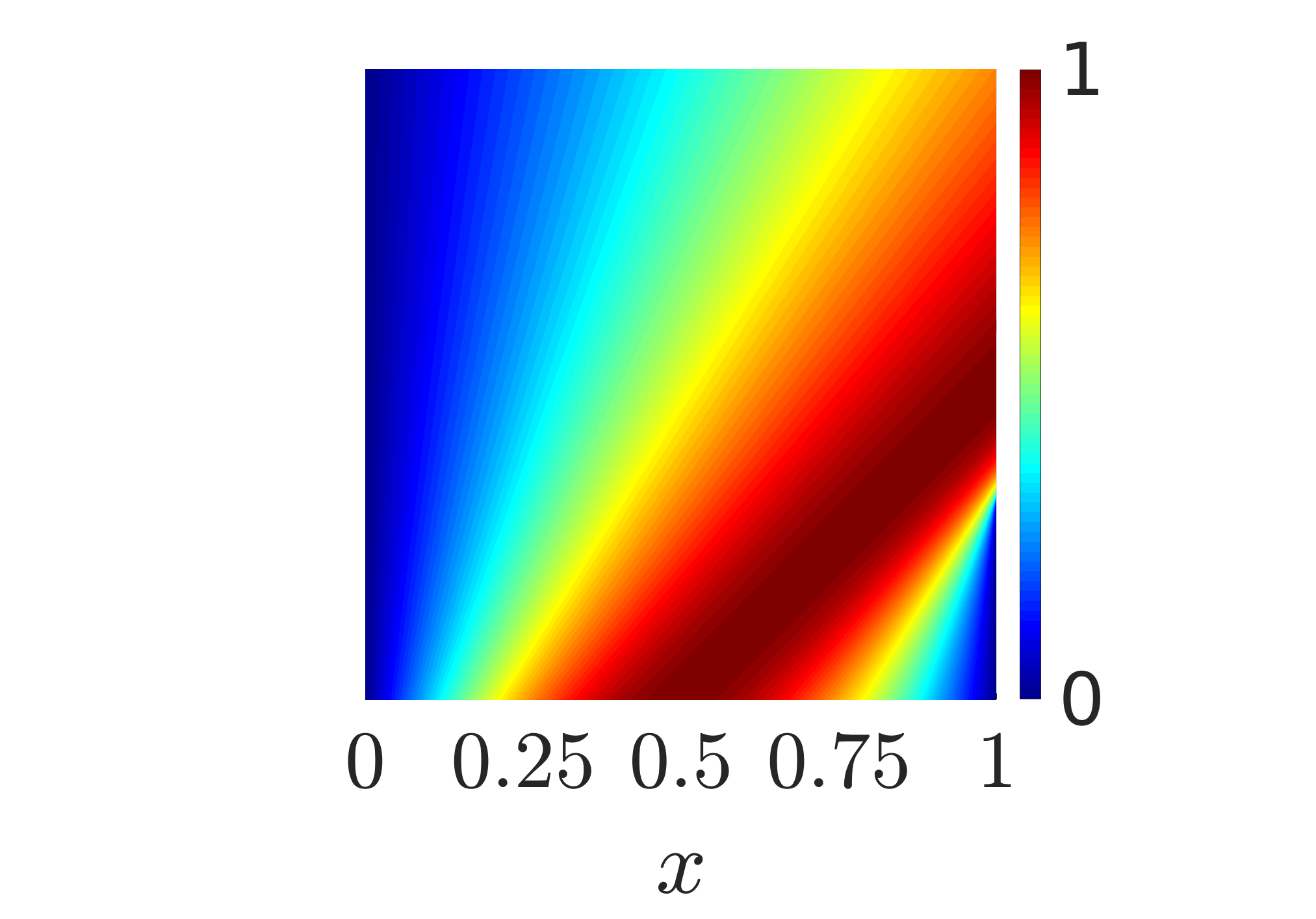}}
\end{subfigmatrix}
\caption{Reference solution contours for the boundary layer problem over a range of $\nu$}
\label{fig:example1_reference_solution}
\end{center}
\end{figure}

The integrals of Eq. \ref{example1_fourier_coefficients} are  convergent for all $\nu \neq 0$. However, for small values of $\nu$ and $t$, the rate of convergence of the series slows down significantly, and results in extremely difficulty computing $u$ using this analytical expression \citep{Dey1983}. A good discussion on this convergence issue is provided by \citep{Ozis1996} and the references within. Since this work concerns solution-tailored numerical solutions with very small viscosities, the poor accuracy of the truncated series may affect convergence rates. Thus, a 5000-element, $p = 1$ FEM solution is used as a reference instead, with the Crank-Nicolson scheme implemented with a step size of $\Delta t = \frac{1}{5000}$. For $\nu = 0$, $p = 1$ FEM is incapable of obtaining a convergent solution (recall the conservative estimate for stability $h \leq \frac{2 \nu}{\max{u_{IC}(x)}} = 0$; or see Figs. \ref{fig:Example1_p1FEM_L2norm} and \ref{fig:Example1_p1FEM_H1norm}). Thus, for $\nu = 0$ the analytical solution to the inviscid problem (Eq. \ref{inviscid_burger_strong_formulation}) is used as a reference, and is found by solving for $u$ in the implicit equation $u = u_{IC}(x - ut)$. Reference solutions for kinematic viscosities $\nu = \Big\{\frac{1}{10}, \frac{1}{50}, \frac{1}{100}, 0\Big\}$ are shown in Fig. \ref{fig:example1_reference_solution}. A boundary layer forms near $x = 1$ with thickness decreasing as $\nu$ decreases. When $\nu = 0$, a discontinuity forms at $x = 1$, starting at time $t_b = \frac{1}{\pi}$, and persisting through later times. 

\subsubsection{Finite element solutions}
\begin{figure}[ht!]
\begin{center}
\begin{subfigmatrix}{6}
\subfigure[$t = 0$]{\includegraphics[width=2.1in]{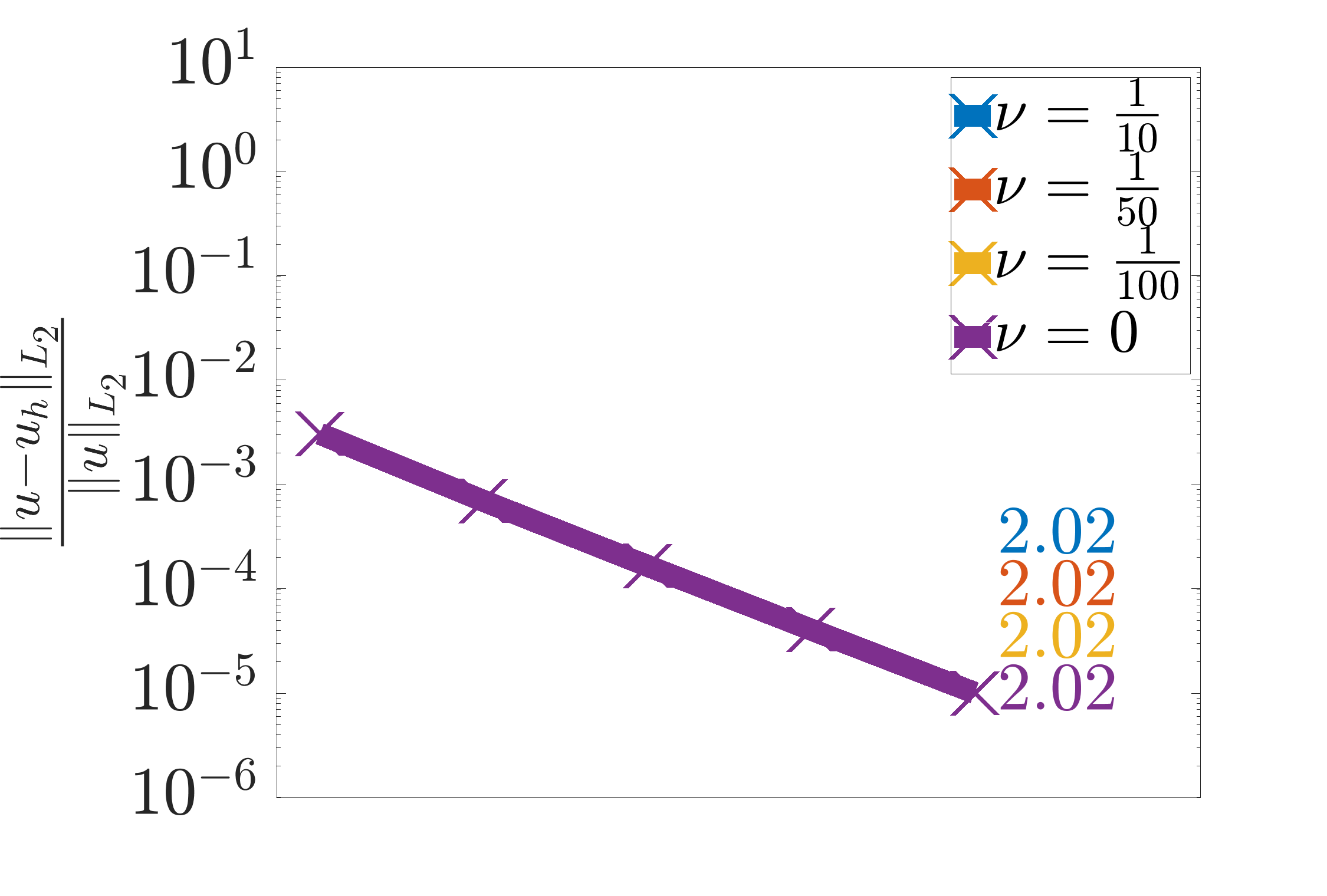}}
\subfigure[$t = 0.25$]{\includegraphics[width=2.1in]{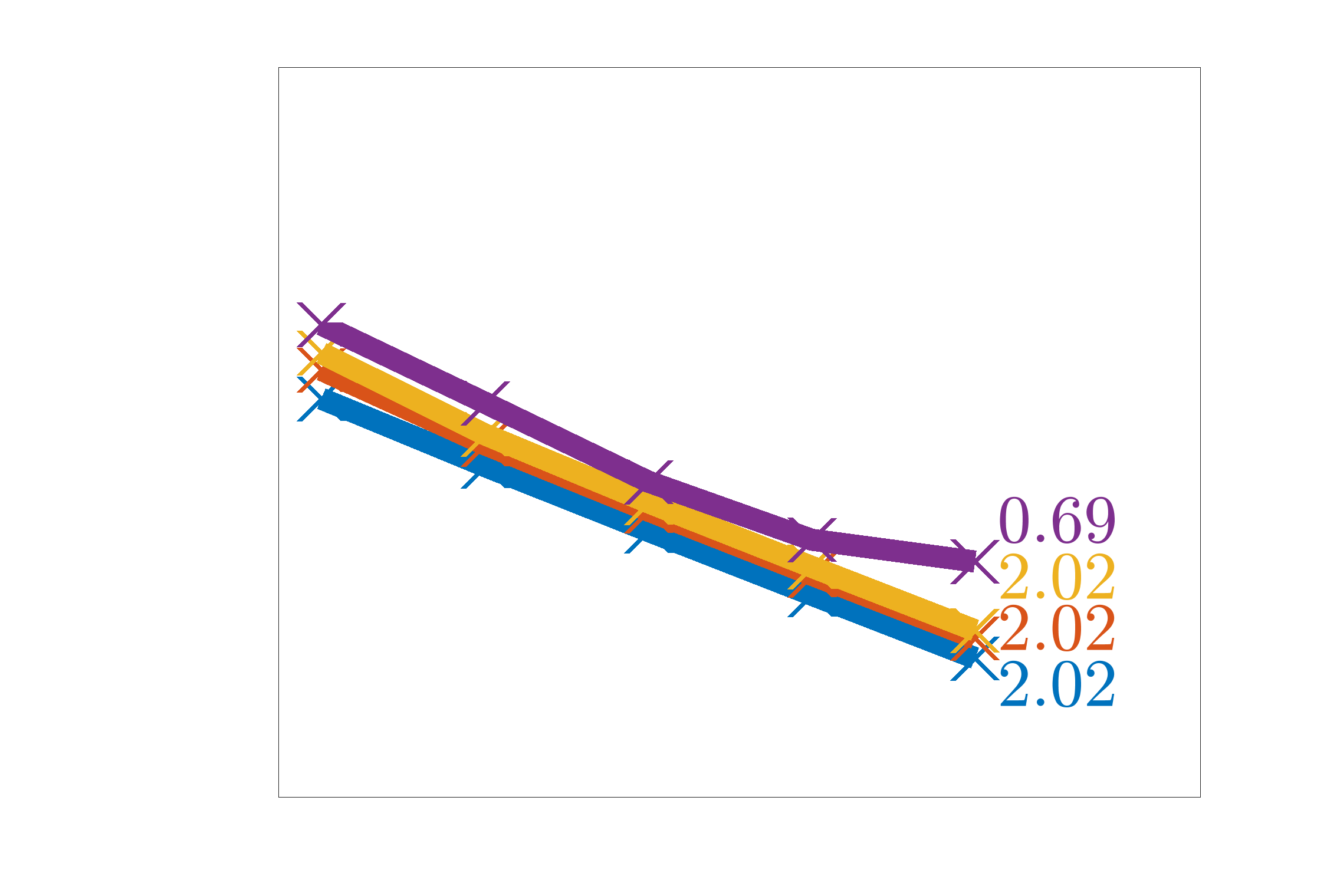}}
\subfigure[$t = 0.318 \approx \frac{1}{\pi}$]{\includegraphics[width=2.1in]{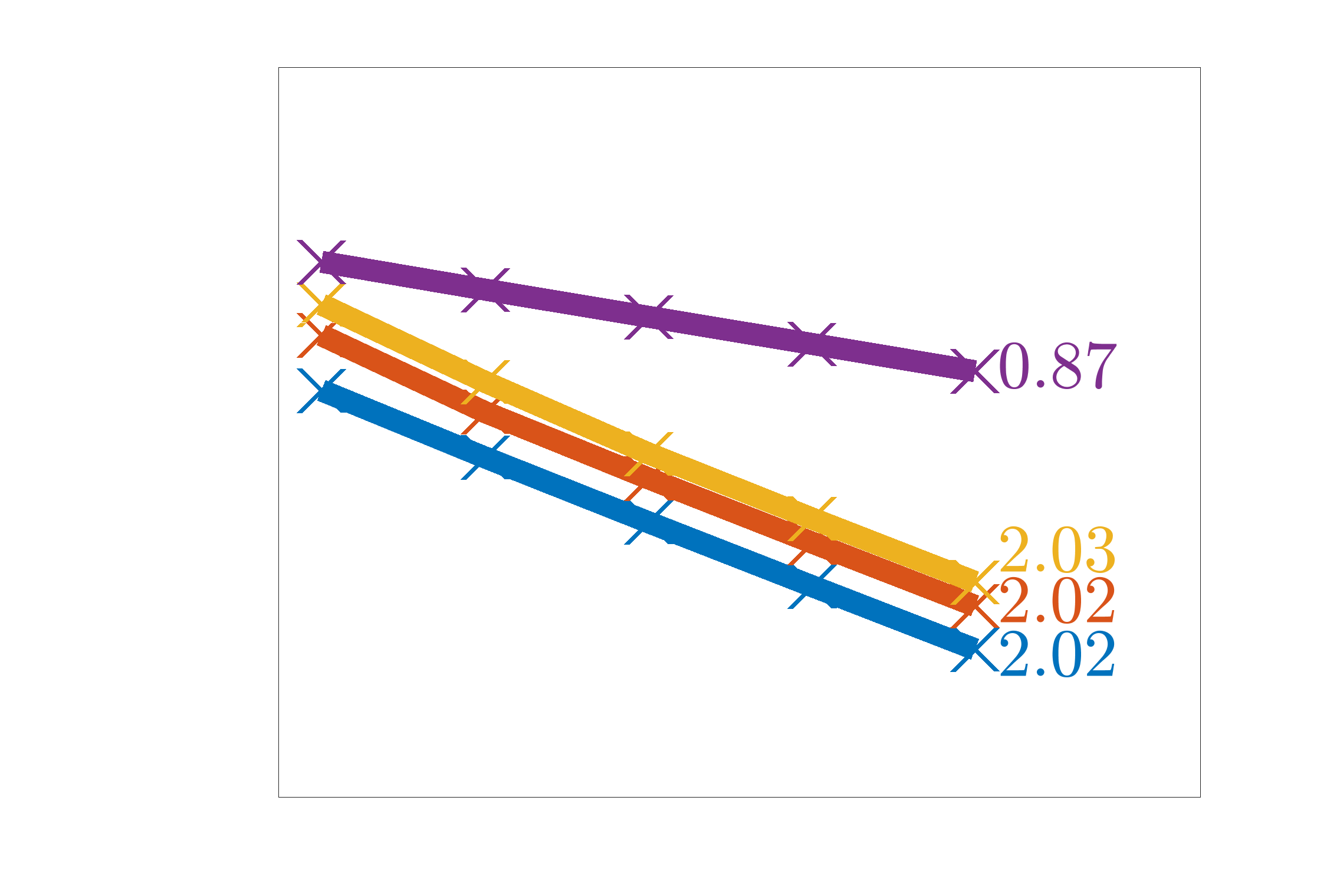}}
\subfigure[$t = 0.5$]{\includegraphics[width=2.1in]{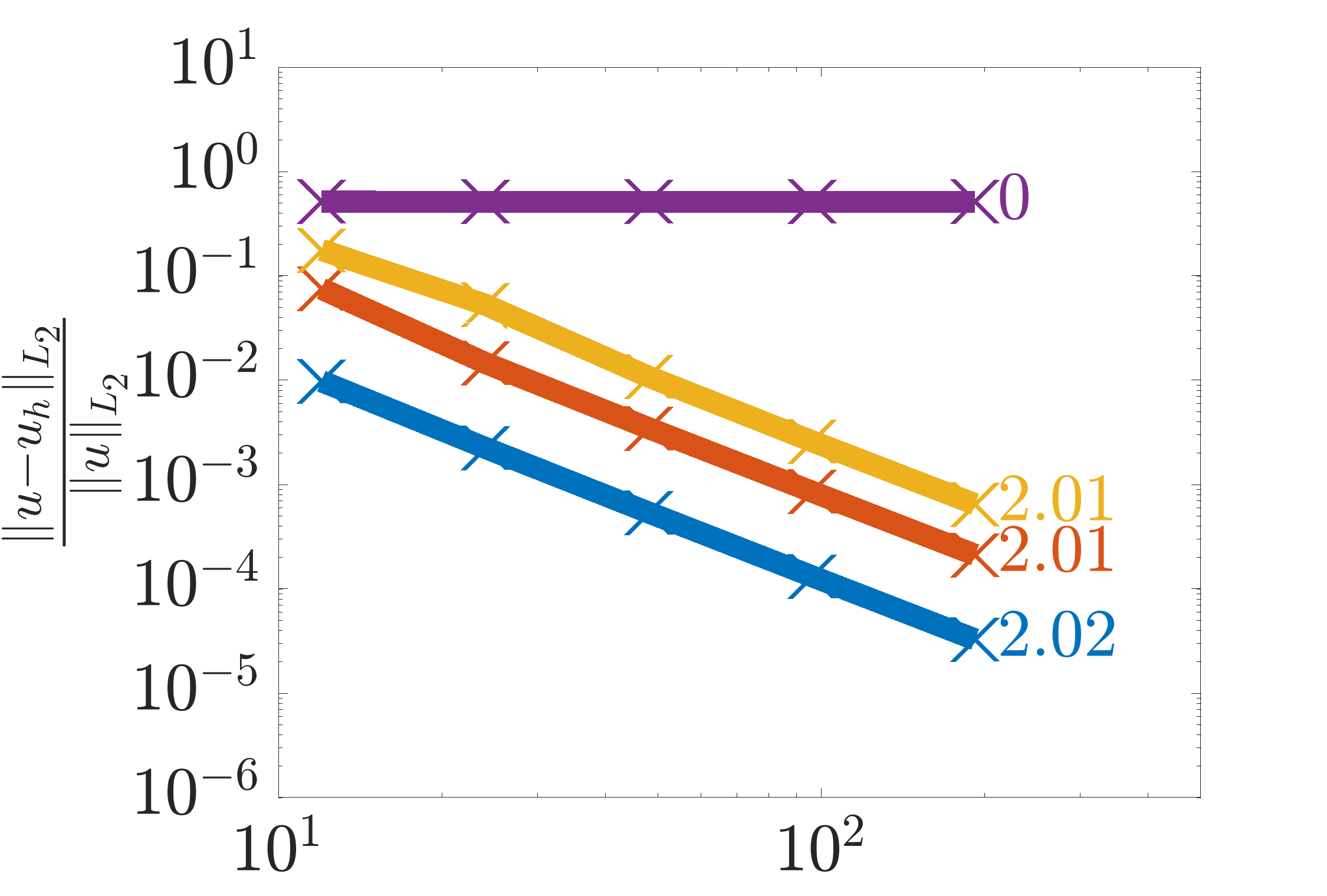}}
\subfigure[$t = 0.75$]{\includegraphics[width=2.1in]{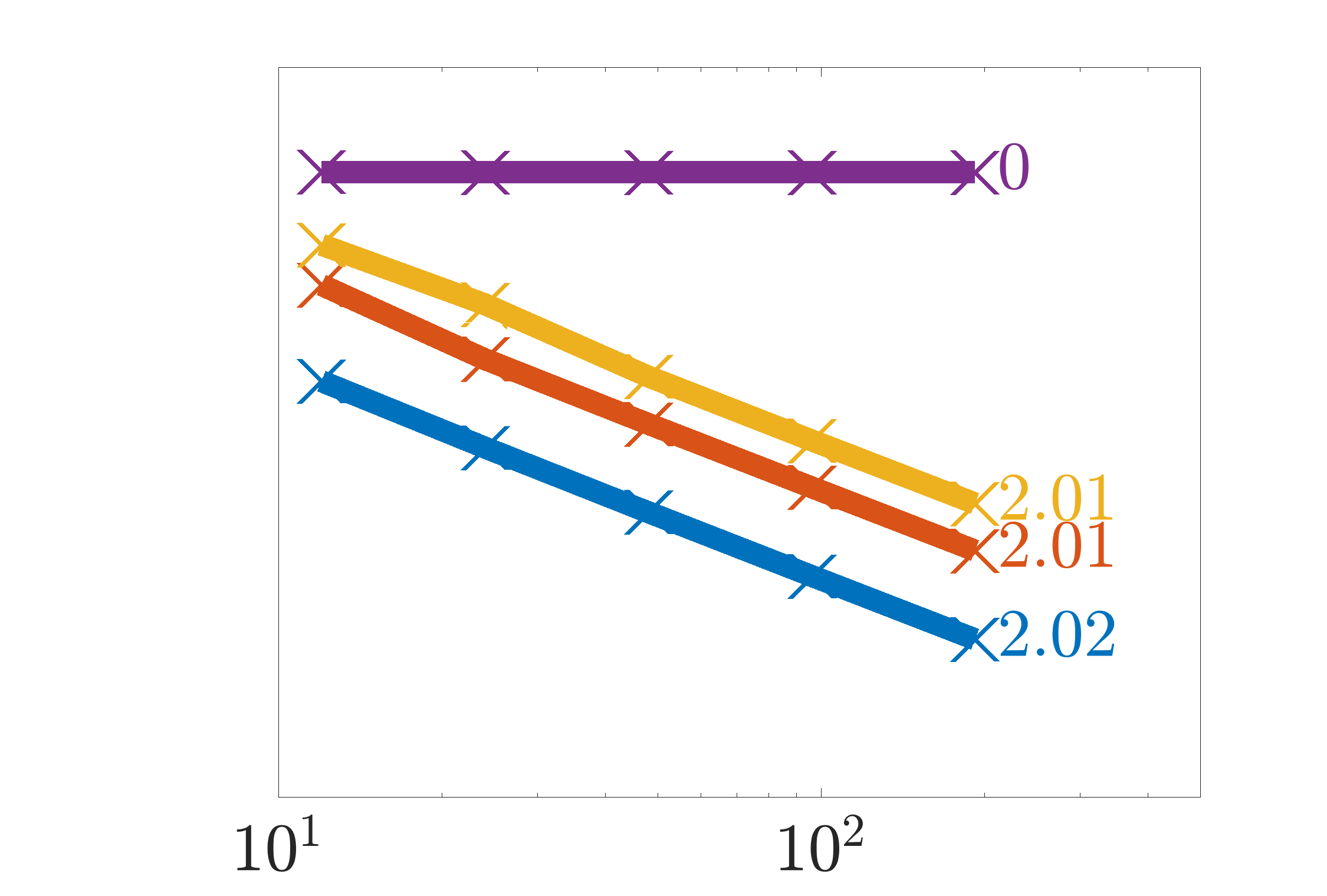}}
\subfigure[$t = 1$]{\includegraphics[width=2.1in]{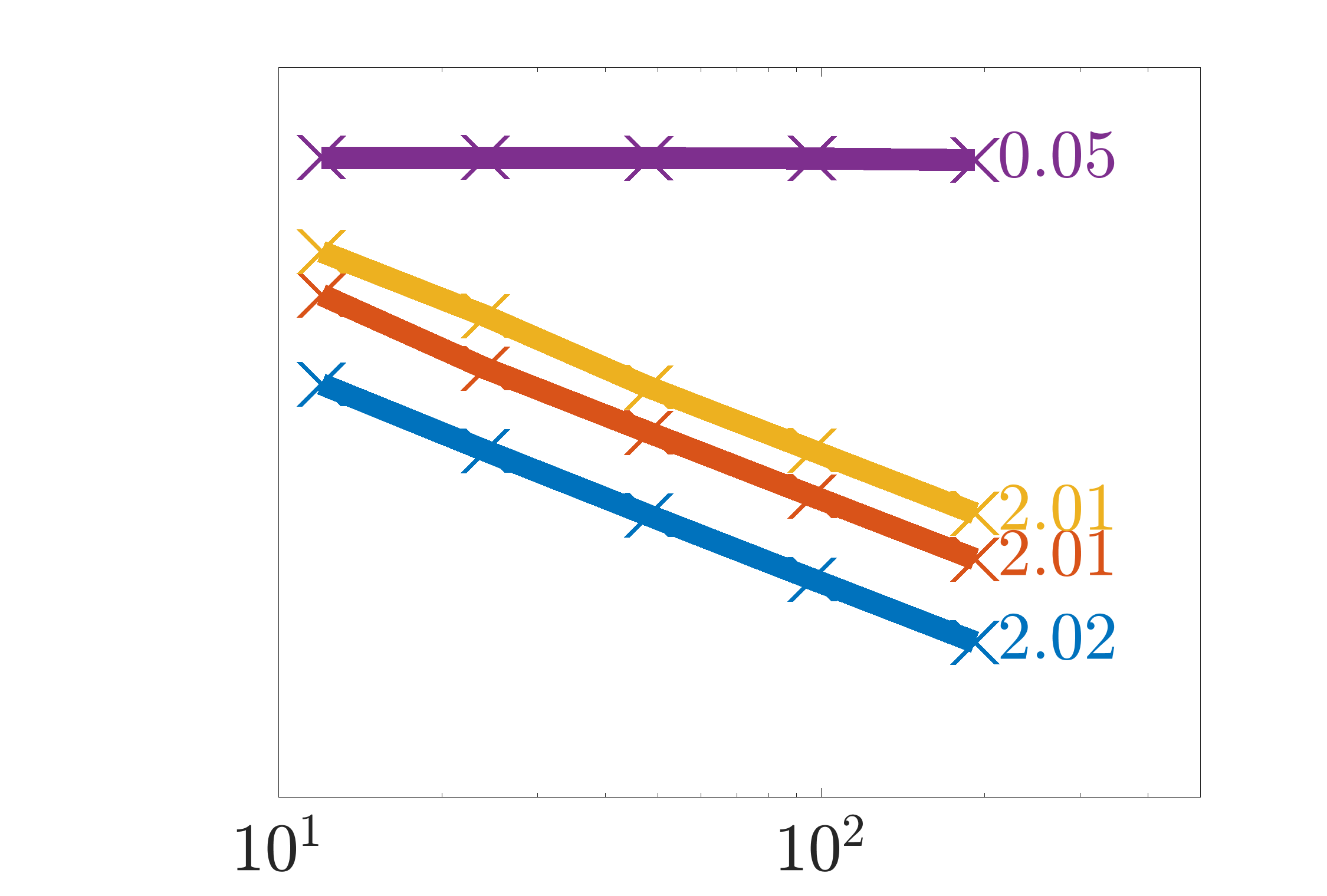}}
\end{subfigmatrix}
\caption{$p = 1$ FEM convergence in the relative $L_2$ norm for the boundary layer problem at various times}
\label{fig:Example1_p1FEM_L2norm}
\end{center}
\end{figure}

\begin{figure}[ht!]
\begin{center}
\begin{subfigmatrix}{6}
\subfigure[$t = 0$]{\includegraphics[width=2.1in]{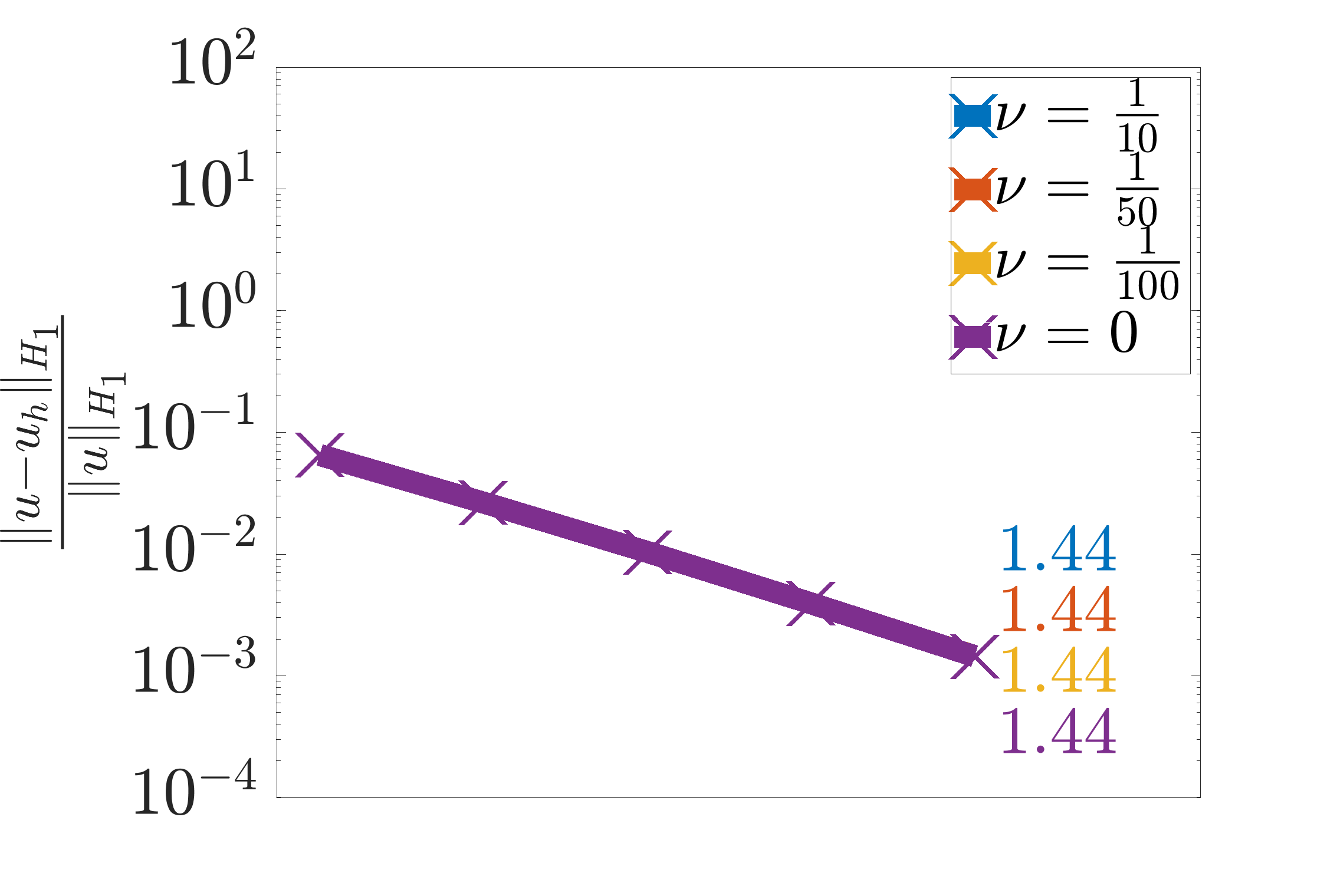}}
\subfigure[$t = 0.25$]{\includegraphics[width=2.1in]{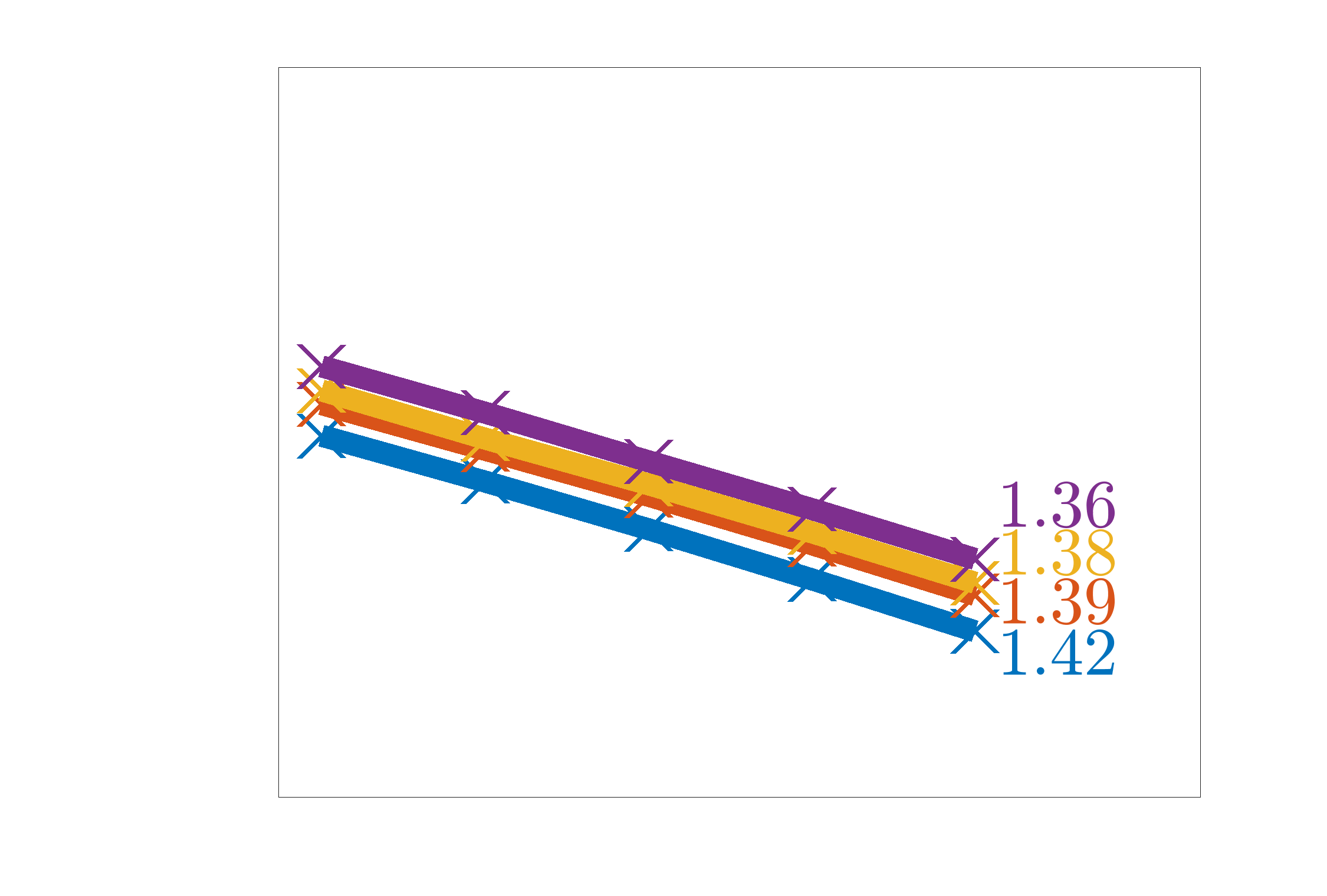}}
\subfigure[$t = 0.318 \approx \frac{1}{\pi}$]{\includegraphics[width=2.1in]{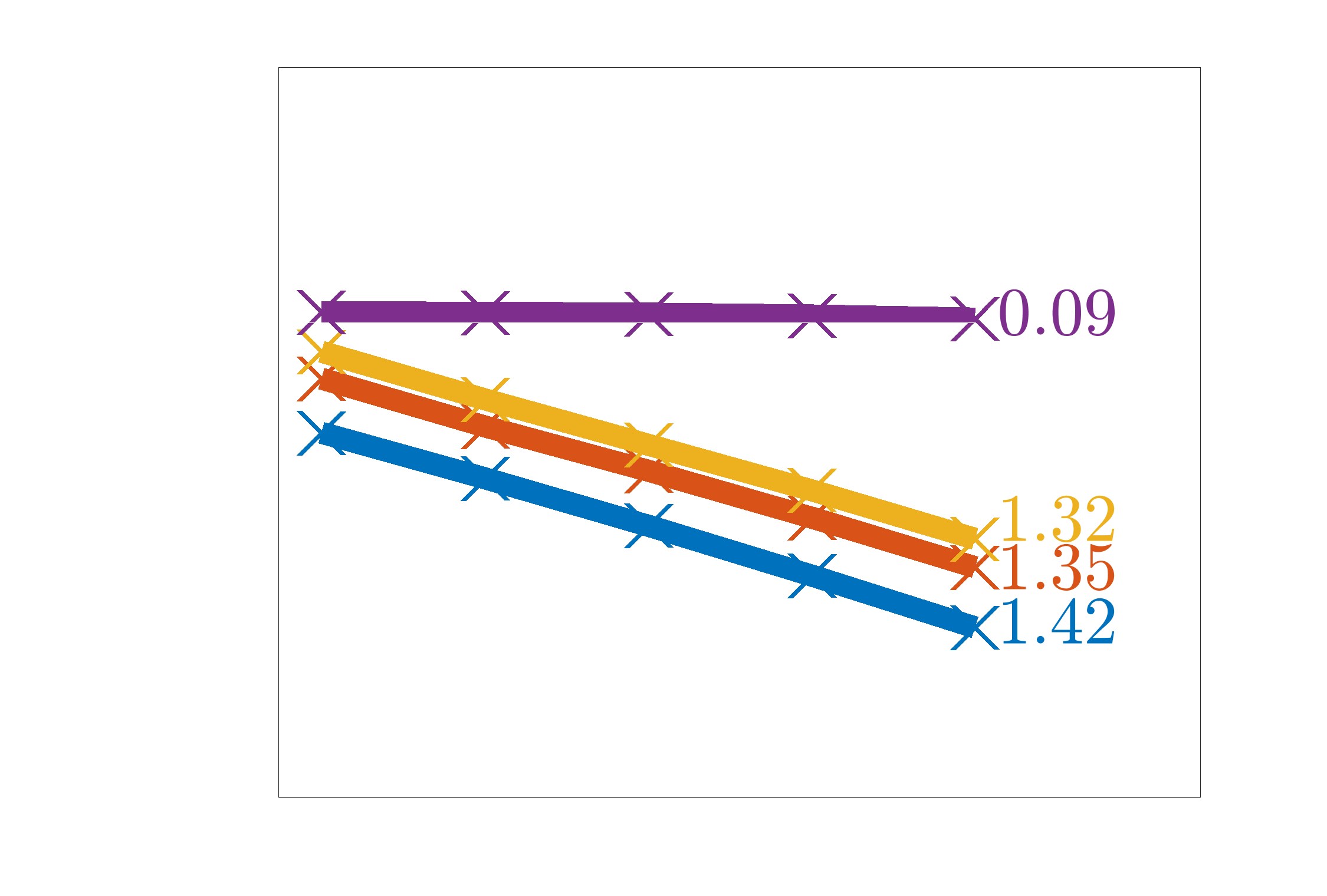}}
\subfigure[$t = 0.5$]{\includegraphics[width=2.1in]{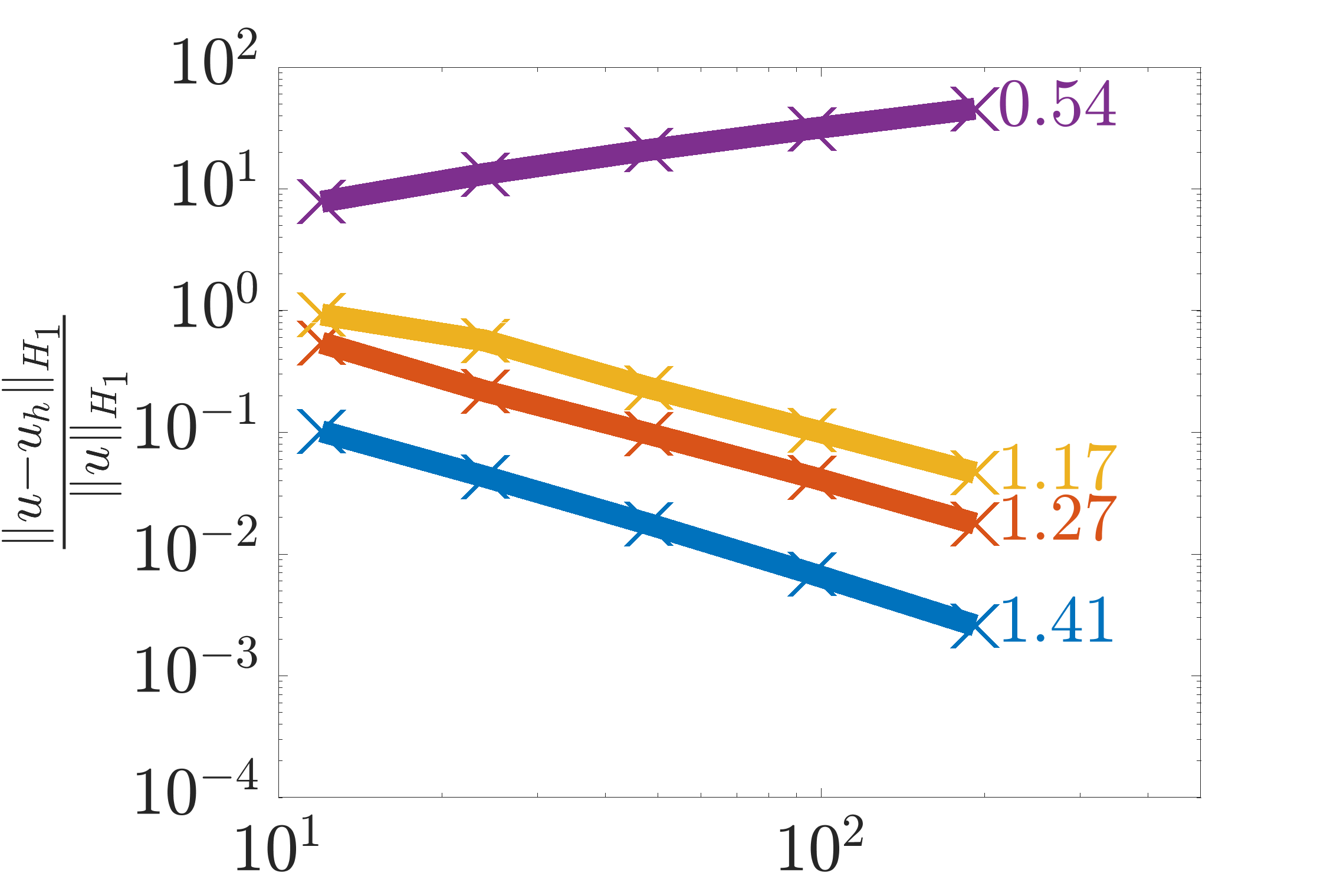}}
\subfigure[$t = 0.75$]{\includegraphics[width=2.1in]{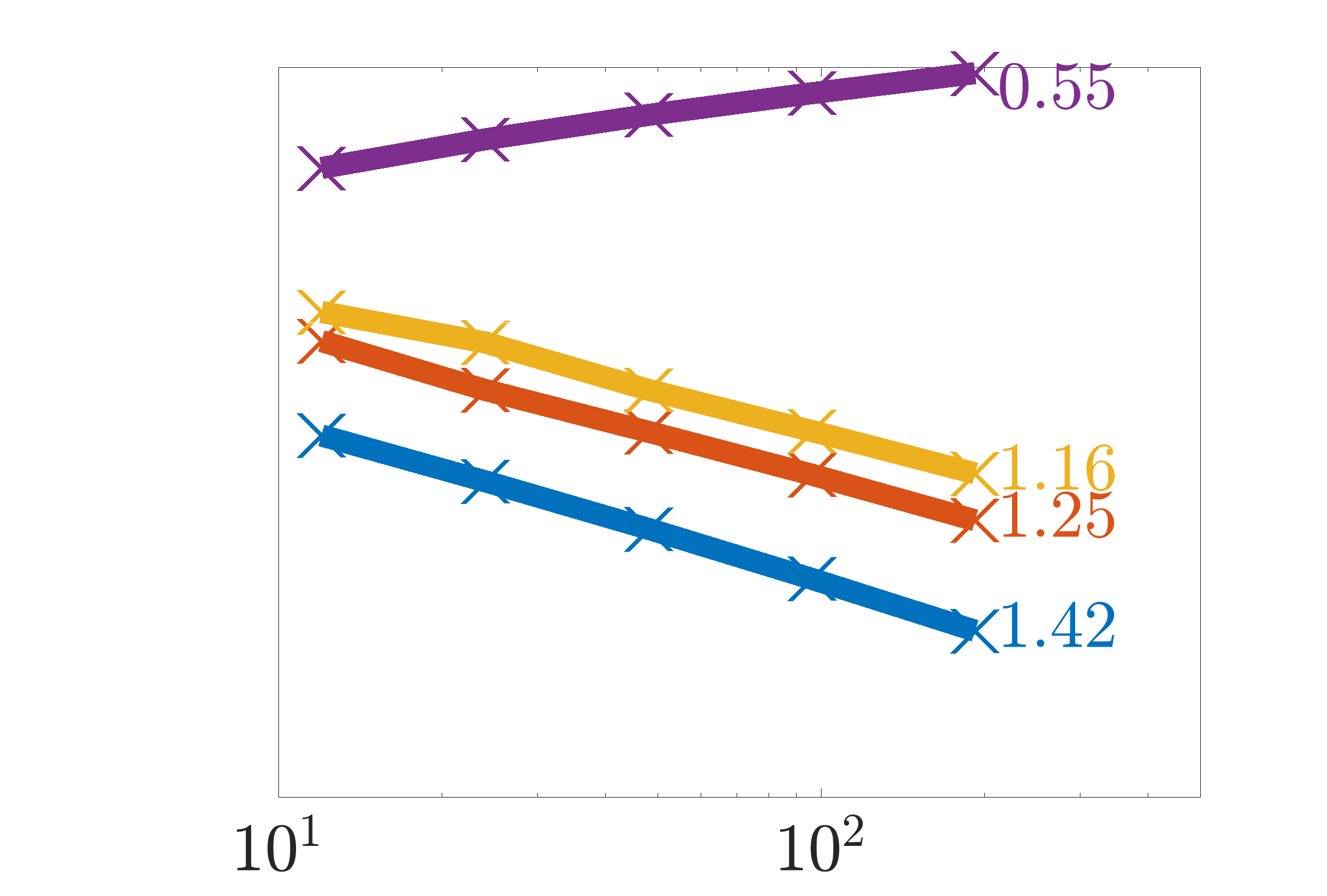}}
\subfigure[$t = 1$]{\includegraphics[width=2.1in]{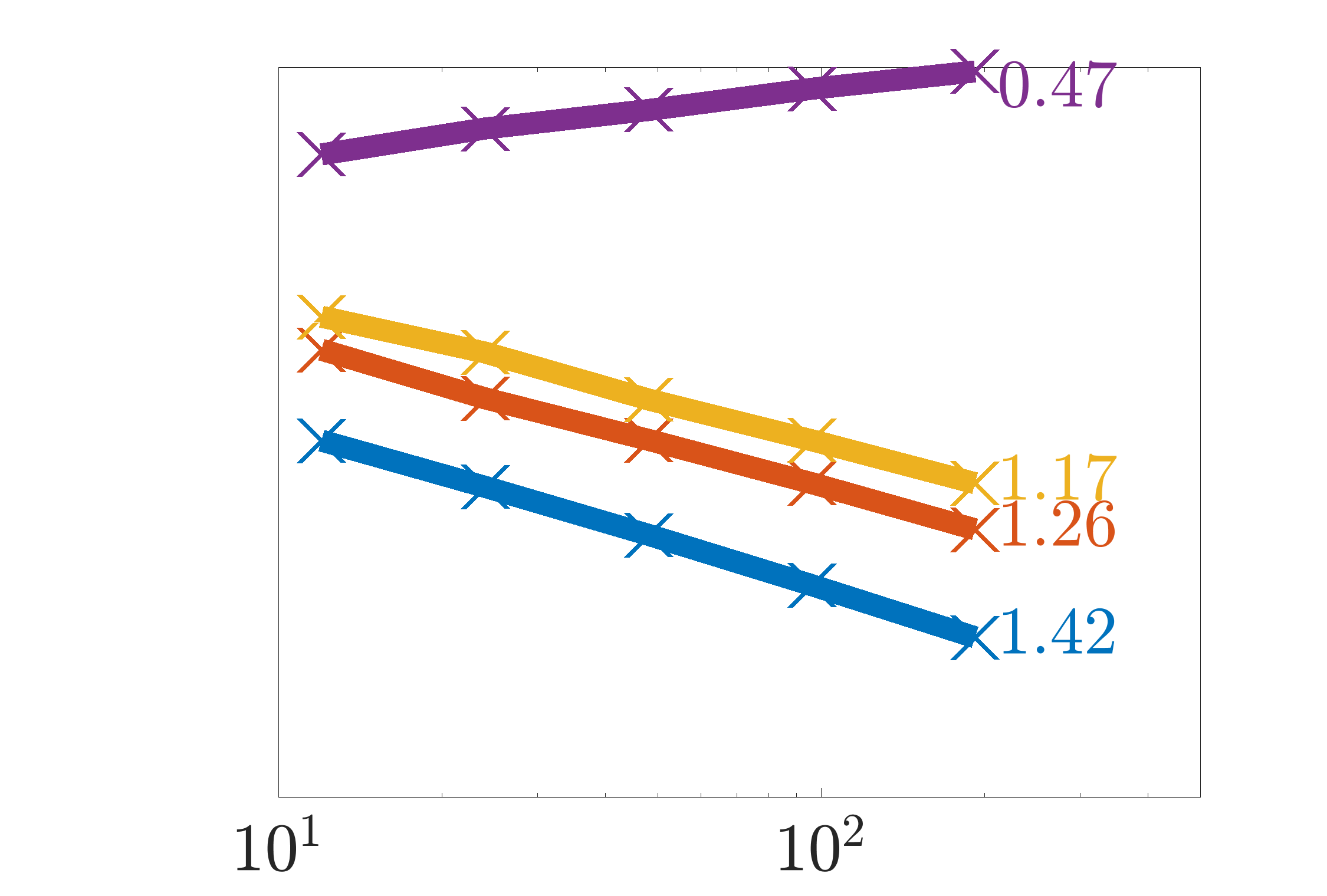}}
\end{subfigmatrix}
\caption{$p = 1$ FEM convergence in the relative $H_1$ norm for the boundary layer problem at various times}
\label{fig:Example1_p1FEM_H1norm}
\end{center}
\end{figure}

\begin{figure}[ht!]
\begin{center}
\begin{subfigmatrix}{4}
\subfigure[$\nu = \frac{1}{10}$]{\includegraphics[width=1.6in]{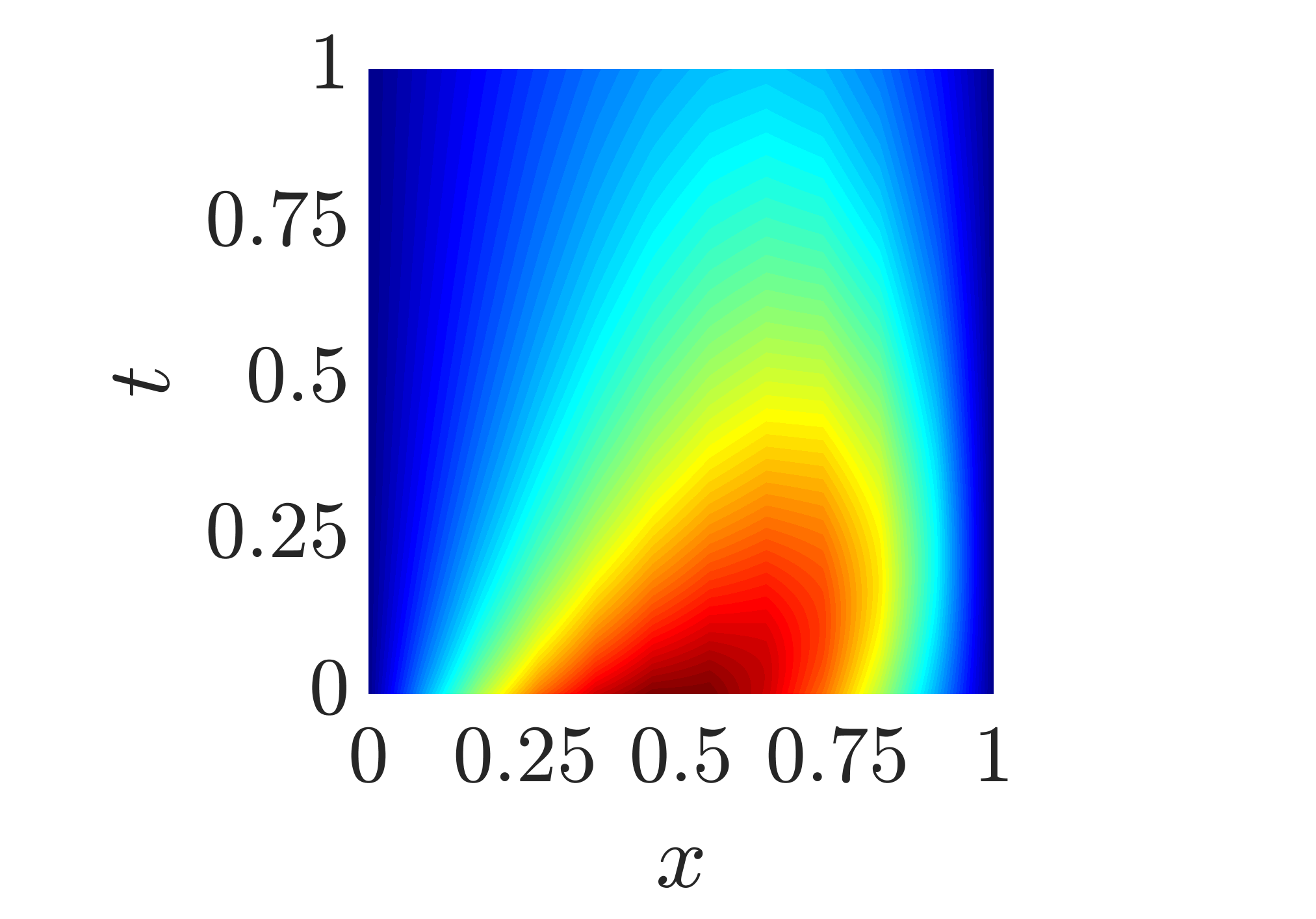}}
\subfigure[$\nu = \frac{1}{50}$]{\includegraphics[width=1.6in]{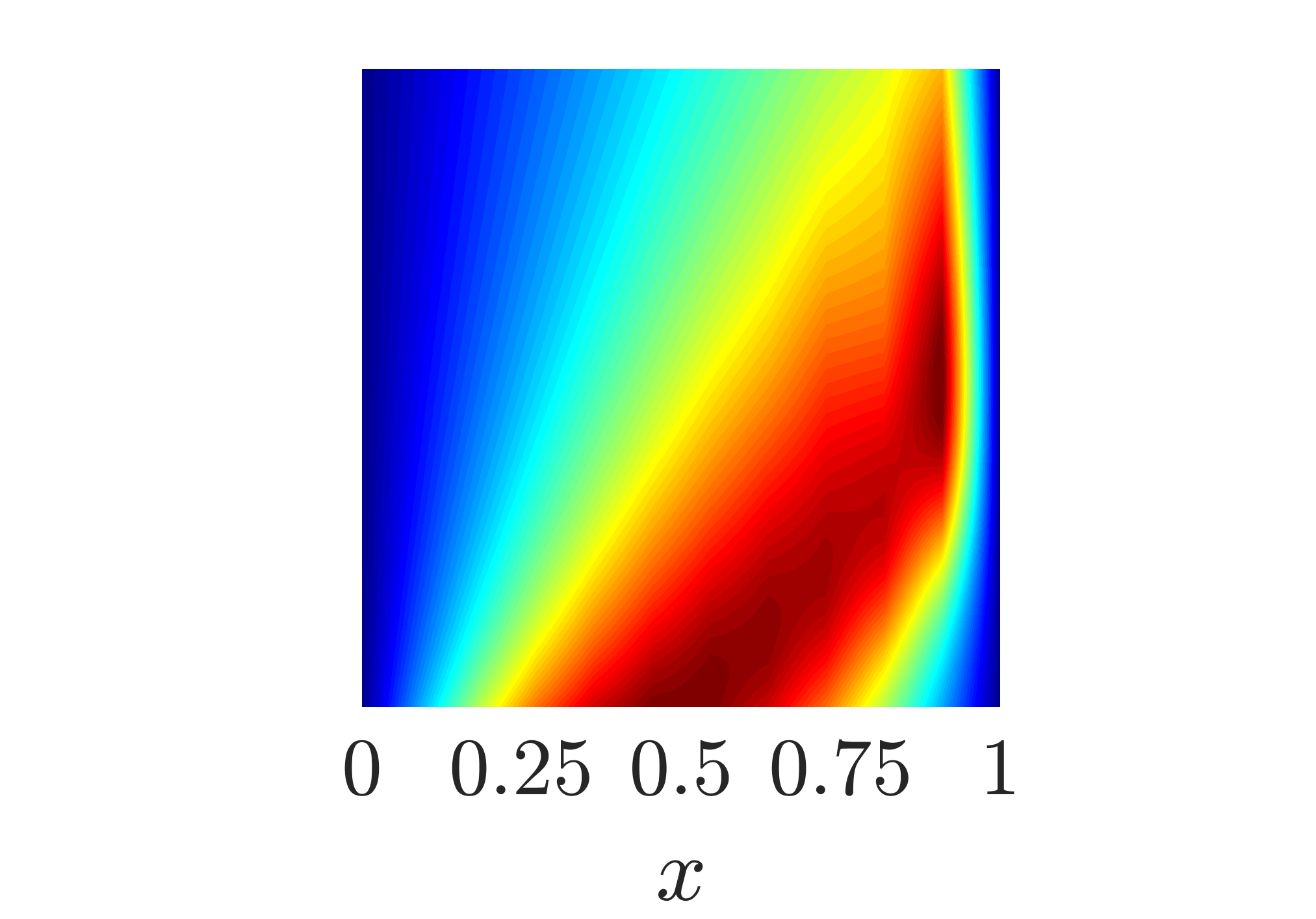}}
\subfigure[$\nu = \frac{1}{100}$]{\includegraphics[width=1.6in]{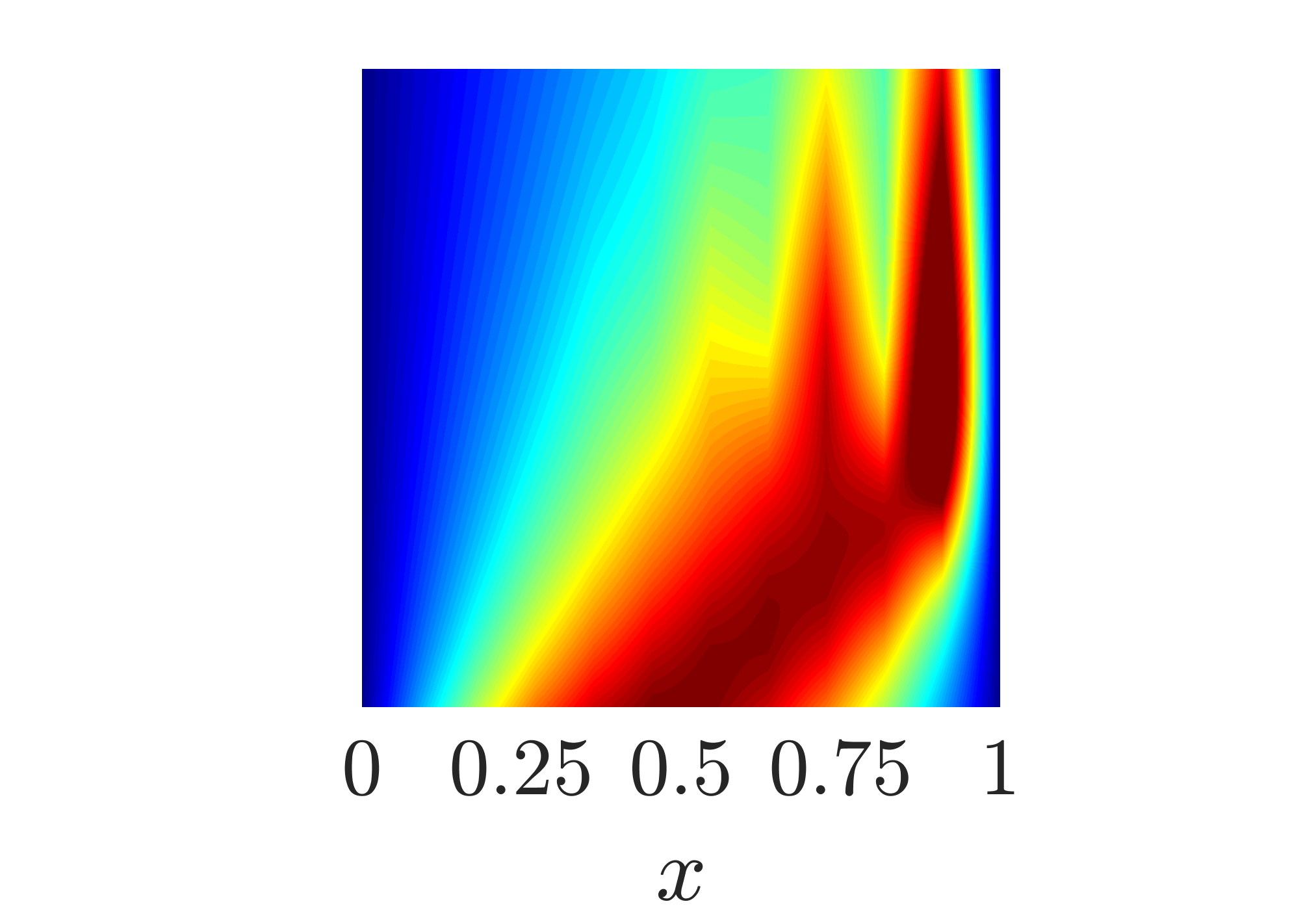}}
\subfigure[$\nu = 0$]{\includegraphics[width=1.6in]{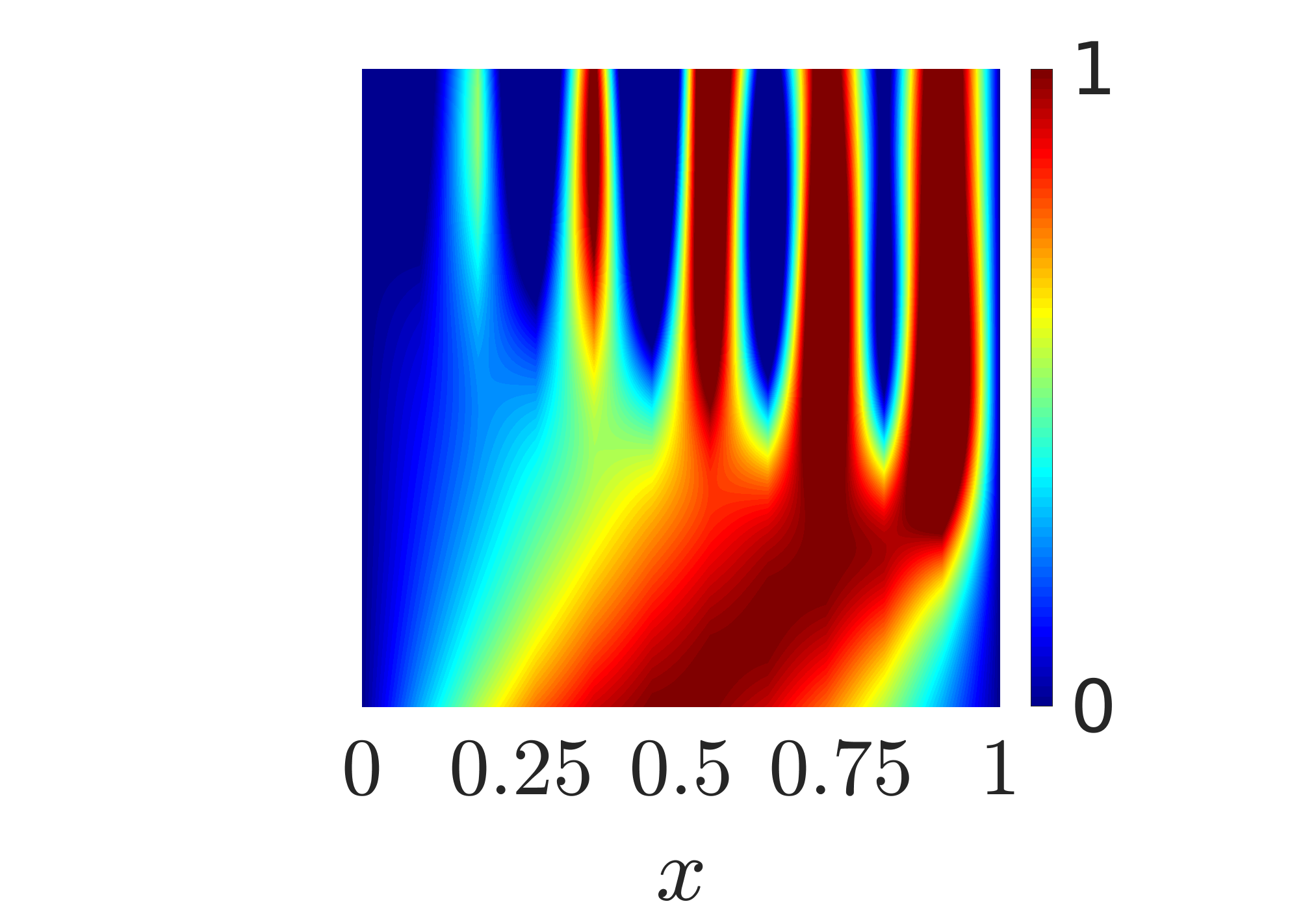}}
\end{subfigmatrix}
\caption{11-element, $p = 1$ FEM solution contours for the boundary layer problem over a range of $\nu$}
\label{fig:Example1_11_element_FEM_soln}
\end{center}
\end{figure}

\begin{figure}[ht!]
\begin{center}
\begin{subfigmatrix}{4}
\subfigure[$\nu = \frac{1}{10}$]{\includegraphics[width=1.6in]{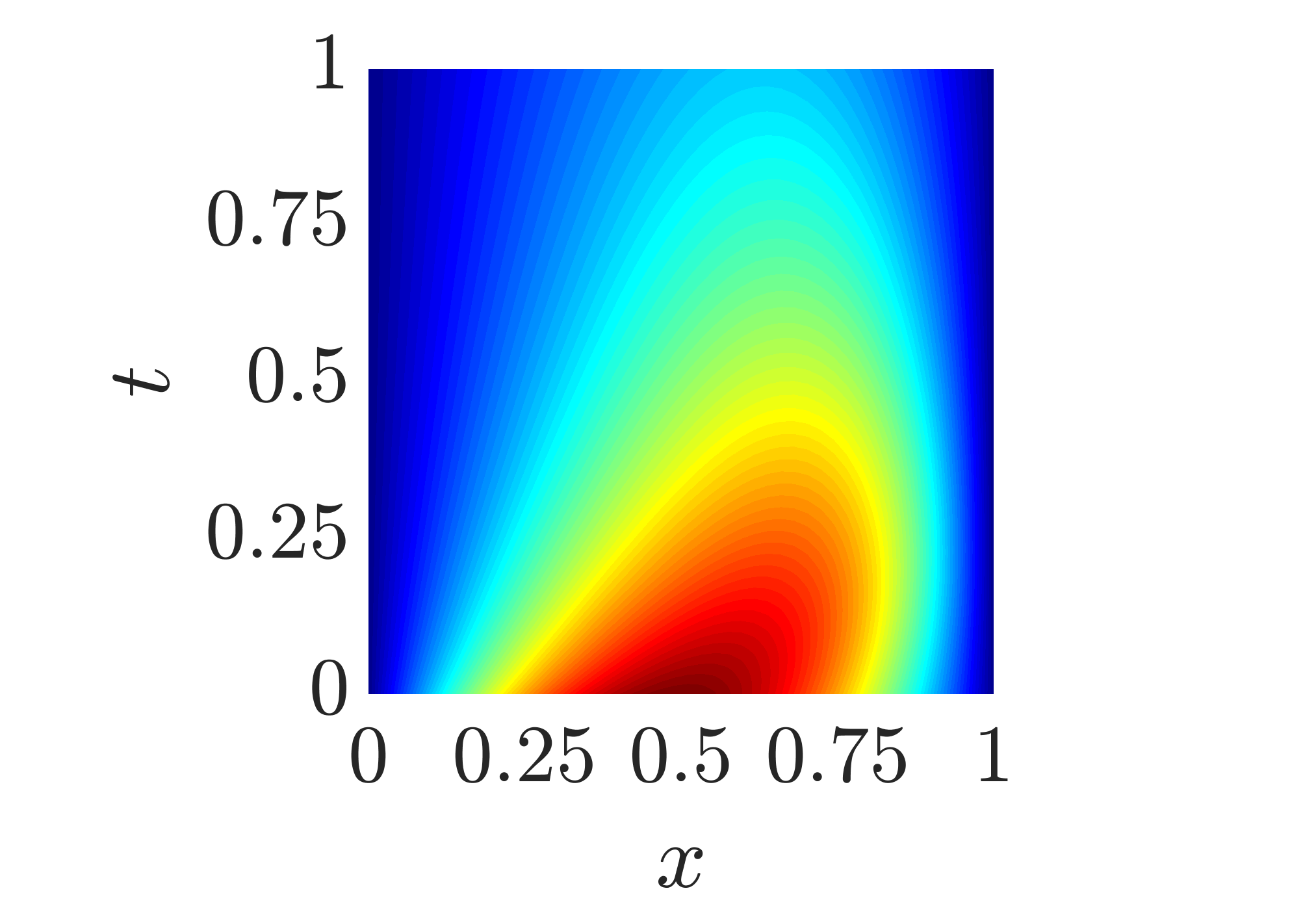}}
\subfigure[$\nu = \frac{1}{50}$]{\includegraphics[width=1.6in]{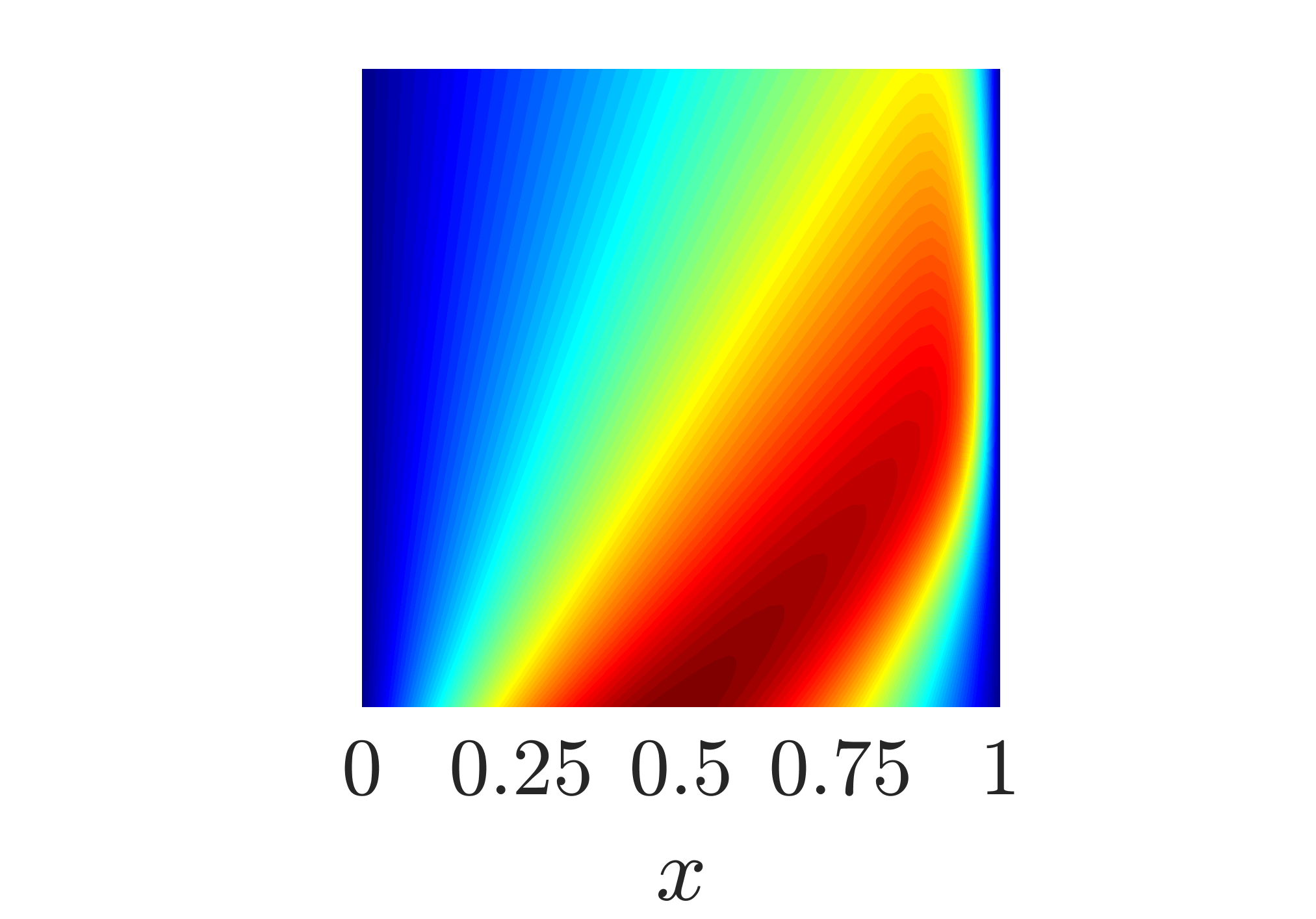}}
\subfigure[$\nu = \frac{1}{100}$]{\includegraphics[width=1.6in]{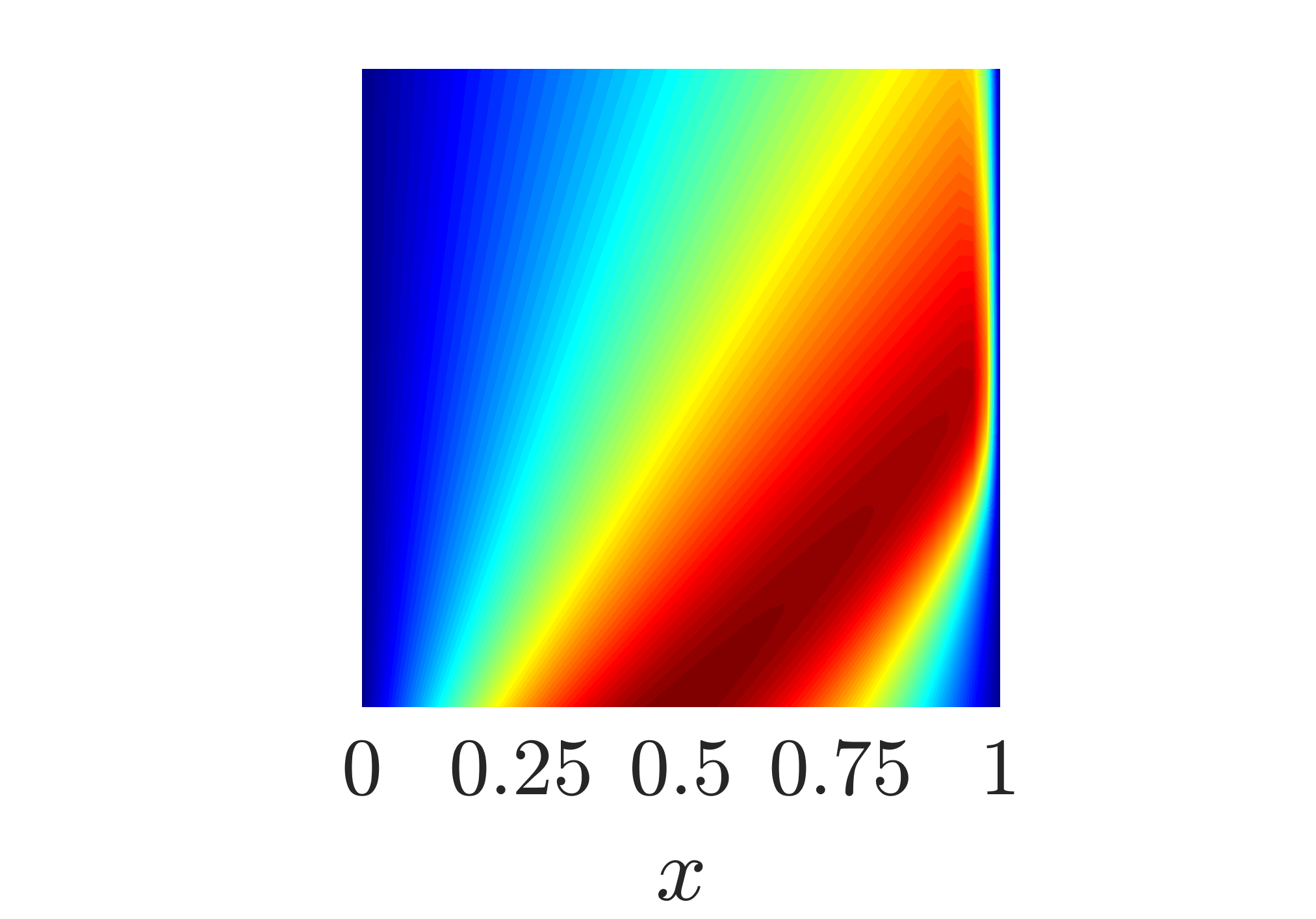}}
\subfigure[$\nu = 0$]{\includegraphics[width=1.6in]{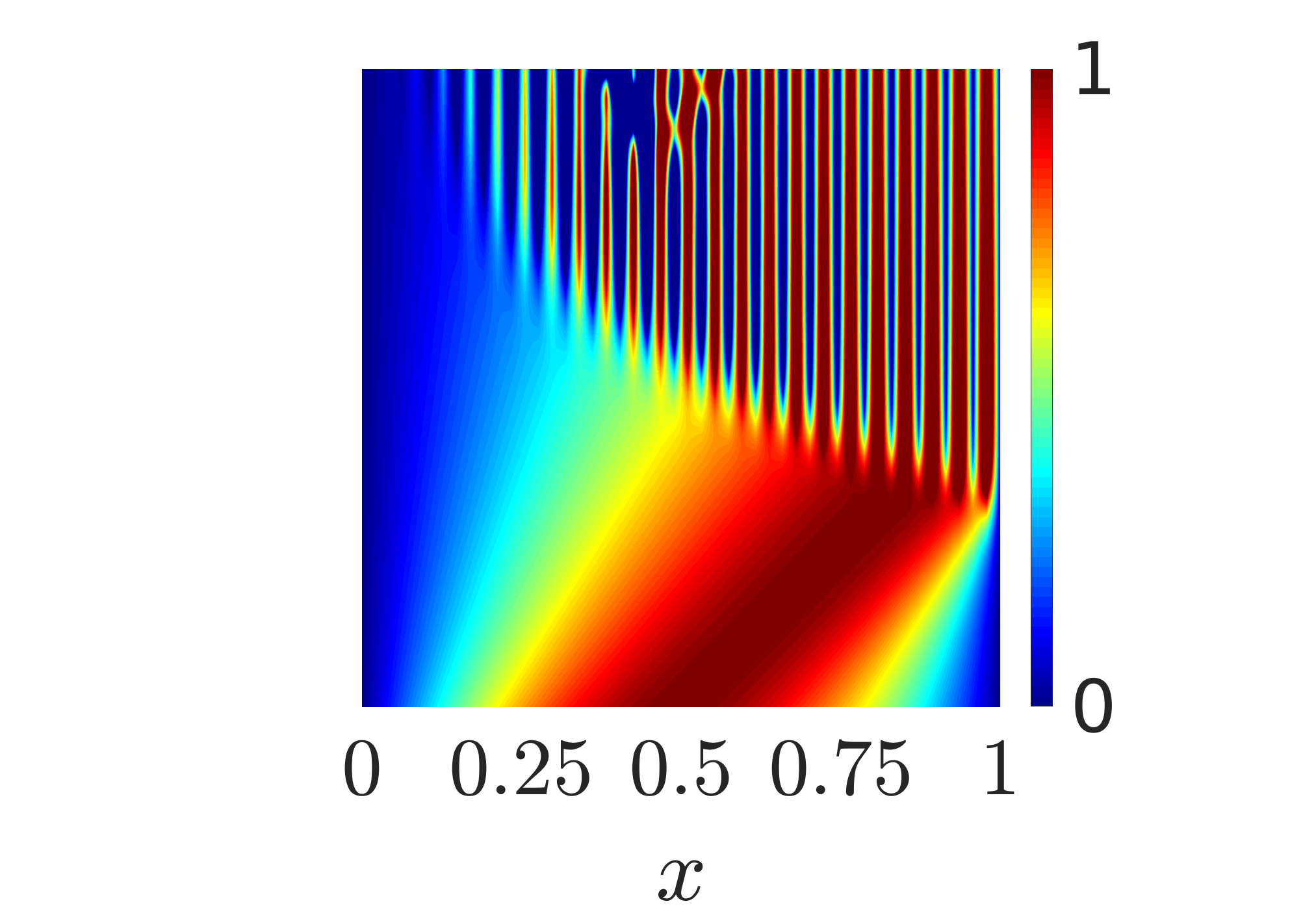}}
\end{subfigmatrix}
\caption{47-element, $p = 1$ FEM solution contours for the boundary layer problem over a range of $\nu$}
\label{fig:Example1_47_element_FEM_soln}
\end{center}
\end{figure}

Equation \ref{example1_viscous_burger_strong_formulation} was initially solved over uniform grids ($h = \Big\{\frac{1}{11}, \frac{1}{23}, \frac{1}{47}, \frac{1}{95}, \frac{1}{191}\Big\}$) using linear ($p = 1$) FEM. The Crank-Nicolson scheme is used for temporal discretization with a step size of $\Delta t = \frac{1}{5000}$. At each time step the Newton-Raphson method is used to iteratively solve the nonlinear set of equations. At times $t = [0, 0.25, 0.318, 0.5, 0.75, 1]$, relative $L_2$ and $H_1$ integral norms versus total degrees of freedom are shown in Figs. \ref{fig:Example1_p1FEM_L2norm} and \ref{fig:Example1_p1FEM_H1norm}, respectively. Convergence rates are computed and presented using $h = \Big\{\frac{1}{95}, \frac{1}{191}\Big\}$. It is observed as $\nu$ decreases, errors in the $L_2$ and $H_1$ norms increase. Specifically, a shift in the relative $L_2$ and $H_1$ norm is observed, with convergence rates remaining relatively unaffected except when $\nu = 0$. This is a result of the increasing difficulty in resolving the steep boundary layer that forms around $x = 1$. For $\nu = 0$, convergence does not occur once a discontinuity arises at $t_b \geq \frac{1}{\pi}$. Spurious oscillations are visually observed in the numerical solutions for small viscosities as shown in Fig. \ref{fig:Example1_11_element_FEM_soln}, which displays 11-element solution contours. For relatively large kinematic vicosities ($\nu = \frac{1}{10}$) no oscillations are visually observed. However, when $\nu = \frac{1}{100}$ and $\nu = 0$, severe nonphysical oscillations arise in the 11-element solutions. With sufficient grid refinement, the boundary layer is captured for $\nu = \frac{1}{100}$, as exemplified by the 47-element solutions contours in Fig. \ref{fig:Example1_47_element_FEM_soln}. However grid refinement does not improve the numerical solution when $\nu = 0$. 

\subsubsection{Generalized finite element solutions}
\begin{figure}[ht!]
\begin{center}
\begin{subfigmatrix}{6}
\subfigure[$t = 0$]{\includegraphics[width=2.1in]{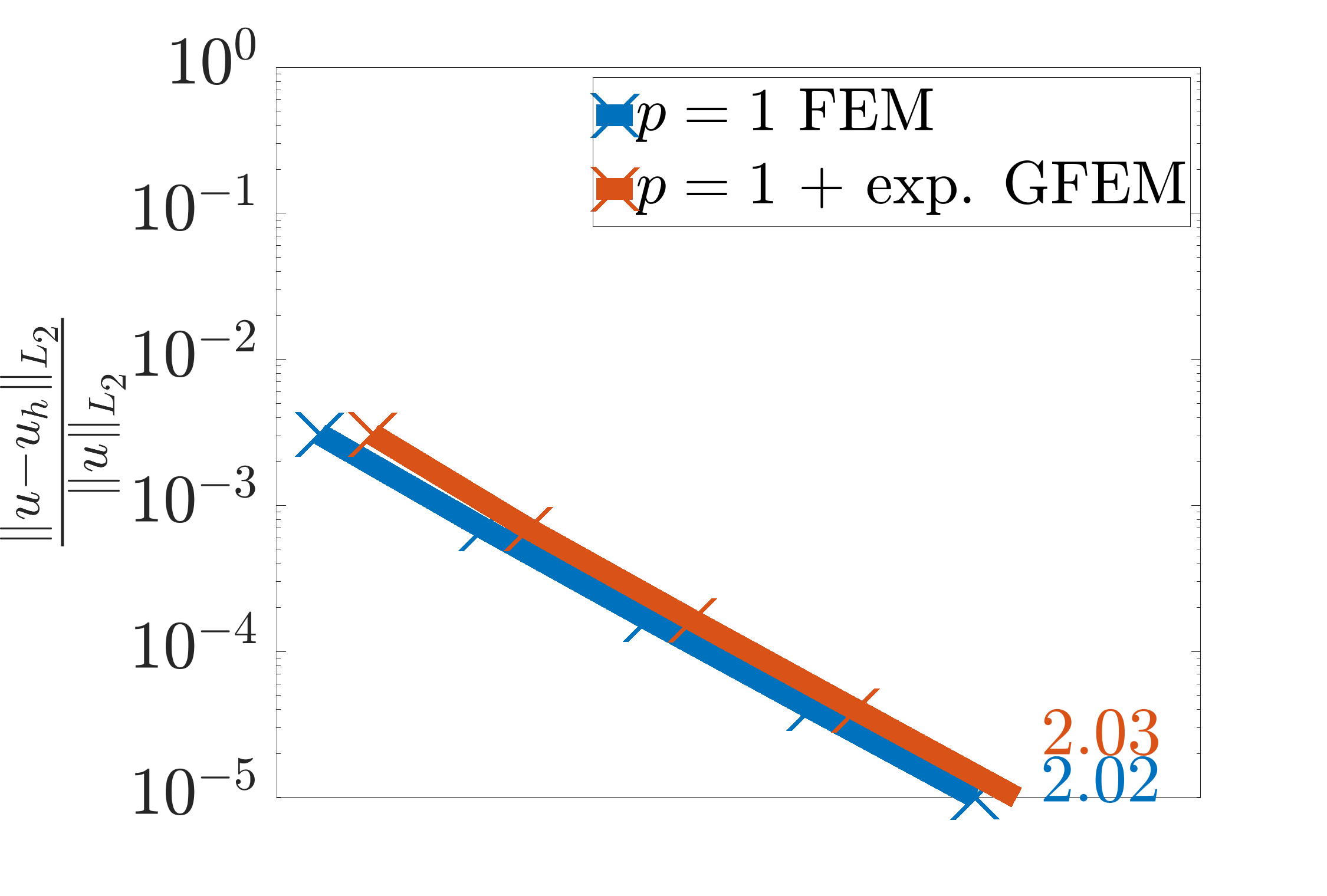}}
\subfigure[$t = 0.25$]{\includegraphics[width=2.1in]{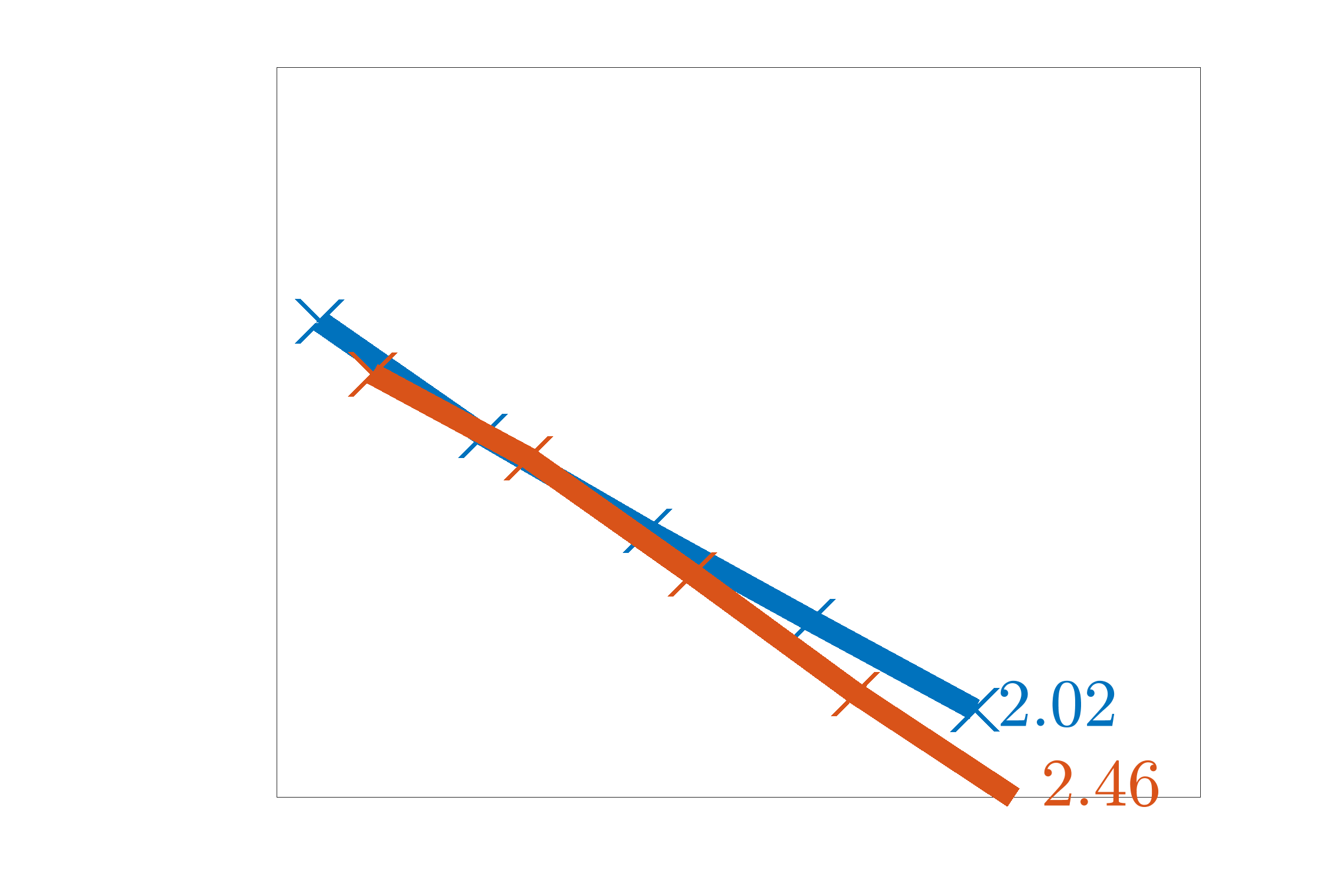}}
\subfigure[$t = 0.318 \approx \frac{1}{\pi}$]{\includegraphics[width=2.1in]{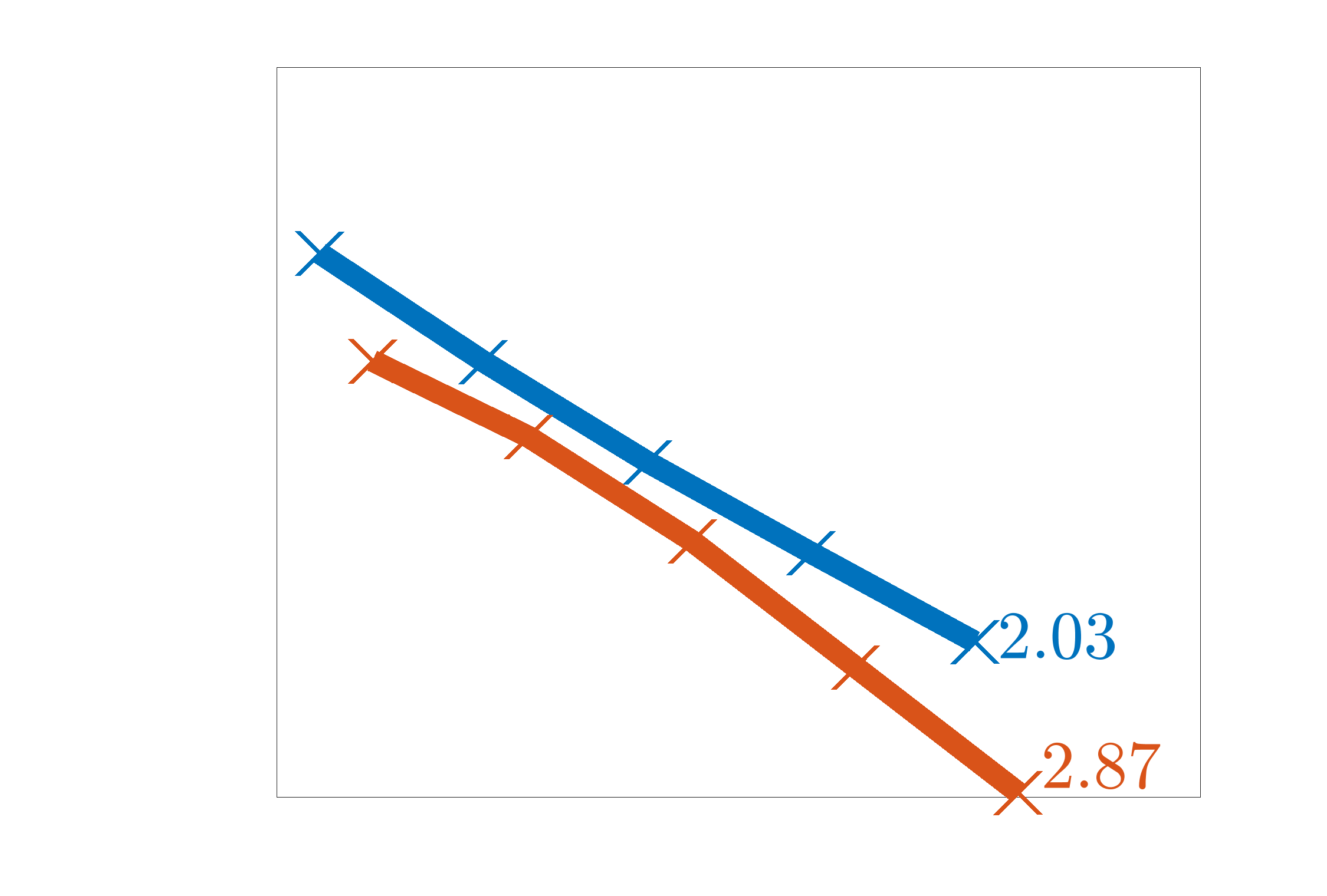}}
\subfigure[$t = 0.5$]{\includegraphics[width=2.1in]{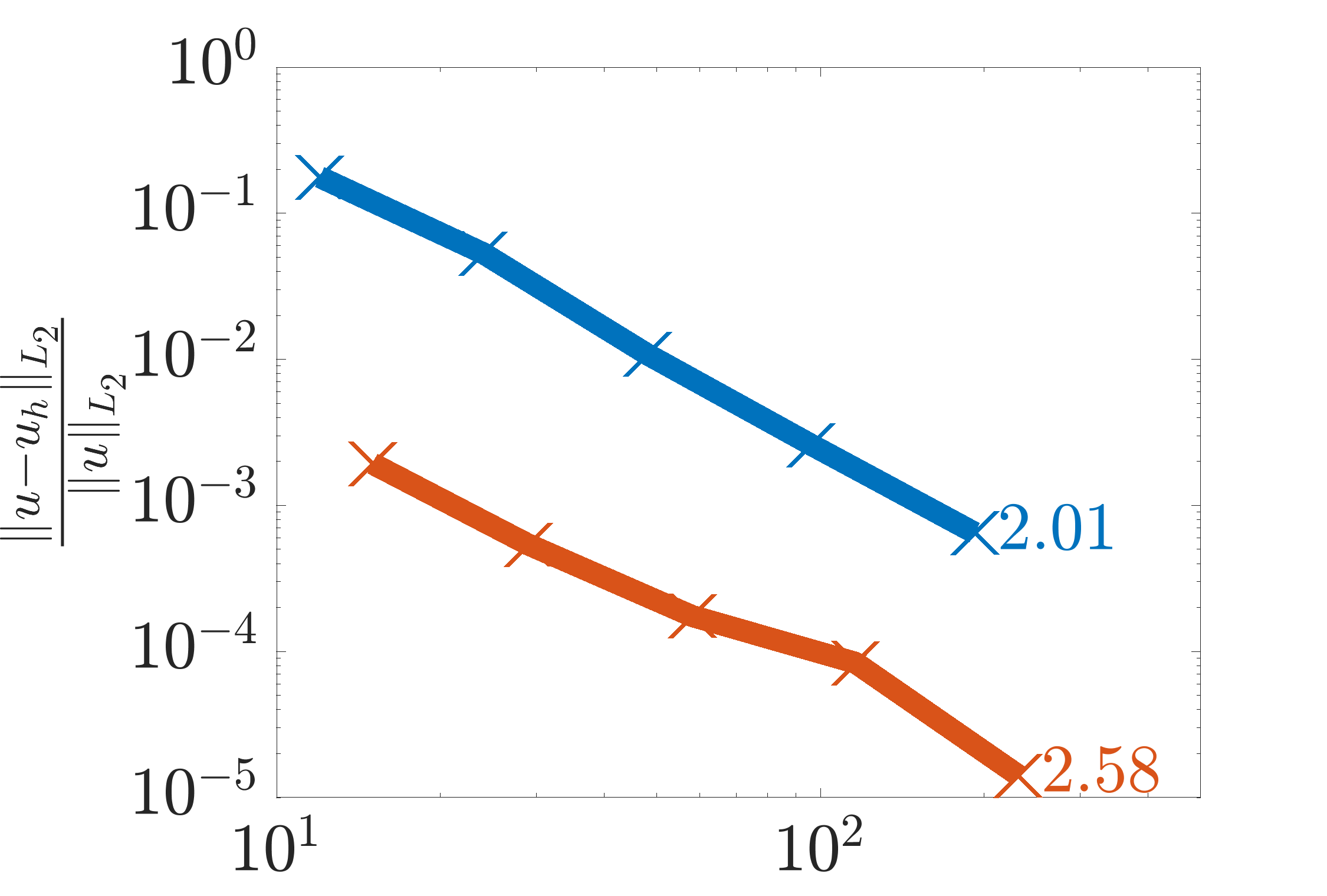}}
\subfigure[$t = 0.75$]{\includegraphics[width=2.1in]{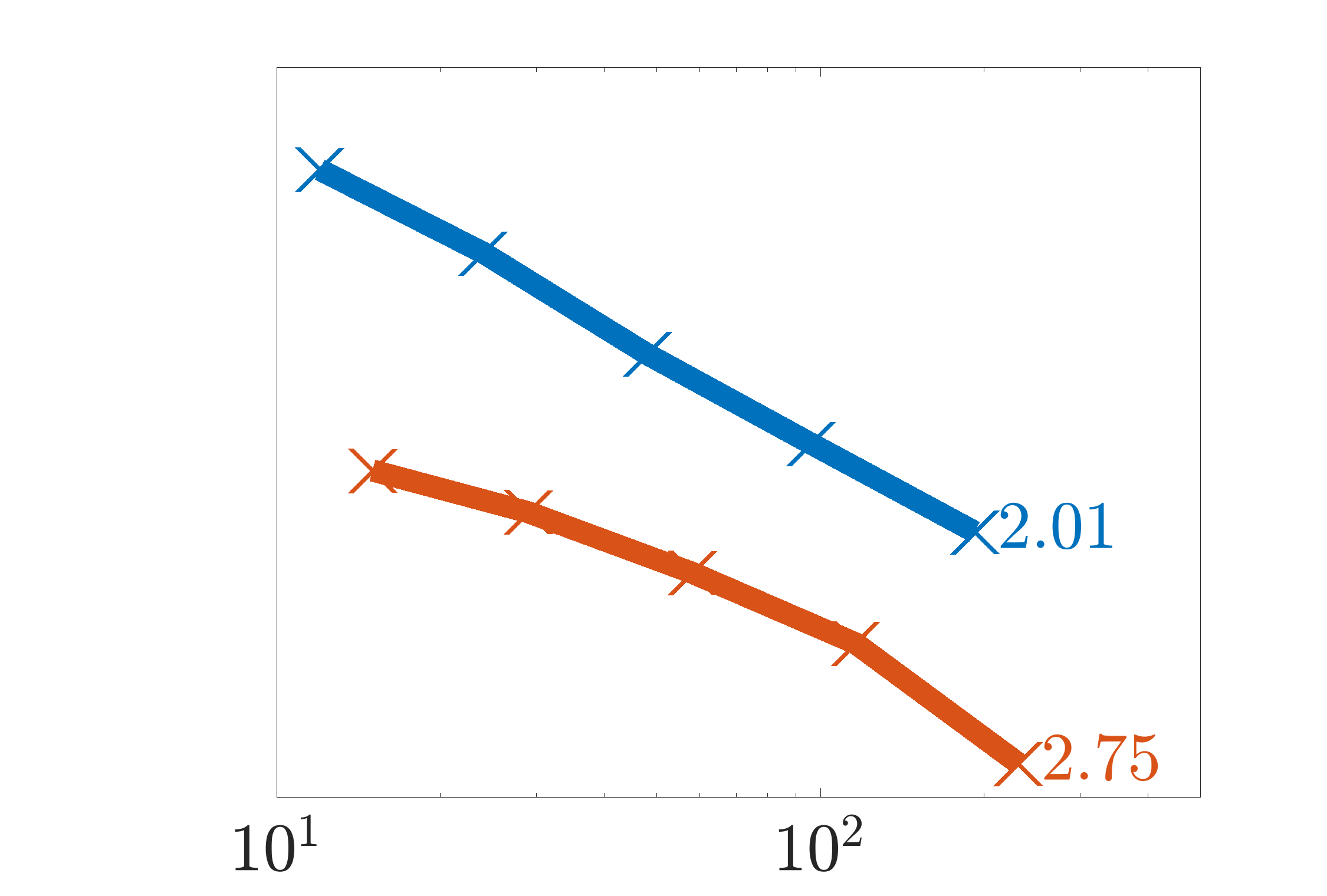}}
\subfigure[$t = 1$]{\includegraphics[width=2.1in]{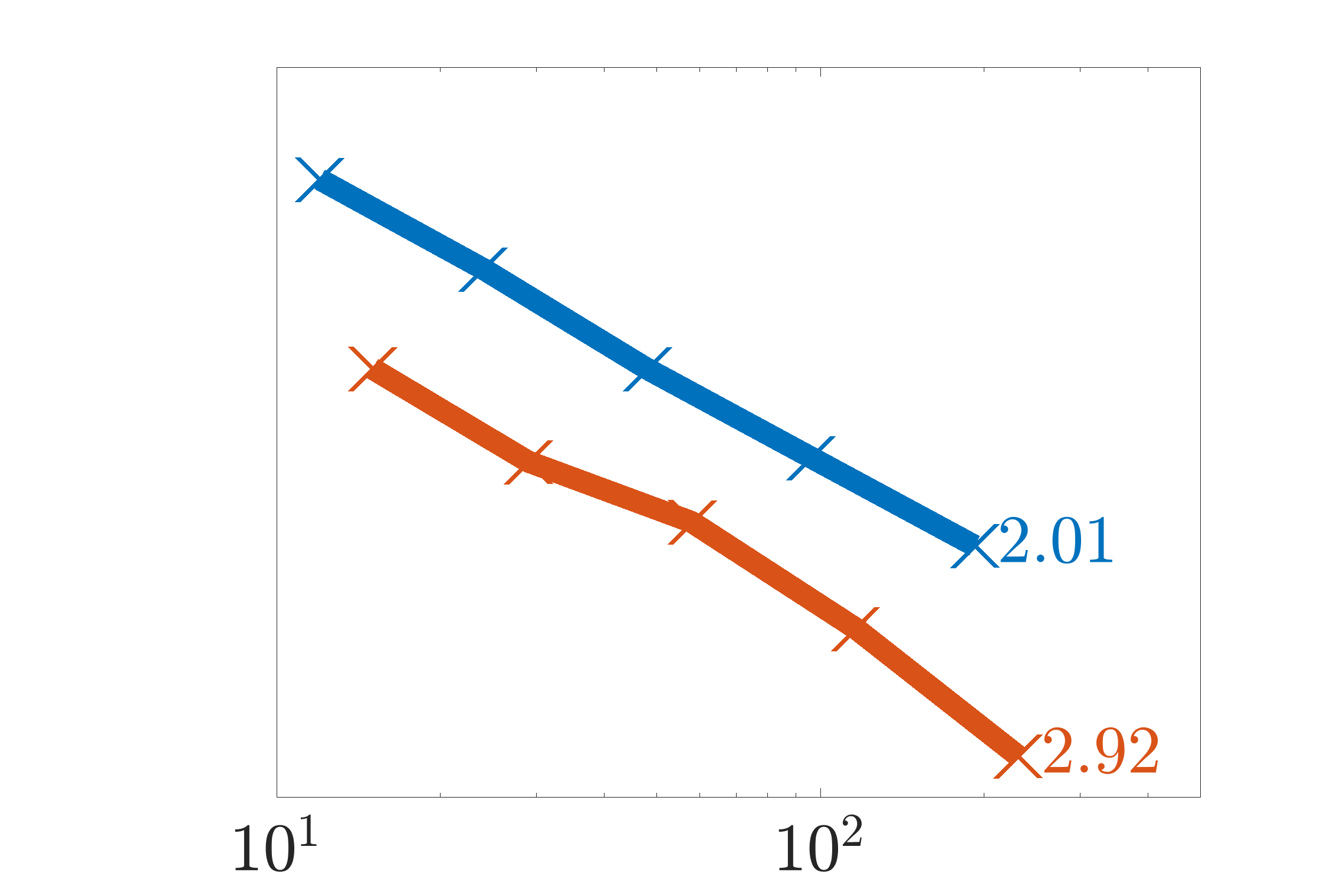}}
\end{subfigmatrix}
\caption{Convergence in the relative $L_2$ integral norms for the boundary layer problem with $\nu = \frac{1}{100}$ using exponential enrichments}
\label{fig:Example1_p1GFEM_nu1over100_L2norm}
\end{center}
\end{figure}

\begin{figure}[ht!]
\begin{center}
\begin{subfigmatrix}{6}
\subfigure[$t = 0$]{\includegraphics[width=2.1in]{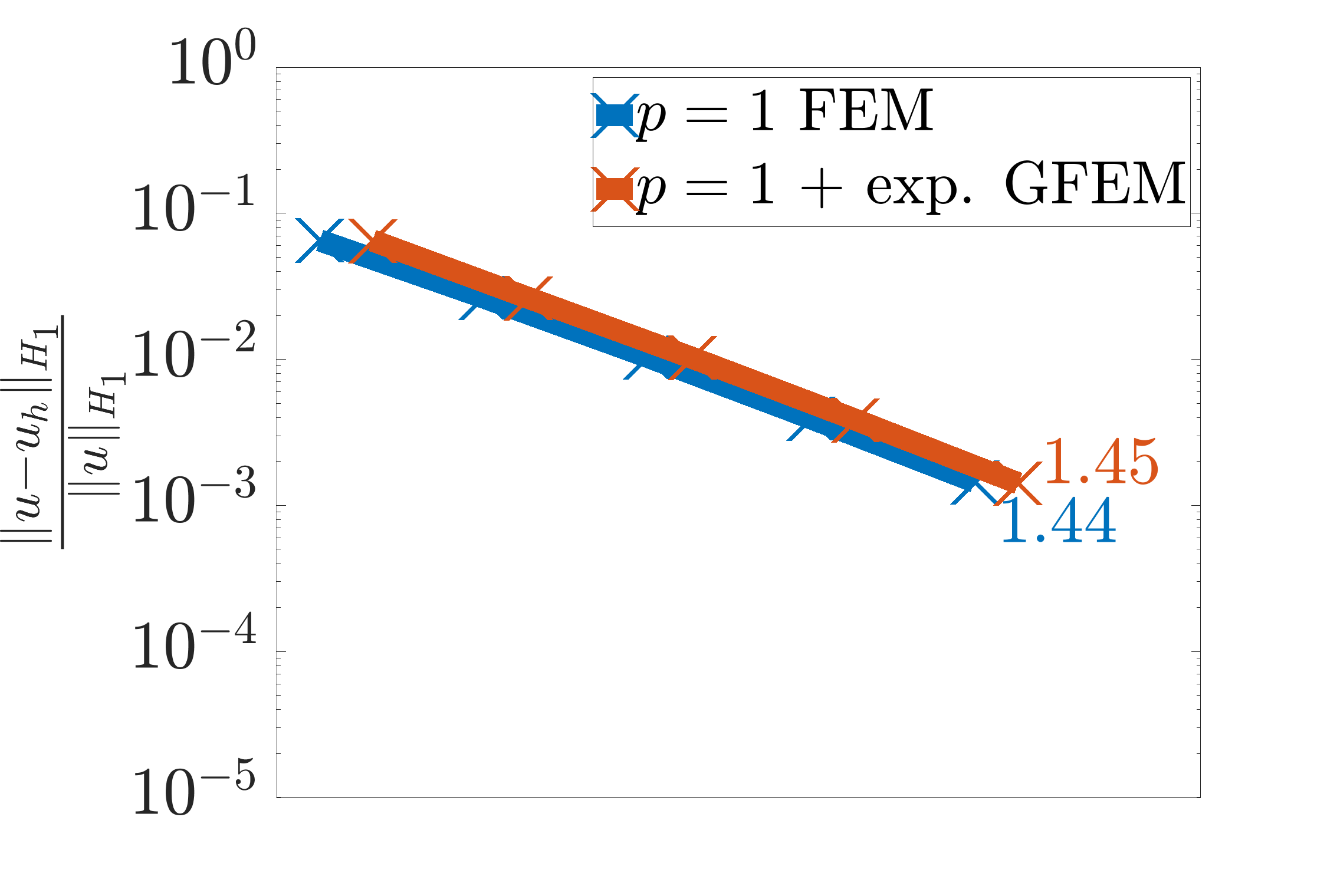}}
\subfigure[$t = 0.25$]{\includegraphics[width=2.1in]{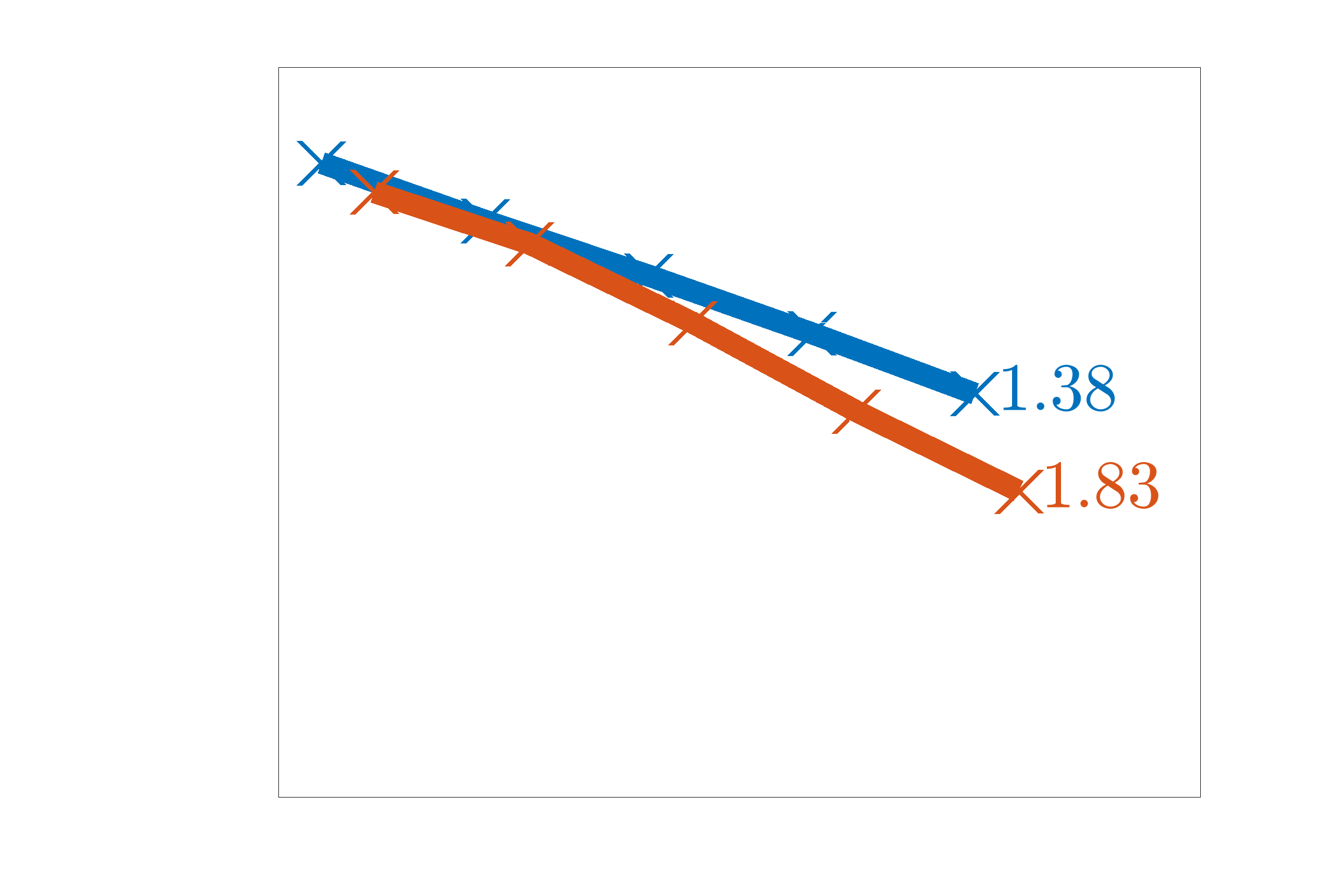}}
\subfigure[$t = 0.318 \approx \frac{1}{\pi}$]{\includegraphics[width=2.1in]{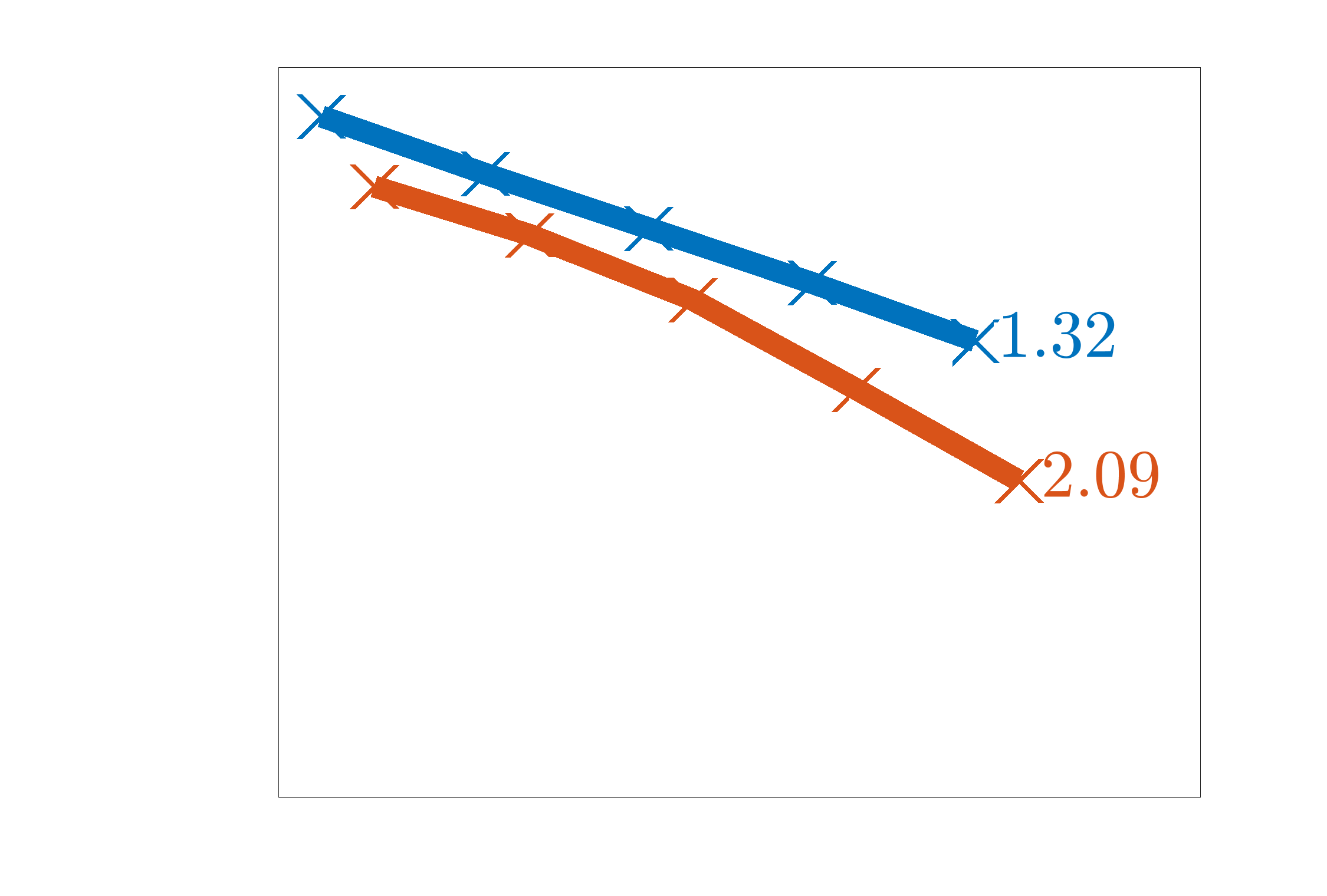}}
\subfigure[$t = 0.5$]{\includegraphics[width=2.1in]{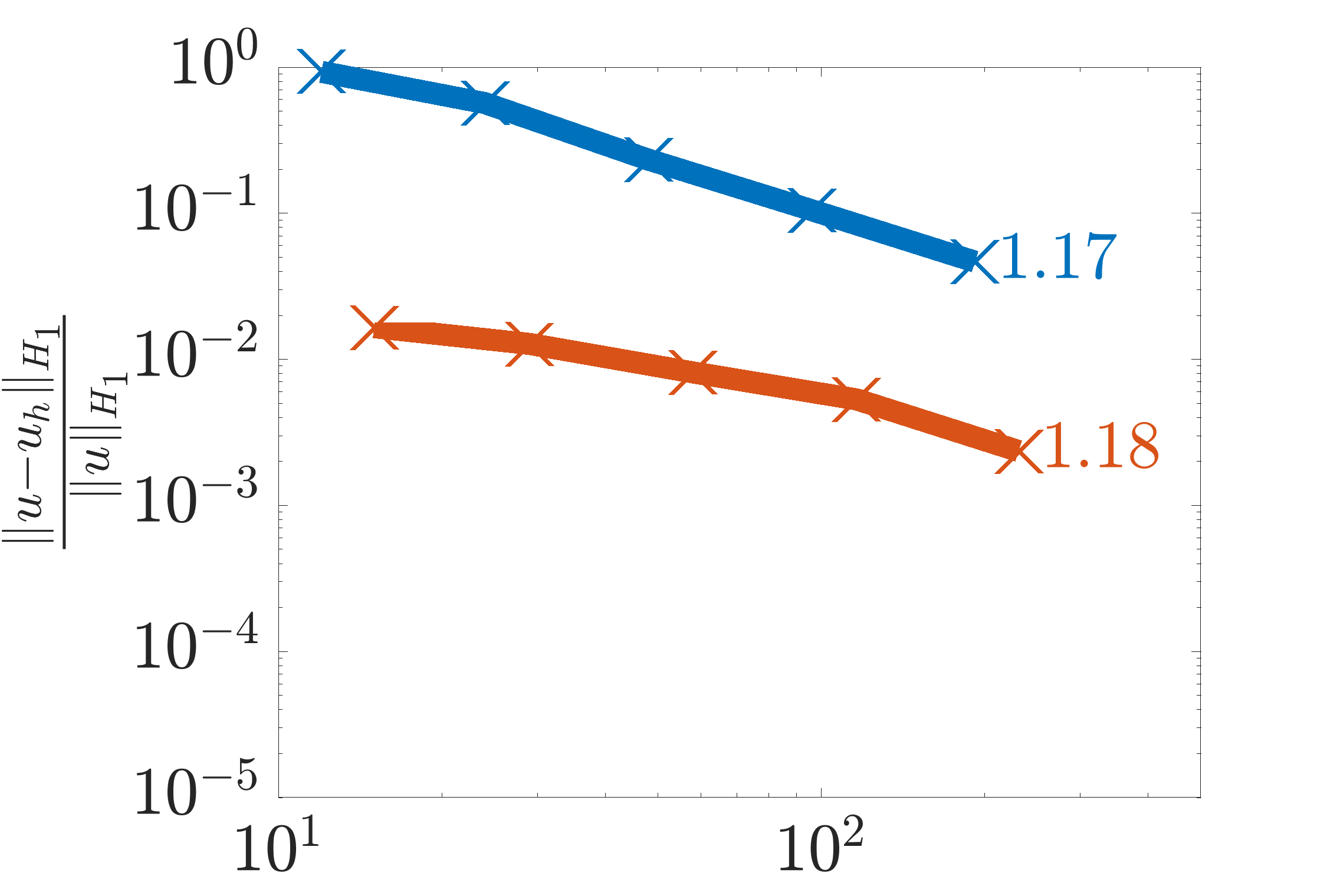}}
\subfigure[$t = 0.75$]{\includegraphics[width=2.1in]{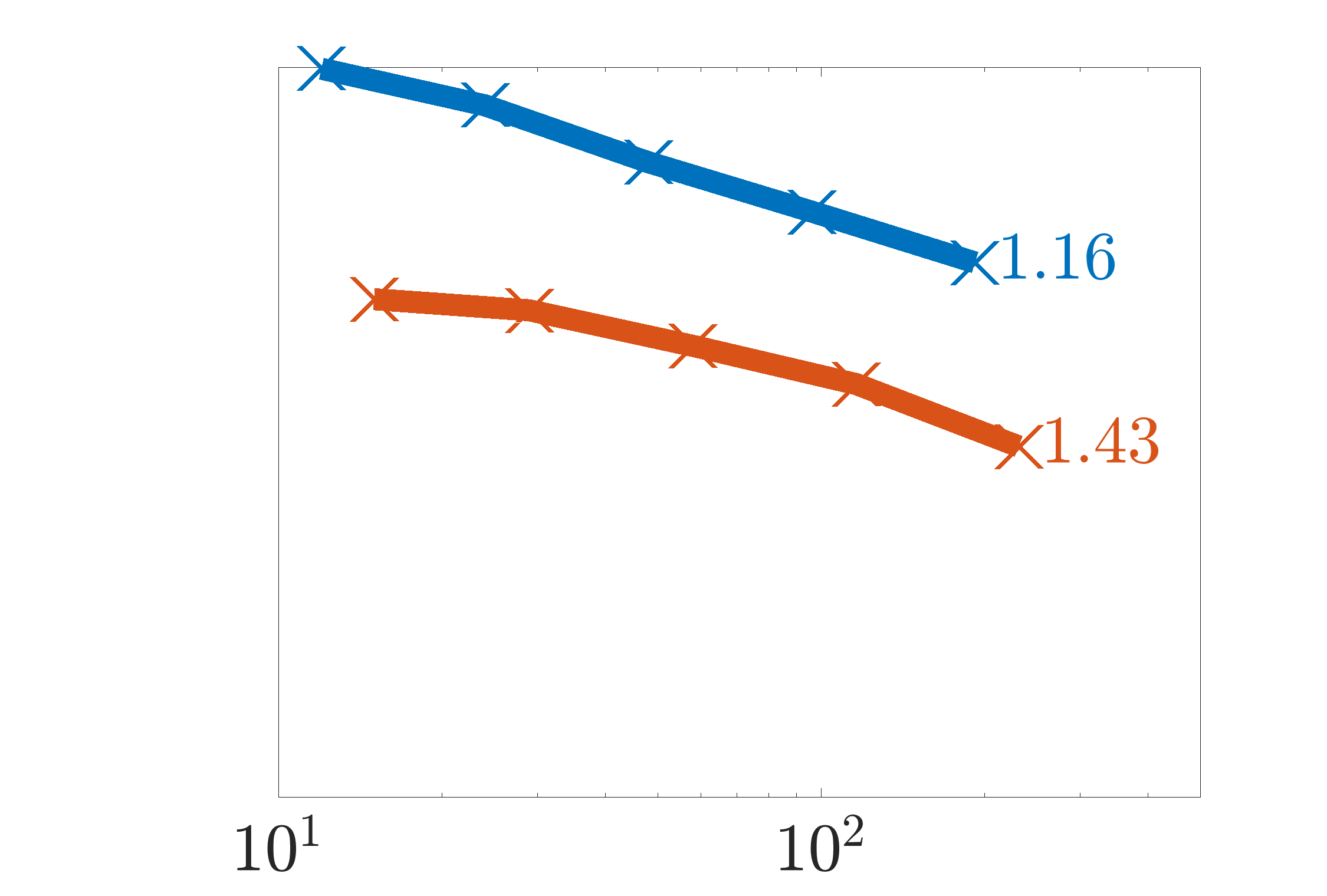}}
\subfigure[$t = 1$]{\includegraphics[width=2.1in]{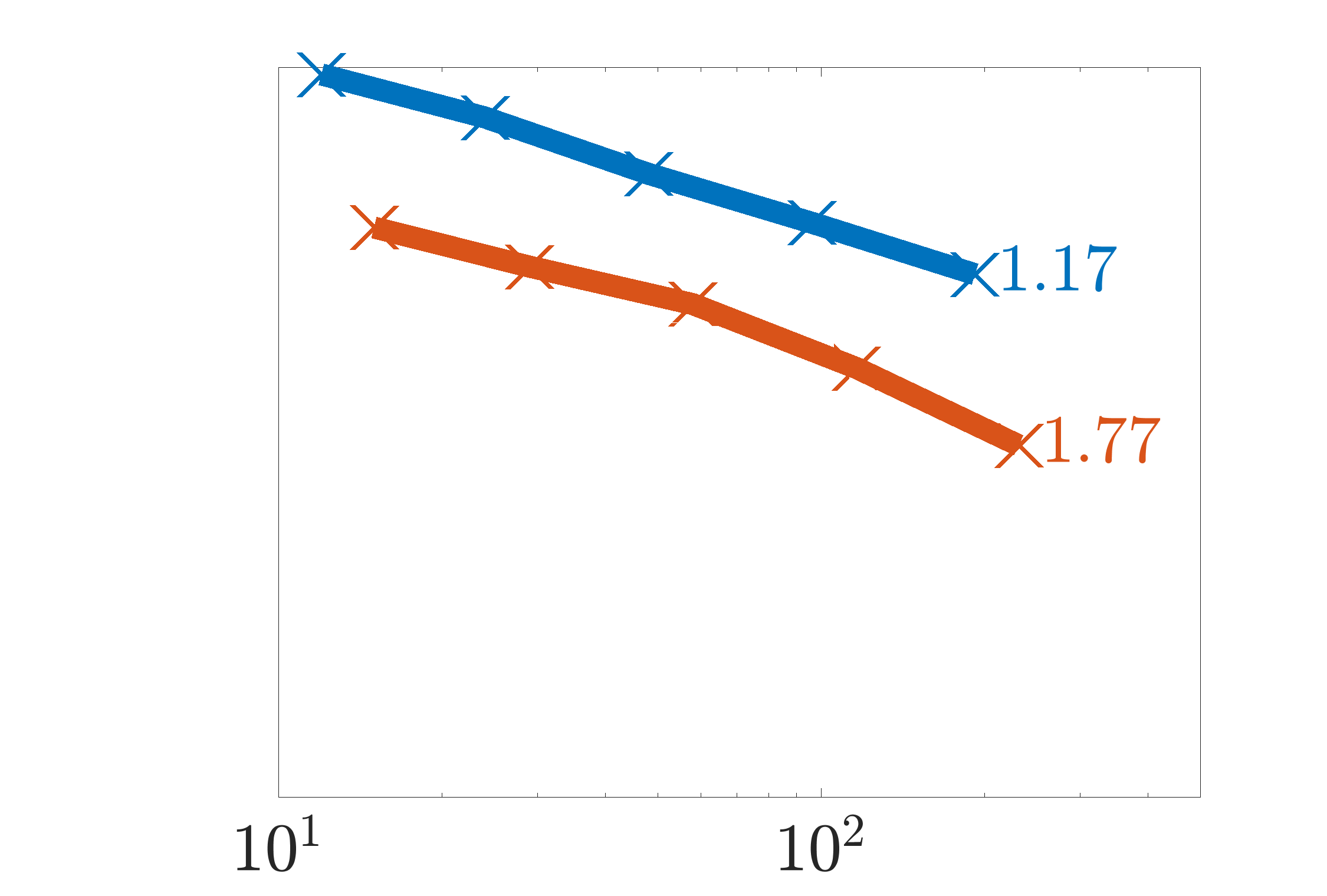}}
\end{subfigmatrix}
\caption{Convergence in the relative $H_1$ integral norms for the boundary layer problem with $\nu = \frac{1}{100}$ using exponential enrichments}
\label{fig:Example1_p1GFEM_nu1over100_H1norm}
\end{center}
\end{figure}

\begin{figure}[ht!]
\begin{center}
\begin{subfigmatrix}{2}
\subfigure[Relative $L_2$ integral norm]{\includegraphics[width=3in]{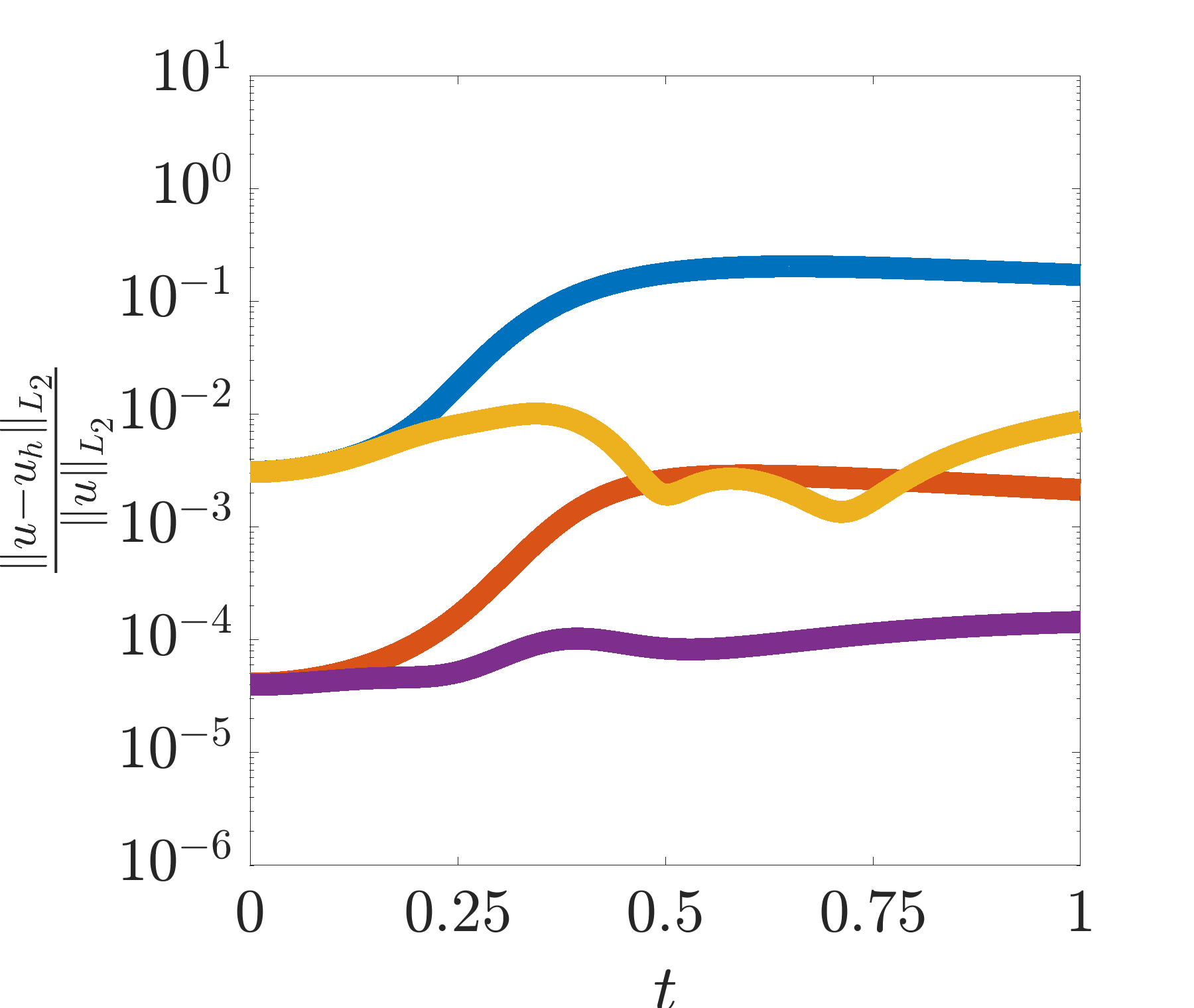}}
\subfigure[Relative $H_1$ integral norm]{\includegraphics[width=3in]{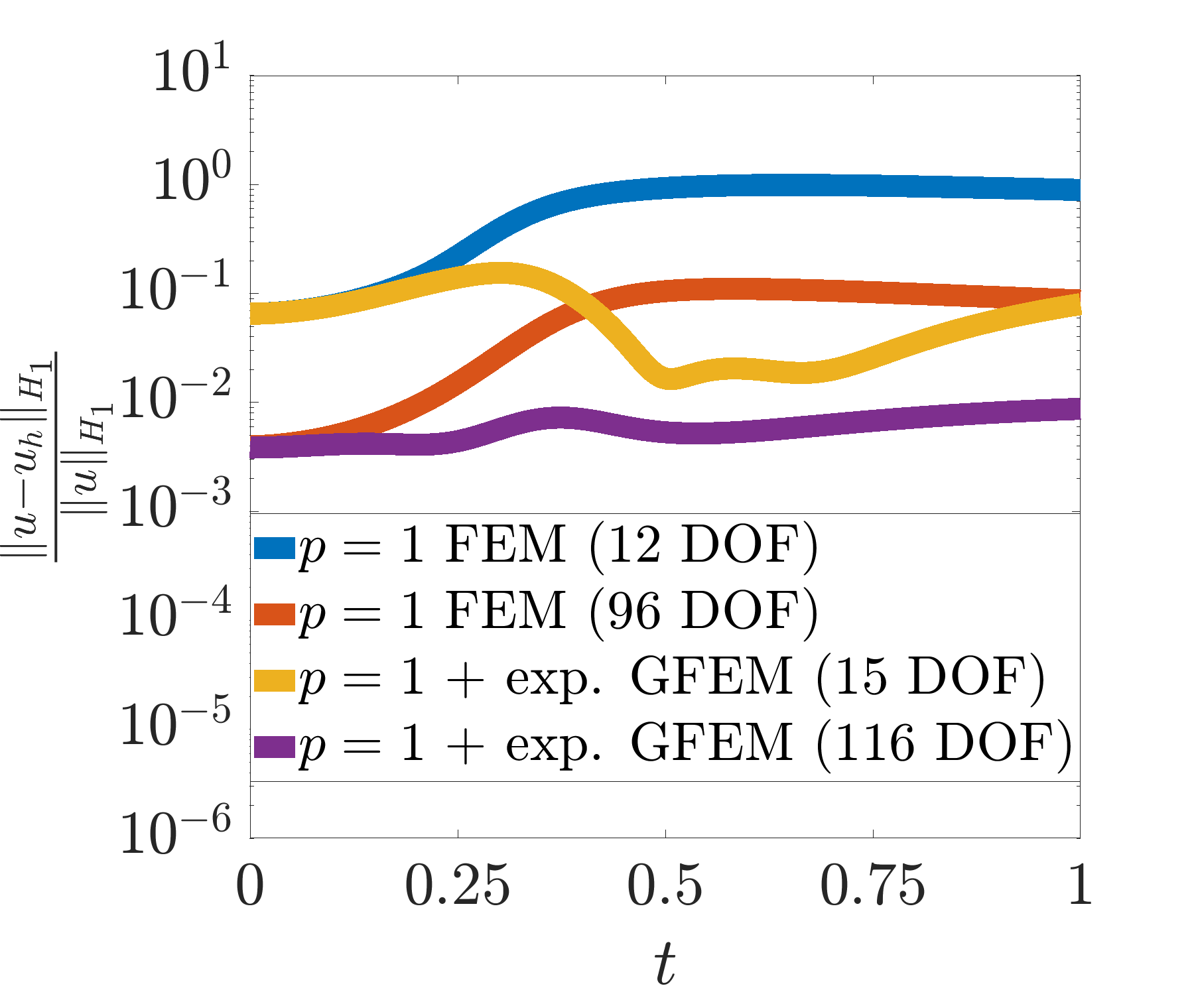}}
\end{subfigmatrix}
\caption{Relative $L_2$ and $H_1$ integral norms versus time for 11-element and 95-element FEM and $p = 1$ + exp. GFEM for the boundary layer problem with $\nu = \frac{1}{100}$}
\label{fig:Example1_p1expGFEM_L2H1vstime_nu1over100}
\end{center}
\end{figure}

\begin{figure}[ht!]
\begin{center}
\begin{subfigmatrix}{3}
\subfigure[Reference]{\includegraphics[width=2.1in]{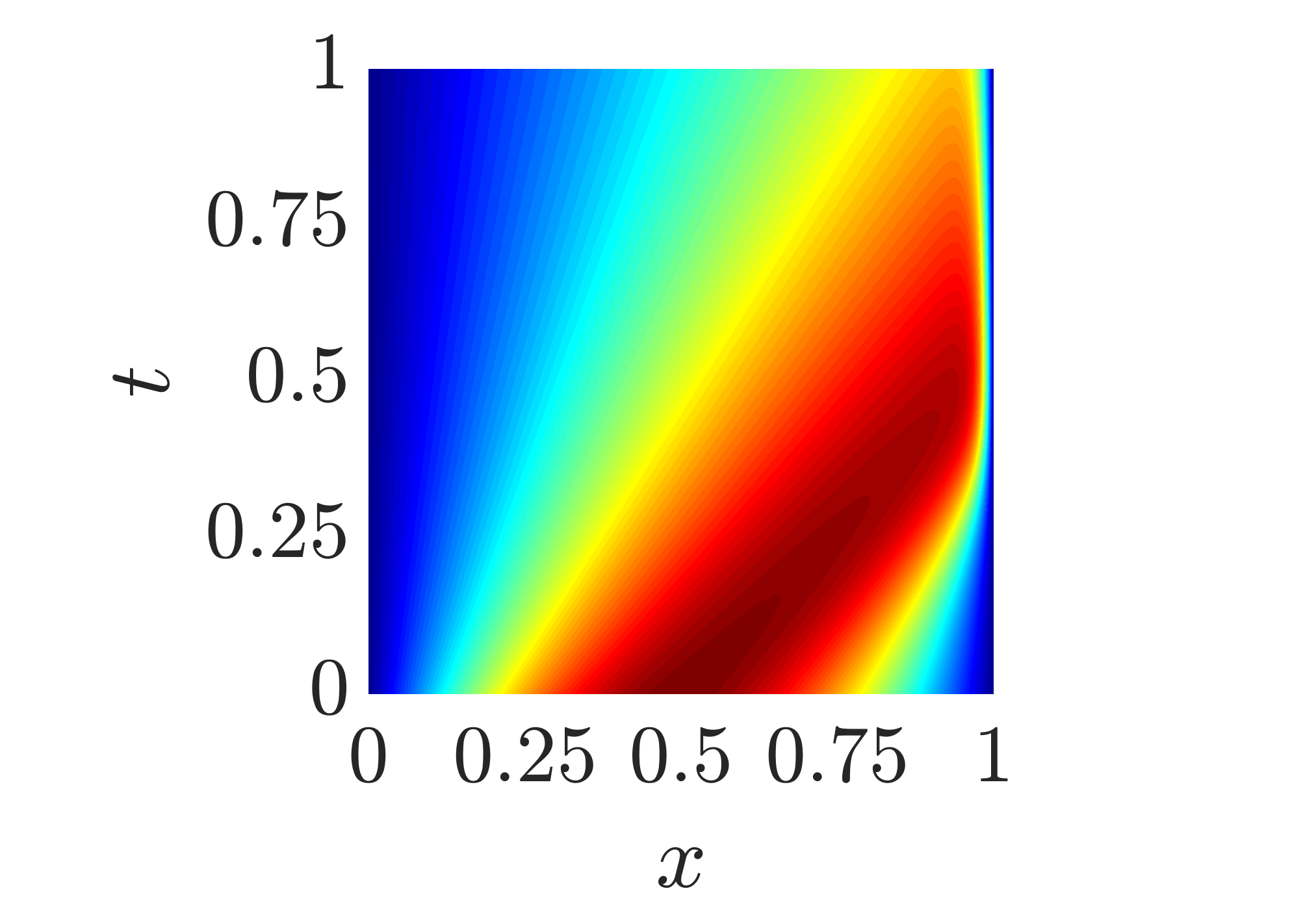}}
\subfigure[12 DOF FEM]{\includegraphics[width=2.1in]{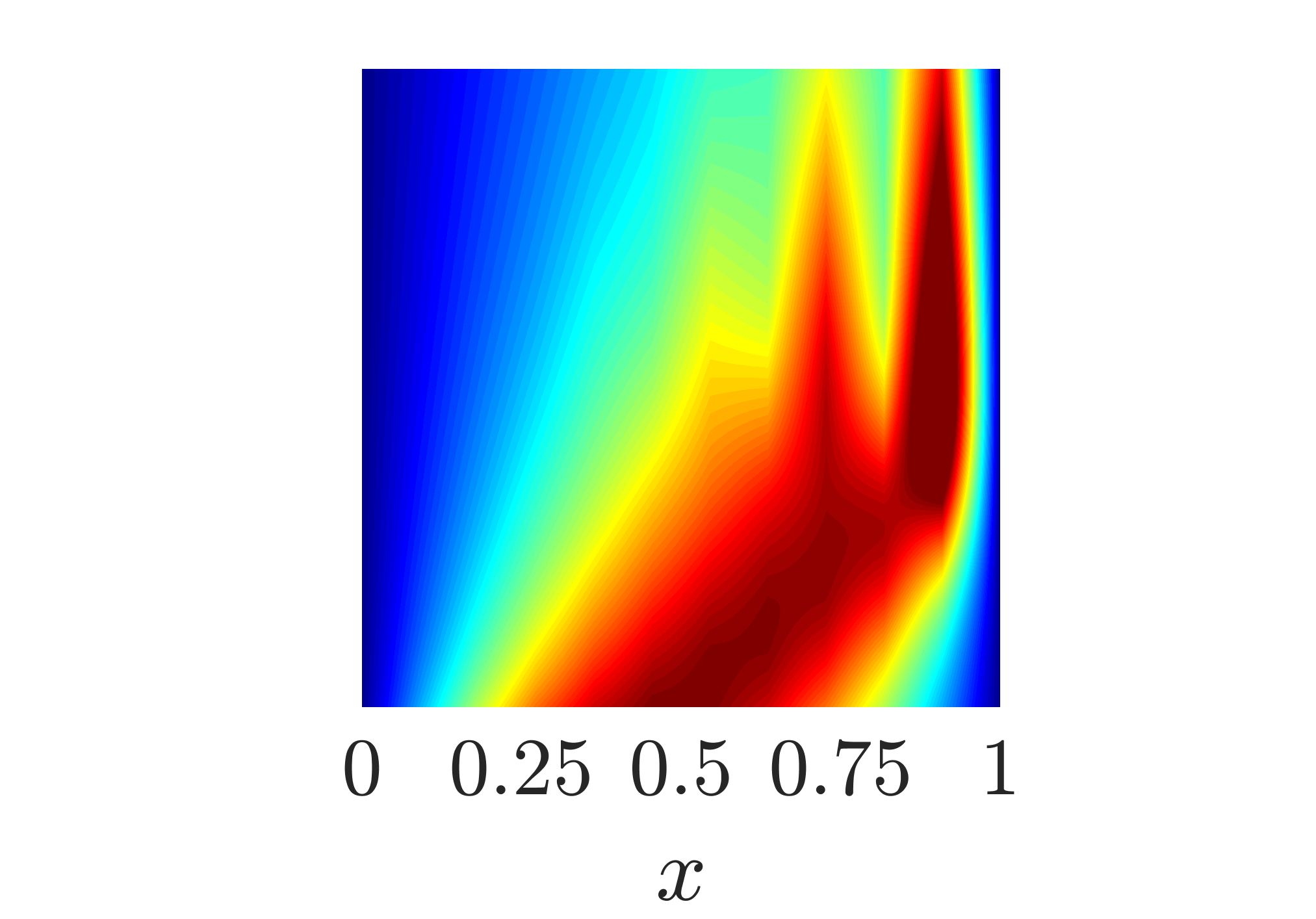}}
\subfigure[15 DOF GFEM]{\includegraphics[width=2.1in]{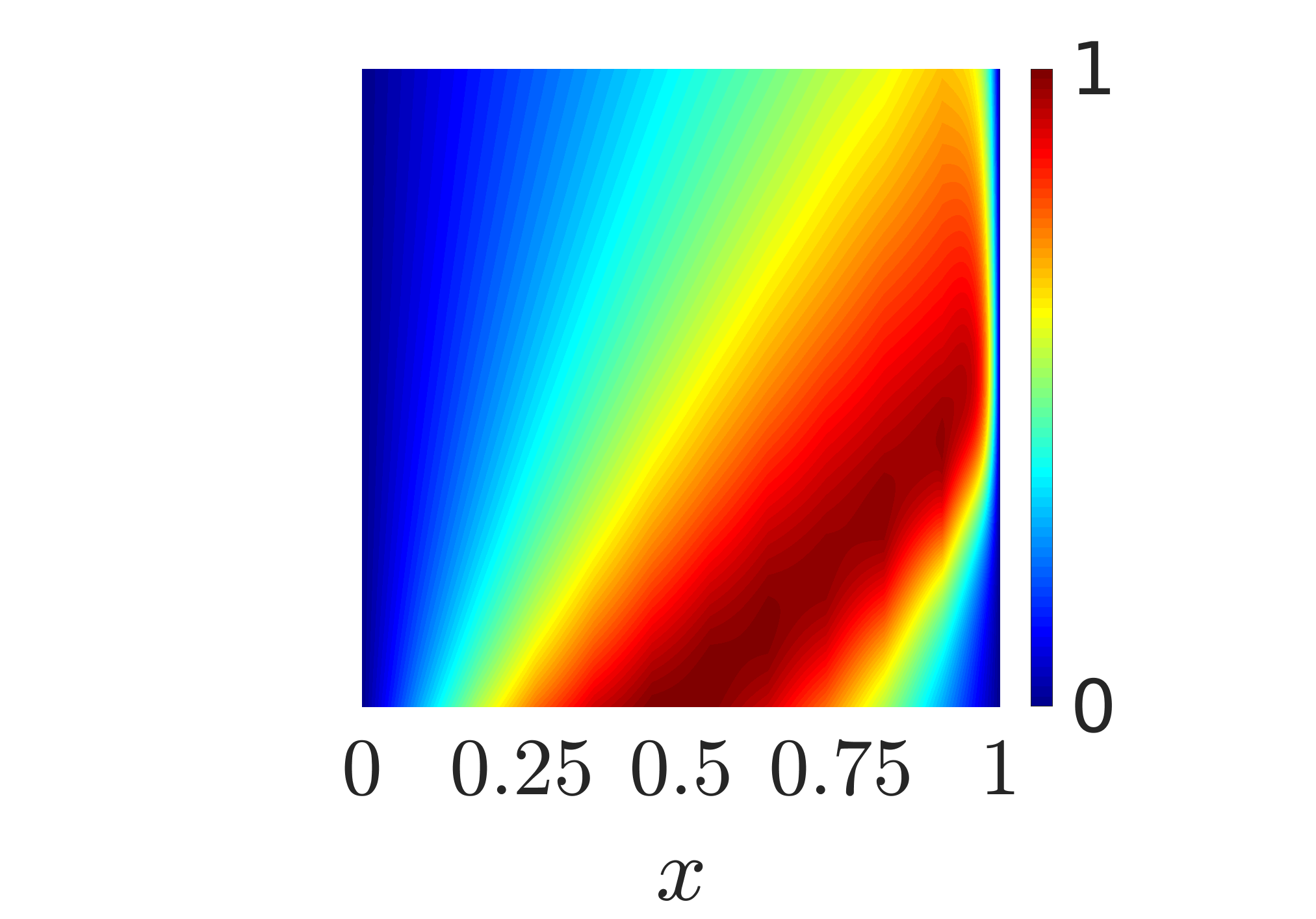}}
\end{subfigmatrix}
\caption{11-element, $p = 1$ FEM (12 DOF) and $p = 1$ + exp. GFEM (15 DOF) solution contours compared to the reference for the boundary layer problem with $\nu = \frac{1}{100}$}
\label{fig:Example1_nu001_FEMGFEM_contour_comparison}
\end{center}
\end{figure}

For $\nu = \frac{1}{100}$, the $p = 1$ FEM solutions are improved upon using GFEM with exponential functions as enrichments. These enrichments are applied to the local domain, $\Omega_{local} = [0.8, 1]$, roughly where the boundary layer forms. Use of these enrichments is motivated by findings in \citep{Shilt2021}, which demonstrate exponential enrichments stabilize the GFEM solution around boundary layers arising in the advection-diffusion equation. Recall the Burgers' equation is of similar form to the advection-diffusion, except the rate of advection is replaced with the solution variable, $u$. Physically, the solution variable $u$ may never exceed the maximum or minimum value provided by the initial conditions $u_{IC}(x)$. Thus, the specific exponential enrichment used  for this problem is chosen to be $E_1 = e^{\frac{|\max{u_{IC}(x)}| x}{\nu}} = e^{100x}$. Results using these exponential enrichments are denoted as $p = 1$ + exp. GFEM solutions. Grid sizes, temporal discretization, and nonlinear iteration are idential to those used for the $p = 1$ FEM solutions. Convergence in the relative $L_2$ and $H_1$ integral norm at various times are shown in Figs. \ref{fig:Example1_p1GFEM_nu1over100_L2norm} and \ref{fig:Example1_p1GFEM_nu1over100_H1norm}, respectively. Convergence rates are computed and presented using $h = \Big\{\frac{1}{95}, \frac{1}{191}\Big\}$. Recall no theoretical convergence rates are formally developed for GFEM using the enrichments in this example. It is observed in Figs \ref{fig:Example1_p1GFEM_nu1over100_L2norm} and \ref{fig:Example1_p1GFEM_nu1over100_H1norm} that use of the exponential enrichments results in a significant reduction of error after the boundary layer forms, and at the same number of DOFs when compared to $p = 1$ FEM. Plots of the relative $L_2$ and $H_1$ integral norms versus time are shown in Fig. \ref{fig:Example1_p1expGFEM_L2H1vstime_nu1over100}  for the 11- and 95-element GFEM solutions. The FEM and GFEM solutions have nearly identical error until $t \approx 0.25$ where boundary layer gradients become larger. The FEM solutions increase in error due to spurious oscillations, whereas the GFEM error levels remain relatively unaffected. The result is approximately 10 times reduction of error in the GFEM solutions at later time steps. Lastly, 11-element solution contours are shown in Fig. \ref{fig:Example1_nu001_FEMGFEM_contour_comparison}. Here, severe nonphysical oscillations in the $p = 1$ FEM solution are observed; whereas the $p = 1$ + exp. GFEM solution successfully captures the steep boundary layer, presenting smooth solution contours at roughly the same number of degrees of freedom. 

\begin{figure}[ht!]
\begin{center}
\begin{subfigmatrix}{6}
\subfigure[$t = 0$]{\includegraphics[width=2.1in]{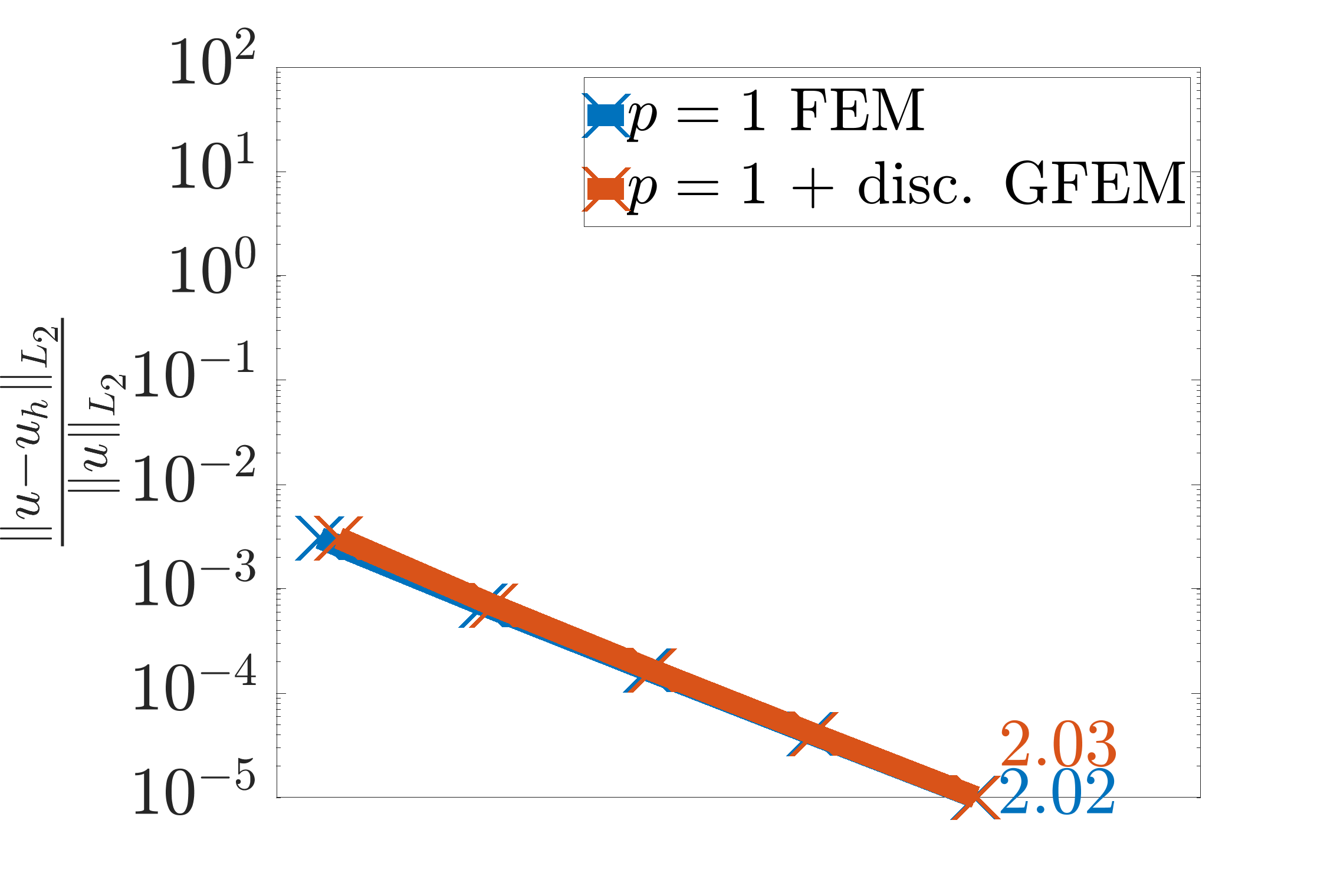}}
\subfigure[$t = 0.25$]{\includegraphics[width=2.1in]{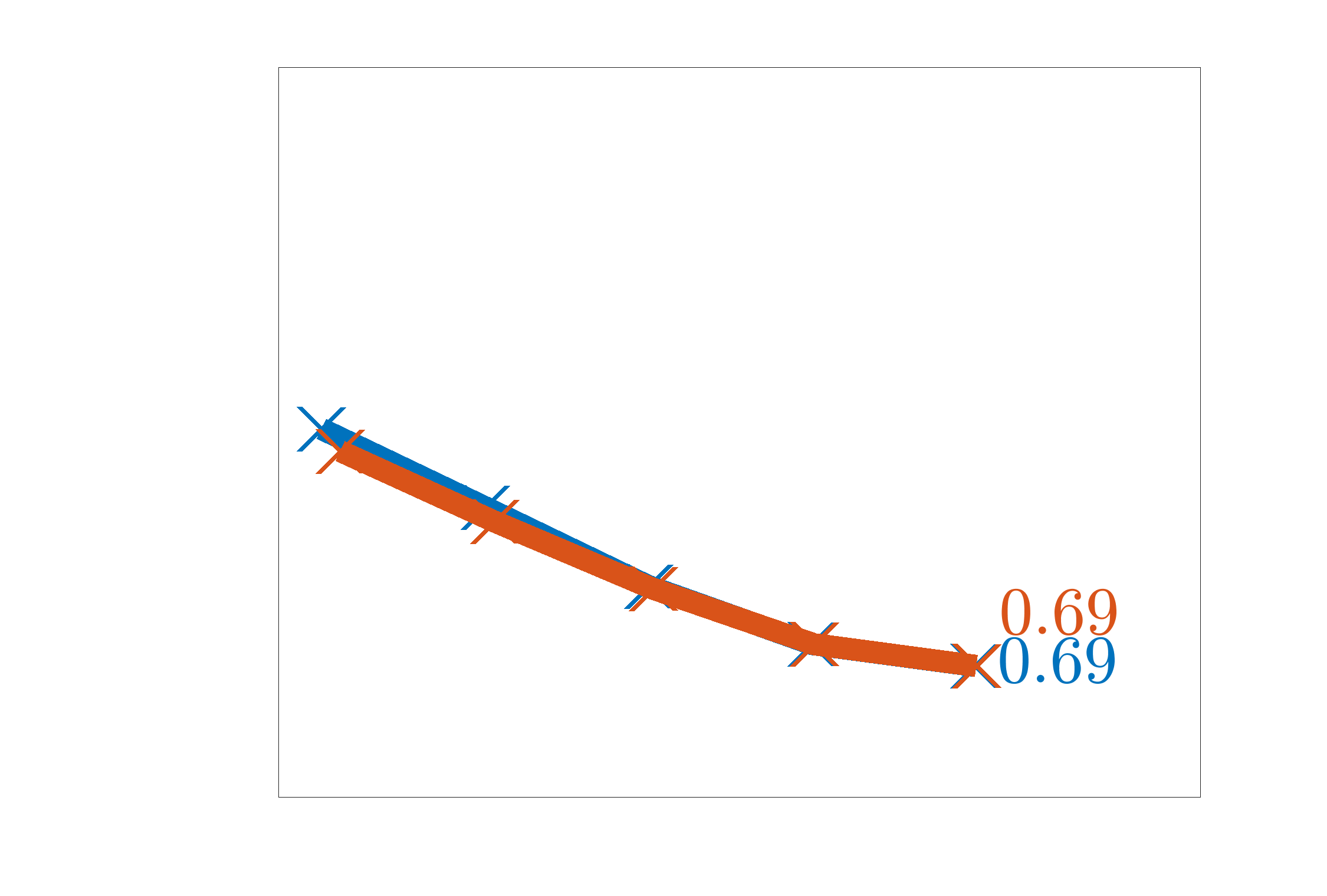}}
\subfigure[$t = 0.318 \approx \frac{1}{\pi}$]{\includegraphics[width=2.1in]{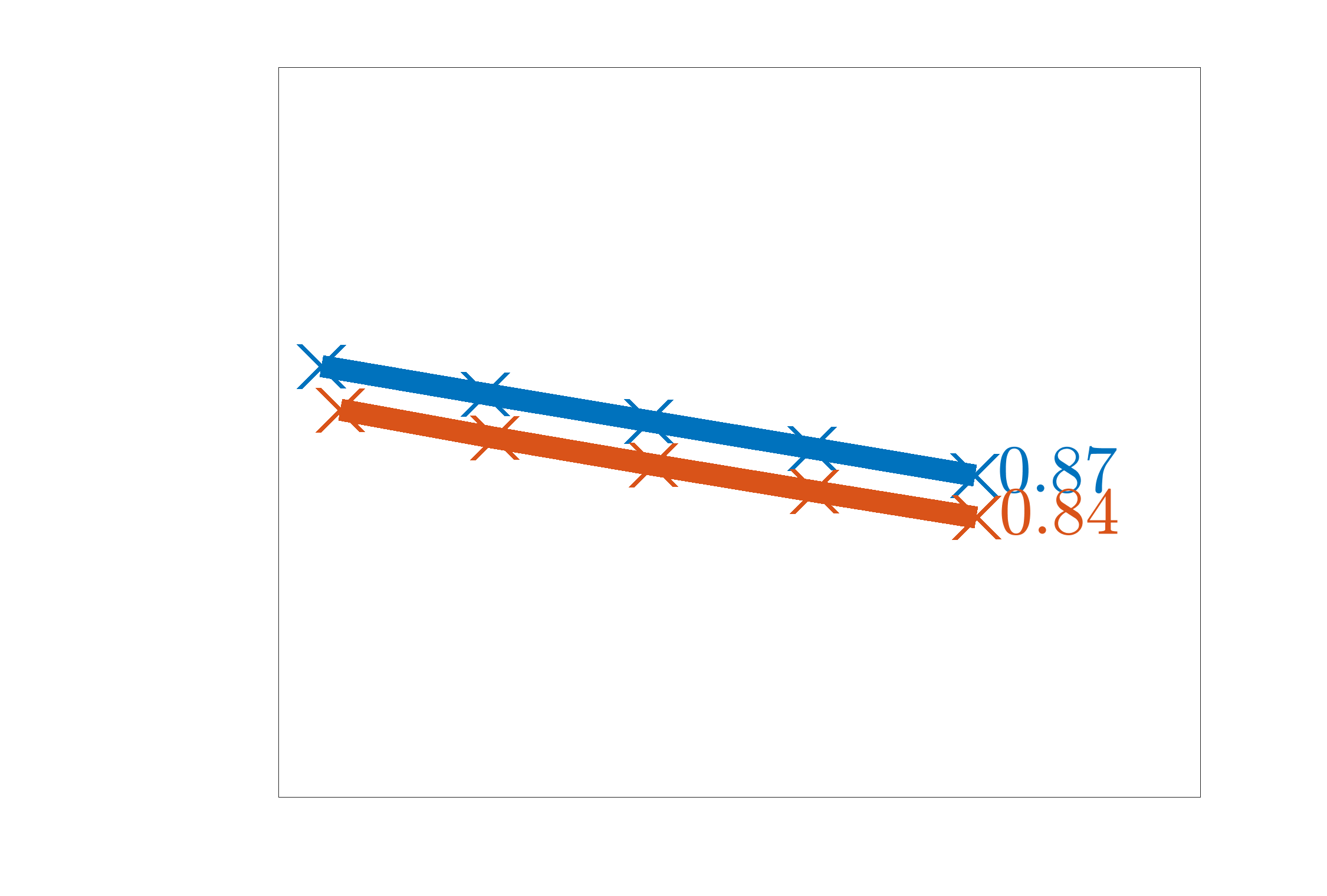}}
\subfigure[$t = 0.5$]{\includegraphics[width=2.1in]{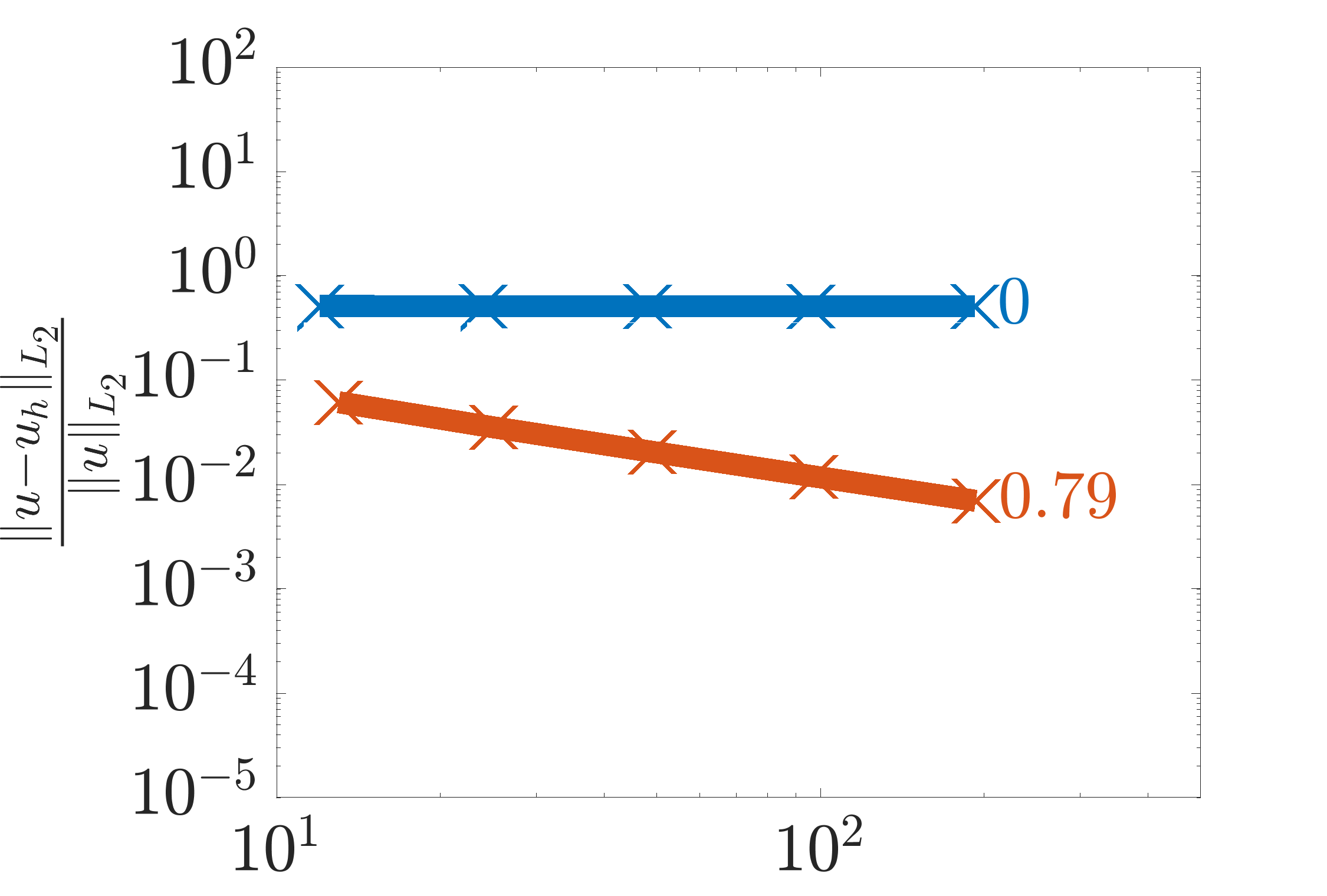}}
\subfigure[$t = 0.75$]{\includegraphics[width=2.1in]{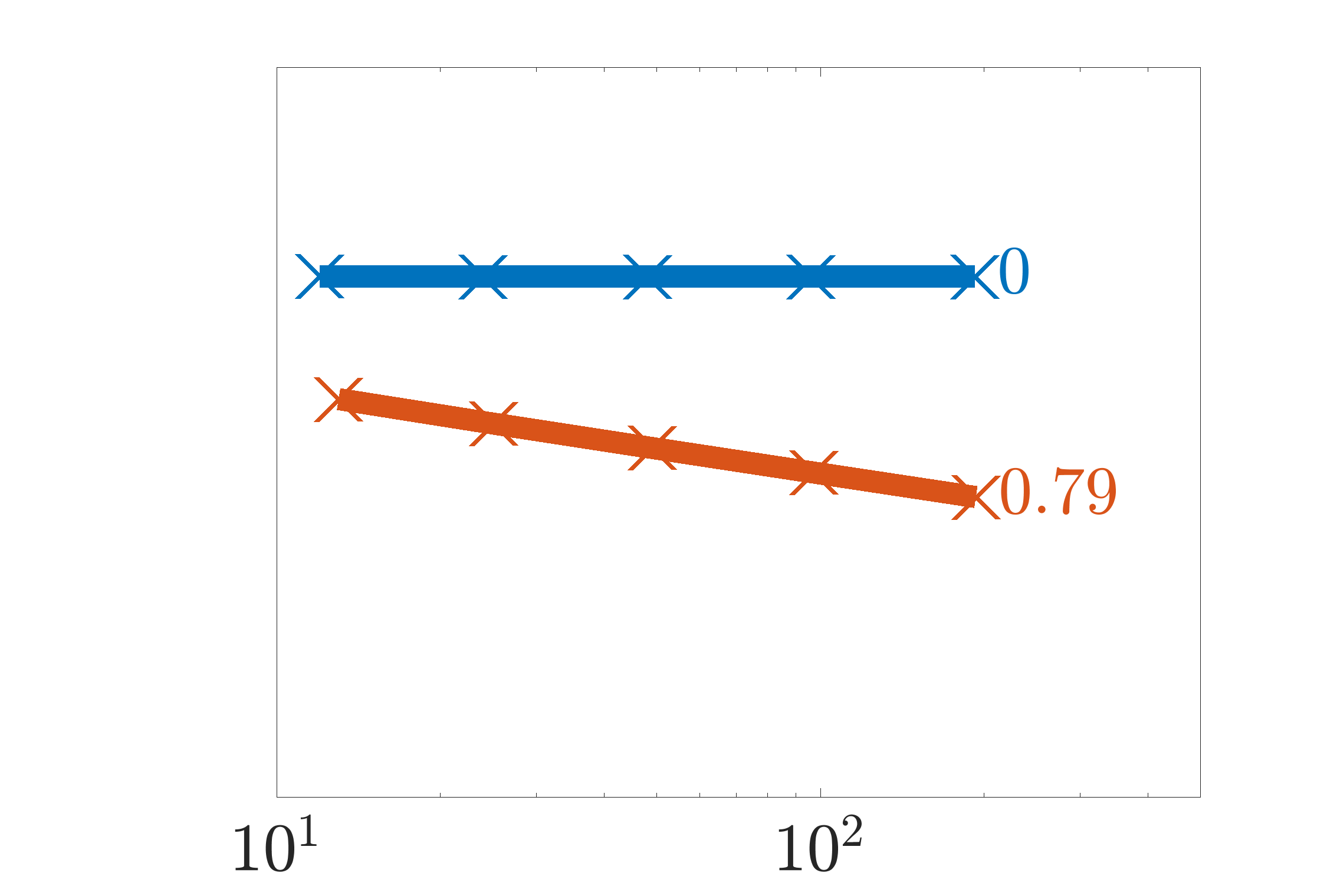}}
\subfigure[$t = 1$]{\includegraphics[width=2.1in]{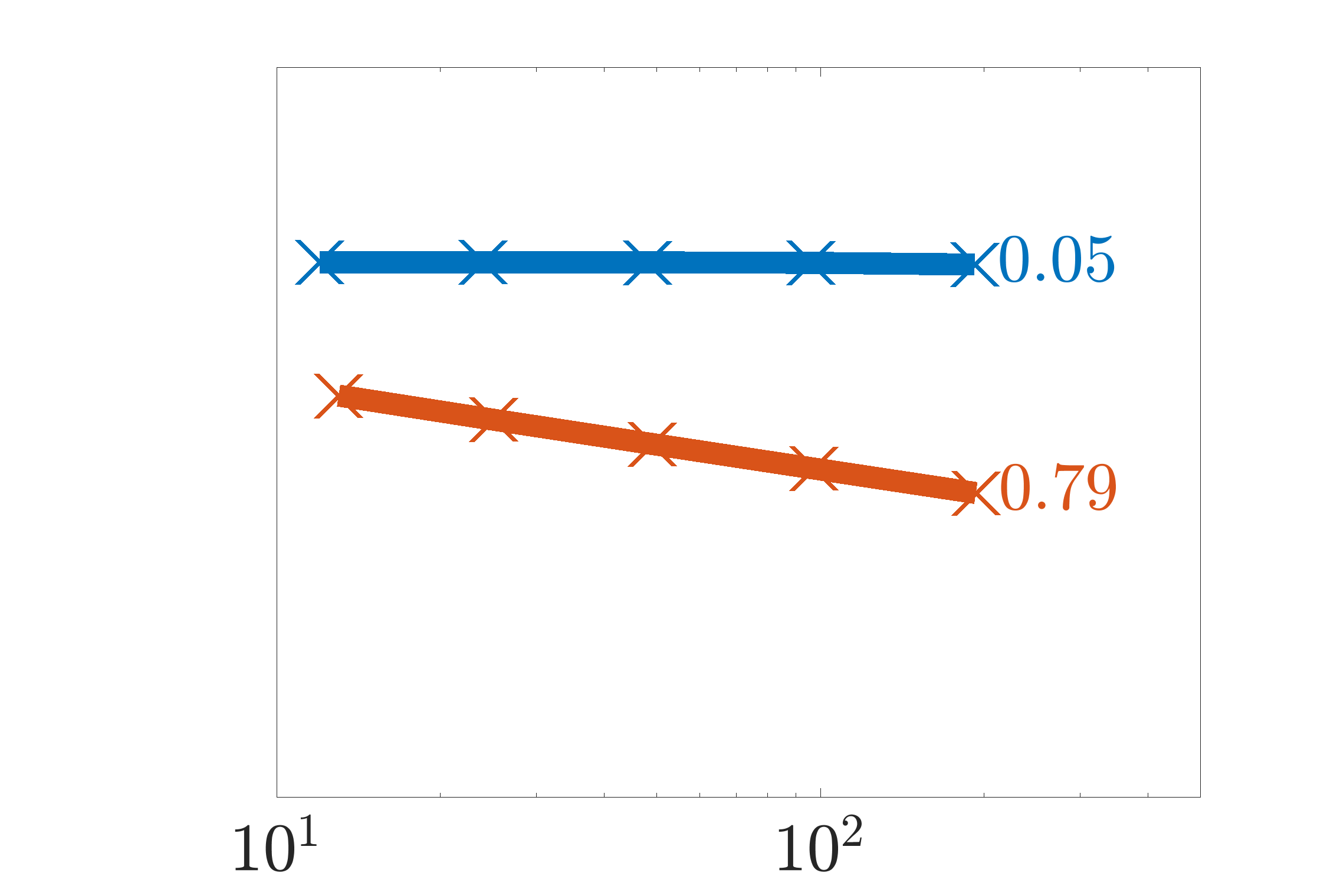}}
\end{subfigmatrix}
\caption{Convergence in the relative $L_2$ integral norms for the boundary layer problem with $\nu = 0$ using Heaviside enrichments}
\label{fig:Example1_p1GFEM_nu0_L2norm}
\end{center}
\end{figure}

\begin{figure}[ht!]
\begin{center}
\begin{subfigmatrix}{6}
\subfigure[$t = 0$]{\includegraphics[width=2.1in]{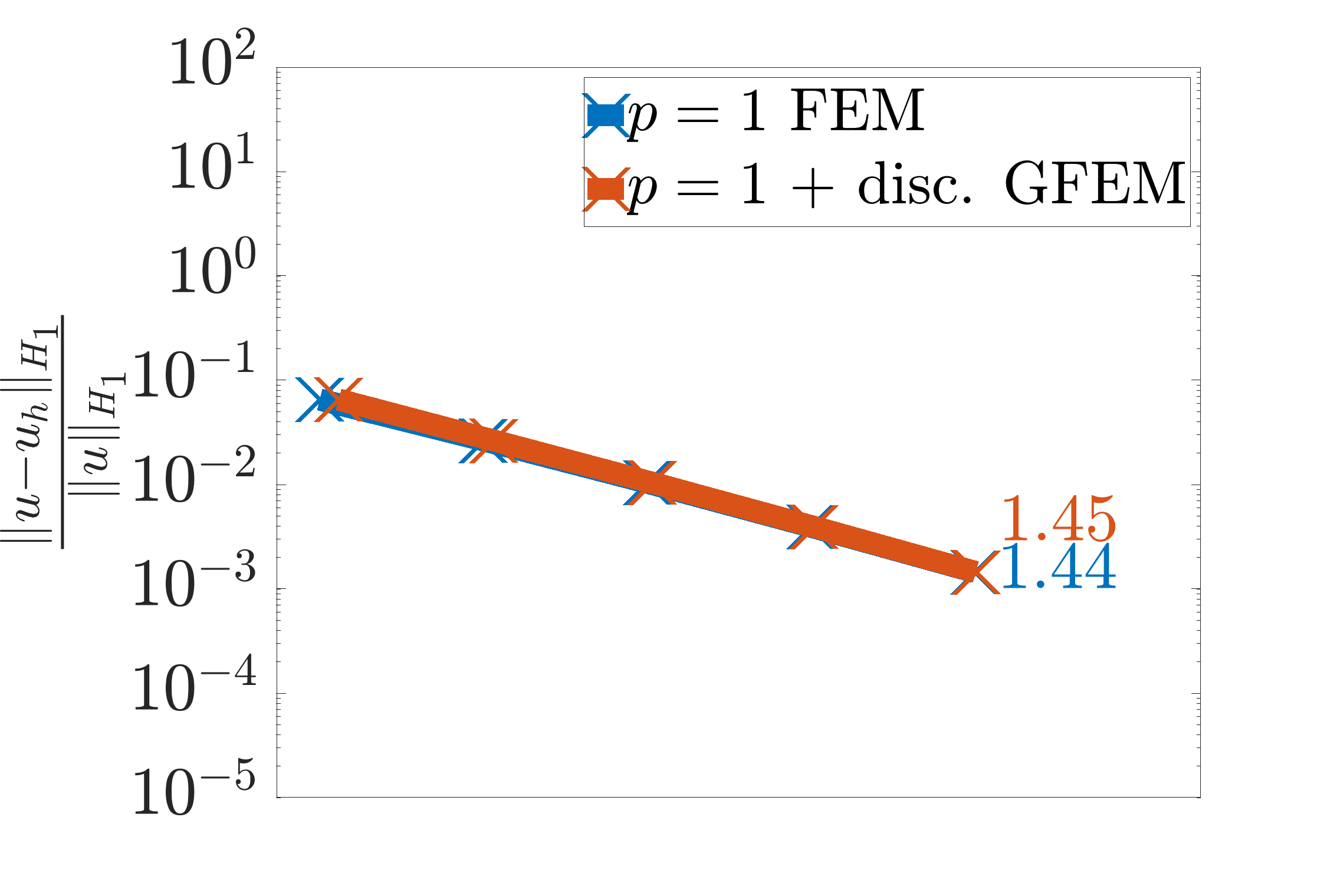}}
\subfigure[$t = 0.25$]{\includegraphics[width=2.1in]{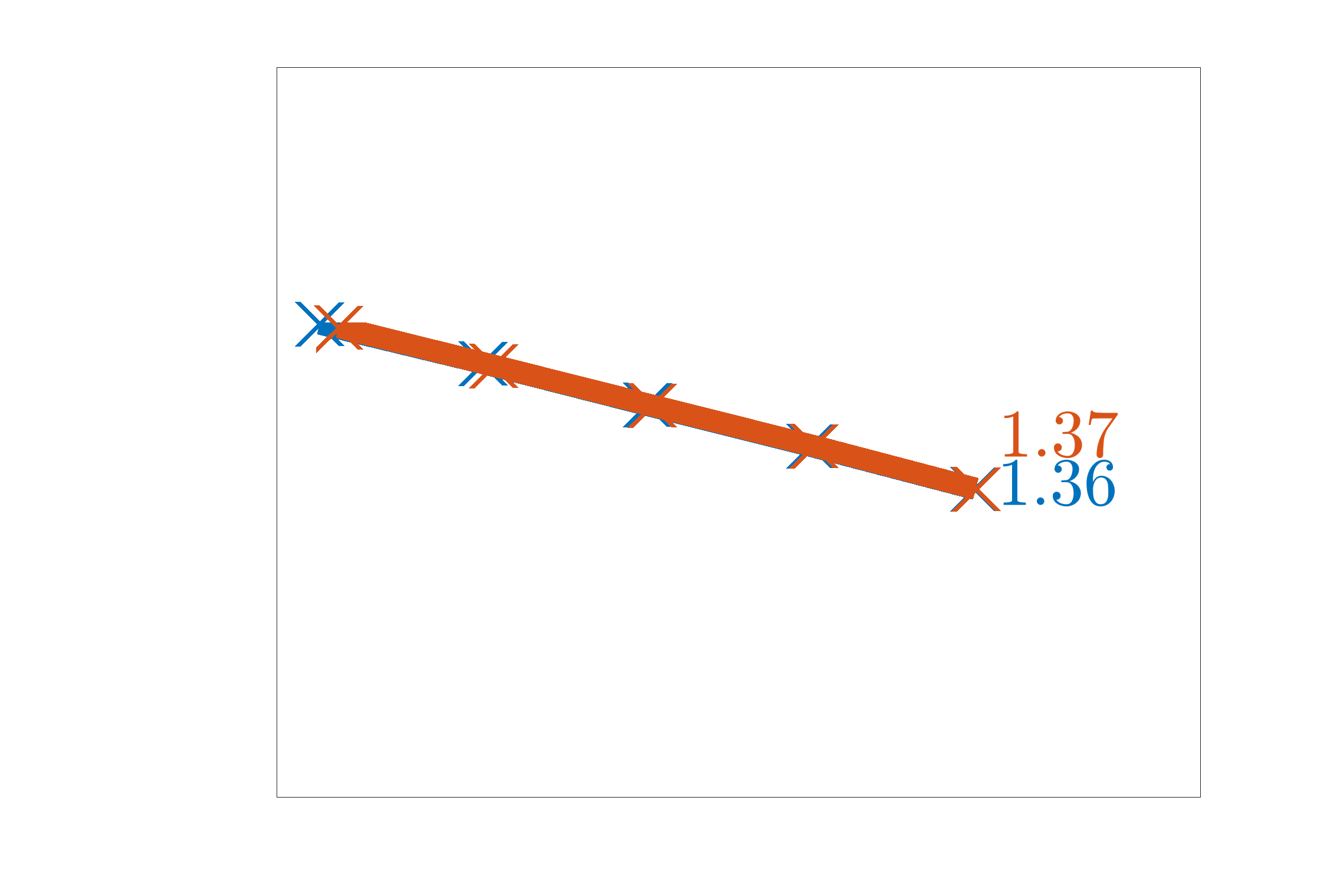}}
\subfigure[$t = 0.318 \approx \frac{1}{\pi}$]{\includegraphics[width=2.1in]{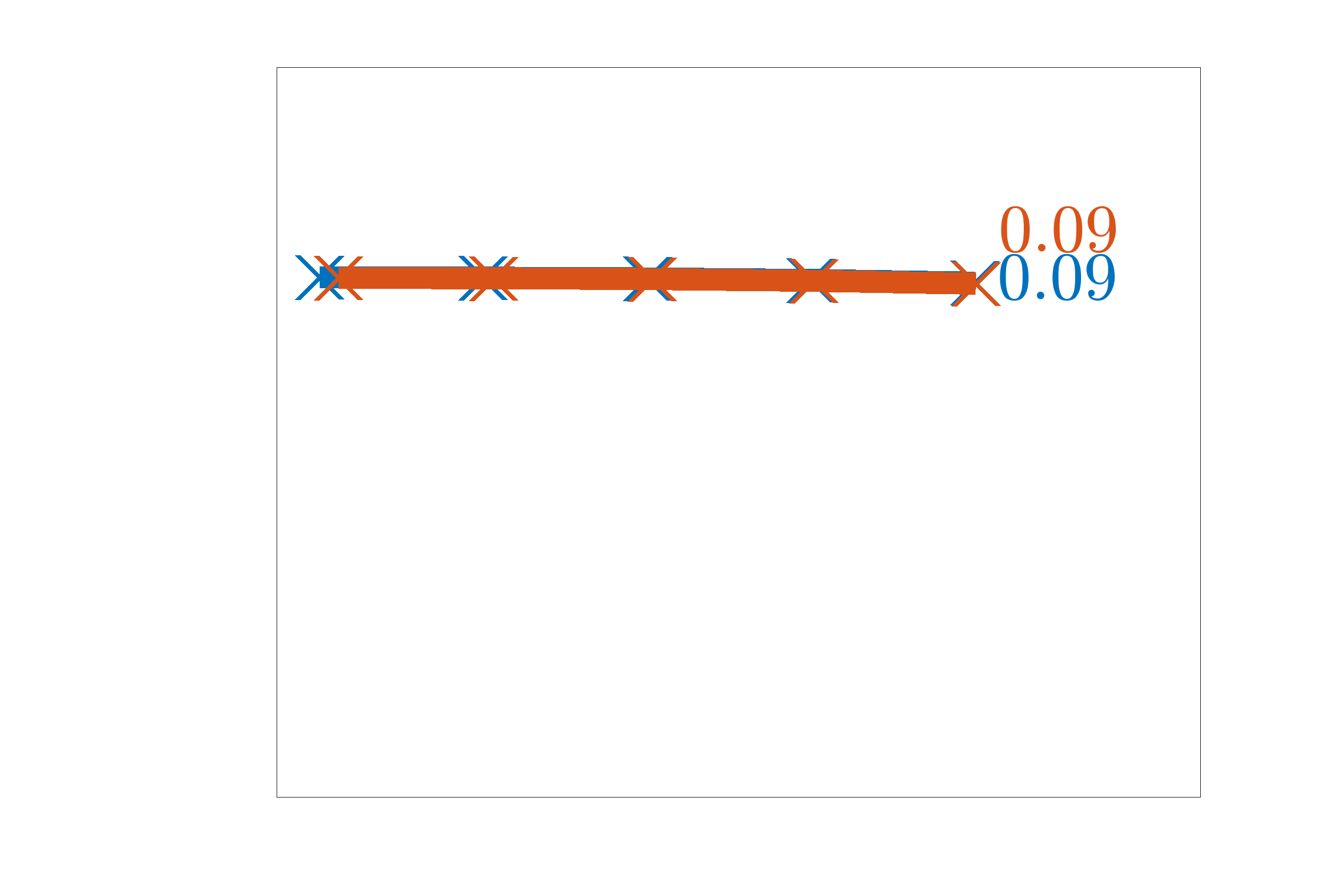}}
\subfigure[$t = 0.5$]{\includegraphics[width=2.1in]{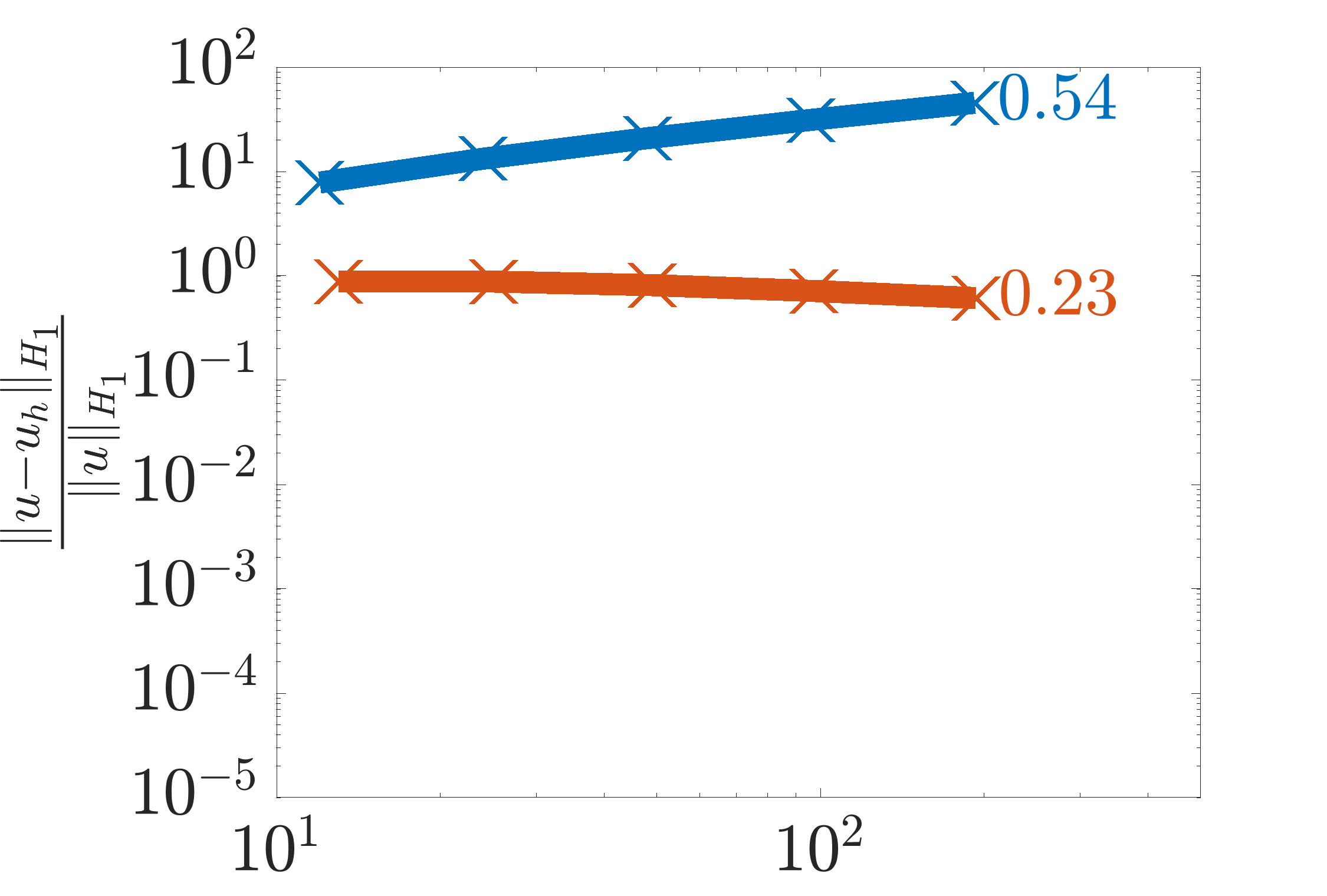}}
\subfigure[$t = 0.75$]{\includegraphics[width=2.1in]{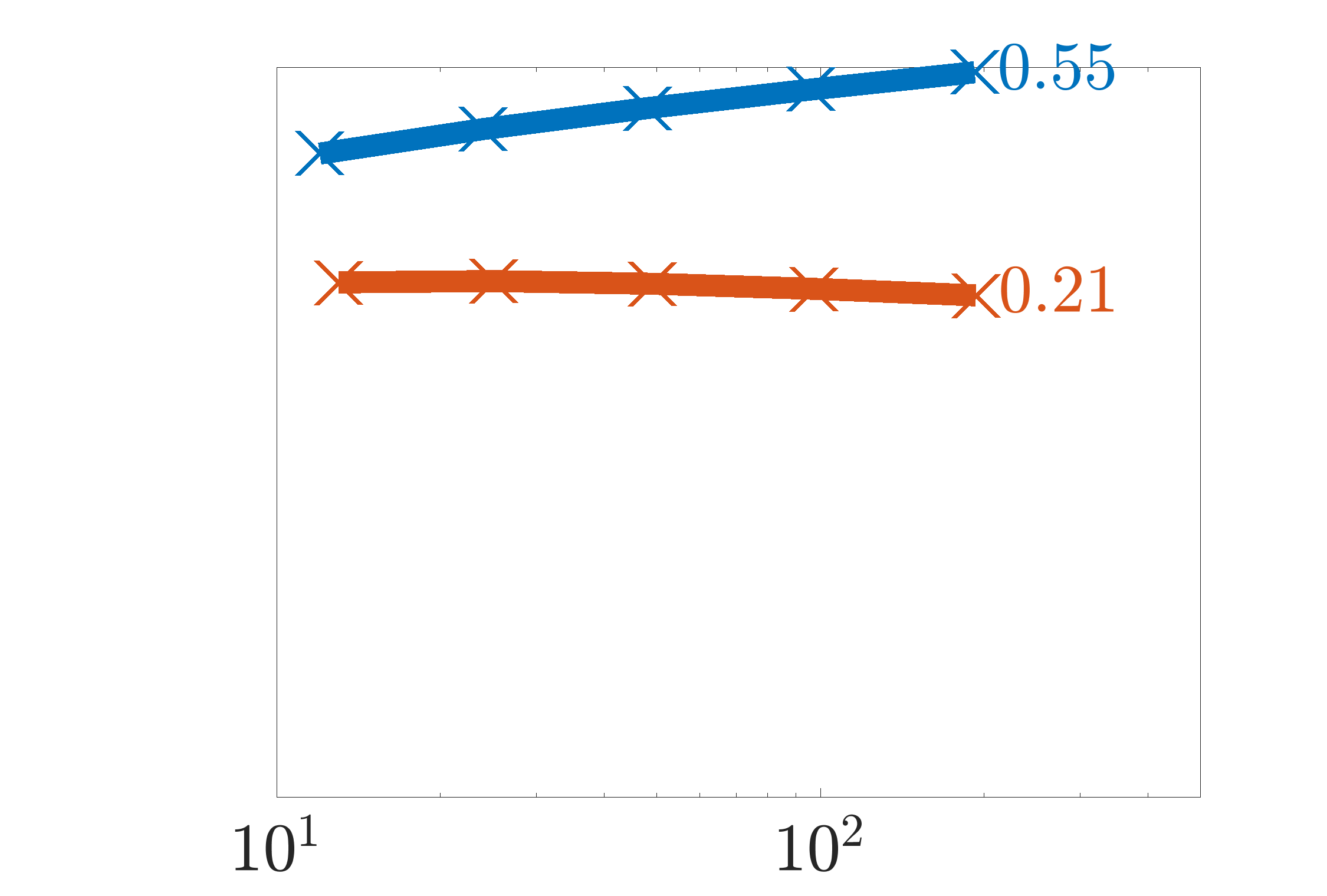}}
\subfigure[$t = 1$]{\includegraphics[width=2.1in]{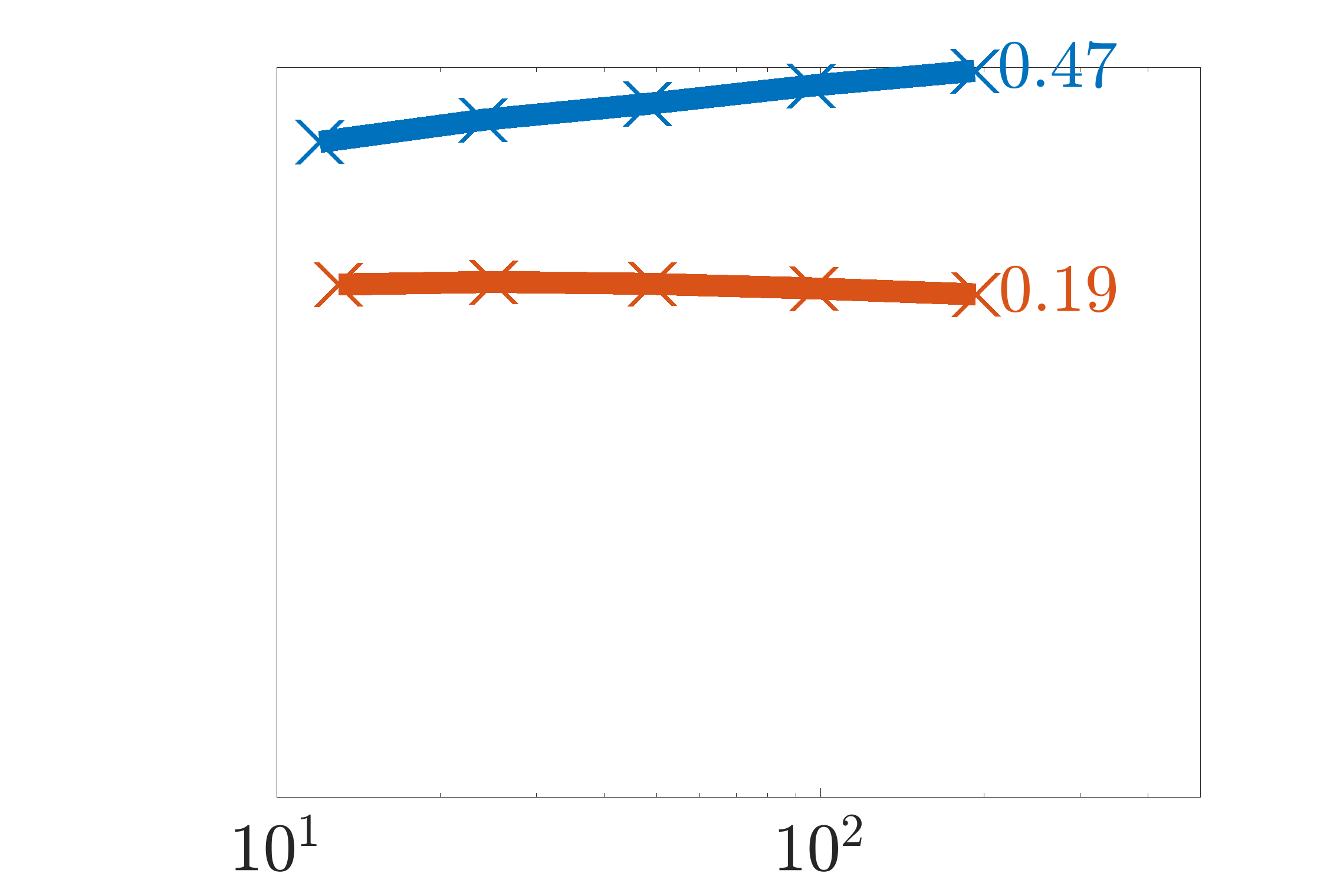}}
\end{subfigmatrix}
\caption{Convergence in the relative $H_1$ integral norms for the boundary layer problem with $\nu = 0$ using Heaviside enrichments}
\label{fig:Example1_p1GFEM_nu0_H1norm}
\end{center}
\end{figure}

\begin{figure}[ht!]
\begin{center}
\begin{subfigmatrix}{2}
\subfigure[Relative $L_2$ integral norm]{\includegraphics[width=3in]{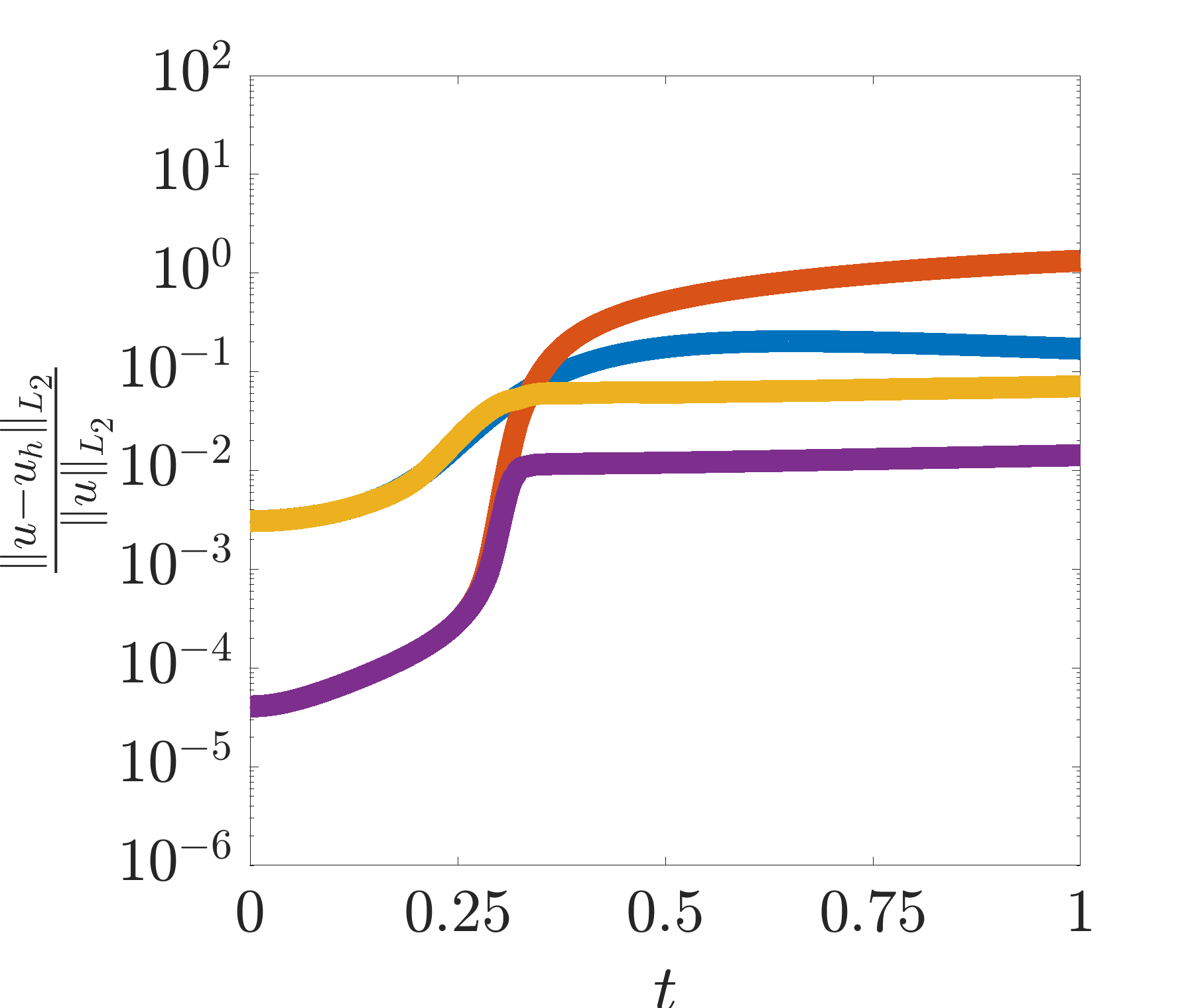}}
\subfigure[Relative $H_1$ integral norm]{\includegraphics[width=3in]{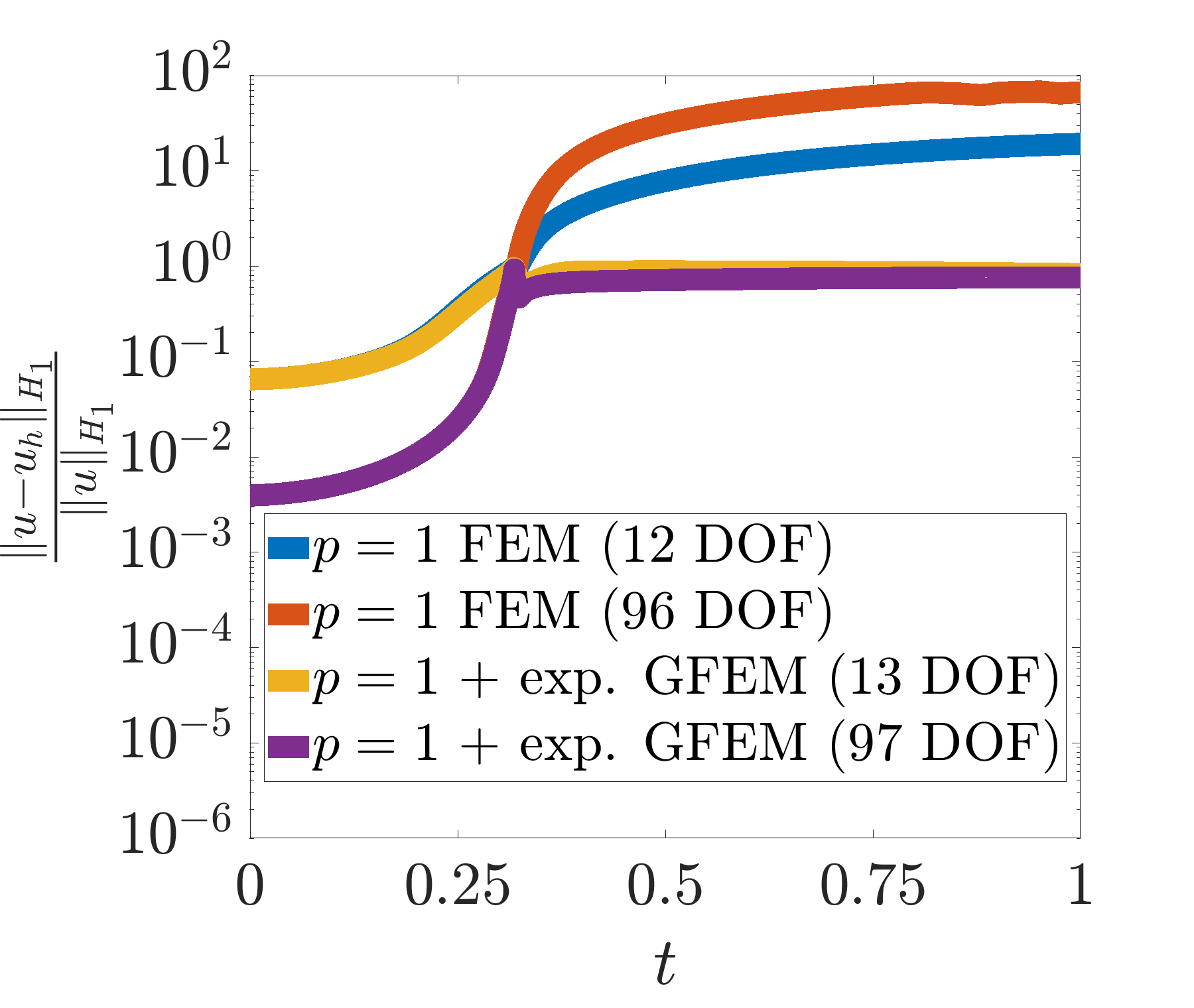}}
\end{subfigmatrix}
\caption{Relative $L_2$ and $H_1$ integral norms versus time for 11-element and 95-element FEM and $p = 1$ + disc. GFEM for the boundary layer problem with $\nu = 0$}
\label{fig:Example1_p1expGFEM_L2H1vstime_nu0}
\end{center}
\end{figure}

\begin{figure}[ht!]
\begin{center}
\begin{subfigmatrix}{3}
\subfigure[Reference]{\includegraphics[width=2.1in]{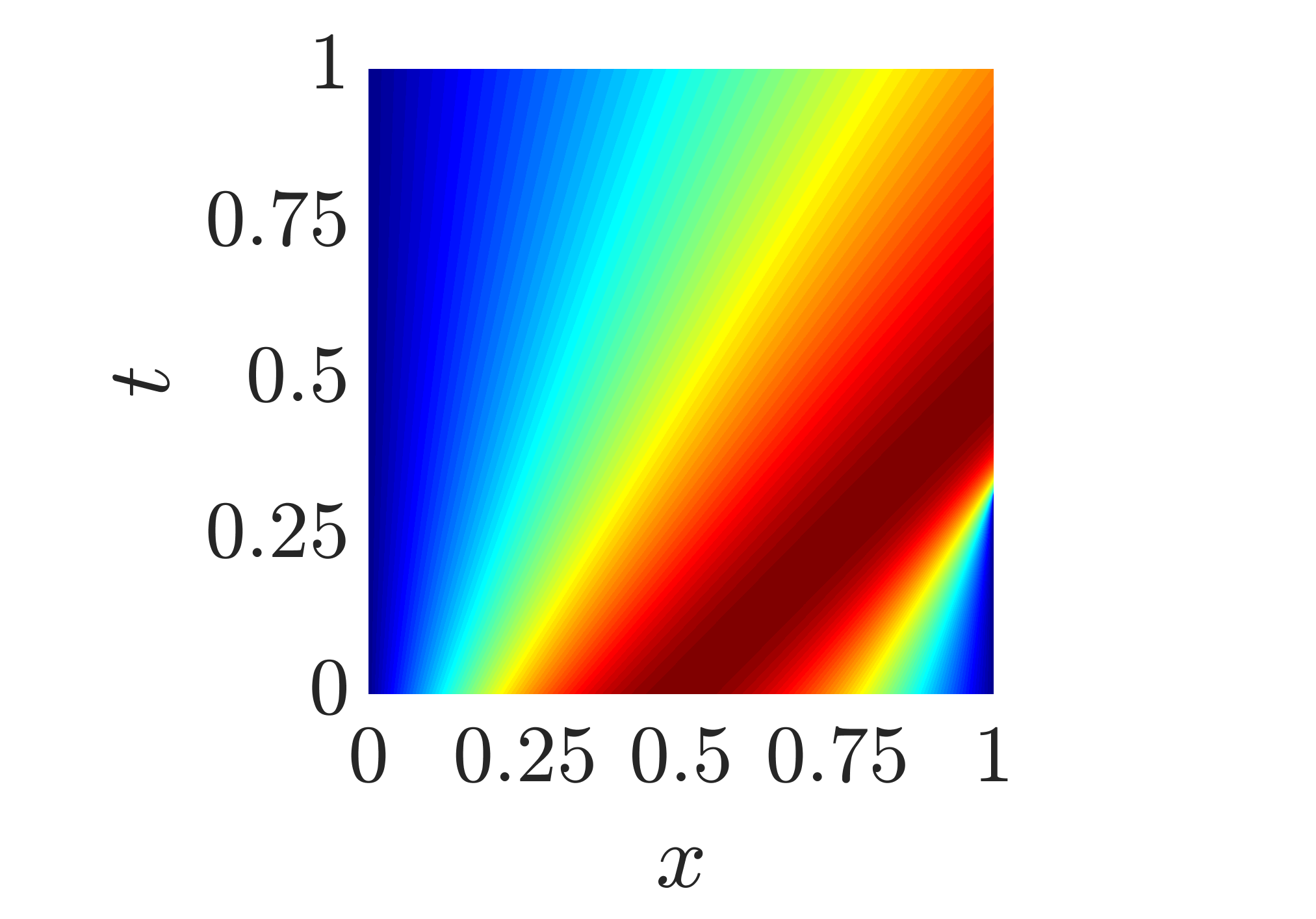}}
\subfigure[12 DOF FEM]{\includegraphics[width=2.1in]{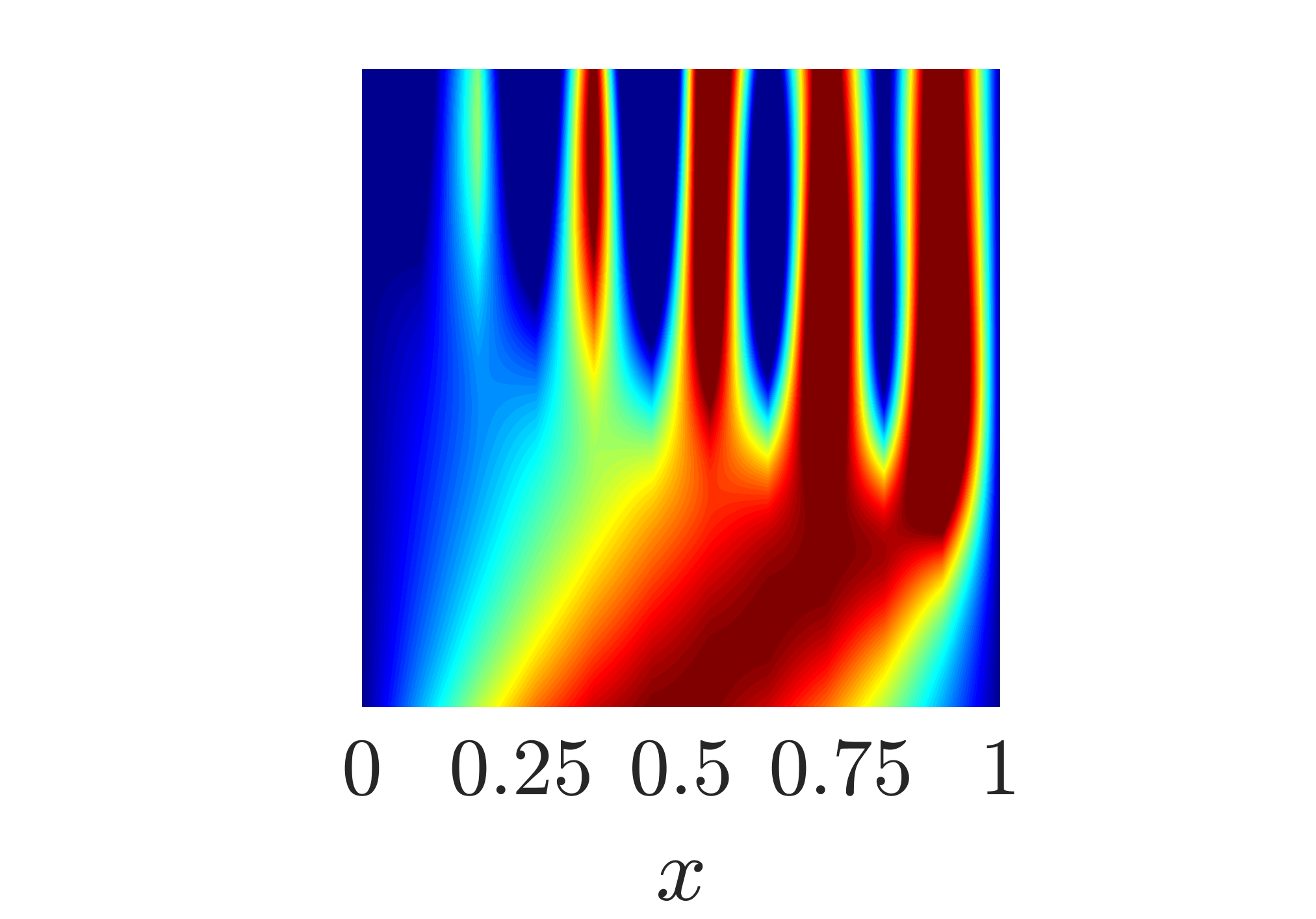}}
\subfigure[13 DOF GFEM]{\includegraphics[width=2.1in]{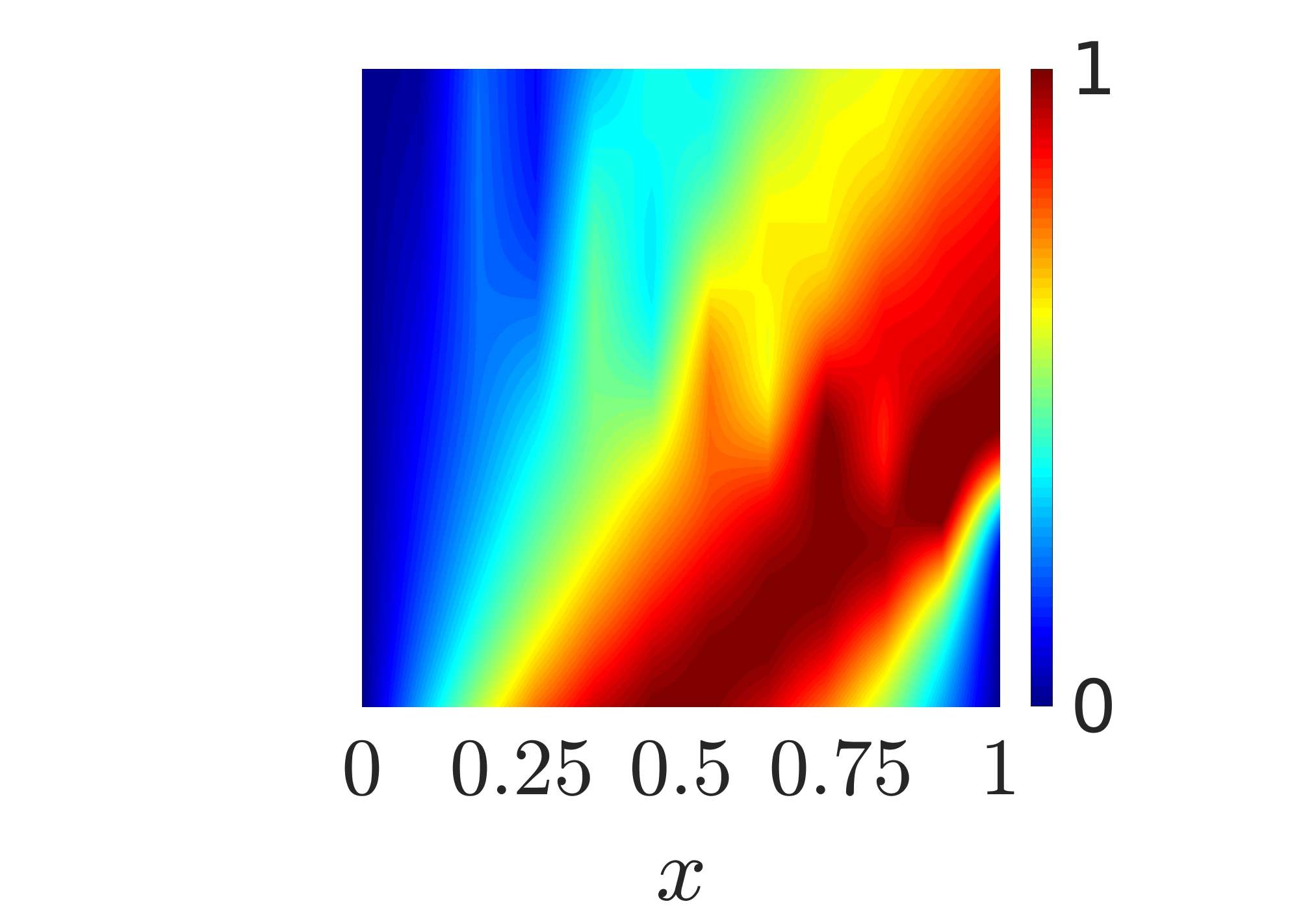}}
\end{subfigmatrix}
\caption{11-element $p = 1$ FEM (12 DOF) $p = 1$ + disc. GFEM (13 DOF) solution contours compared to the reference for the boundary layer problem with $\nu = 0$}
\label{fig:Example1_nu0_FEMGFEM_contour_comparison_1}
\end{center}
\end{figure}

\begin{figure}[ht!]
\begin{center}
\begin{subfigmatrix}{3}
\subfigure[Reference]{\includegraphics[width=2.1in]{figures/Example1_reference_nu0_v2.png}}
\subfigure[48 DOF FEM]{\includegraphics[width=2.1in]{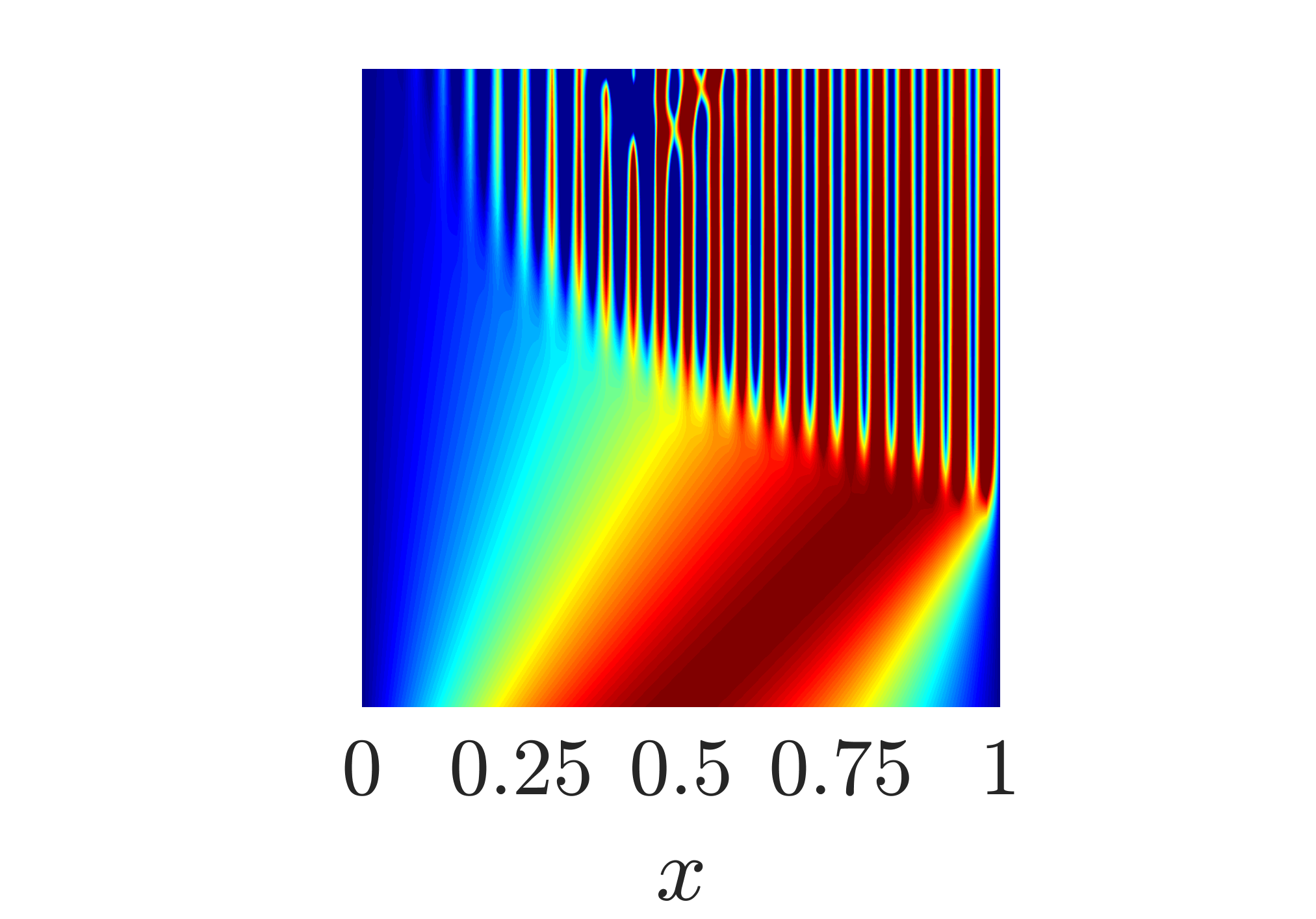}}
\subfigure[49 DOF GFEM]{\includegraphics[width=2.1in]{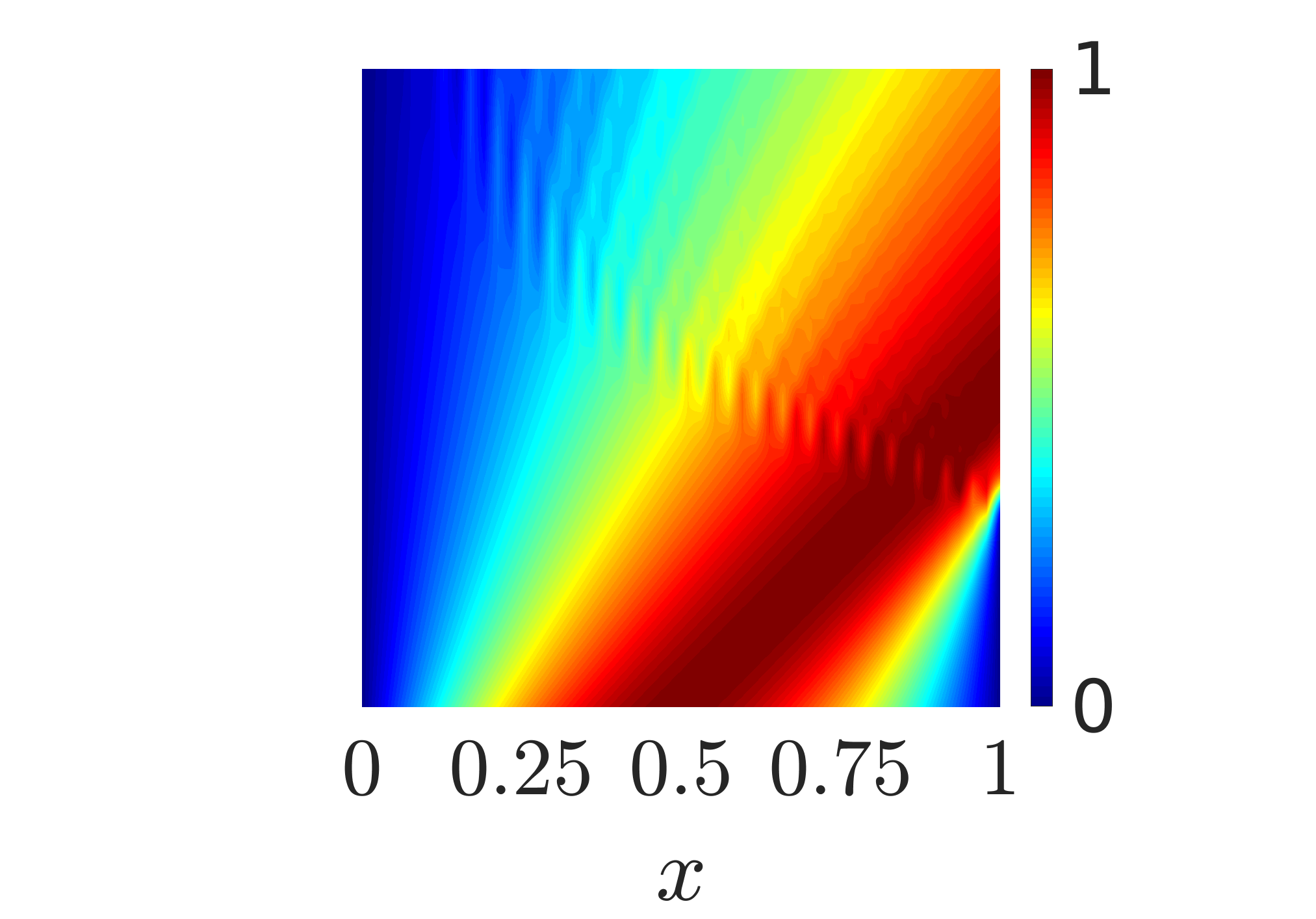}}
\end{subfigmatrix}
\caption{47-element $p = 1$ FEM (48 DOF) $p = 1$ + disc. GFEM (49 DOF) solution contours compared to the reference for the boundary layer problem with $\nu = 0$}
\label{fig:Example1_nu0_FEMGFEM_contour_comparison_2}
\end{center}
\end{figure}

When $\nu = 0$, the presented $p = 1$ FEM is ill-suited for two reasons:  1) enforcement of the Dirichlet boundary condition $u(1,t) = 0$ inhibits the FEM solution from capturing the discontinuity occuring at $x = 1$, and 2) increasingly steep gradients leading up to $t = \frac{1}{\pi}$ (infinitely steep at $t = \frac{1}{\pi}$) results in spurious oscillations which propagate with time. To address the first challenge, the last node in the domain at $x = 1$ is enriched with a Heaviside function which is 0 everywhere except the element containing the node. This may be thought of as a linear correction, which allows the GFEM to satisfy $u(1,t) = 0$, but also account for the discontinuity which arises. Results using the Heaviside enrichment are denoted as $p = 1$ + disc. GFEM solutions. Convergence in the relative $L_2$ and $H_1$ integral norm at various times are shown in Figs. \ref{fig:Example1_p1GFEM_nu0_L2norm} and \ref{fig:Example1_p1GFEM_nu0_H1norm}, respectively. Convergence rates are computed and presented using $h = \Big\{\frac{1}{95}, \frac{1}{191}\Big\}$. Recall no theoretical convergence rates are formally developed for GFEM using the enrichments in this example. Here convergence in the $L_2$ norm is significantly improved in the GFEM solutions with respect to linear FEM, however, the GFEM using Heaviside enrichment still yields poor convergence in the $H_1$ norm. Plots of the relative $L_2$ and $H_1$ integral norms versus time are shown in Fig. \ref{fig:Example1_p1expGFEM_L2H1vstime_nu0} for 11- and 95-element GFEM solutions. Similar to the $\nu = \frac{1}{100}$ results, the FEM and GFEM solutions have similar error levels until $t \approx \frac{1}{\pi}$ where boundary layer gradients become larger. Around $t \approx \frac{1}{\pi}$, both the FEM and GFEM solutions rise in error. This is explained as both the linear interpolation and Heaviside function are ill-suited for capturing the increasingly steep gradients leading up to the discontinuity at $t = \frac{1}{\pi}$. This is observed in Fig. \ref{fig:Example1_nu0_FEMGFEM_contour_comparison_1}, which provides 11-element $p = 1$ FEM and $p = 1$ + disc. GFEM solution contours. Severe oscillations are observed in the $p = 1$ FEM solutions, whereas comparatively the GFEM solution is significantly better. However, oscillations still persist in the GFEM solutions, starting around $t = \frac{1}{\pi}$, and propagate with time. With grid refinement as shown in Fig. \ref{fig:Example1_nu0_FEMGFEM_contour_comparison_2}, which provides 47-element $p = 1$ FEM and $p = 1$ + disc. GFEM solution contours, the FEM solution does not improve. However, while oscillations persist with the GFEM solution, they are significantly muted. These results demonstrate the importance of the intermediate, transitional solution features on stability of the GFEM solution for the Burgers' equation. Specific to this example, the infinitely steep boundary formed in the limit as $t \rightarrow \frac{1}{\pi}$ results in spurious oscillations even though the discontinuity is captured. To better explore the effect the intermediate solution features have on the stability of the GFEM, a second example is examined in which a shock forms with a known steady state solution.

\subsection{Example 2: Shock formulation in the domain}

\subsubsection{Problem statement and reference solutions}
Consider the viscous Burgers' equation (Eq. \ref{viscous_burger_strong_formulation}) defined over a unit domain ($\Omega = [0, 1]$) and subject to Dirichlet boundary conditions everywhere ($\Gamma = \Gamma_D$). The problem formulation is as follows: For $t \in [0, 1]$, find $u$ such that:

\begin{equation}
\label{example2_viscous_burger_strong_formulation}
\begin{aligned}
\frac{\partial u}{\partial t}  + u \frac{\partial u}{\partial x} - \nu \frac{\partial^2 u}{\partial x^2} = 0& \quad \mbox{on} \quad \Omega\\
u(x,0) = \cos{\pi x}& \quad \mbox{on} \quad \Omega\\
u(0,t) = 1; u(1,t) = -1& \quad \mbox{on} \quad \Gamma\\
\end{aligned}
\end{equation}

\begin{figure}[ht!]
\begin{center}
\begin{subfigmatrix}{4}
\subfigure[$\nu = \frac{1}{50}$]{\includegraphics[width=1.6in]{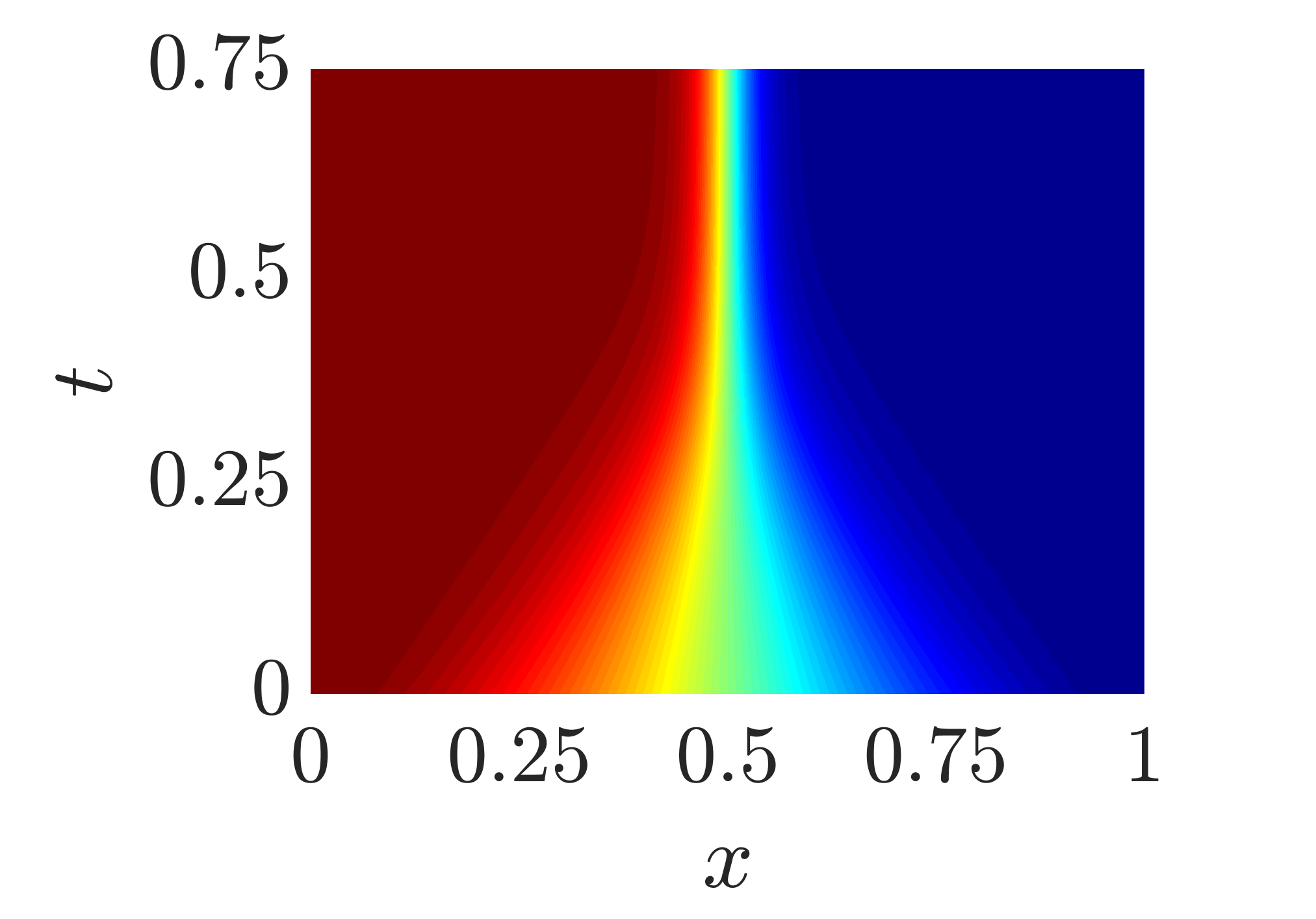}}
\subfigure[$\nu = \frac{1}{100}$]{\includegraphics[width=1.6in]{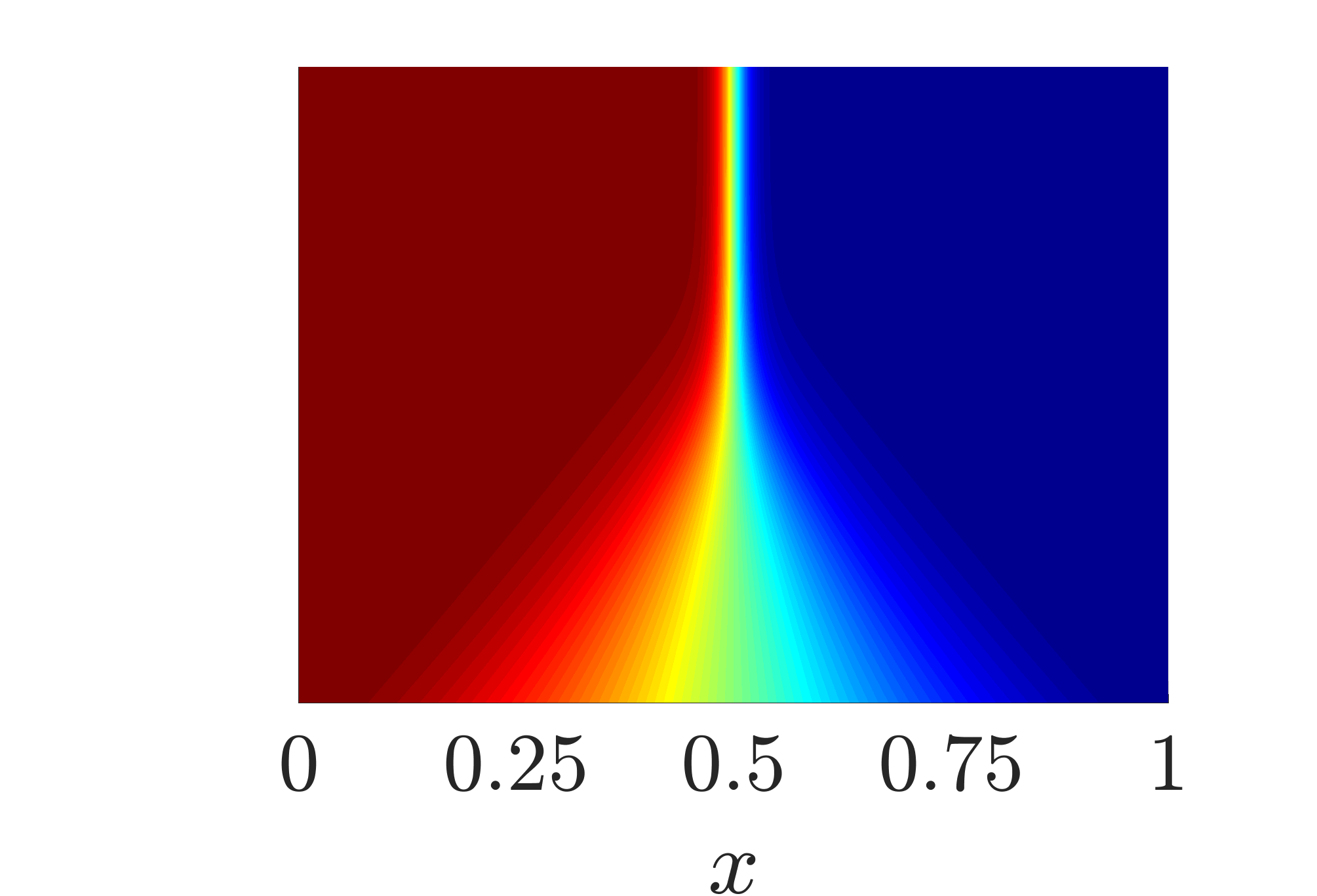}}
\subfigure[$\nu = \frac{1}{500}$]{\includegraphics[width=1.6in]{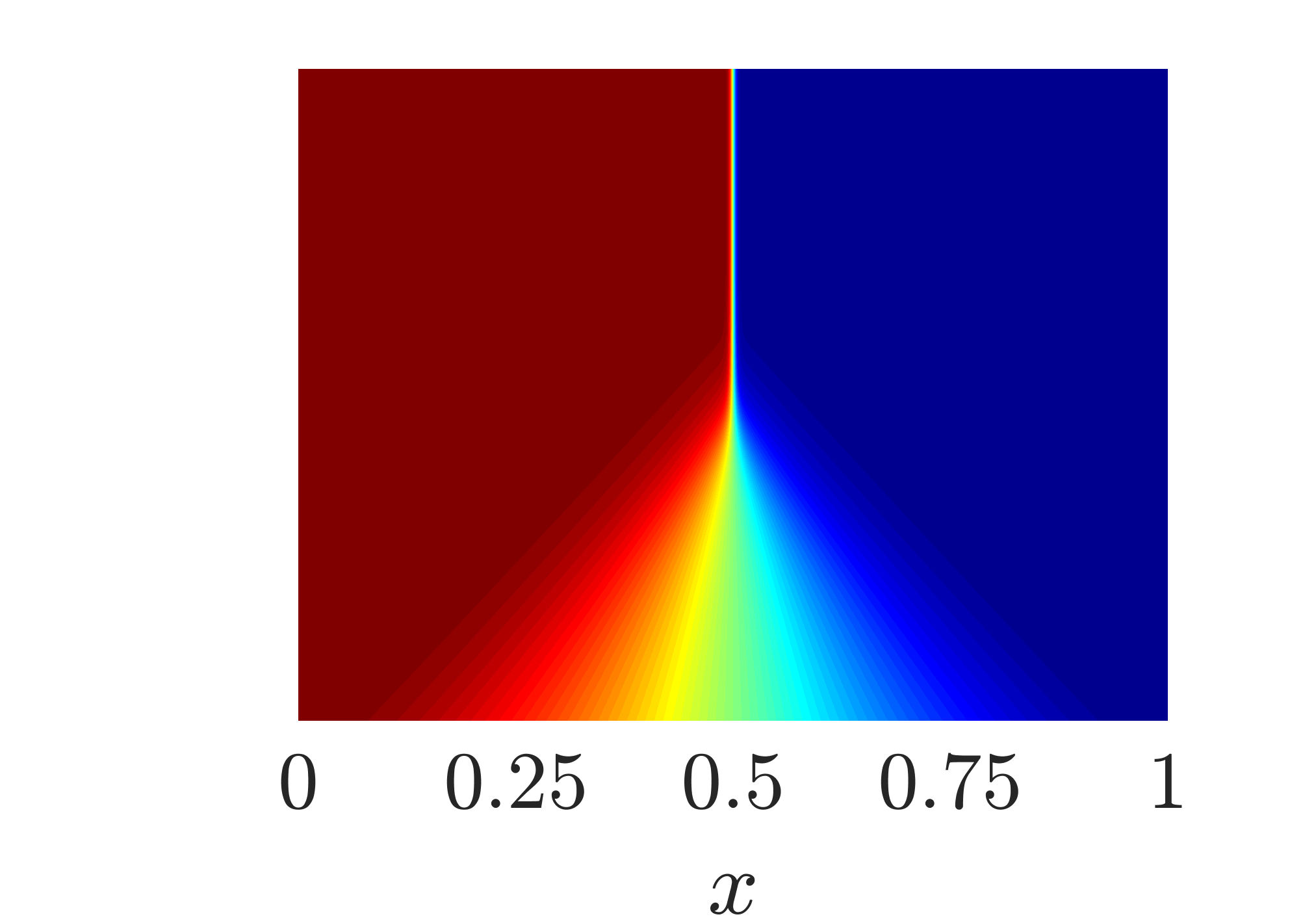}}
\subfigure[$\nu = \frac{1}{1000}$]{\includegraphics[width=1.6in]{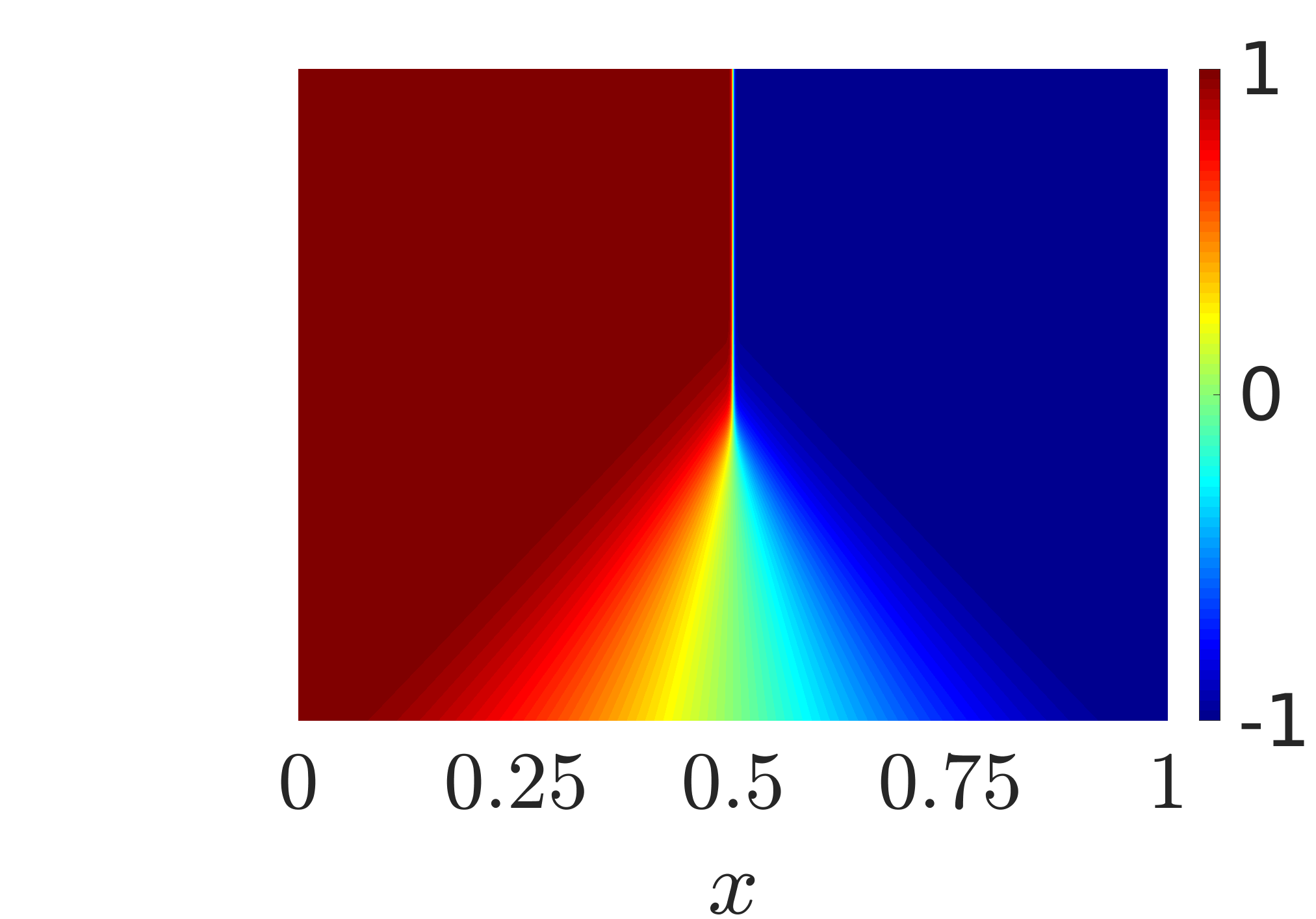}}
\end{subfigmatrix}
\caption{5000-element FEM reference solutions for various kinematic viscosities for the shock problem}
\label{fig:Example2_reference}
\end{center}
\end{figure}

An analytical solution to this problem is unknown, however a steady state solution is provided by $u_{ss} = \sqrt{2 k} \tanh{ \Big[\sqrt{\frac{k}{2 \nu^2}} \Big( \frac{1}{2} - x\Big)} \Big]$, where $k$ is a constant solvable from the nonlinear equation $\sqrt{2k} \tanh \sqrt{\frac{k}{8 \nu^2}} - 1 = 0$. For sufficiently small kinematic visocity $\nu$, the constant  $k \approx \frac{1}{2}$, simplifying the steady state solution to $u_{ss} \approx \tanh{ \Big[\frac{1}{2 \nu} \Big( \frac{1}{2} - x\Big) \Big]}$. The steady state solution represents the instance when the shock has formed entirely, with the shock thickness decreasing as $\nu$ decreases. Shock location is at $x = 0.5$. The 5000-element, $p = 1$ FEM reference solution is shown in Fig. \ref{fig:Example2_reference} for $\nu = \Big[\frac{1}{50}, \frac{1}{100}, \frac{1}{500}, \frac{1}{1000} \Big]$. Note the temporal term for these references were solved over $t = [0, 0.75]$  using the Crank-Nicolson scheme with a step size of $\Delta t = \frac{1}{5000}$. 

\subsubsection{Numerical solutions}
Equation \ref{example2_viscous_burger_strong_formulation} was initially solved over uniform grids $\Big( h = \Big[\frac{1}{11}, \frac{1}{23}, \frac{1}{47}, \frac{1}{95}, \frac{1}{191}\Big] \Big)$ using $p = 1$ FEM and GFEM enriched with the steady state solution for $\nu = \Big[\frac{1}{50}, \frac{1}{100}, \frac{1}{500}, \frac{1}{1000} \Big]$. The local domain the steady state enrichment is applied is given by $\Omega_{local} = \Big[\frac{1}{2} - 2\nu \tanh^{-1}{0.99}-h_{e}, \frac{1}{2} + 2\nu \tanh^{-1}{0.99} + h_{e}\Big]$, which includes all nodes around the shock location $x = \frac{1}{2}$ where the steady state solution $|u_{ss}| \leq 0.99$. GFEM solutions using the steady state solution as an enrichment are denoted by $p = 1$ + ss GFEM solutions. The Crank-Nicolson scheme is used for temporal discretization with a step size of $\Delta t = \frac{1}{5000}$. At each time step the Newton-Raphson method is used to iteratively solve the nonlinear set of equations. 

\begin{figure}[ht!]
\begin{center}
\begin{subfigmatrix}{8}
\subfigure[FEM; $\nu = \frac{1}{50}$]{\includegraphics[width=1.6in]{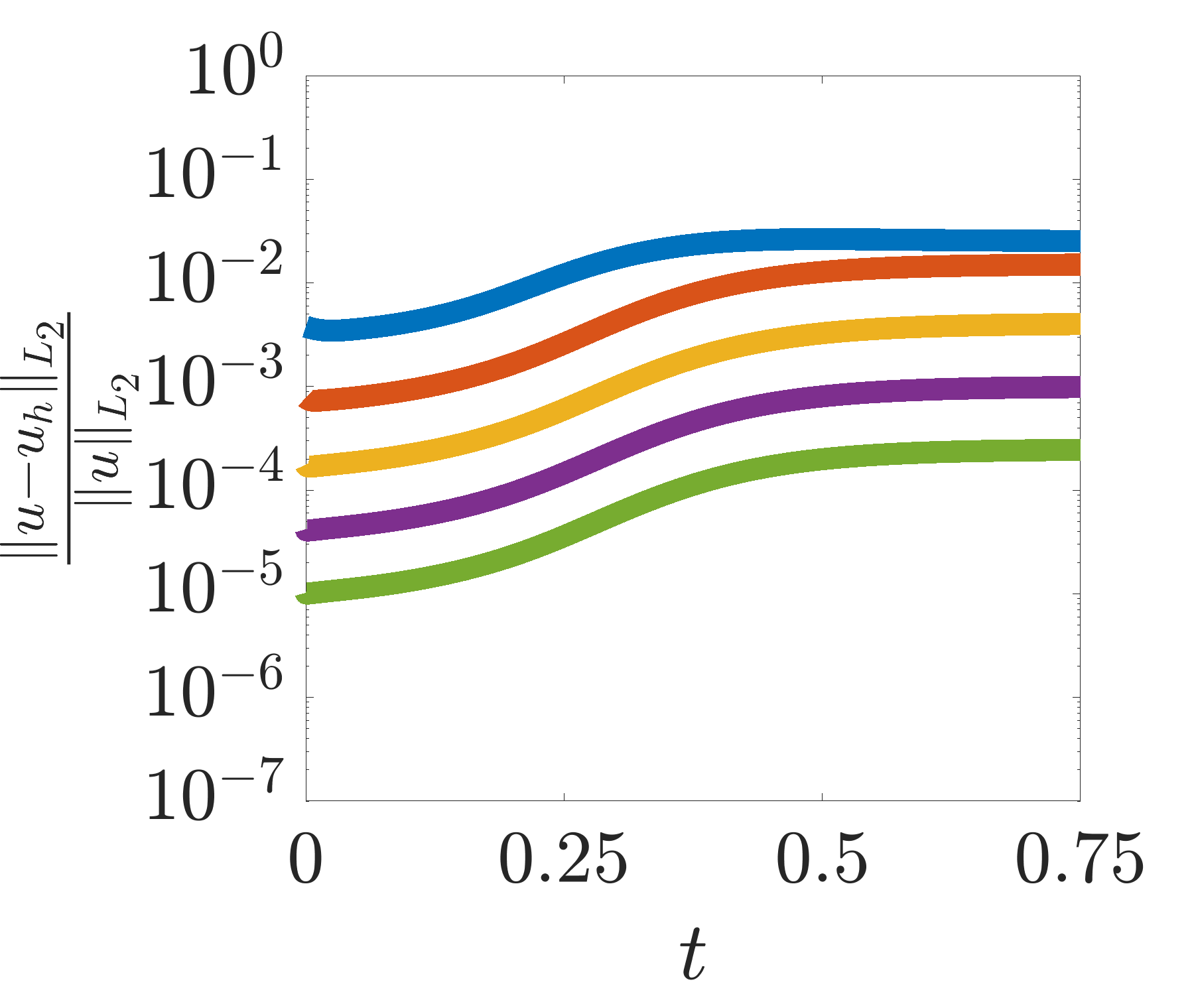}}
\subfigure[FEM; $\nu = \frac{1}{100}$]{\includegraphics[width=1.6in]{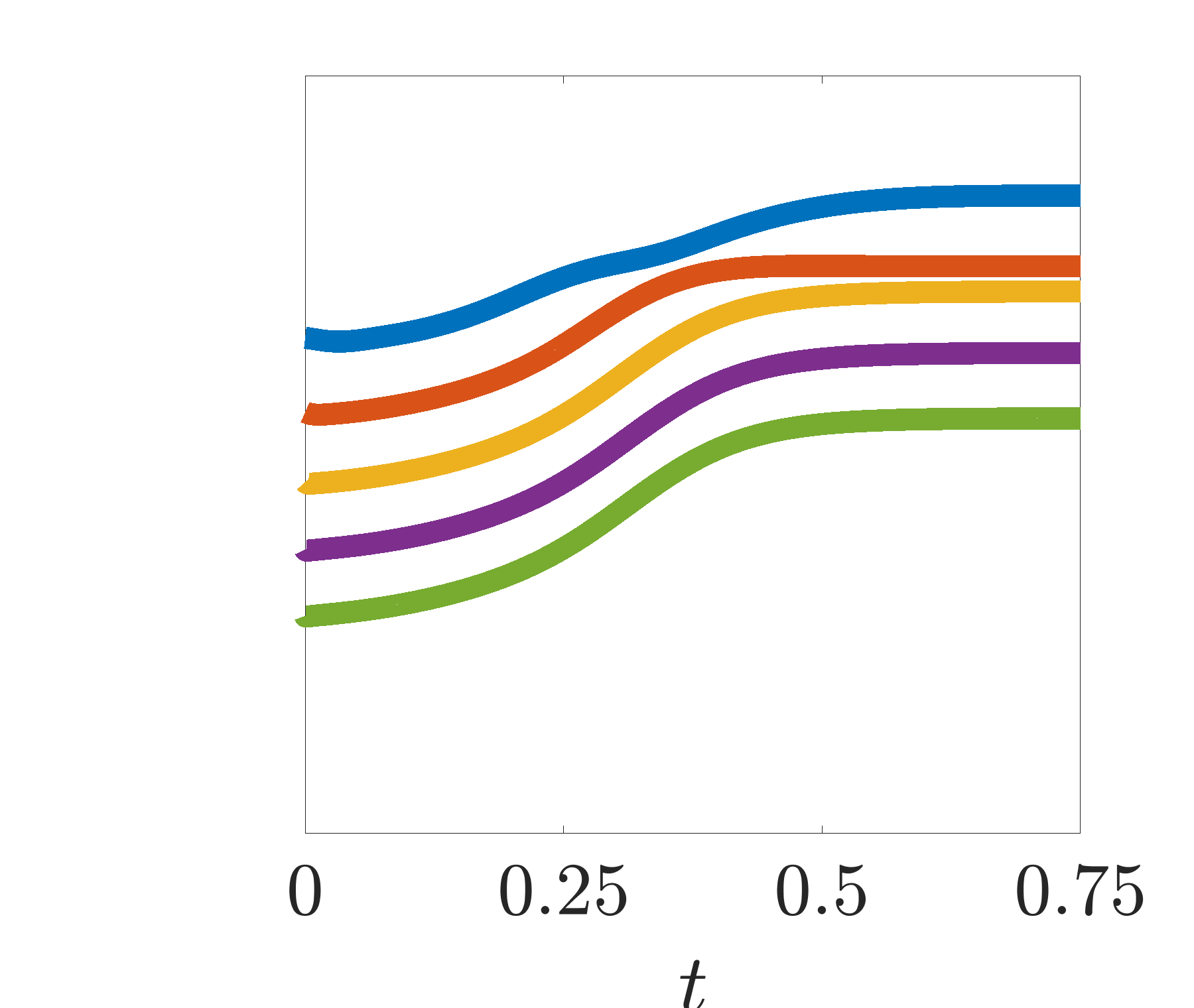}}
\subfigure[FEM; $\nu = \frac{1}{500}$]{\includegraphics[width=1.6in]{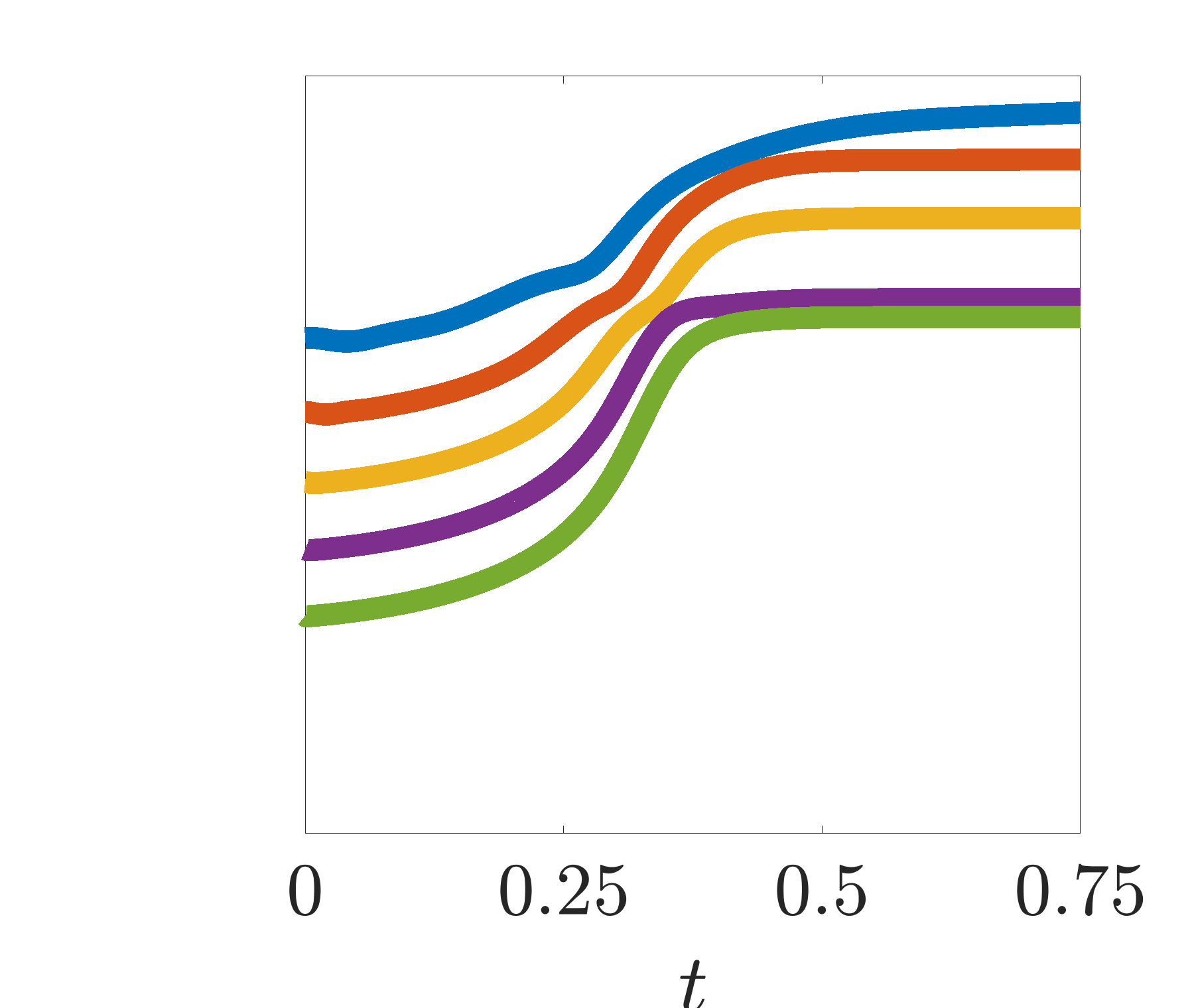}}
\subfigure[FEM; $\nu = \frac{1}{1000}$]{\includegraphics[width=1.6in]{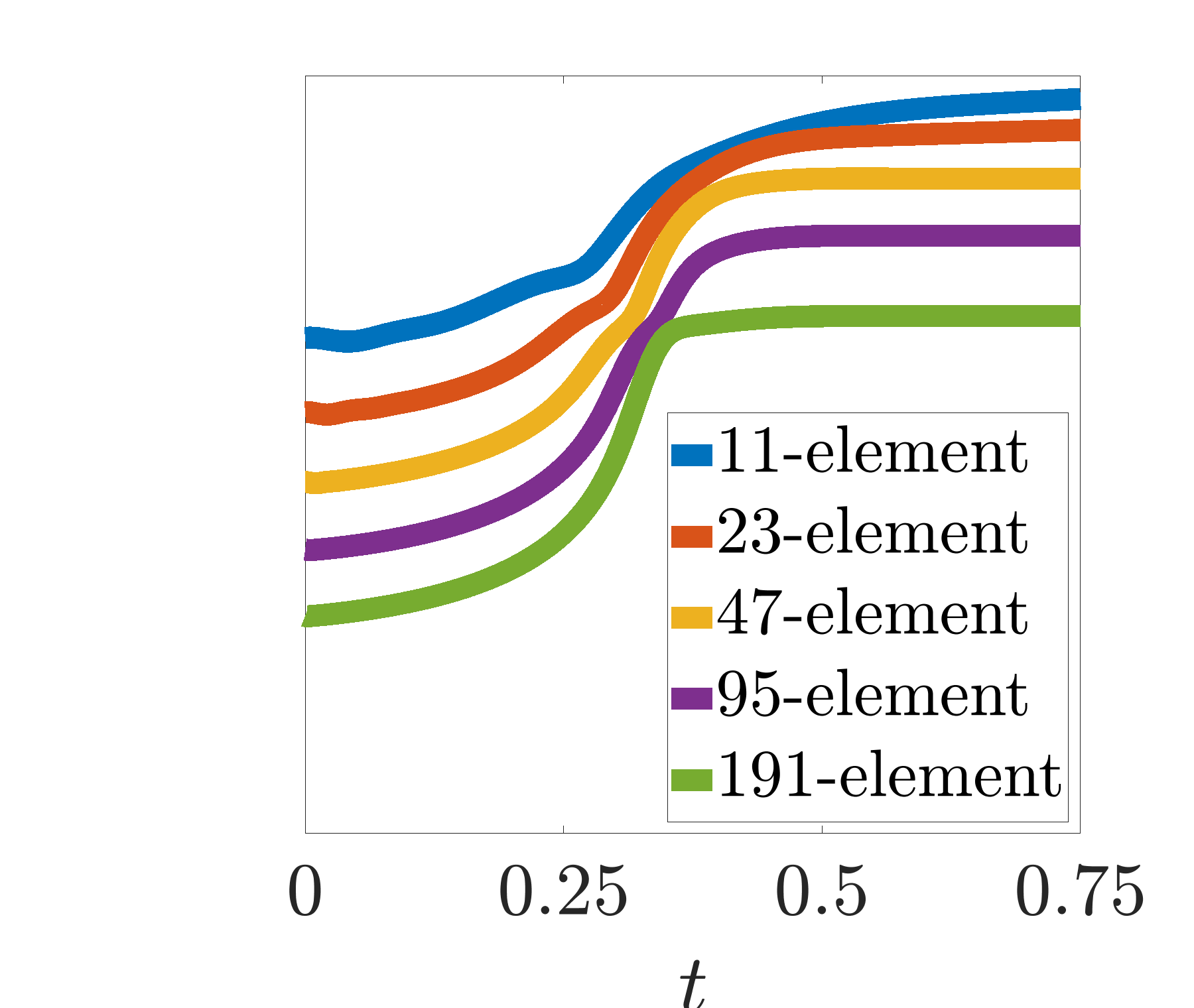}}

\subfigure[GFEM; $\nu = \frac{1}{50}$]{\includegraphics[width=1.6in]{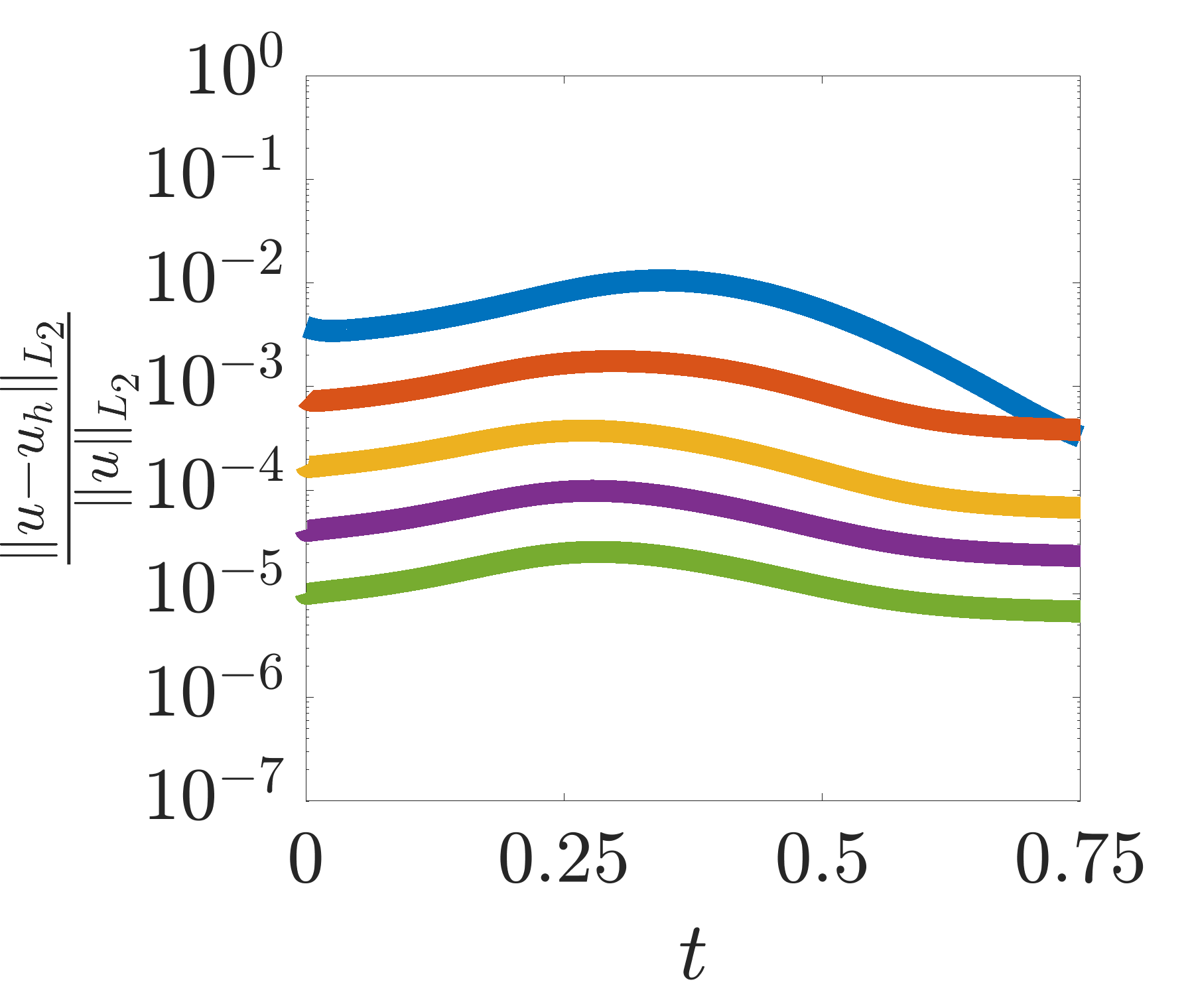}}
\subfigure[GFEM; $\nu = \frac{1}{100}$]{\includegraphics[width=1.6in]{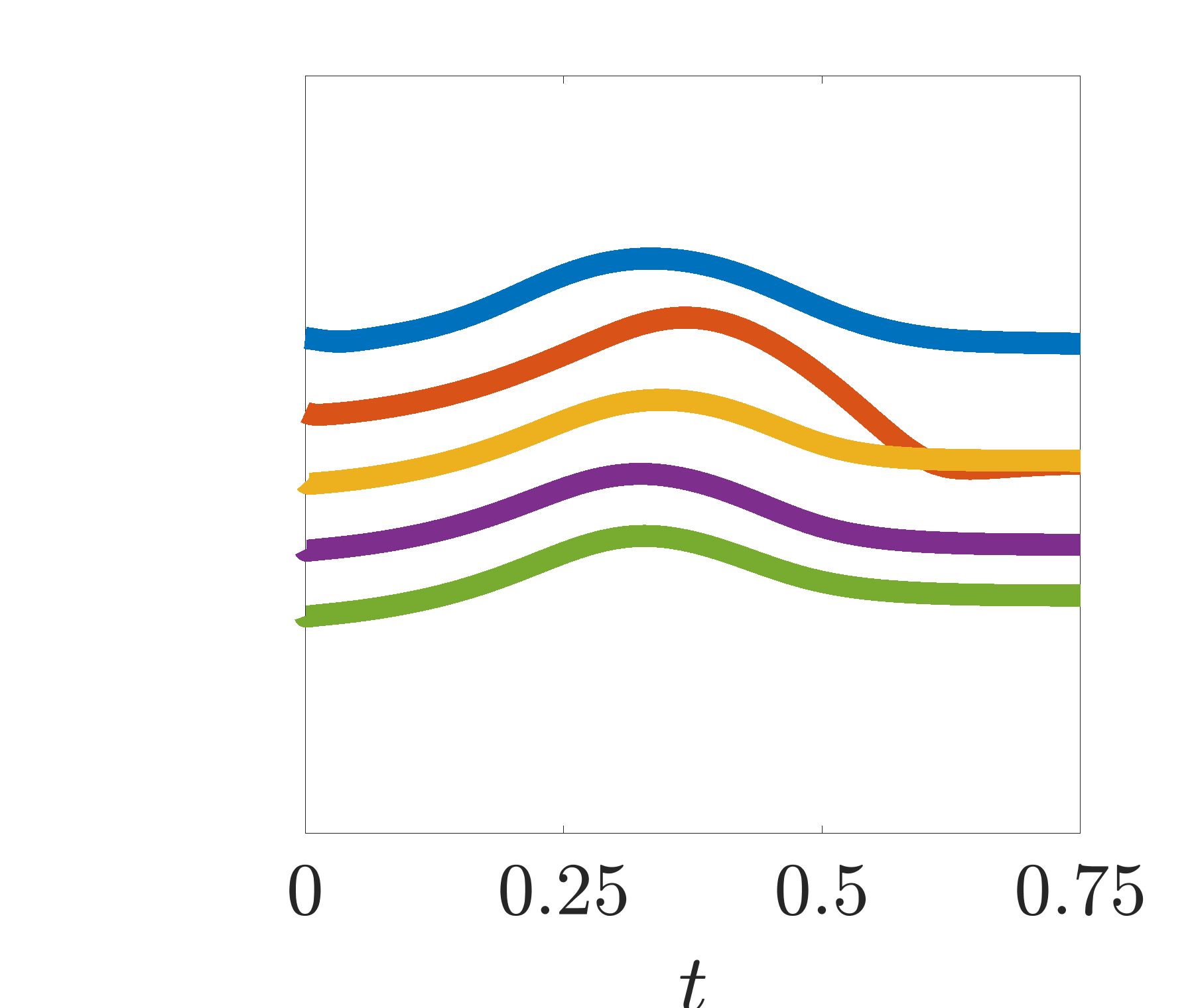}}
\subfigure[GFEM; $\nu = \frac{1}{500}$]{\includegraphics[width=1.6in]{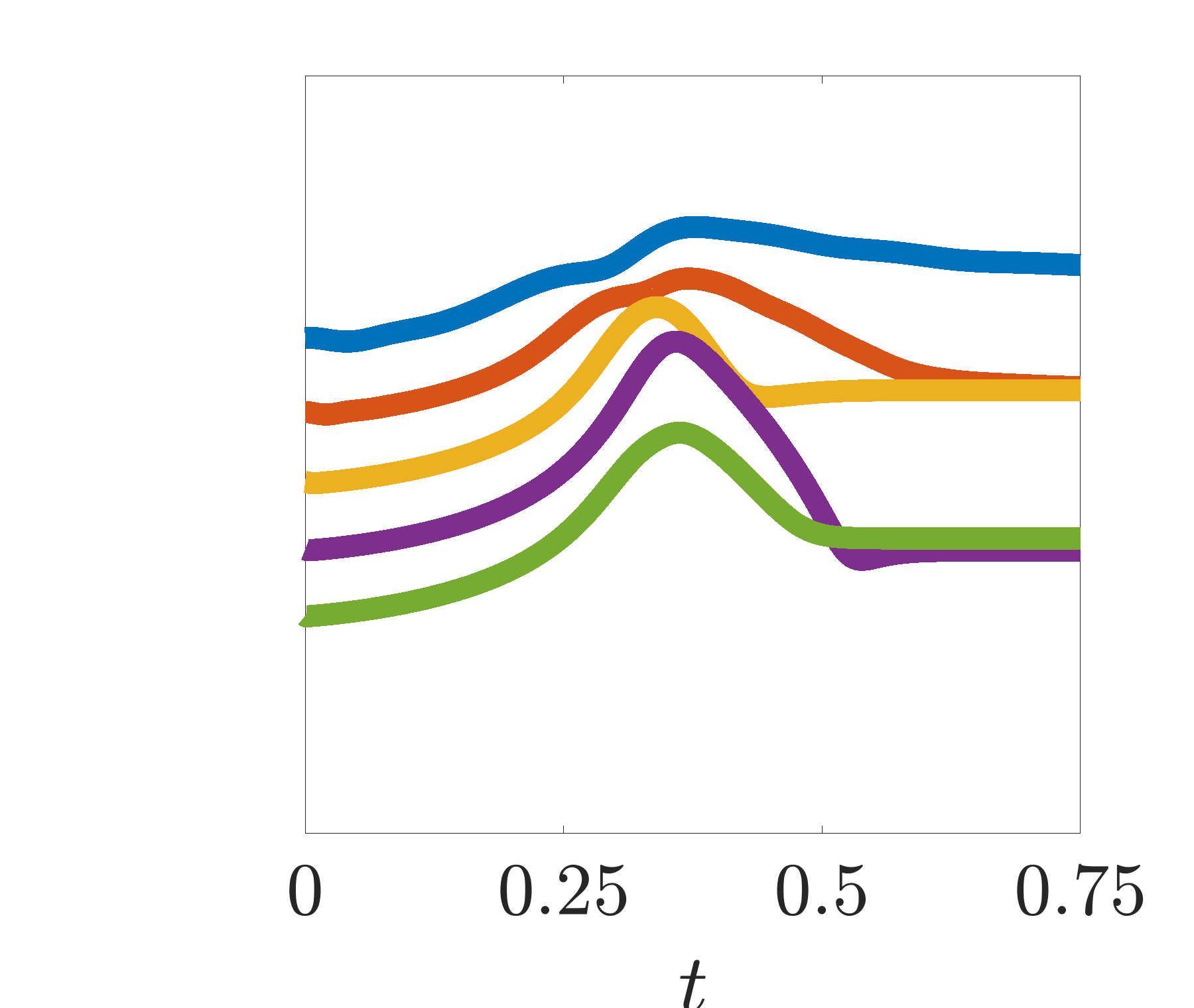}}
\subfigure[GFEM; $\nu = \frac{1}{1000}$]{\includegraphics[width=1.6in]{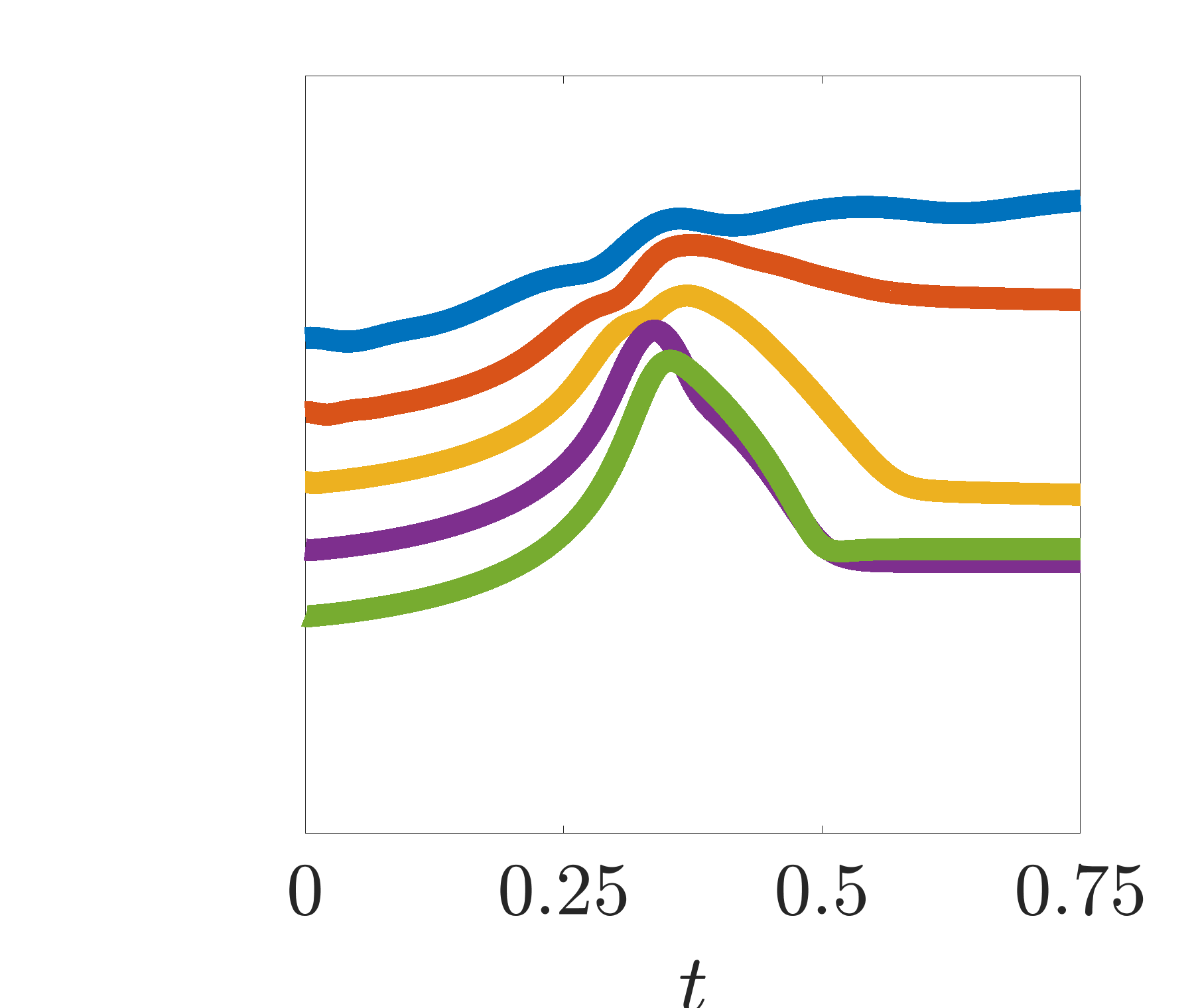}}
\end{subfigmatrix}
\caption{Relative $L_2$ norm versus time for
11-, 23-, 47-, 95-element, and 191-element $p = 1$ FEM and $p = 1$ + ss GFEM solutions over various kinematic viscosities for the shock problem}
\label{fig:Example2_L2vstime}
\end{center}
\end{figure}

\begin{figure}[ht!]
\begin{center}
\begin{subfigmatrix}{8}
\subfigure[FEM; $\nu = \frac{1}{50}$]{\includegraphics[width=1.6in]{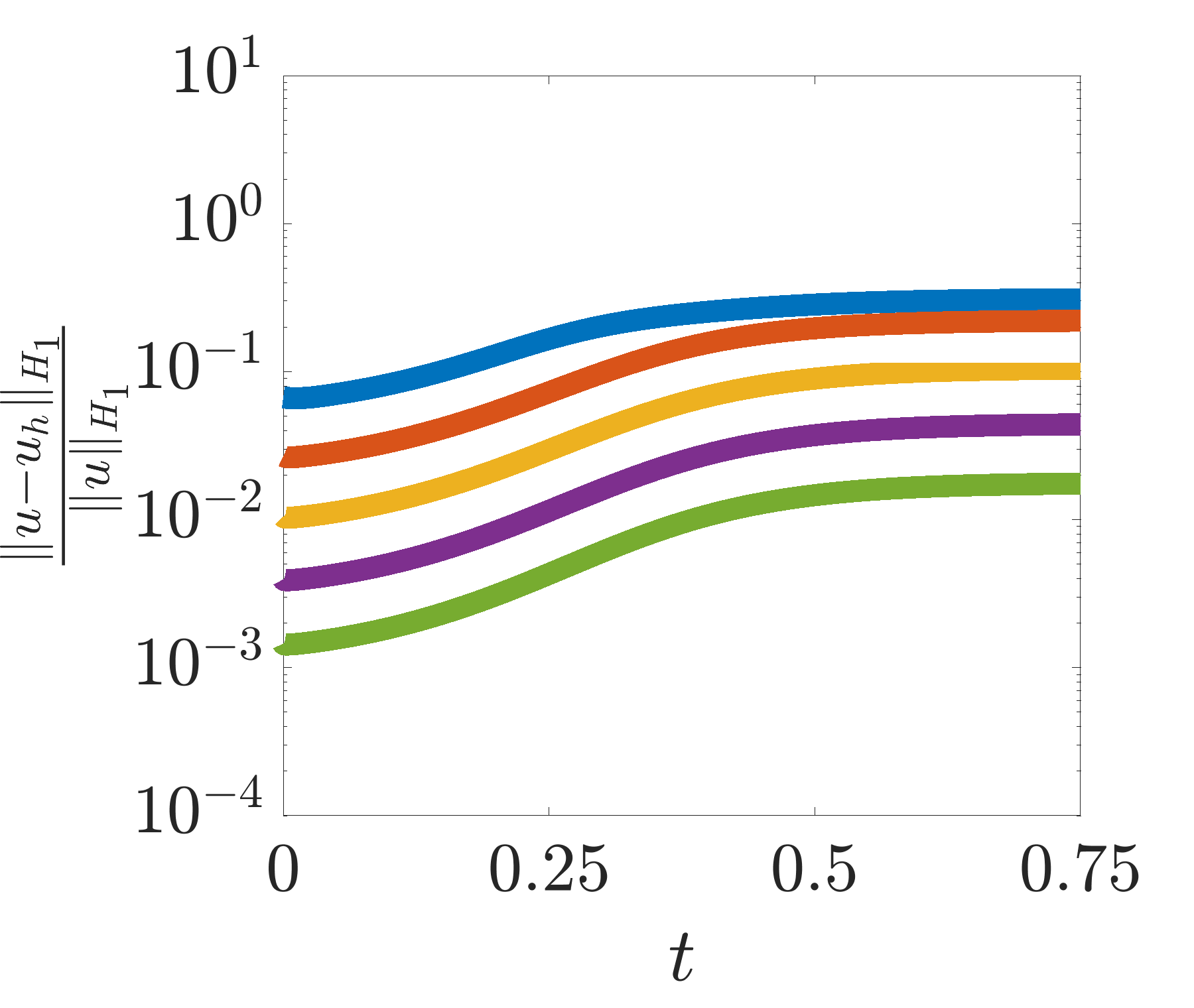}}
\subfigure[FEM; $\nu = \frac{1}{100}$]{\includegraphics[width=1.6in]{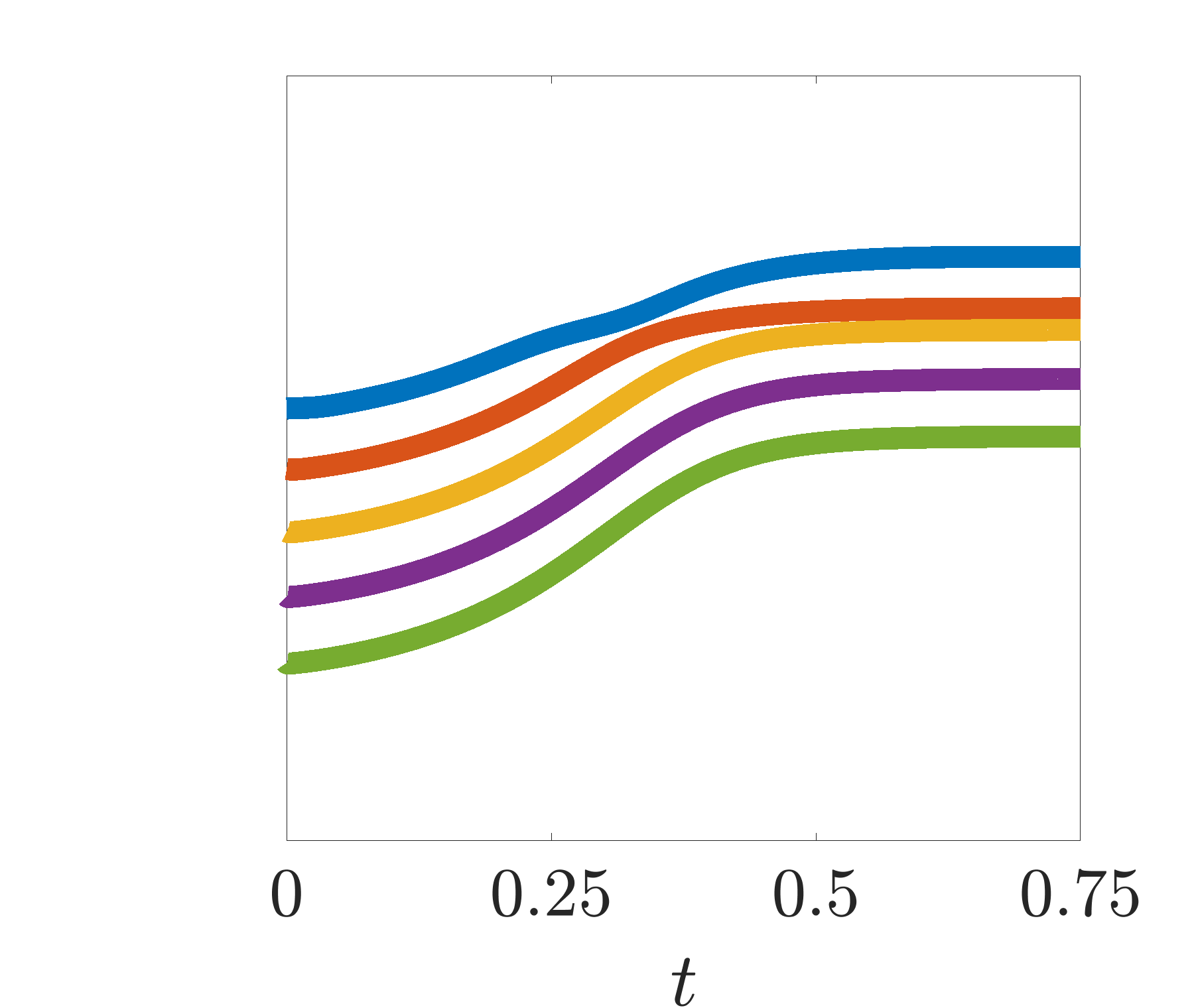}}
\subfigure[FEM; $\nu = \frac{1}{500}$]{\includegraphics[width=1.6in]{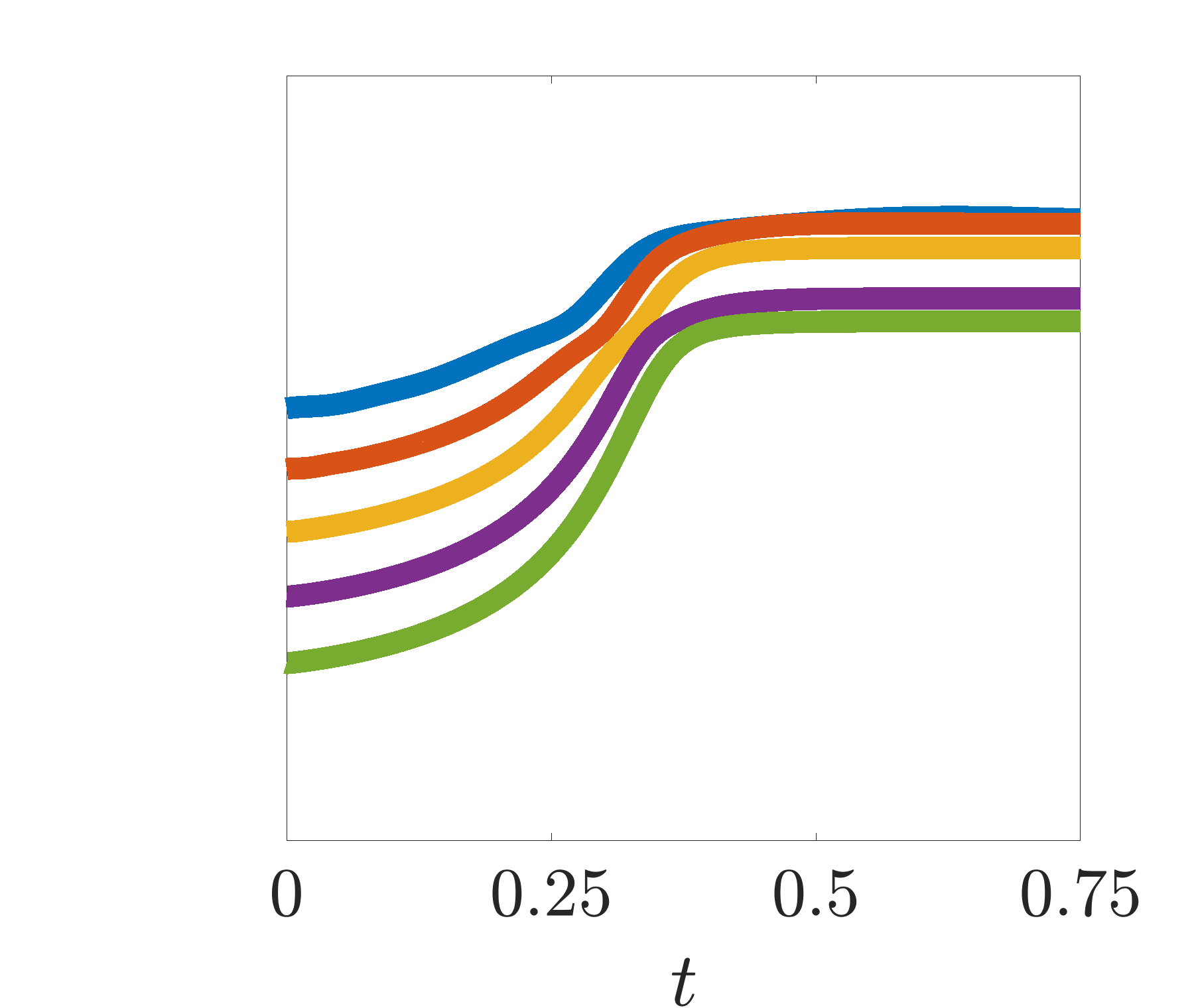}}
\subfigure[FEM; $\nu = \frac{1}{1000}$]{\includegraphics[width=1.6in]{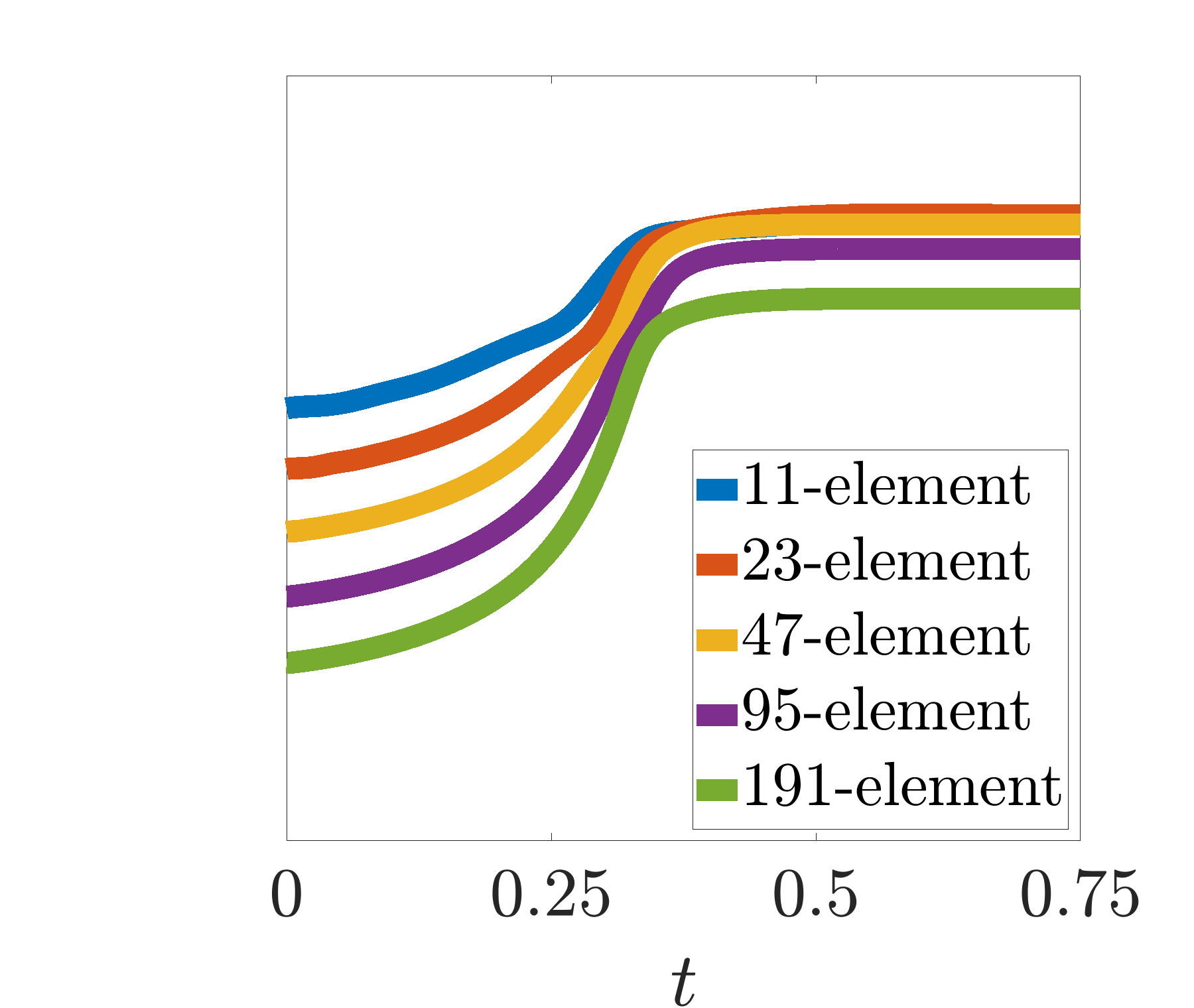}}

\subfigure[GFEM; $\nu = \frac{1}{50}$]{\includegraphics[width=1.6in]{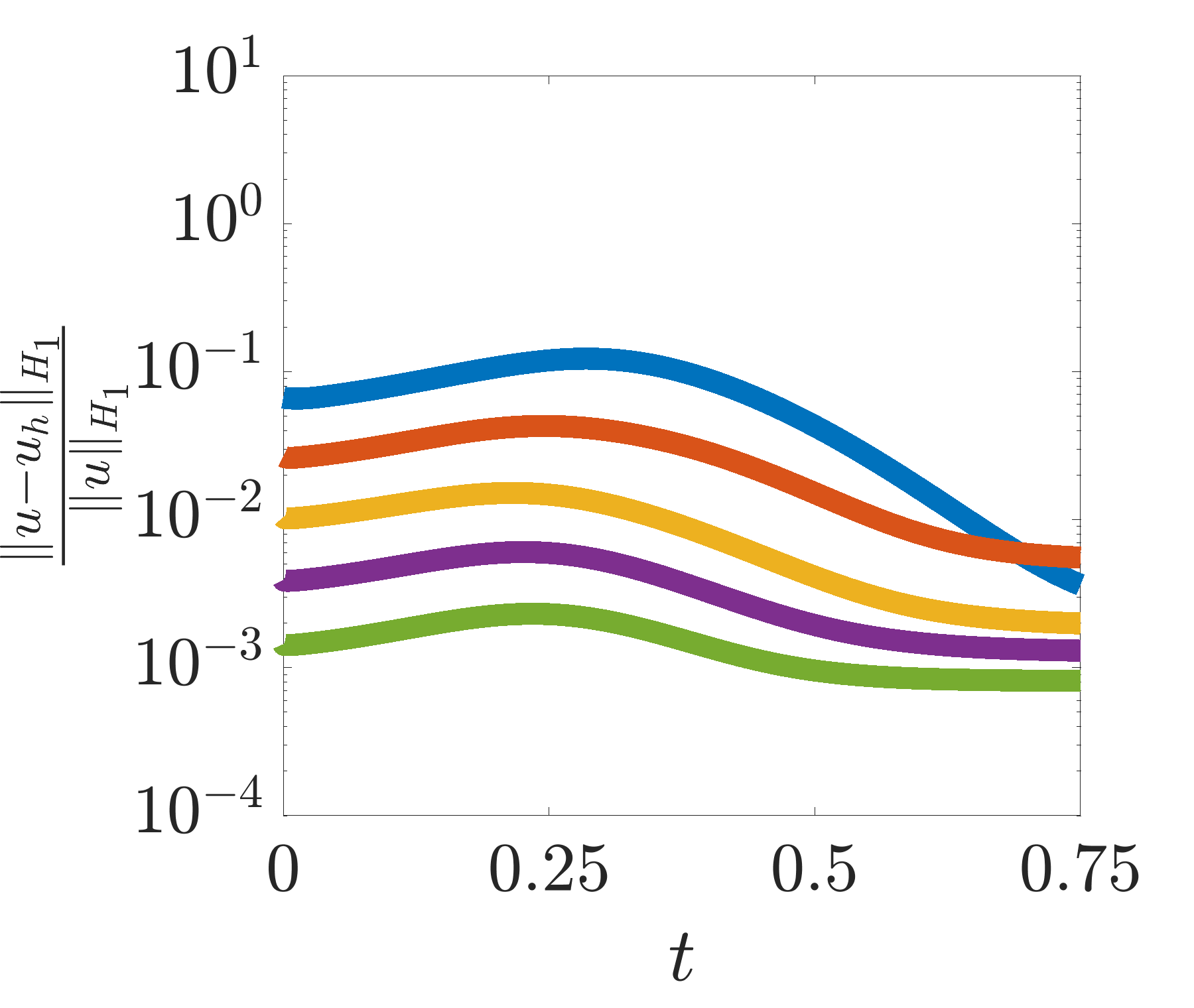}}
\subfigure[GFEM; $\nu = \frac{1}{100}$]{\includegraphics[width=1.6in]{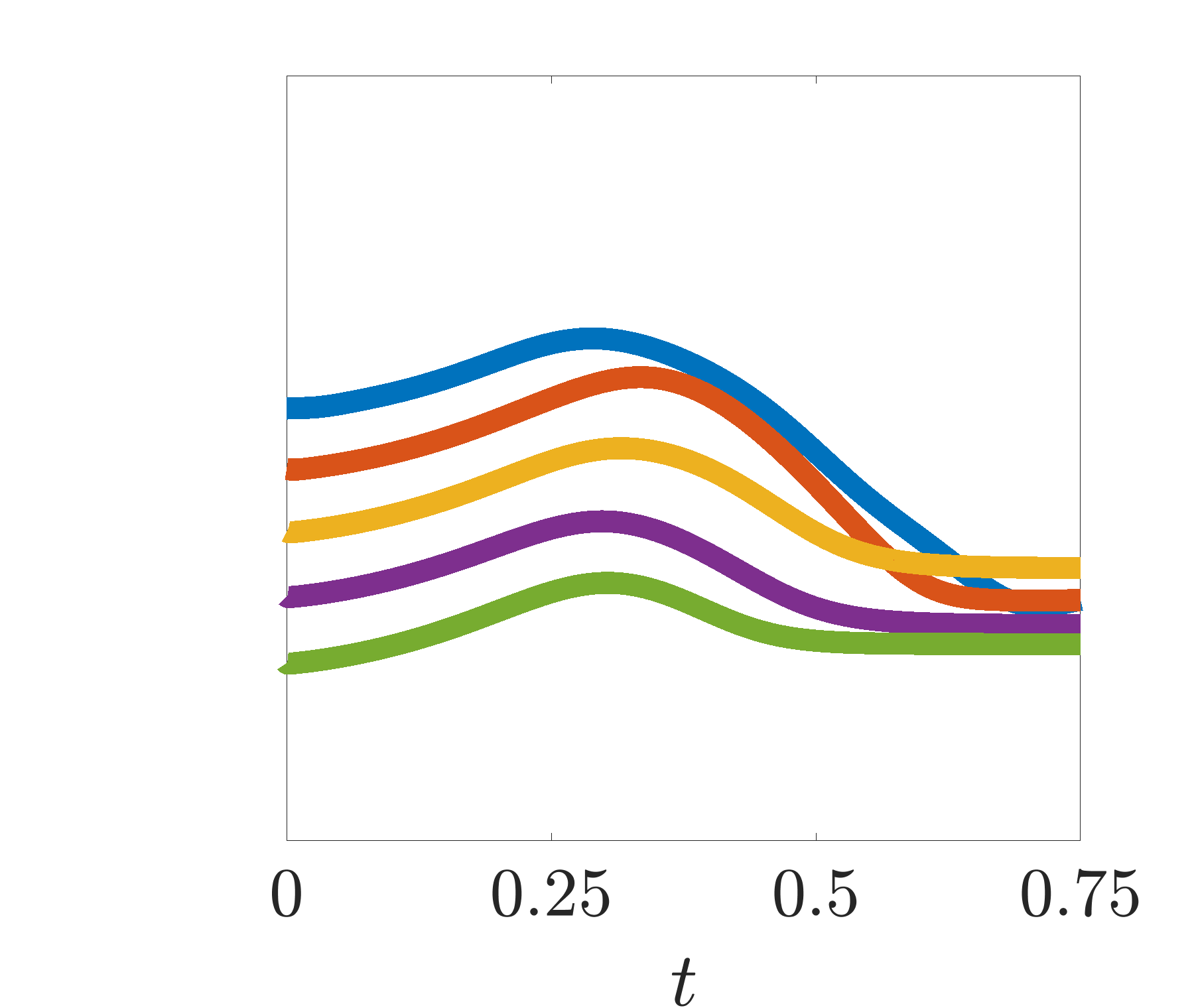}}
\subfigure[GFEM; $\nu = \frac{1}{500}$]{\includegraphics[width=1.6in]{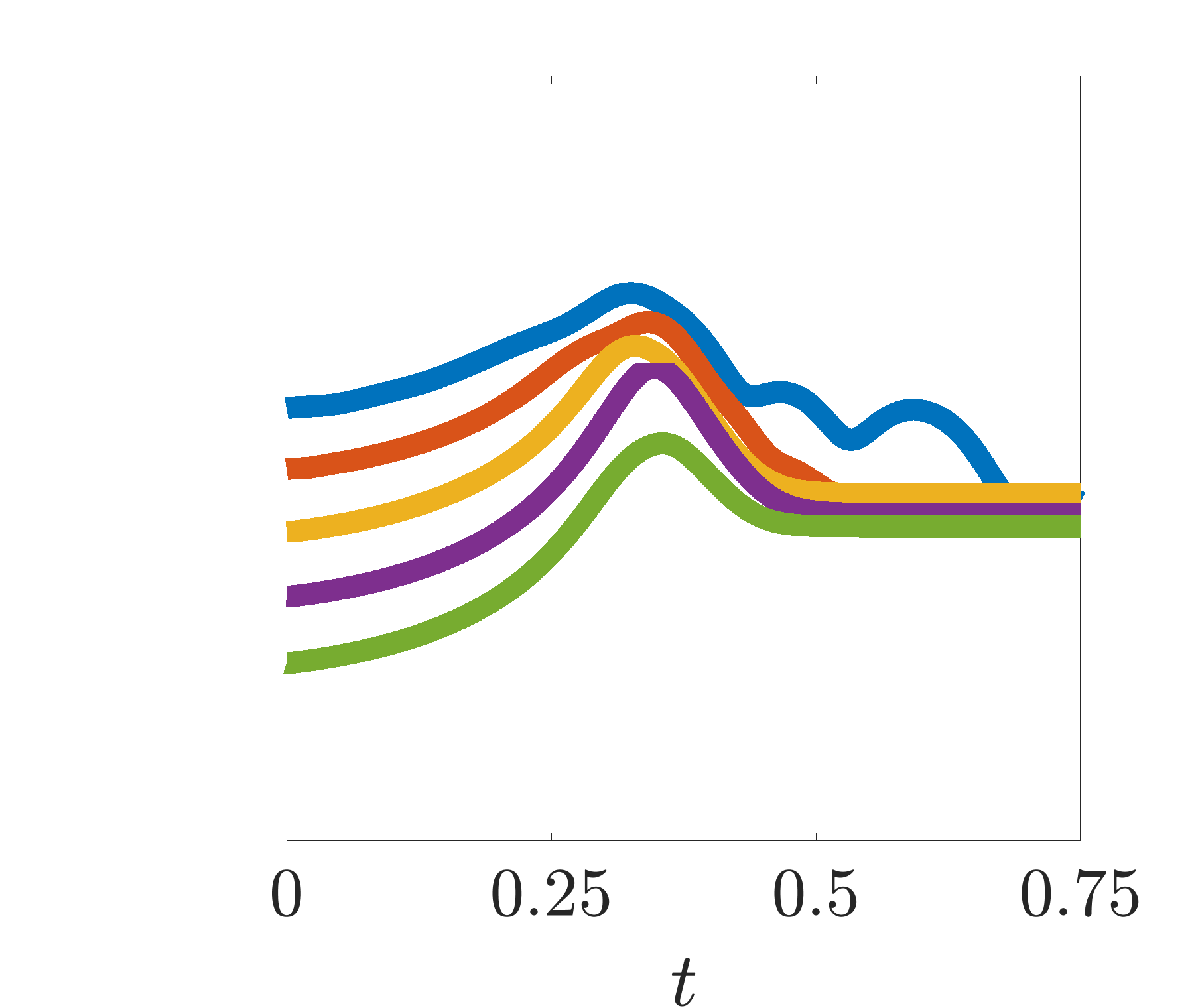}}
\subfigure[GFEM; $\nu = \frac{1}{1000}$]{\includegraphics[width=1.6in]{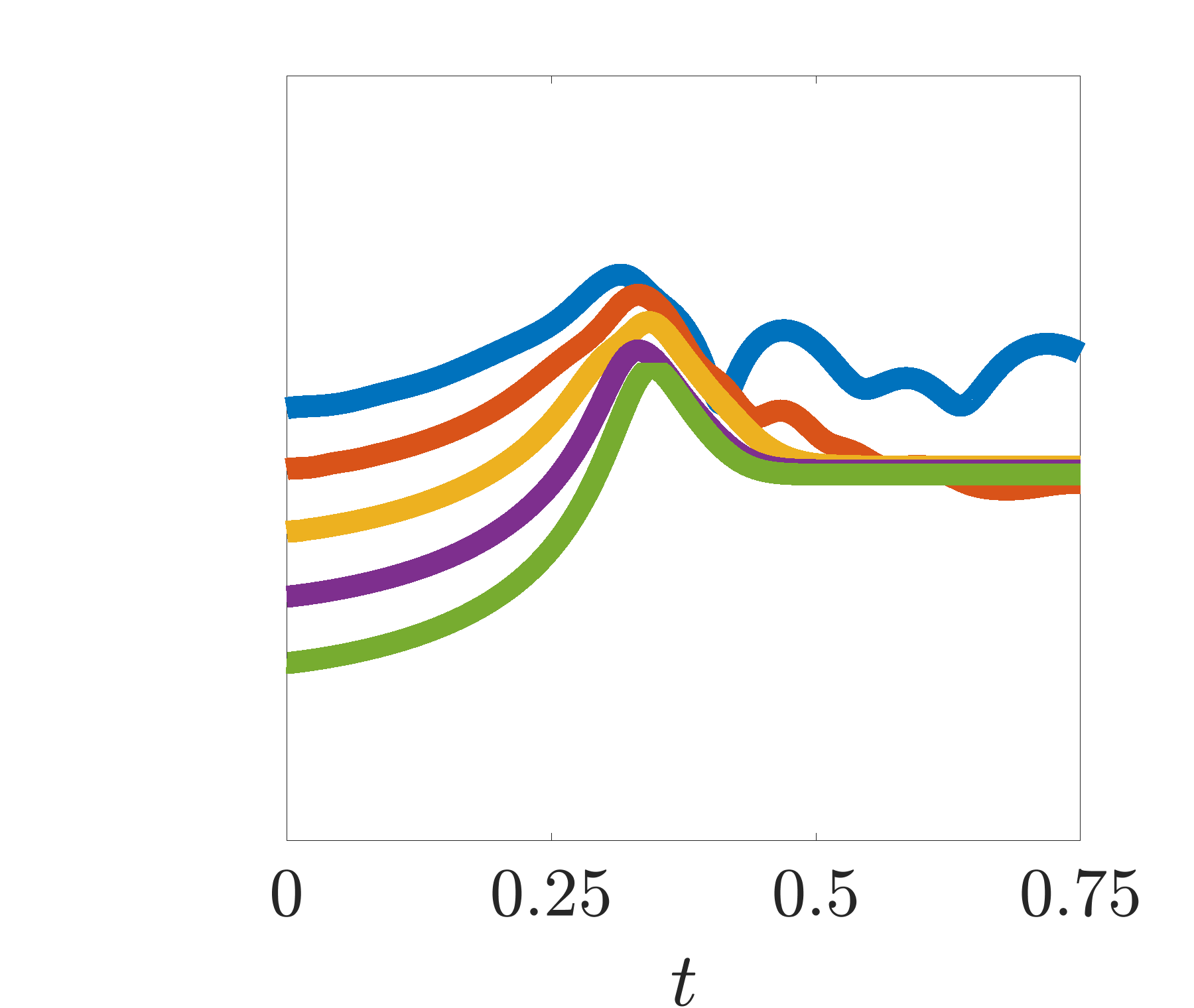}}
\end{subfigmatrix}
\caption{Relative $H_1$ norm versus time for
11-, 23-, 47-, 95-element, and 191-element $p = 1$ FEM and $p = 1$ + ss GFEM solutions over various kinematic viscosities for the shock problem}
\label{fig:Example2_H1vstime}
\end{center}
\end{figure}

\begin{figure}[ht!]
\begin{center}
\begin{subfigmatrix}{8}
\subfigure[12 DOF FEM; $\nu = \frac{1}{50}$]{\includegraphics[width=1.6in]{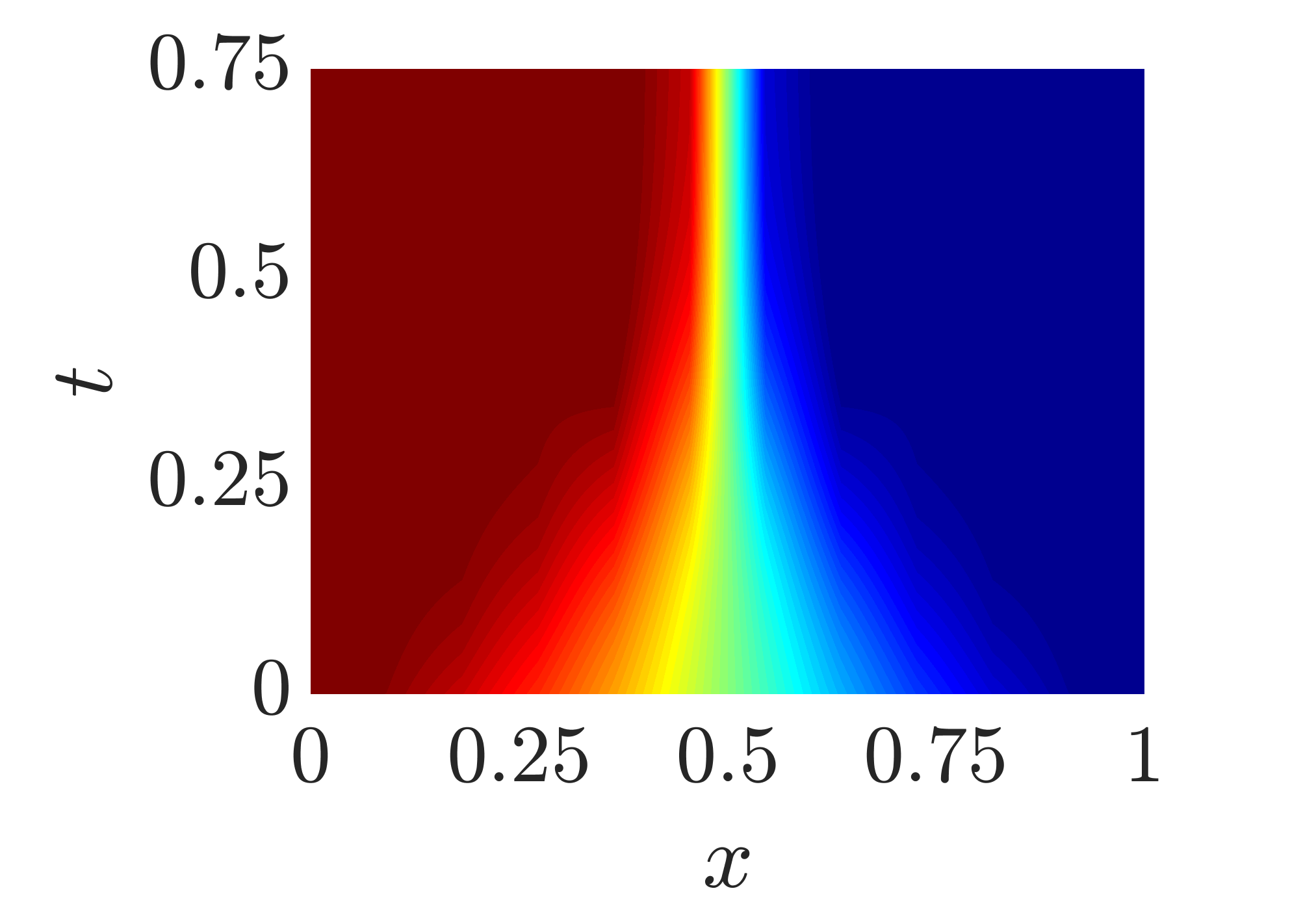}}
\subfigure[12 DOF FEM; $\nu = \frac{1}{100}$]{\includegraphics[width=1.6in]{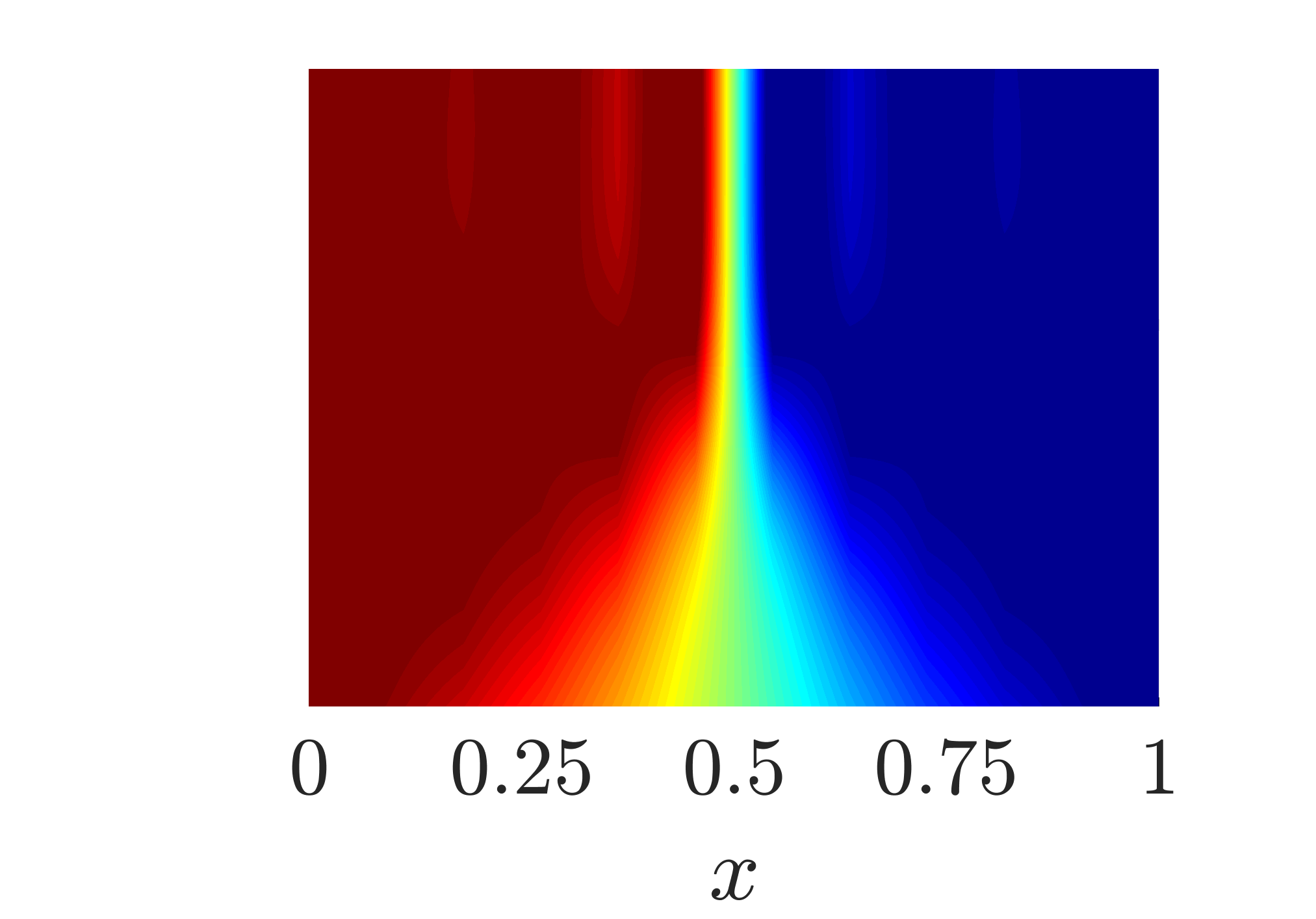}}
\subfigure[12 DOF FEM; $\nu = \frac{1}{500}$]{\includegraphics[width=1.6in]{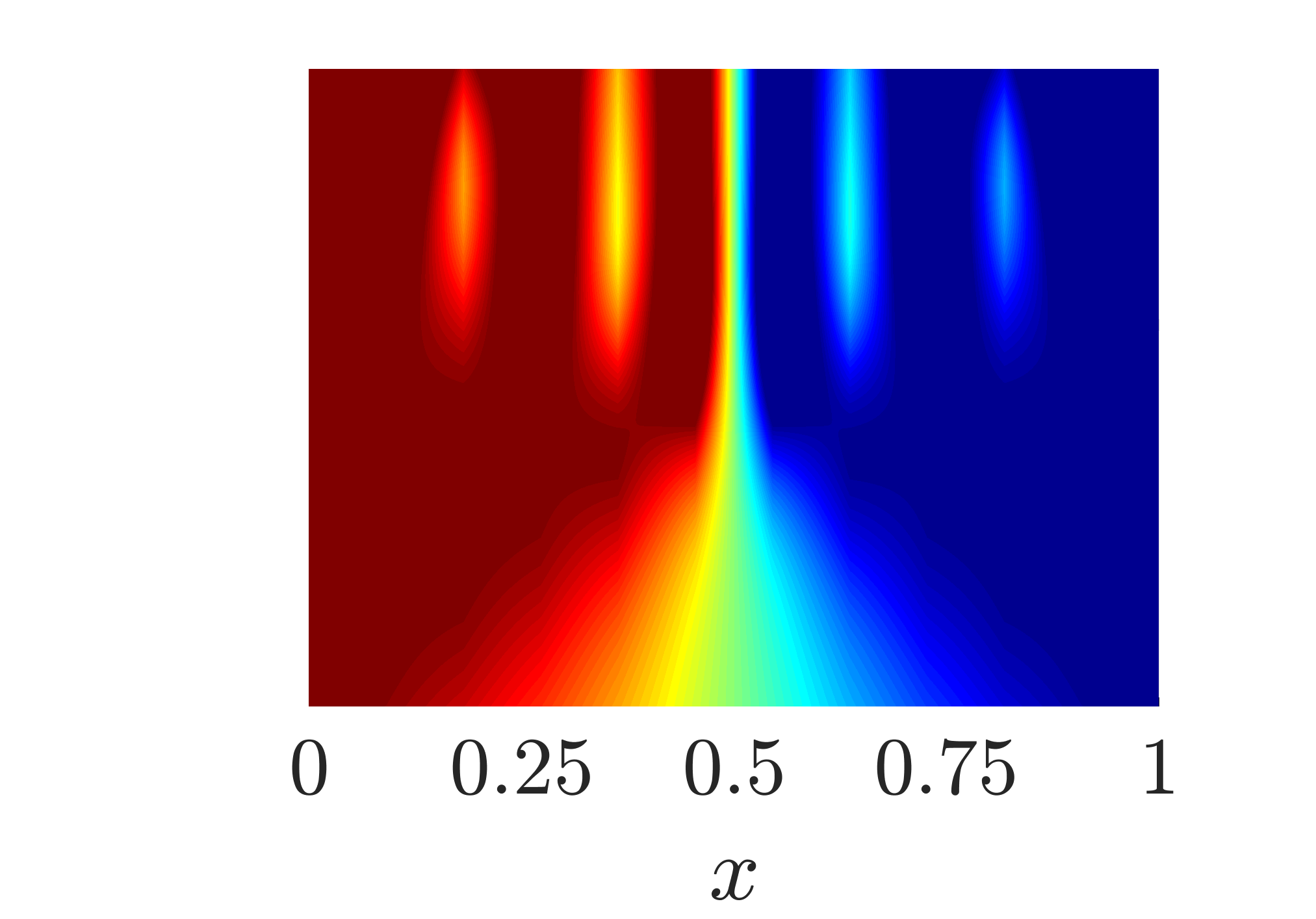}}
\subfigure[12 DOF FEM; $\nu = \frac{1}{1000}$]{\includegraphics[width=1.6in]{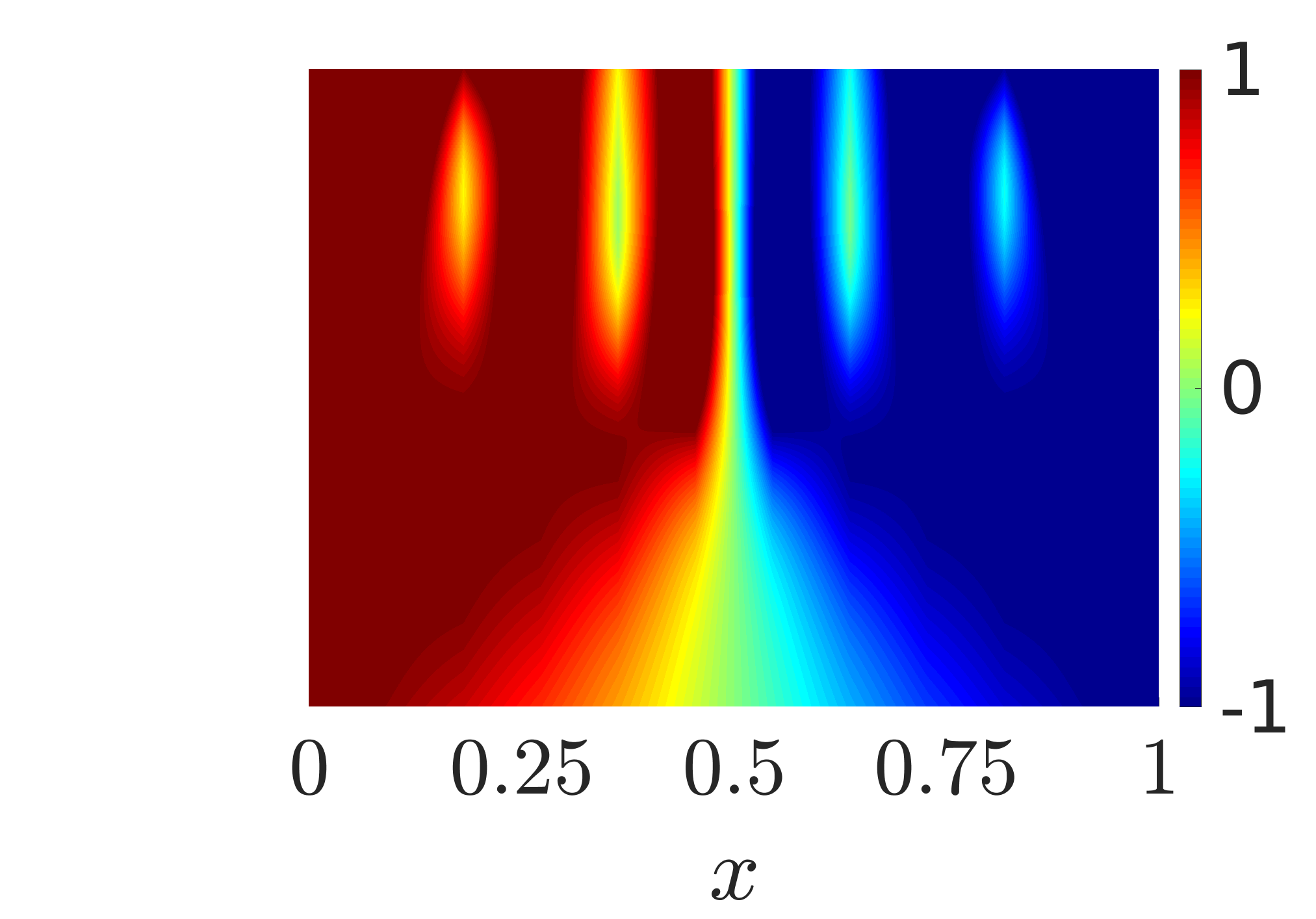}}
\subfigure[16 DOF GFEM; $\nu = \frac{1}{50}$]{\includegraphics[width=1.6in]{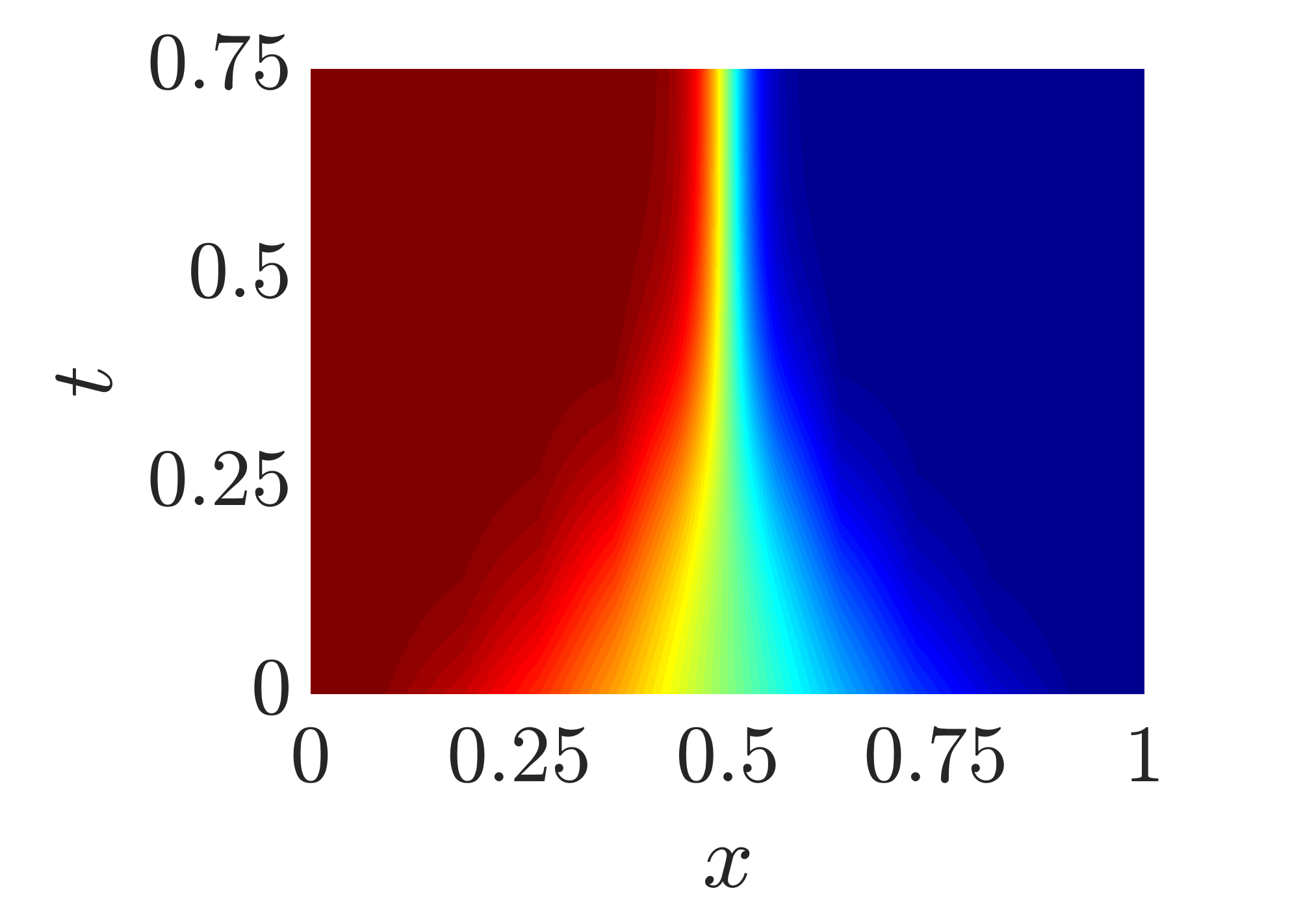}}
\subfigure[16 DOF GFEM; $\nu = \frac{1}{100}$]{\includegraphics[width=1.6in]{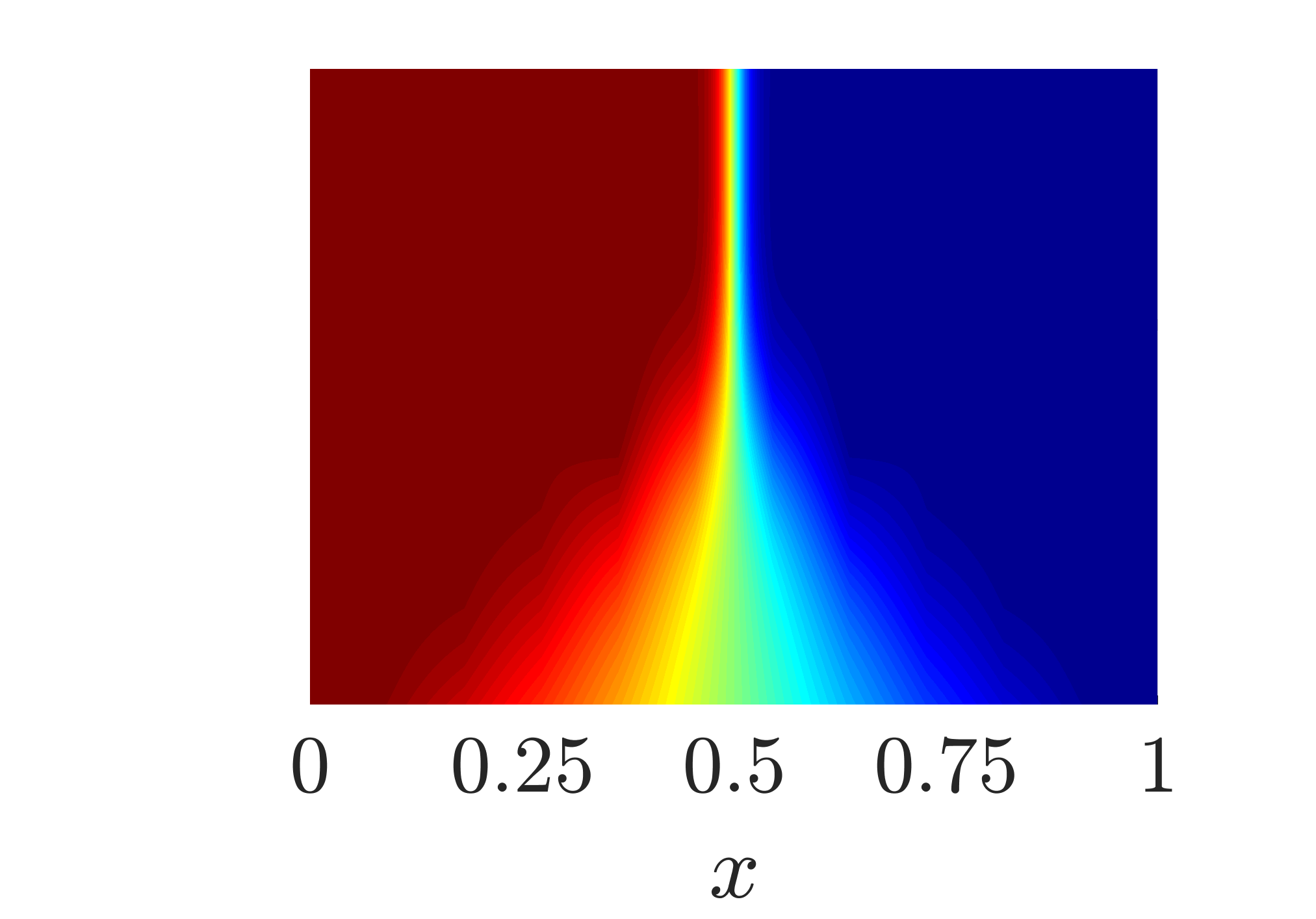}}
\subfigure[14 DOF GFEM; $\nu = \frac{1}{500}$]{\includegraphics[width=1.6in]{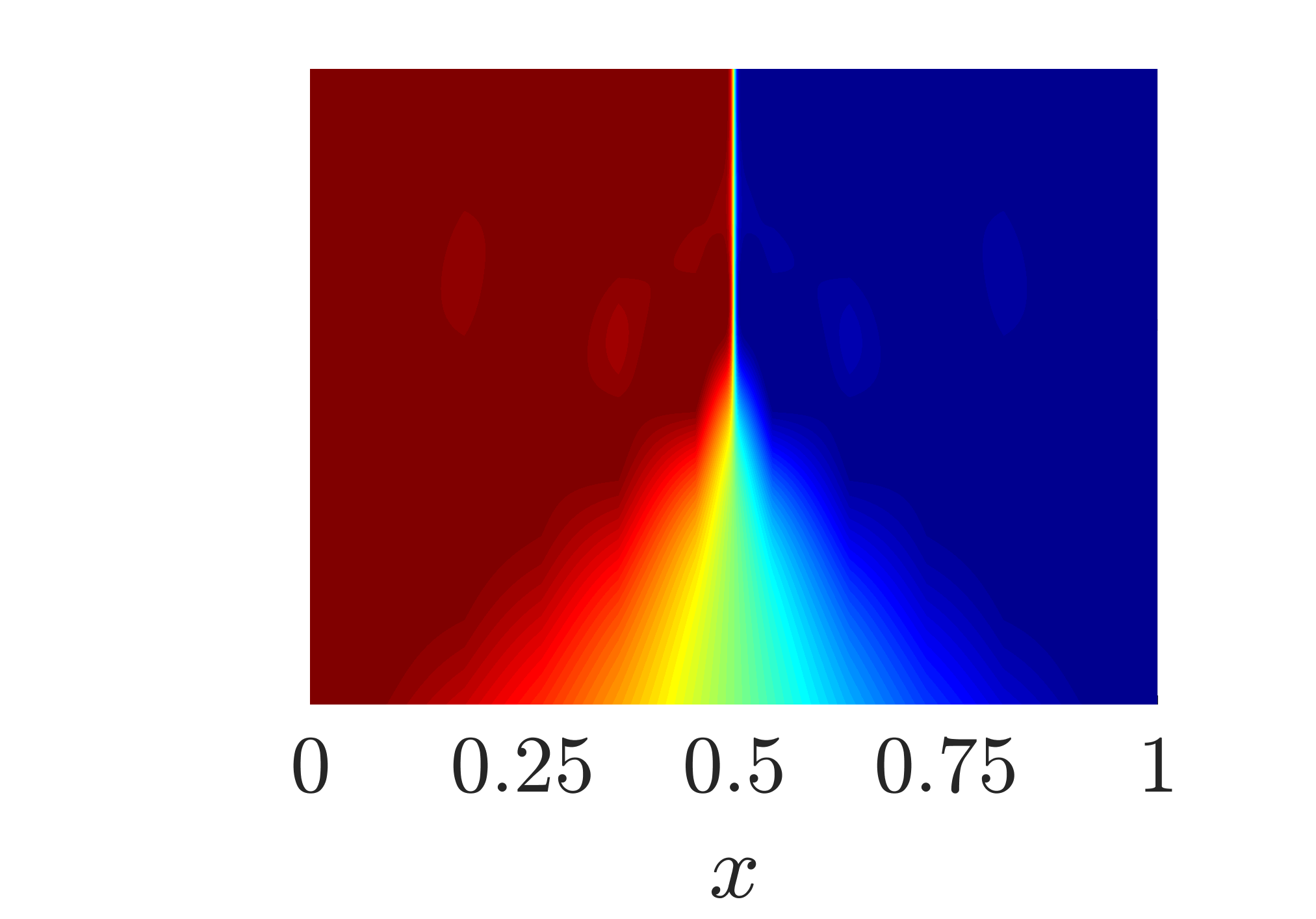}}
\subfigure[14 DOF GFEM; $\nu = \frac{1}{1000}$]{\includegraphics[width=1.6in]{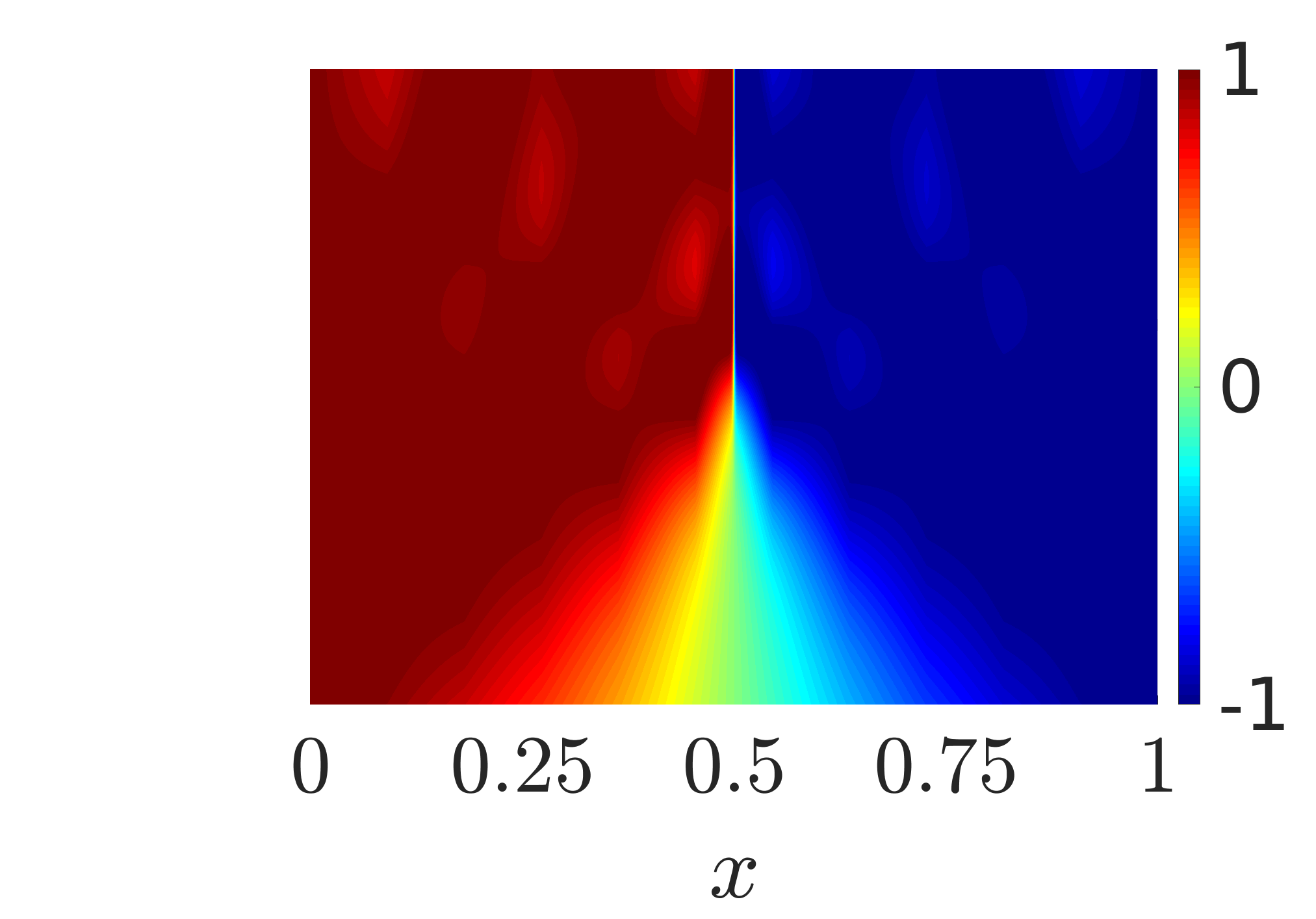}}
\end{subfigmatrix}
\caption{11-element $p = 1$ FEM and $p = 1$ + ss GFEM solution contours over various kinematic viscosities for the shock problem}
\label{fig:Example2_11element_contours}
\end{center}
\end{figure}

\begin{figure}[ht!]
\begin{center}
\begin{subfigmatrix}{8}
\subfigure[48 DOF FEM; $\nu = \frac{1}{50}$]{\includegraphics[width=1.6in]{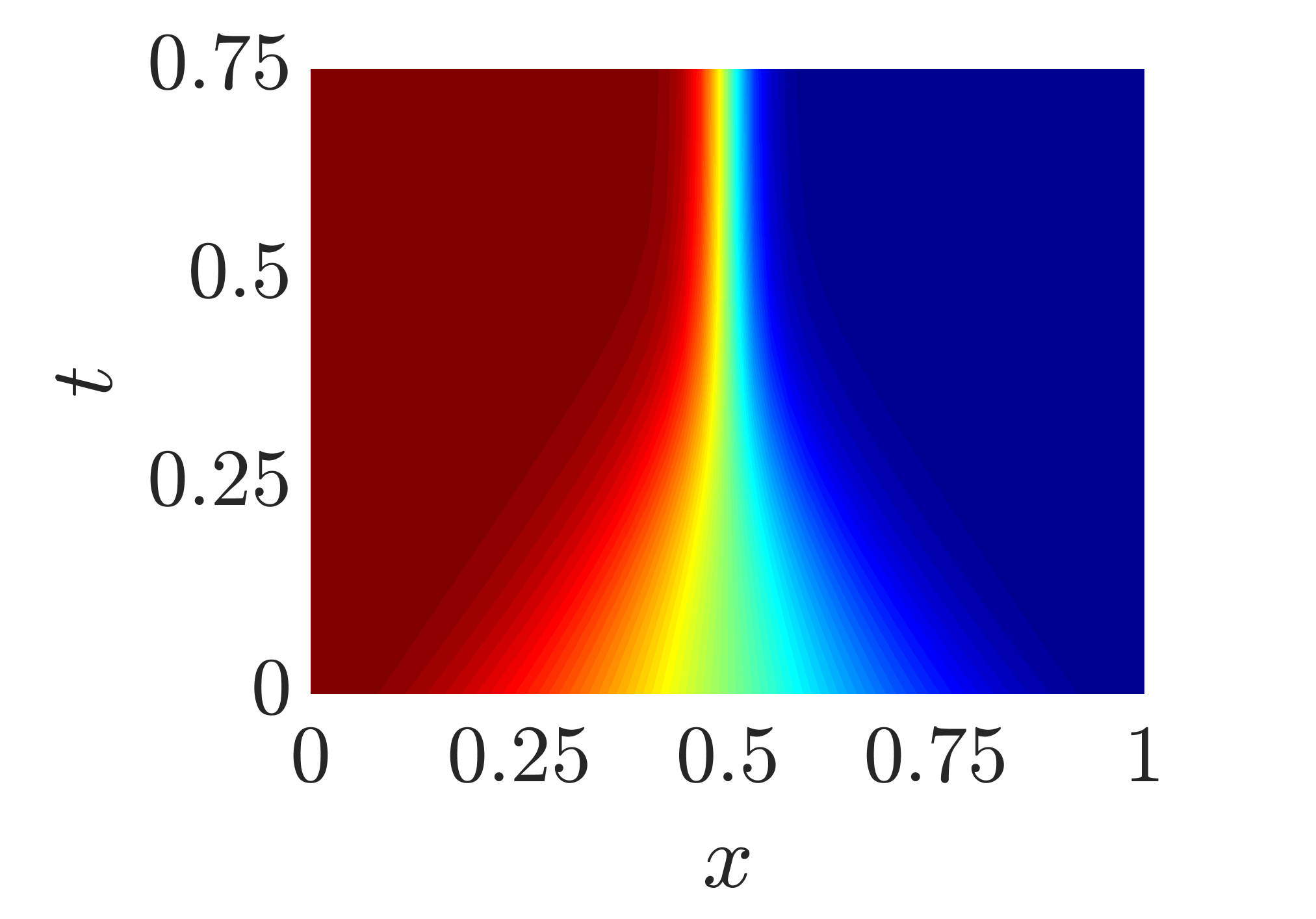}}
\subfigure[48 DOF FEM; $\nu = \frac{1}{100}$]{\includegraphics[width=1.6in]{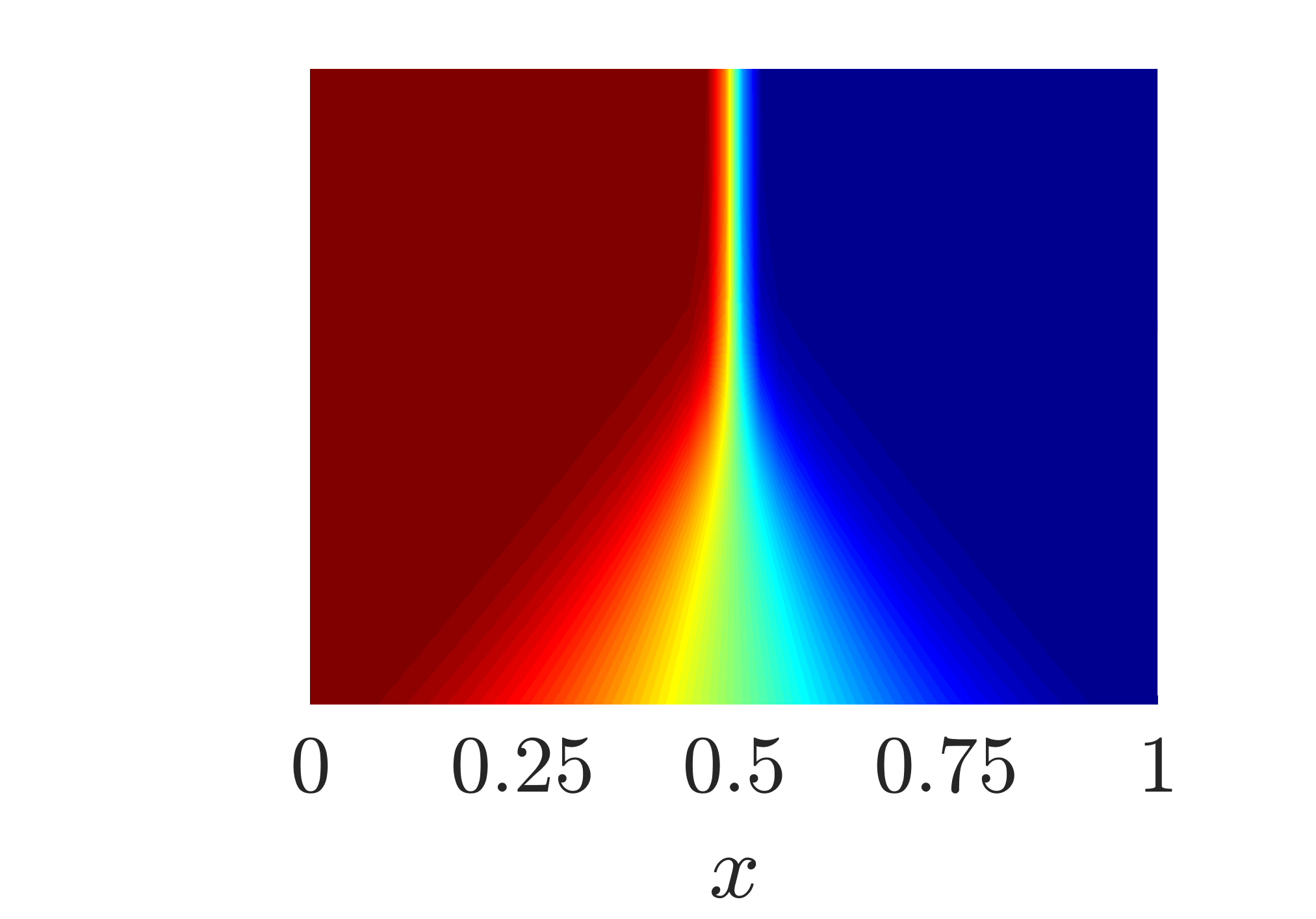}}
\subfigure[48 DOF FEM; $\nu = \frac{1}{500}$]{\includegraphics[width=1.6in]{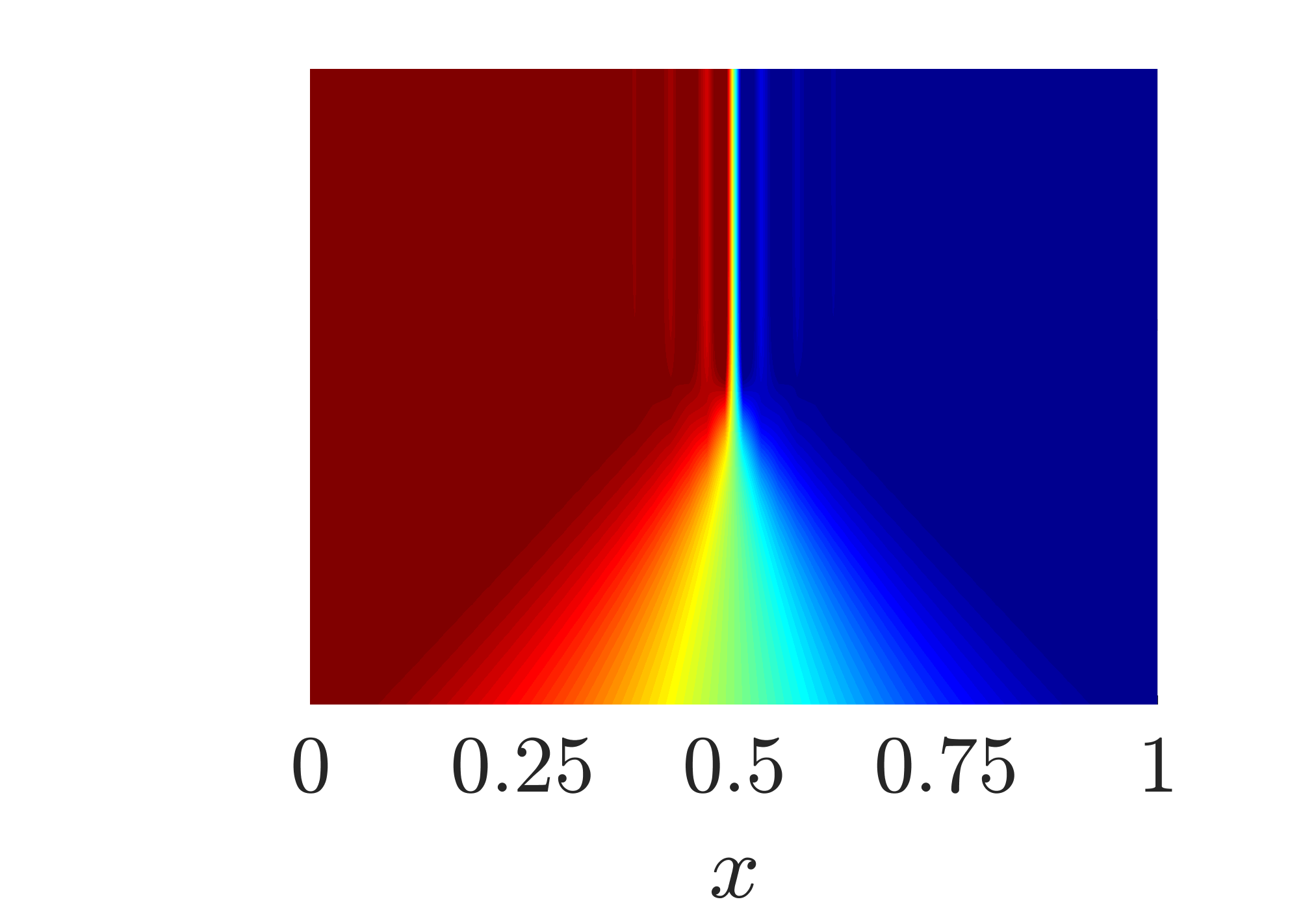}}
\subfigure[48 DOF FEM; $\nu = \frac{1}{1000}$]{\includegraphics[width=1.6in]{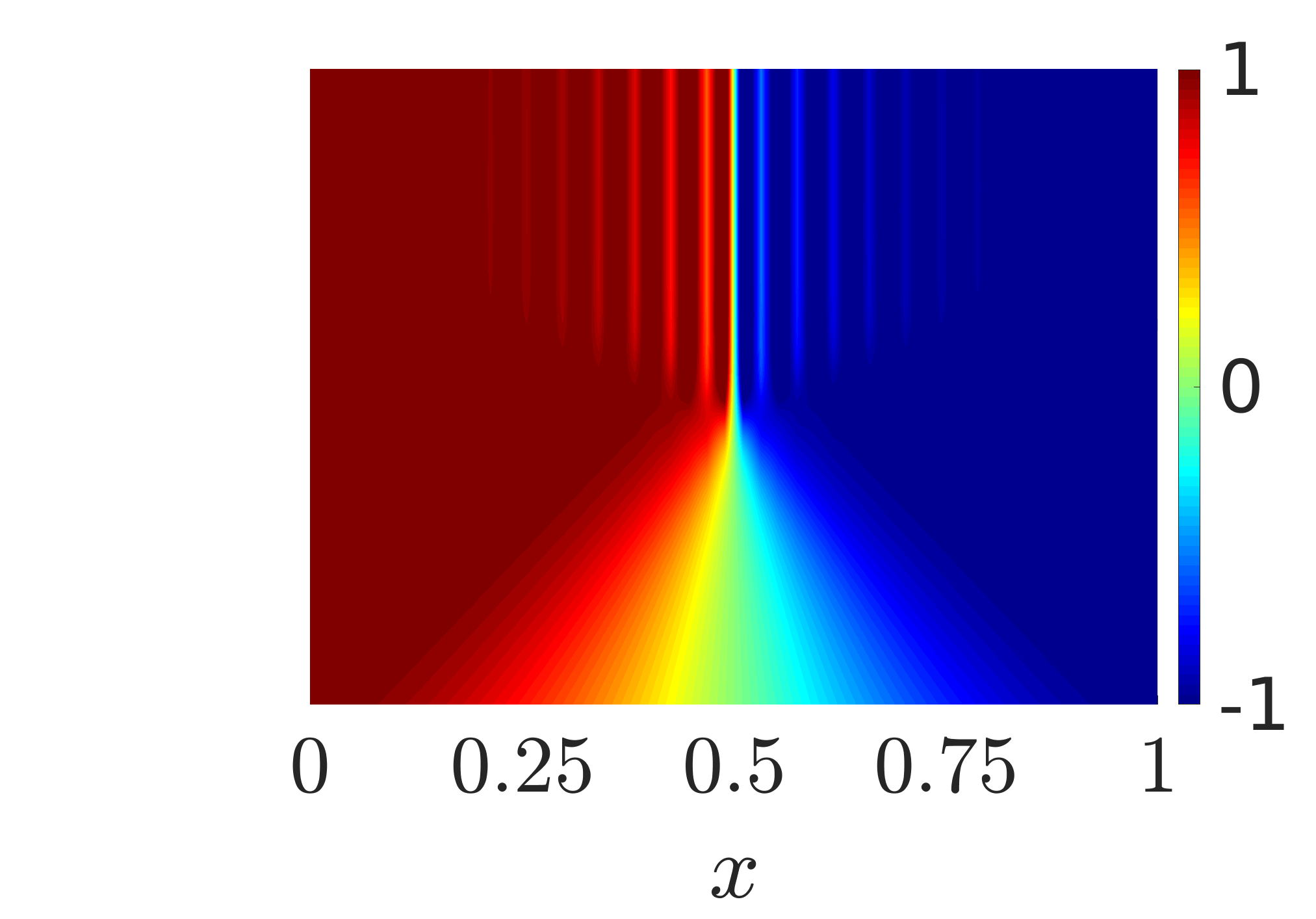}}
\subfigure[60 DOF GFEM; $\nu = \frac{1}{50}$]{\includegraphics[width=1.6in]{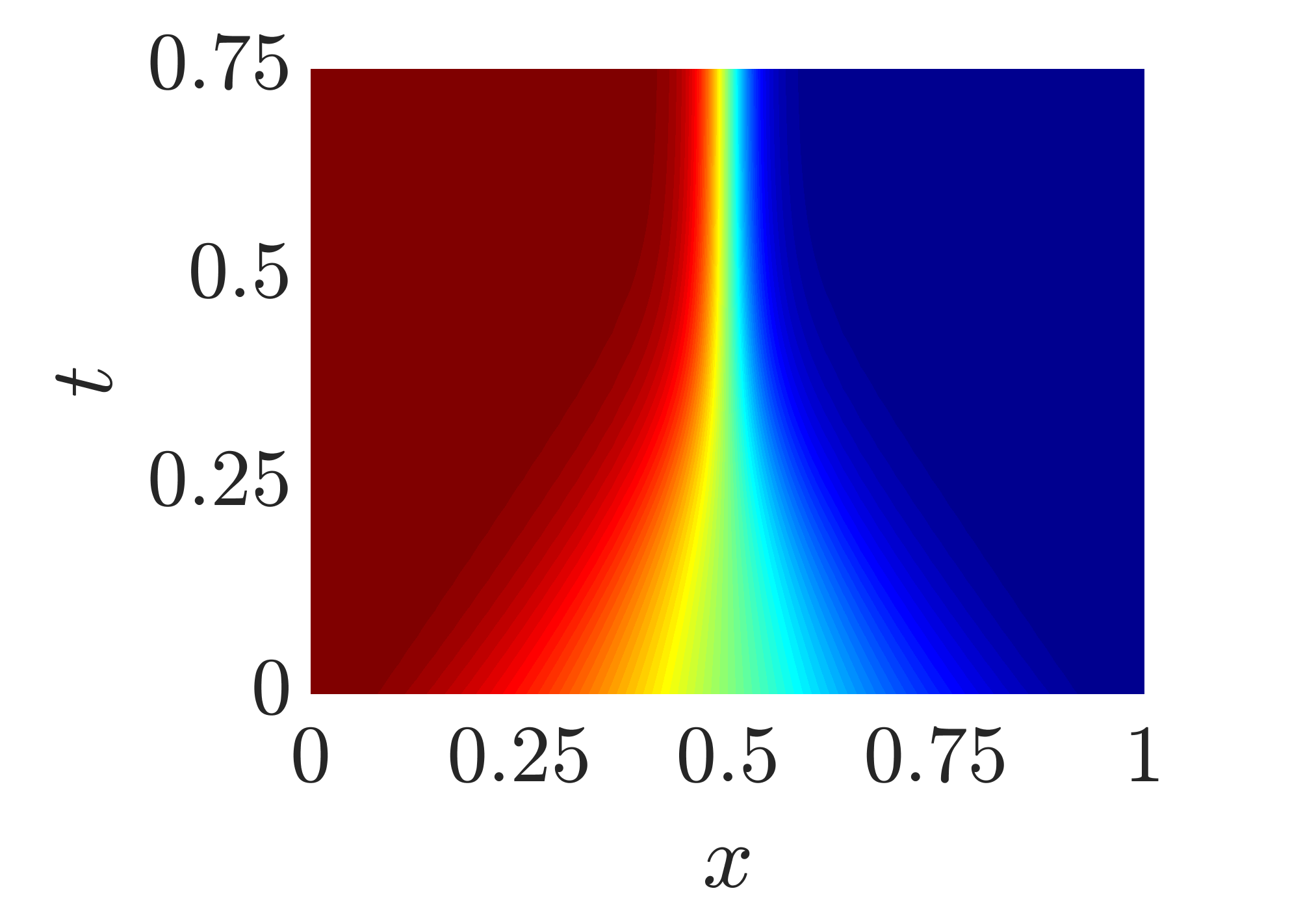}}
\subfigure[54 DOF GFEM; $\nu = \frac{1}{100}$]{\includegraphics[width=1.6in]{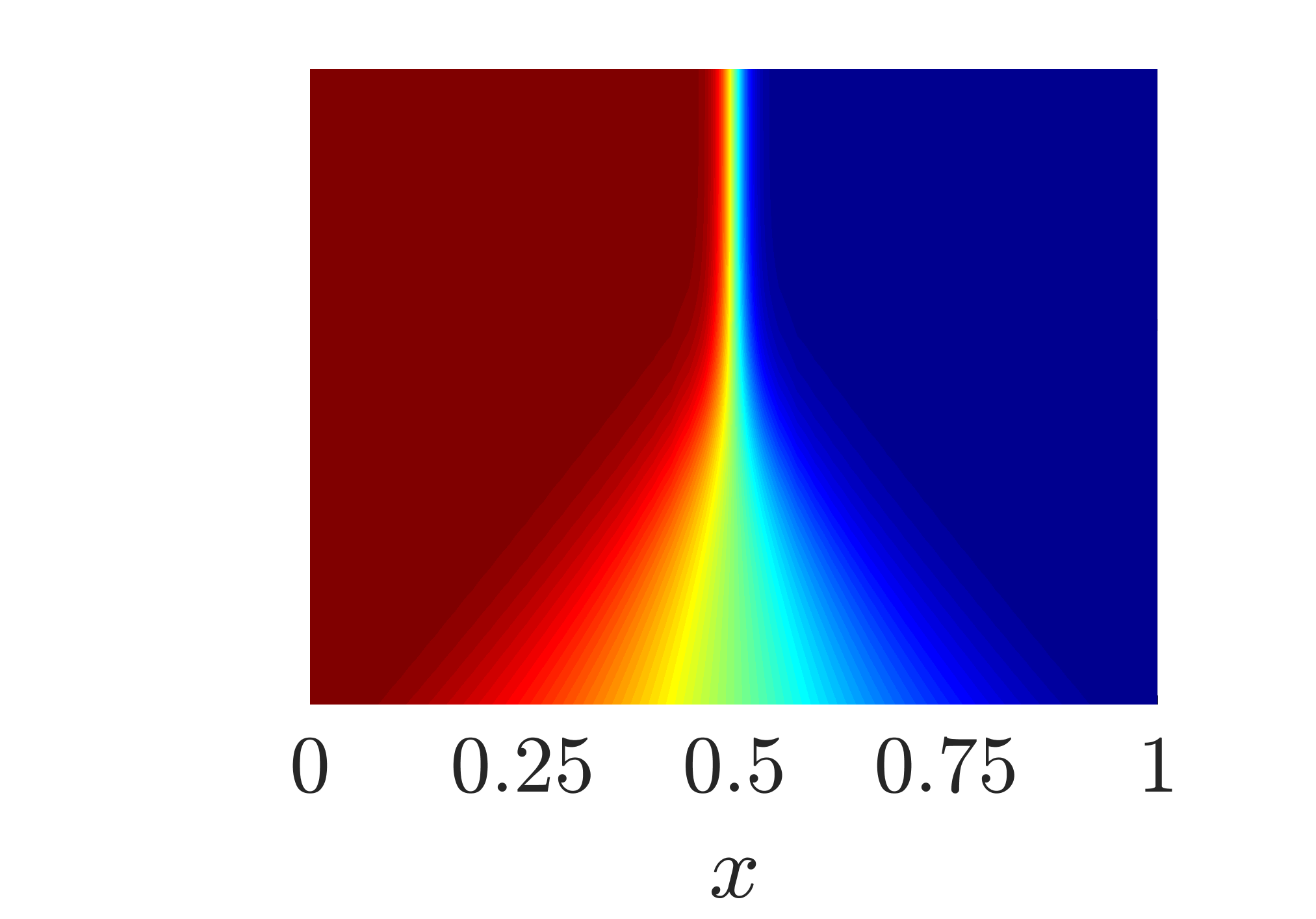}}
\subfigure[50 DOF GFEM; $\nu = \frac{1}{500}$]{\includegraphics[width=1.6in]{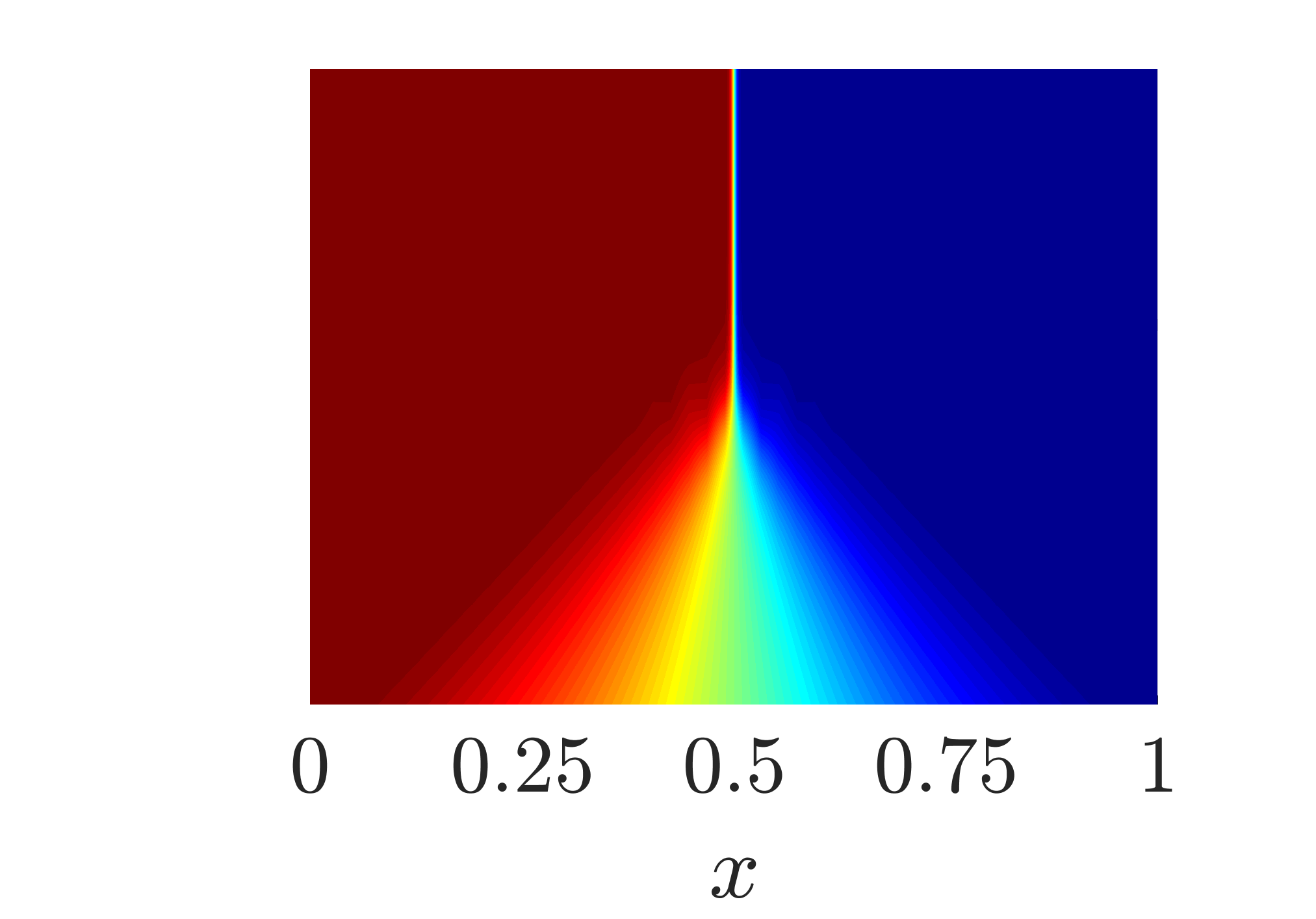}}
\subfigure[50 DOF GFEM; $\nu = \frac{1}{1000}$]{\includegraphics[width=1.6in]{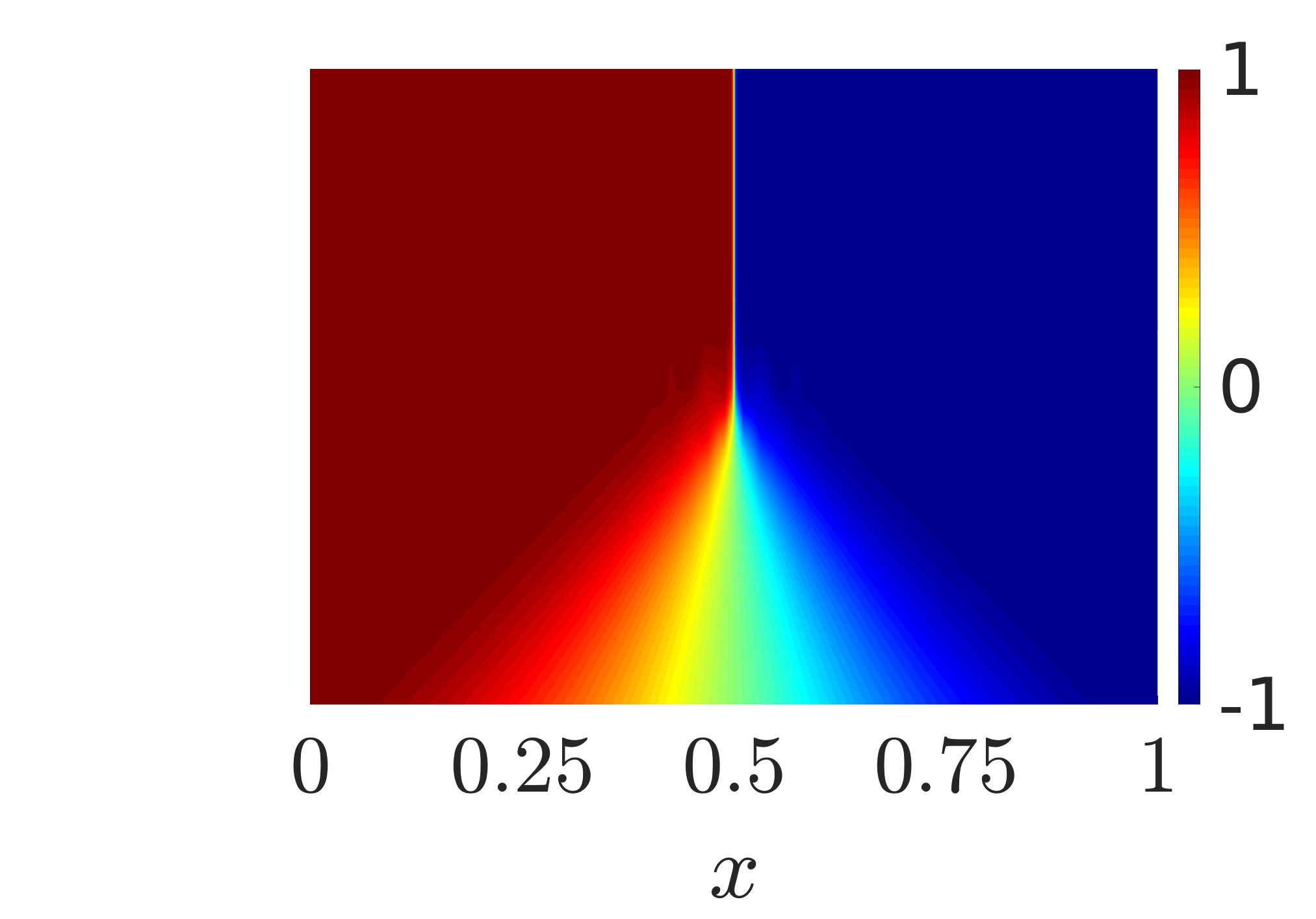}}
\end{subfigmatrix}
\caption{47-element $p = 1$ FEM and $p = 1$ + ss GFEM solution contours over various kinematic viscosities for the shock problem}
\label{fig:Example2_47element_contours}
\end{center}
\end{figure}

Plots of the relative $L_2$ and $H_1$ norm versus time are shown in Figs. \ref{fig:Example2_L2vstime} and \ref{fig:Example2_H1vstime}, respectively. The $p = 1$ FEM and $p = 1$ + ss GFEM solutions return similar error levels up until $t \approx \frac{1}{\pi}$ where shock gradients increase. This is expected since the steady state solution is not closely correlated with the initial transient solution. Around $t = \frac{1}{\pi}$, error levels rise in the $p = 1$ FEM solutions, with the rise increasing as $\nu$ becomes smaller. At $t = 0.75$, 95-element $p = 1$ FEM solutions return errors in the relative $L_2$ are 0.10\%, 0.27\%, 0.86\%, and 3.3\% for $\nu = \Big[\frac{1}{50}, \frac{1}{100}, \frac{1}{500}, \frac{1}{1000} \Big]$, respectively. Similarly in the relative $H_1$ norm the errors are 4.4\%, 10.4\%, 35.0\%, and 73.5\%, respectively. Use of the steady state as an enrichment in GFEM results in a significant reduction of error in both the $L_2$ and $H_1$ norm. At $t = 0.75$ the 95-element GFEM solutions enriched with the steady state solution have errors in the relative $L_2$ norm of 0.0027\%, 0.0047\%, 0.0041\%, and 0.0032\% for $\nu = \Big[\frac{1}{50}, \frac{1}{100}, \frac{1}{500}, \frac{1}{1000} \Big]$, respectively. Similarly in the relative $H_1$ norm the errors are 0.13\%, 0.26\%, 1.36\%, and 2.79\%, respectively. Although the $p = 1$ + ss GFEM solutions significantly reduce relative errors with respect to $p = 1$ FEM, around $t = \frac{1}{\pi}$ errors peak to high values in both the integral norms. The peak value of the 95-element GFEM solution error in the $H_1$ norm for $\nu = \Big[\frac{1}{50}, \frac{1}{100}, \frac{1}{500}, \frac{1}{1000} \Big]$ are 0.62\%, 1.3\%, 12.9\%, and 17.0\%, respectively. These large errors arising in the GFEM solutions while the shock is forming is due to the linear interpolation and the steady state enrichment being ill-suited at capturing the intermediate solution features, similar to the results of example 1 when $\nu = 0$. Visually this is explained in Fig. \ref{fig:Example2_11element_contours} which shows 11-element $p = 1$ FEM and $p = 1$ + ss GFEM solution contours. As expected for $p = 1$ FEM solutions, oscillations arise in the numerical solutions for small $\nu = \Bigg[\frac{1}{100}, \frac{1}{500}, \frac{1}{1000}\Bigg]$. For the $p = 1$ + ss GFEM solutions, muted oscillations are observed for $\nu = \Bigg[\frac{1}{500}, \frac{1}{1000} \Bigg]$, which arise during the formation of the shock around $t = \frac{1}{\pi}$ and which propagate with time. With grid refinement as shown in Fig. \ref{fig:Example2_47element_contours}, which provides 47-element $p = 1$ FEM and $p = 1$ + ss GFEM solution contours, oscillations visually improve. 

\begin{figure}[ht!]
\begin{center}
\begin{subfigmatrix}{6}
\subfigure[$t = 0$]{\includegraphics[width=2.1in]{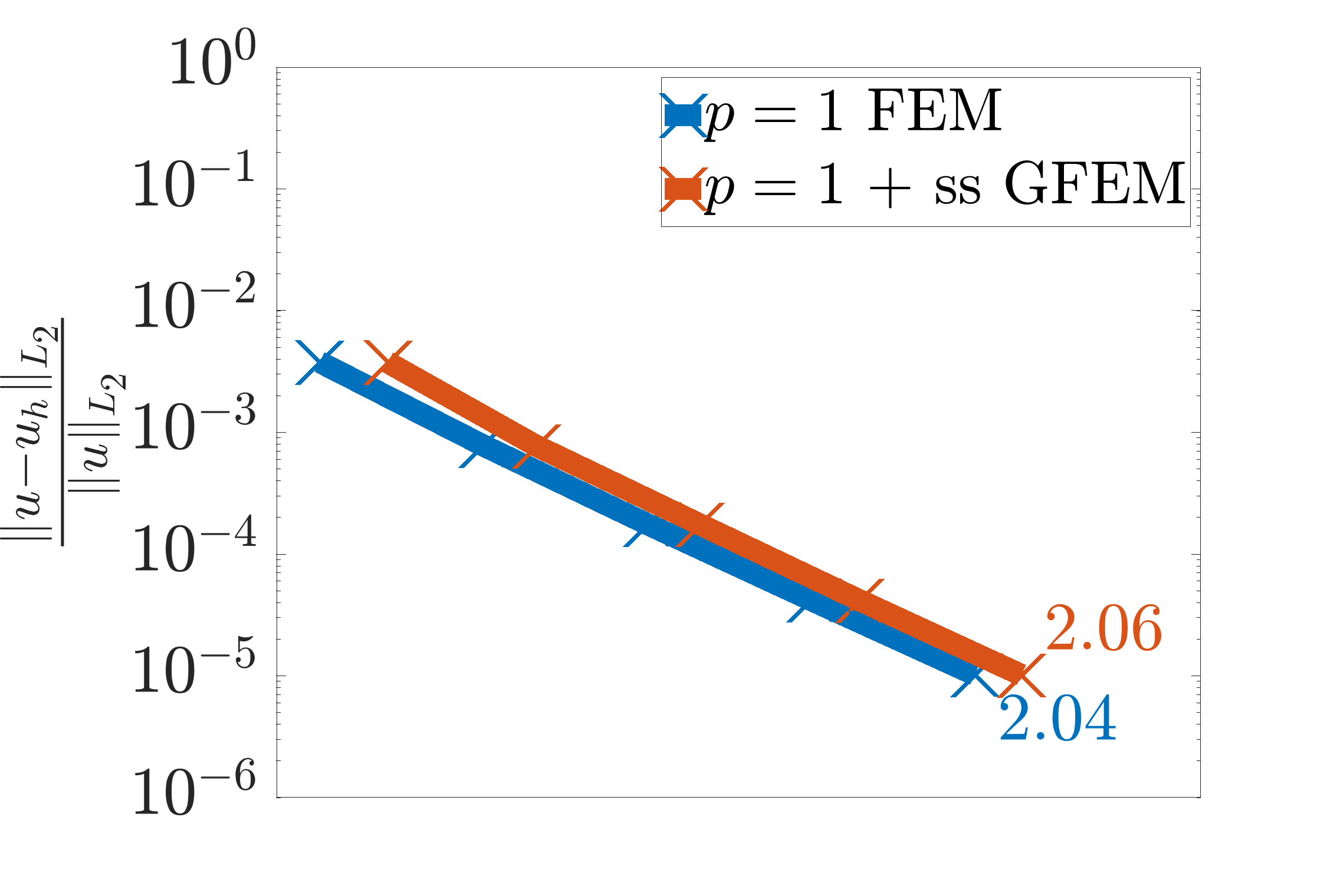}}
\subfigure[$t = 0.25$]{\includegraphics[width=2.1in]{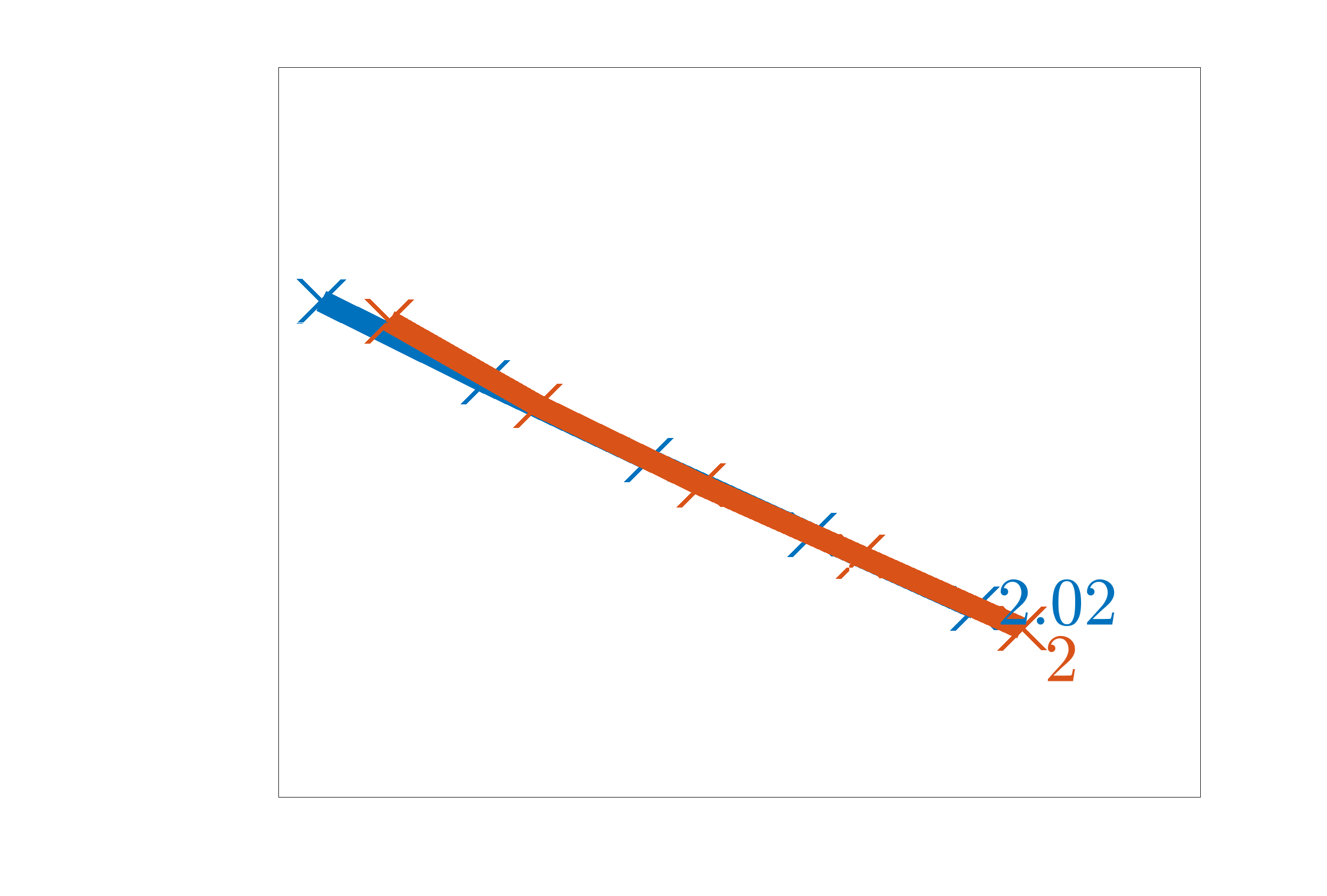}}
\subfigure[$t = 0.3$]{\includegraphics[width=2.1in]{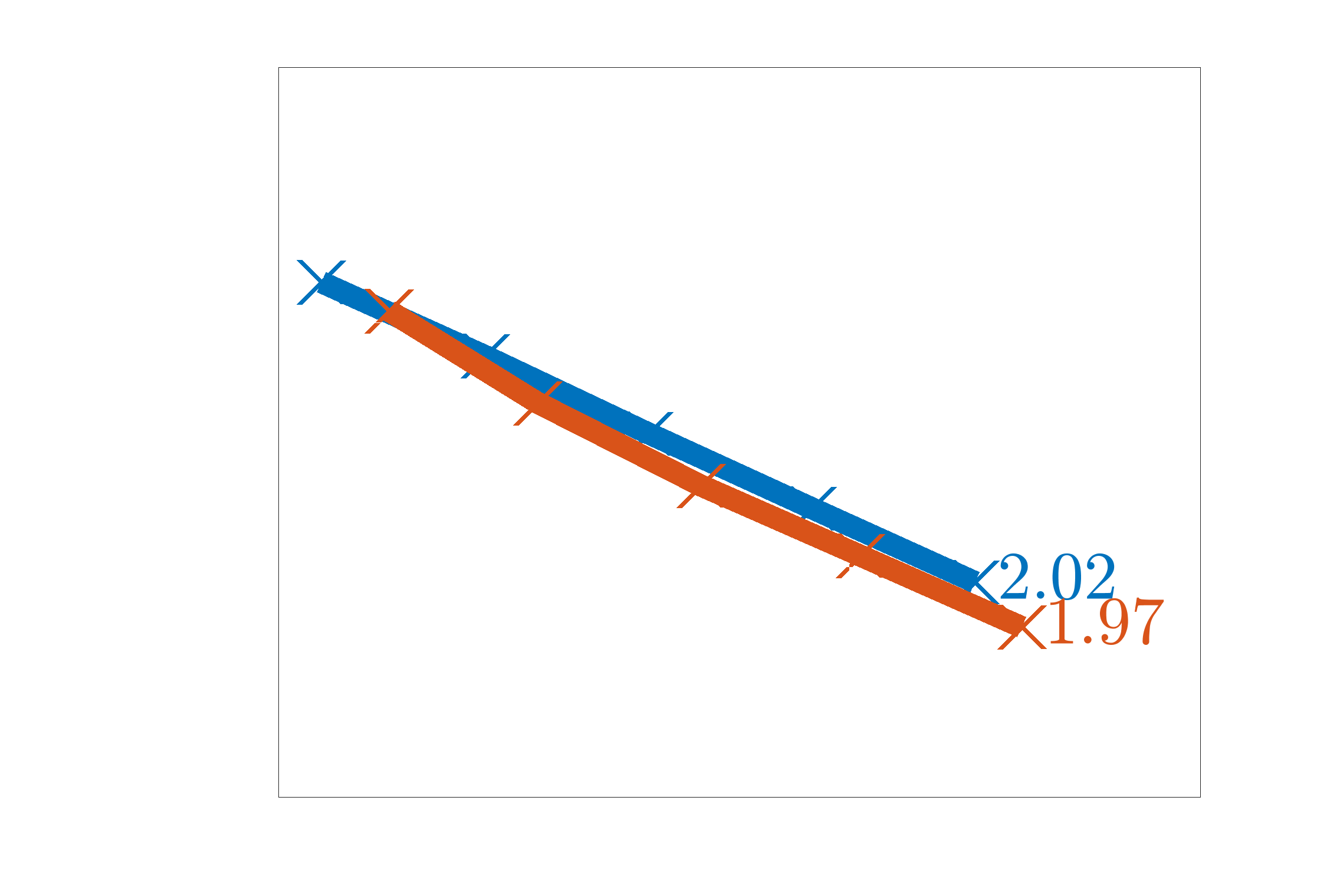}}
\subfigure[$t = 0.35$]{\includegraphics[width=2.1in]{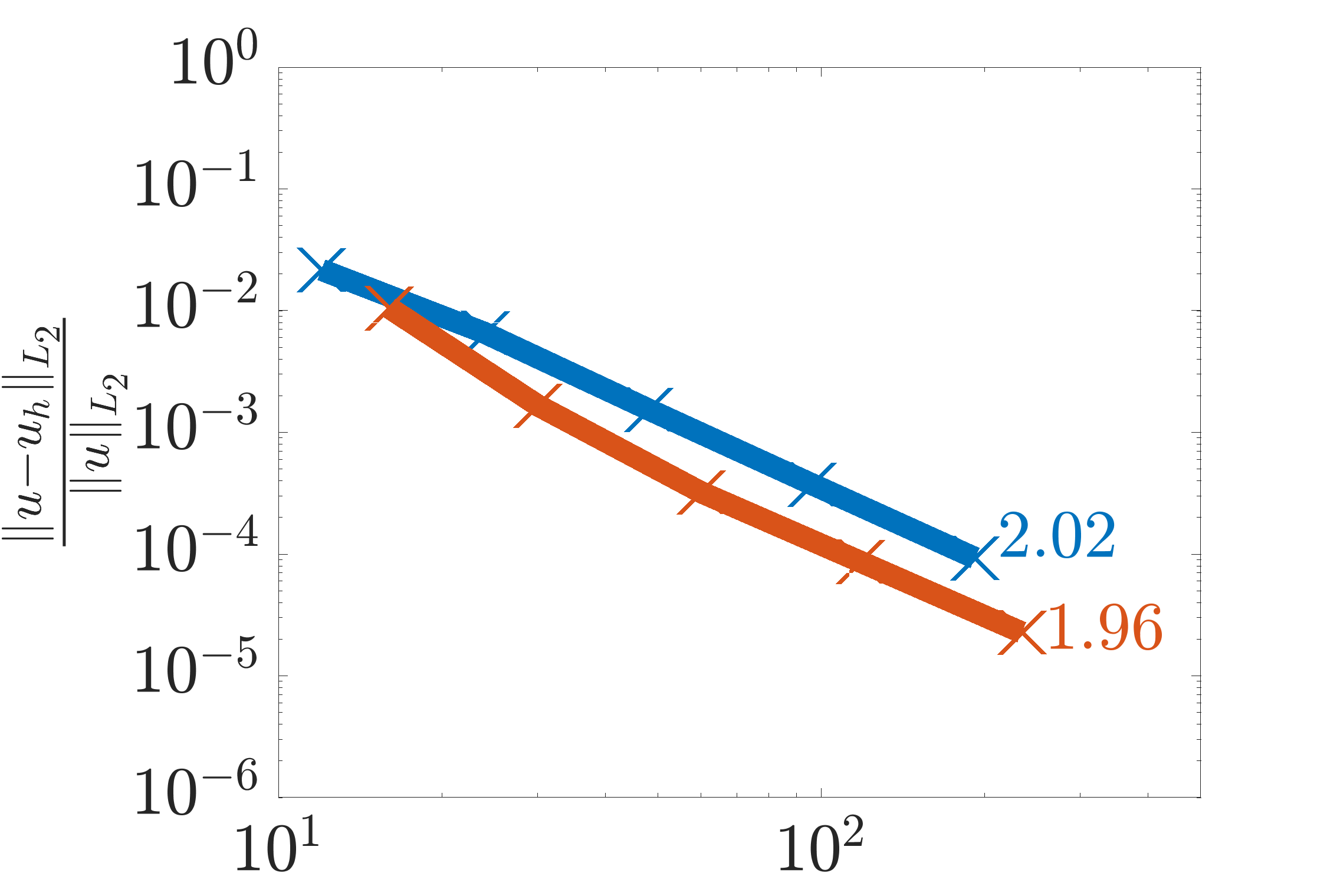}}
\subfigure[$t = 0.5$]{\includegraphics[width=2.1in]{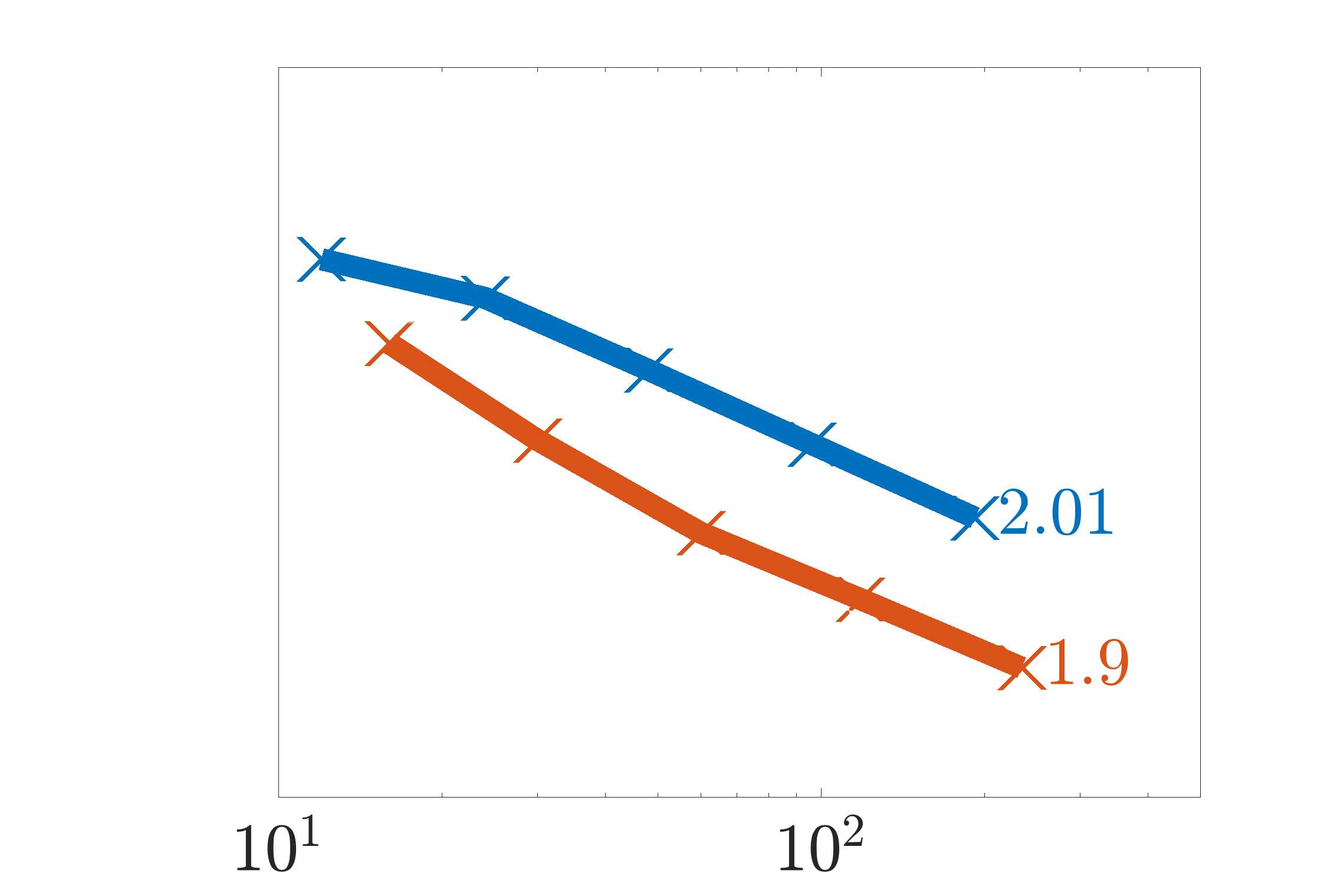}}
\subfigure[$t = 0.75$]{\includegraphics[width=2.1in]{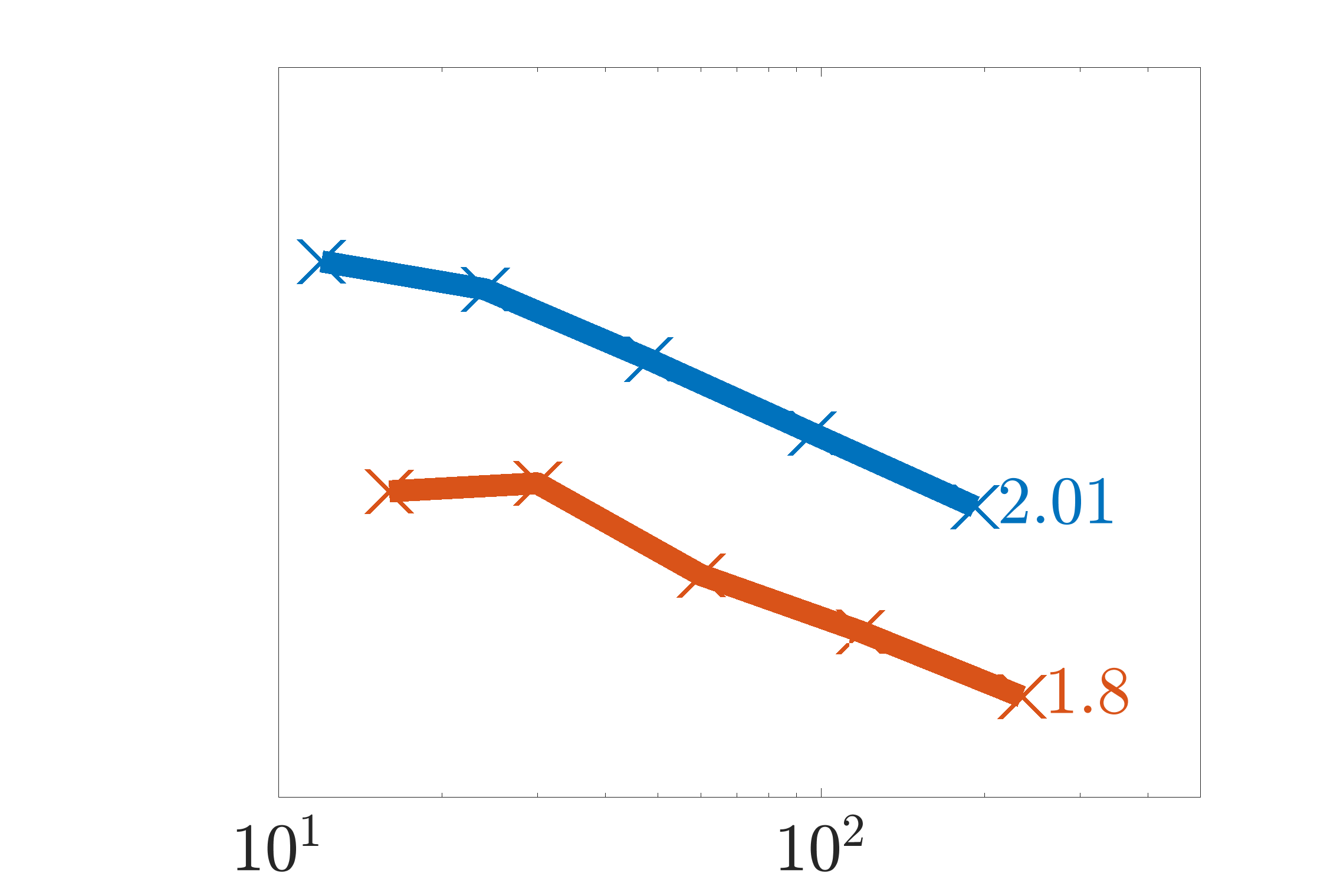}}
\end{subfigmatrix}
\caption{Convergence the relative $L_2$ integral norm for the shock problem with $\nu = \frac{1}{50}$}
\label{fig:Example2_L2_vsdofs_nu1over50}
\end{center}
\end{figure}

\begin{figure}[ht!]
\begin{center}
\begin{subfigmatrix}{6}
\subfigure[$t = 0$]{\includegraphics[width=2.1in]{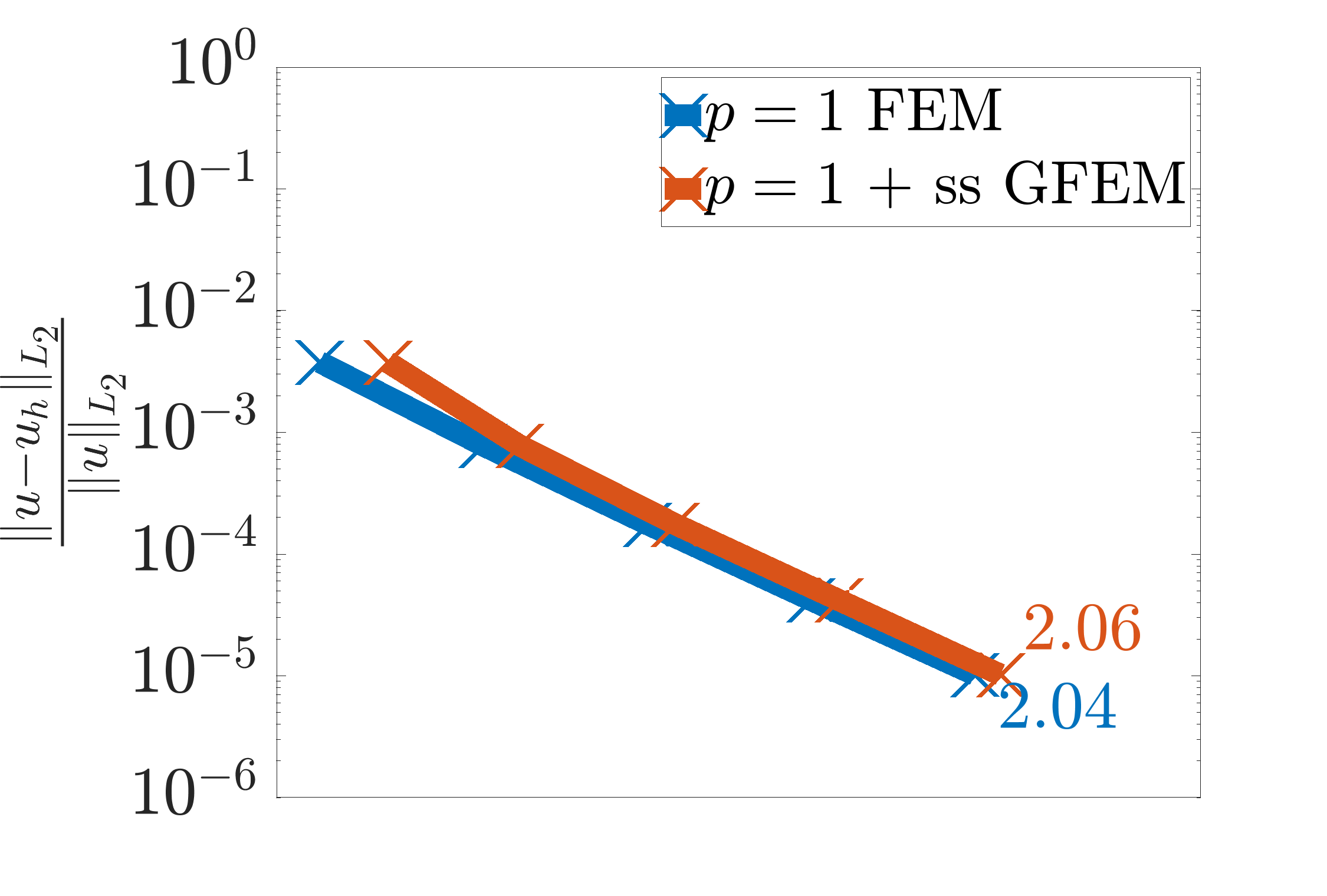}}
\subfigure[$t = 0.25$]{\includegraphics[width=2.1in]{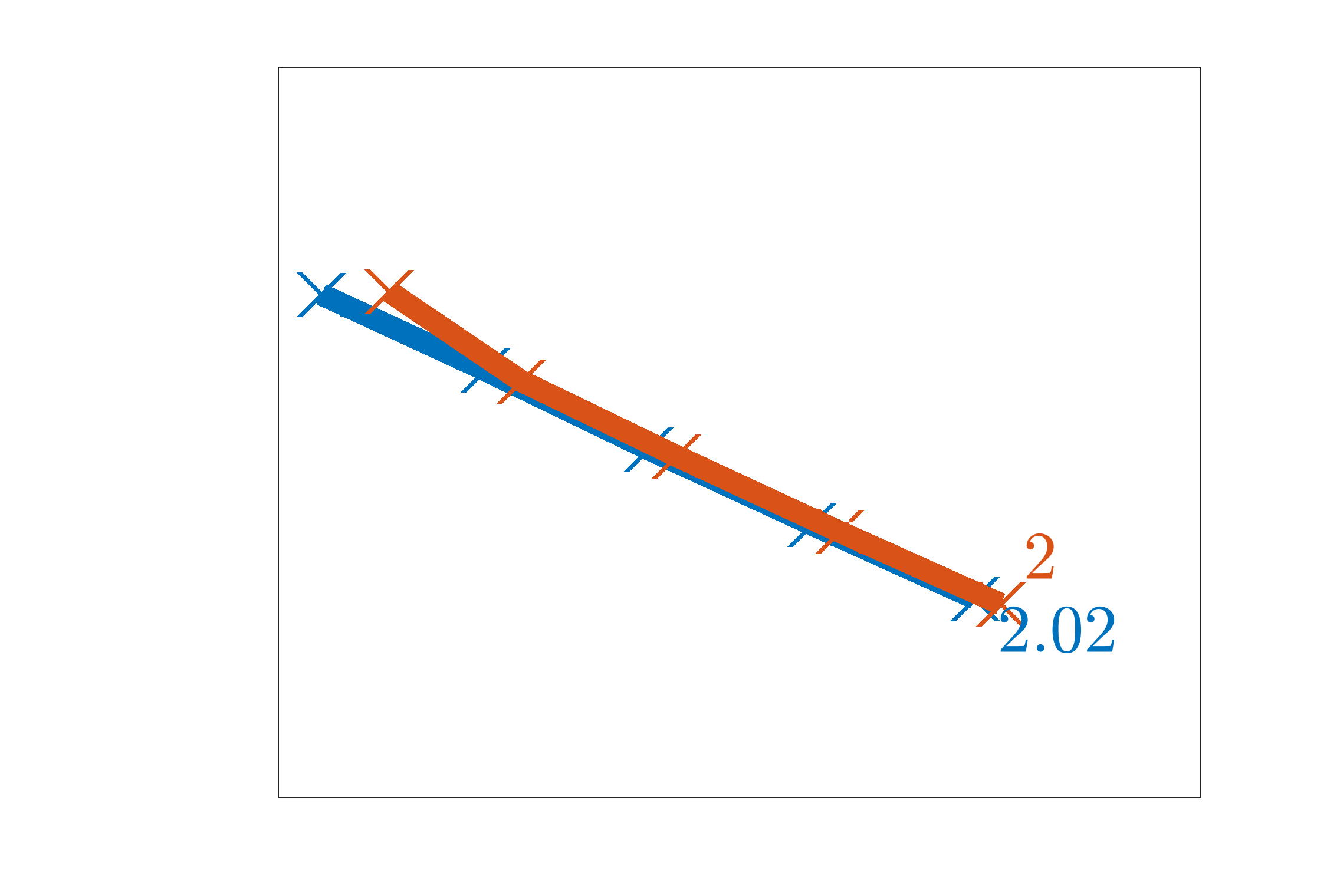}}
\subfigure[$t = 0.3$]{\includegraphics[width=2.1in]{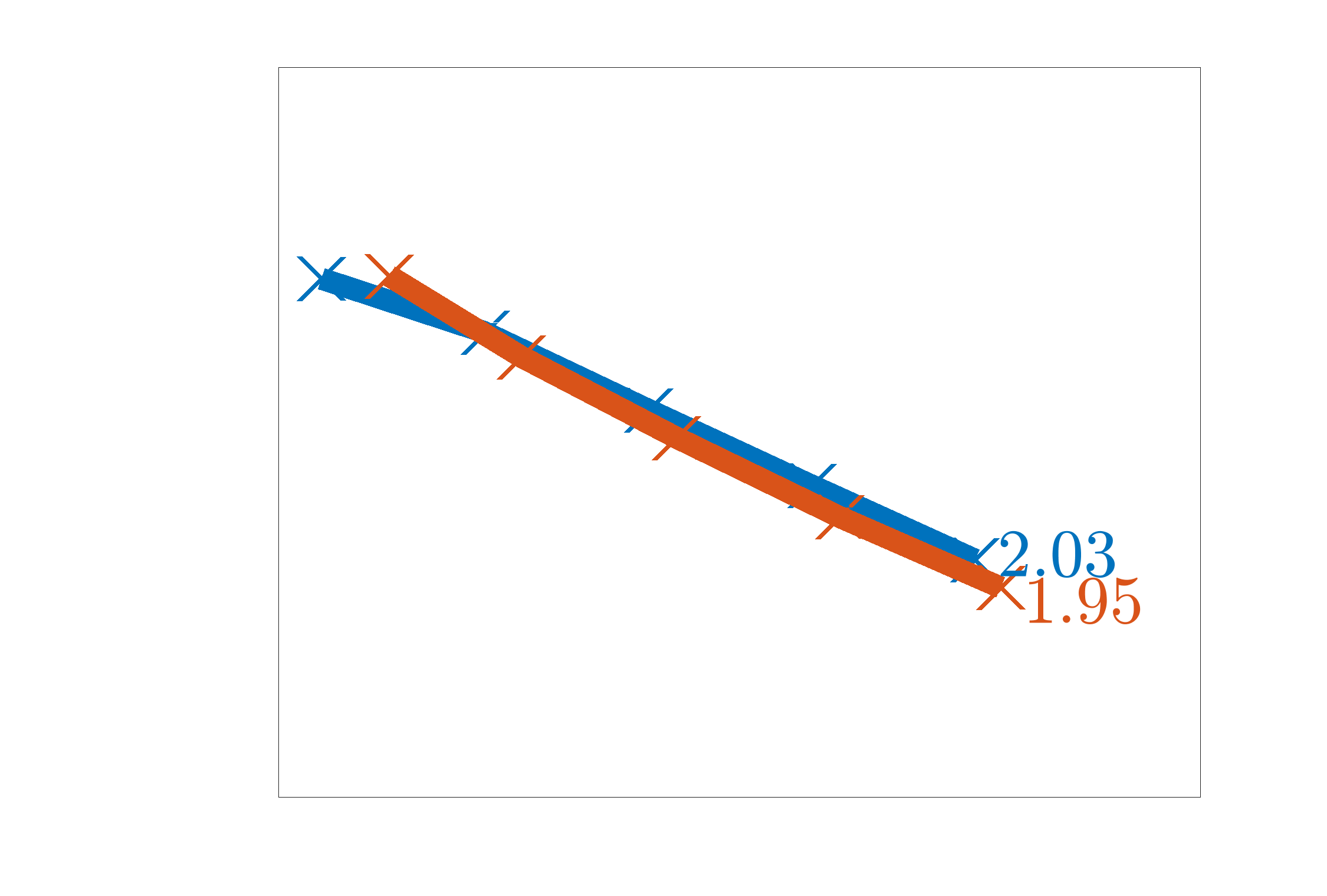}}
\subfigure[$t = 0.35$]{\includegraphics[width=2.1in]{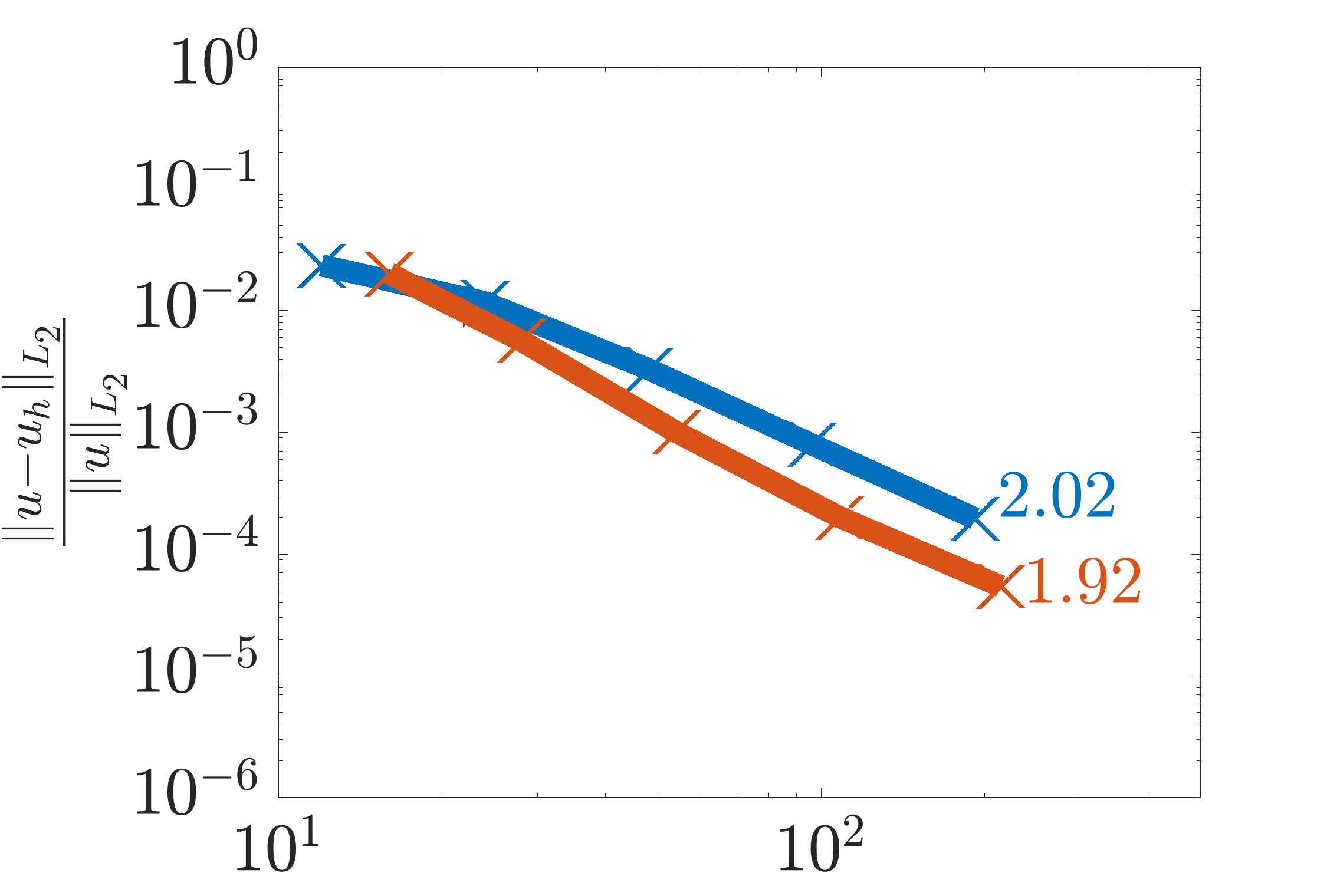}}
\subfigure[$t = 0.5$]{\includegraphics[width=2.1in]{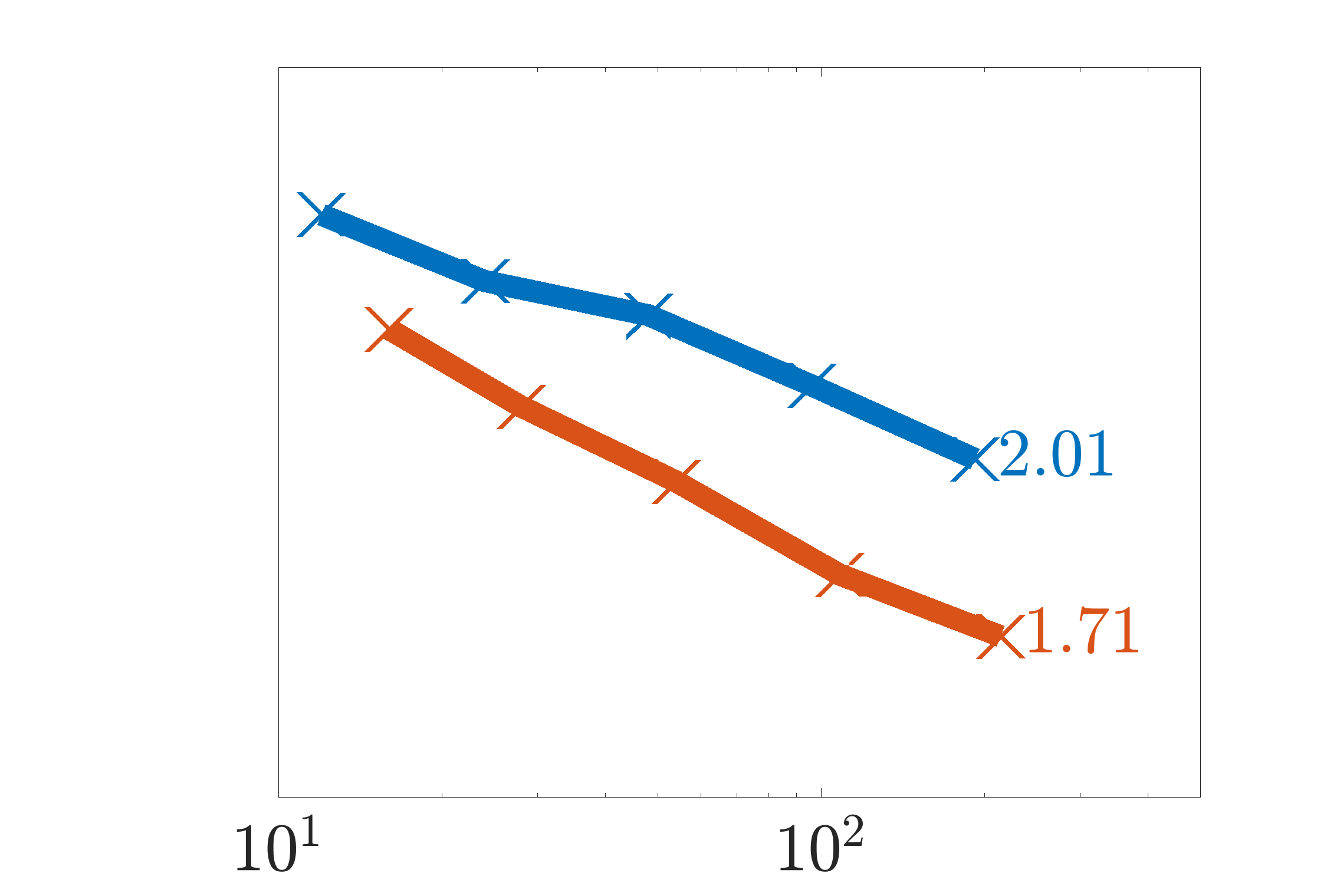}}
\subfigure[$t = 0.75$]{\includegraphics[width=2.1in]{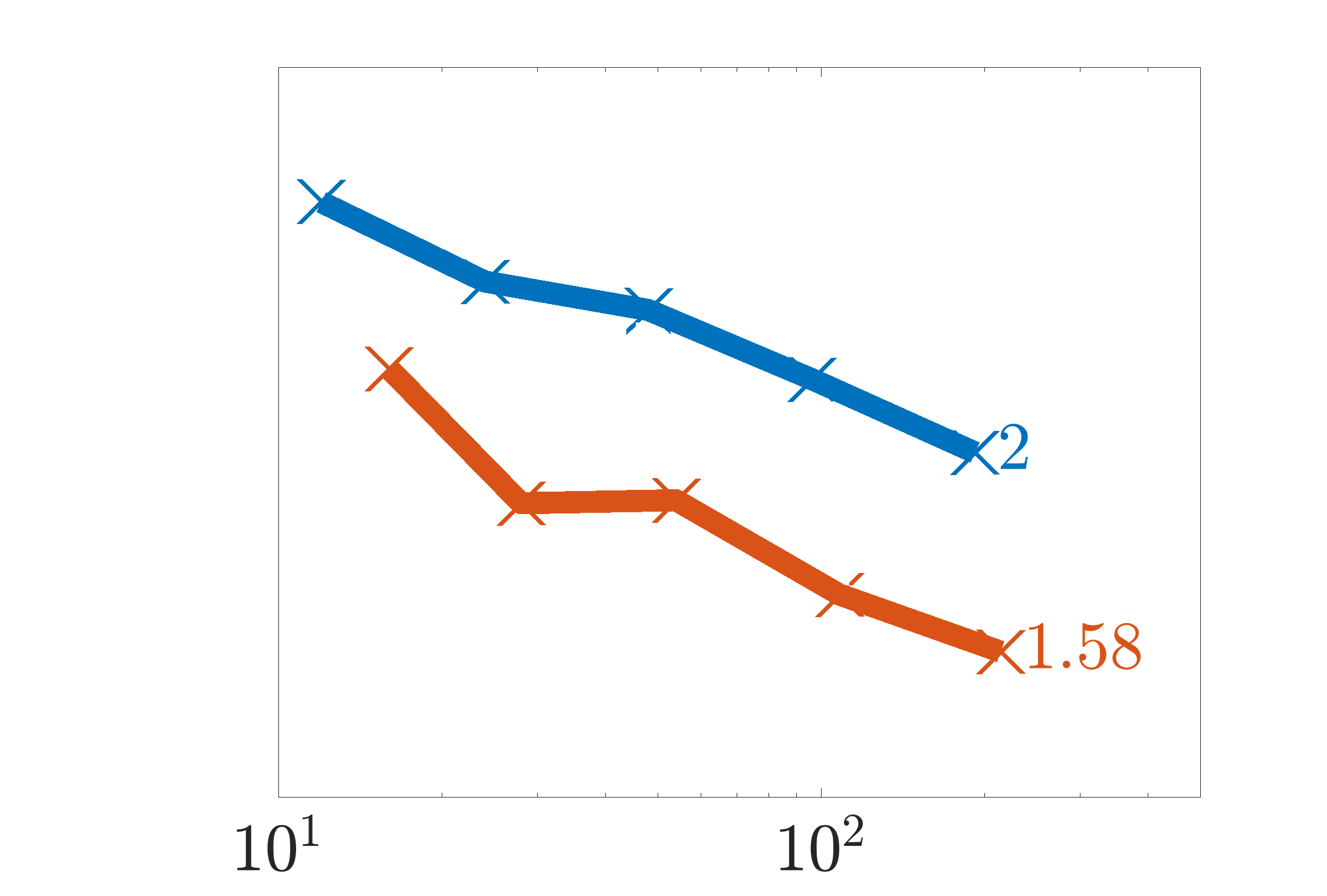}}
\end{subfigmatrix}
\caption{Convergence the relative $L_2$ integral norm for the shock problem with $\nu = \frac{1}{100}$}
\label{fig:Example2_L2_vsdofs_nu1over100}
\end{center}
\end{figure}

\begin{figure}[ht!]
\begin{center}
\begin{subfigmatrix}{6}
\subfigure[$t = 0$]{\includegraphics[width=2.1in]{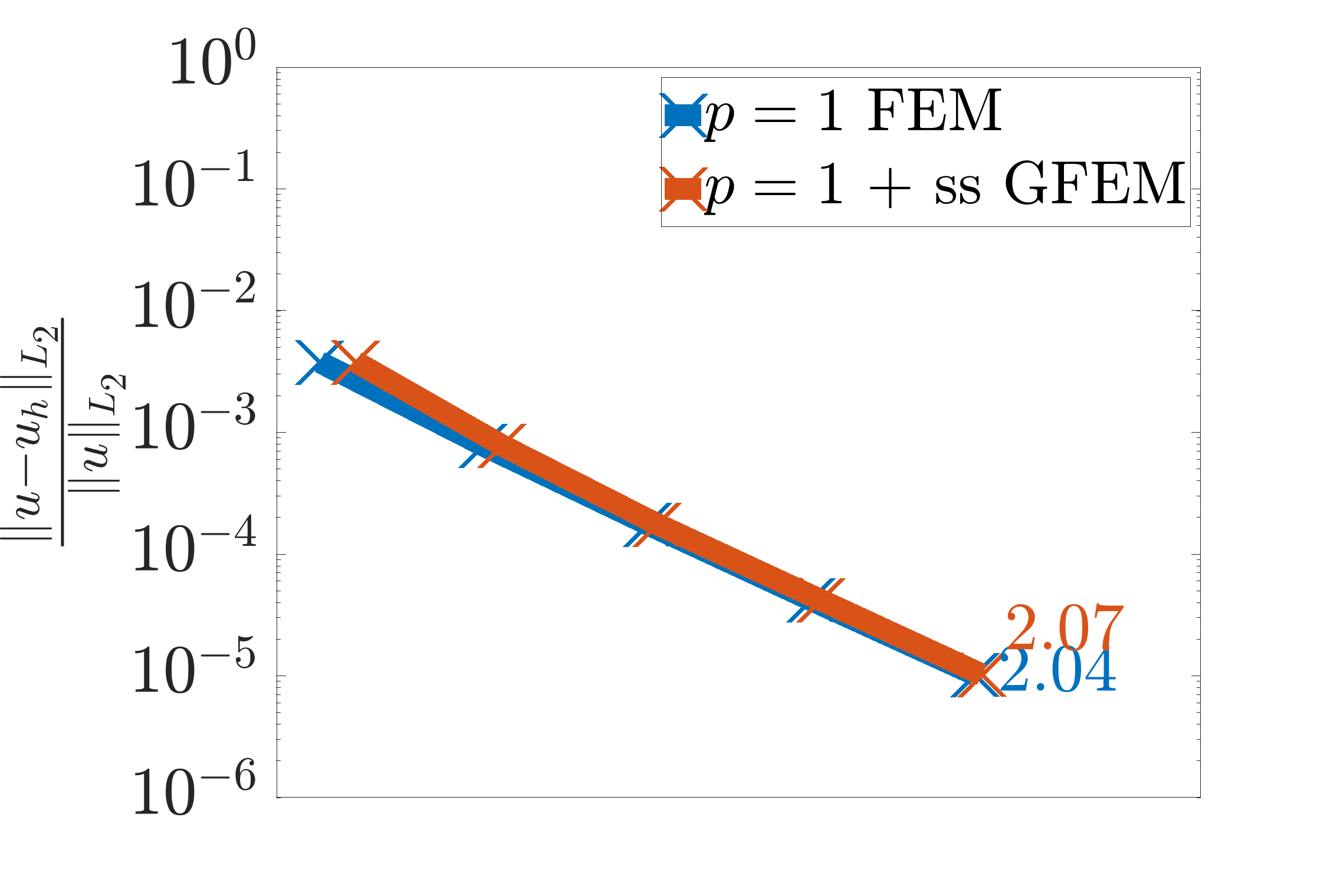}}
\subfigure[$t = 0.25$]{\includegraphics[width=2.1in]{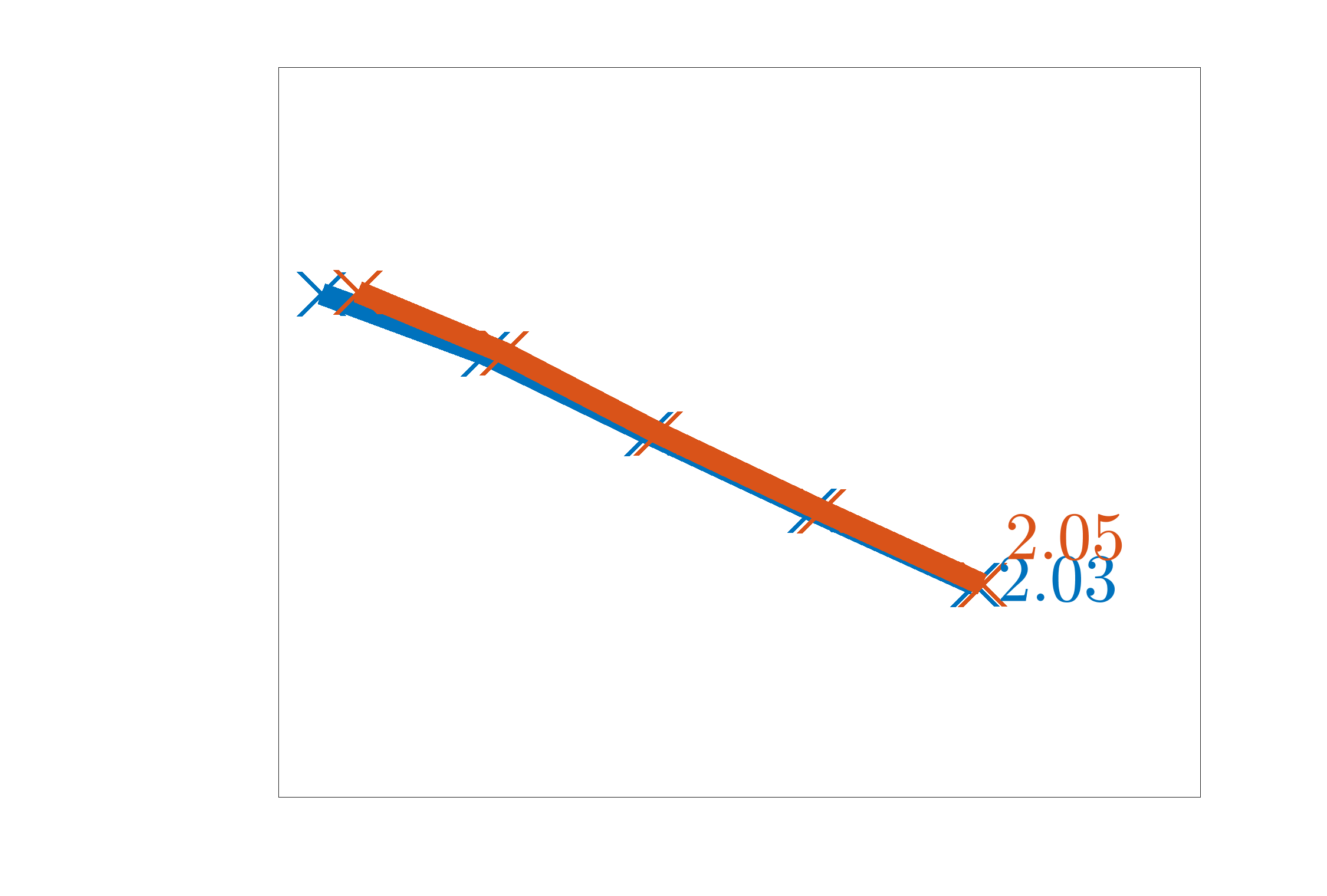}}
\subfigure[$t = 0.3$]{\includegraphics[width=2.1in]{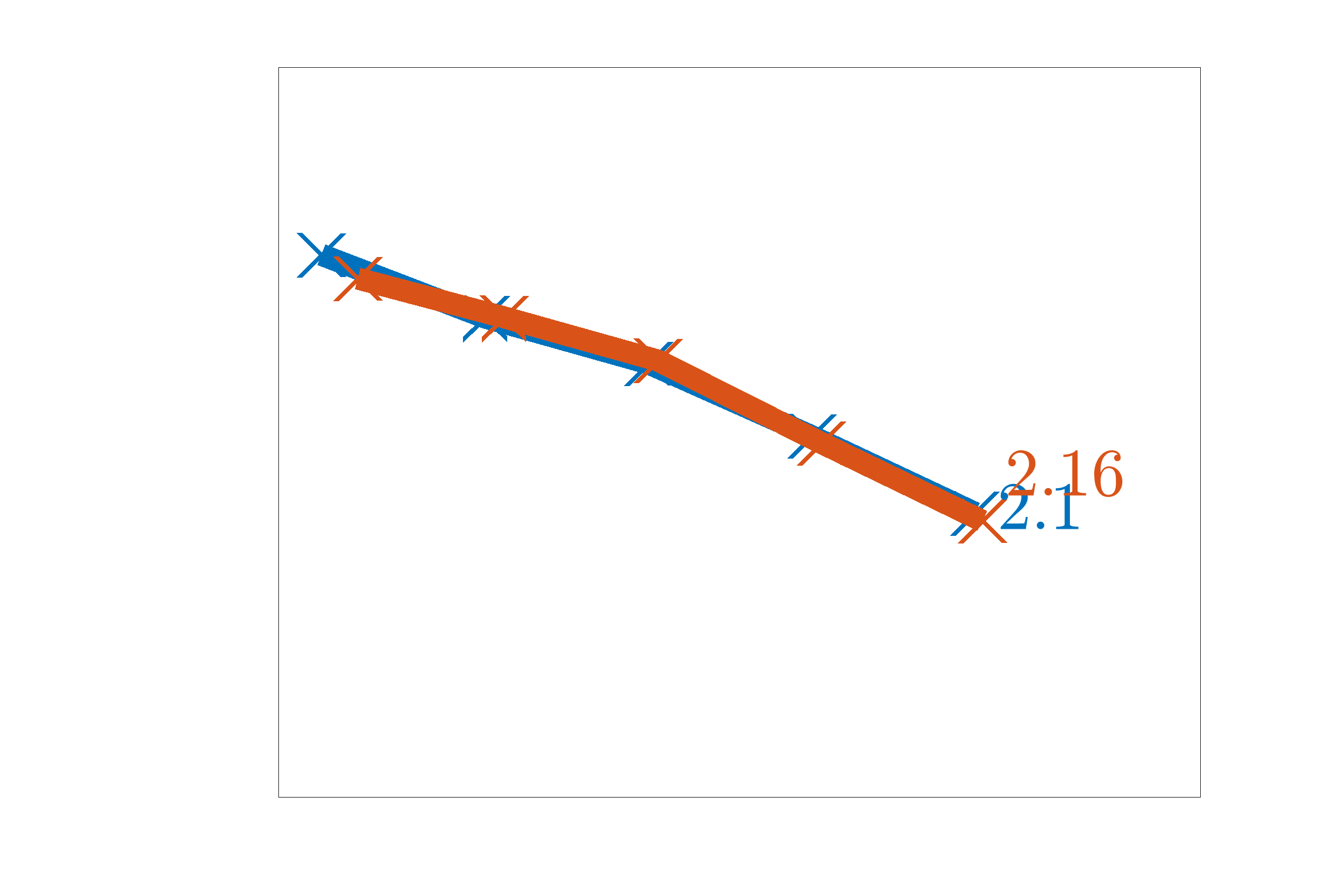}}
\subfigure[$t = 0.35$]{\includegraphics[width=2.1in]{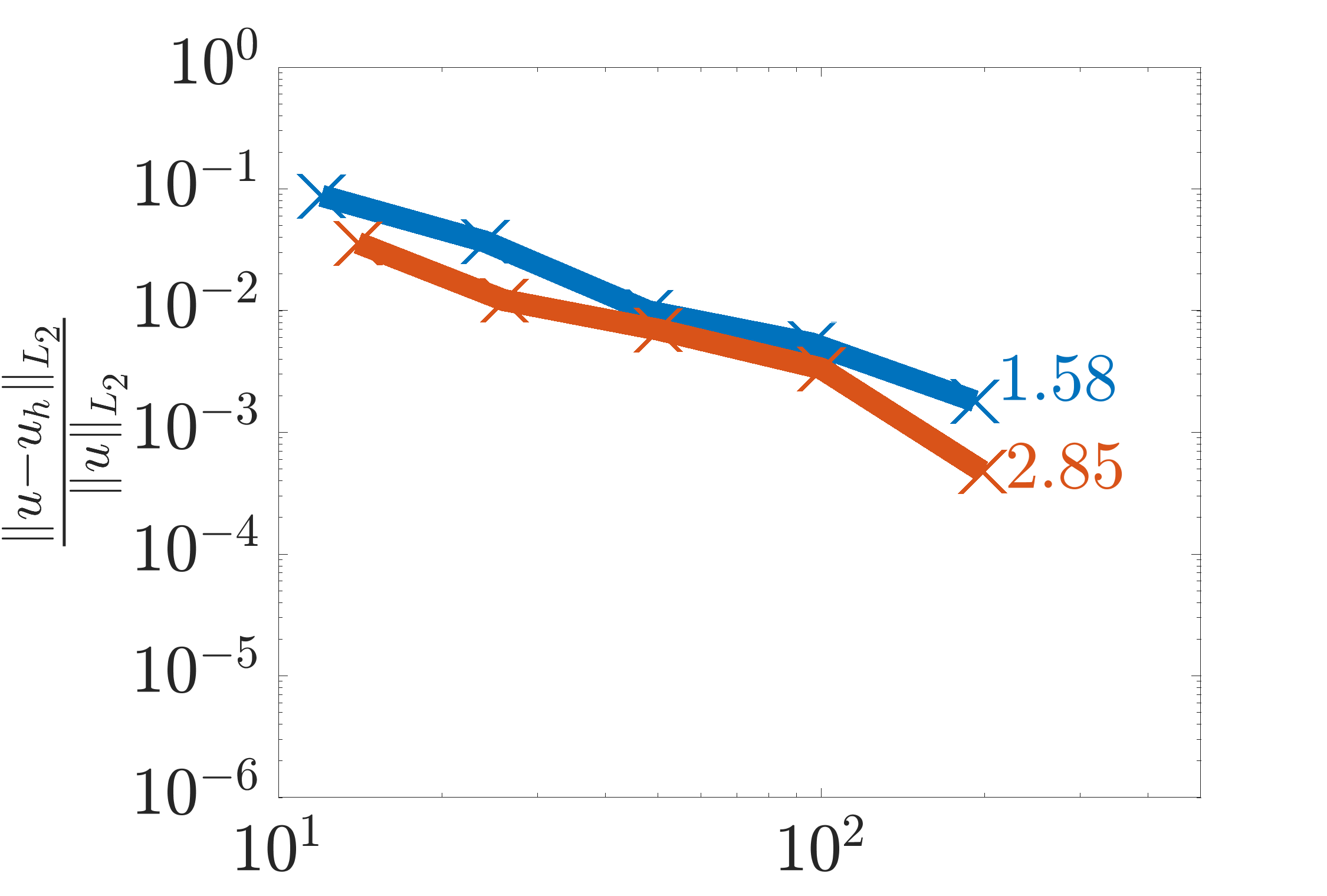}}
\subfigure[$t = 0.5$]{\includegraphics[width=2.1in]{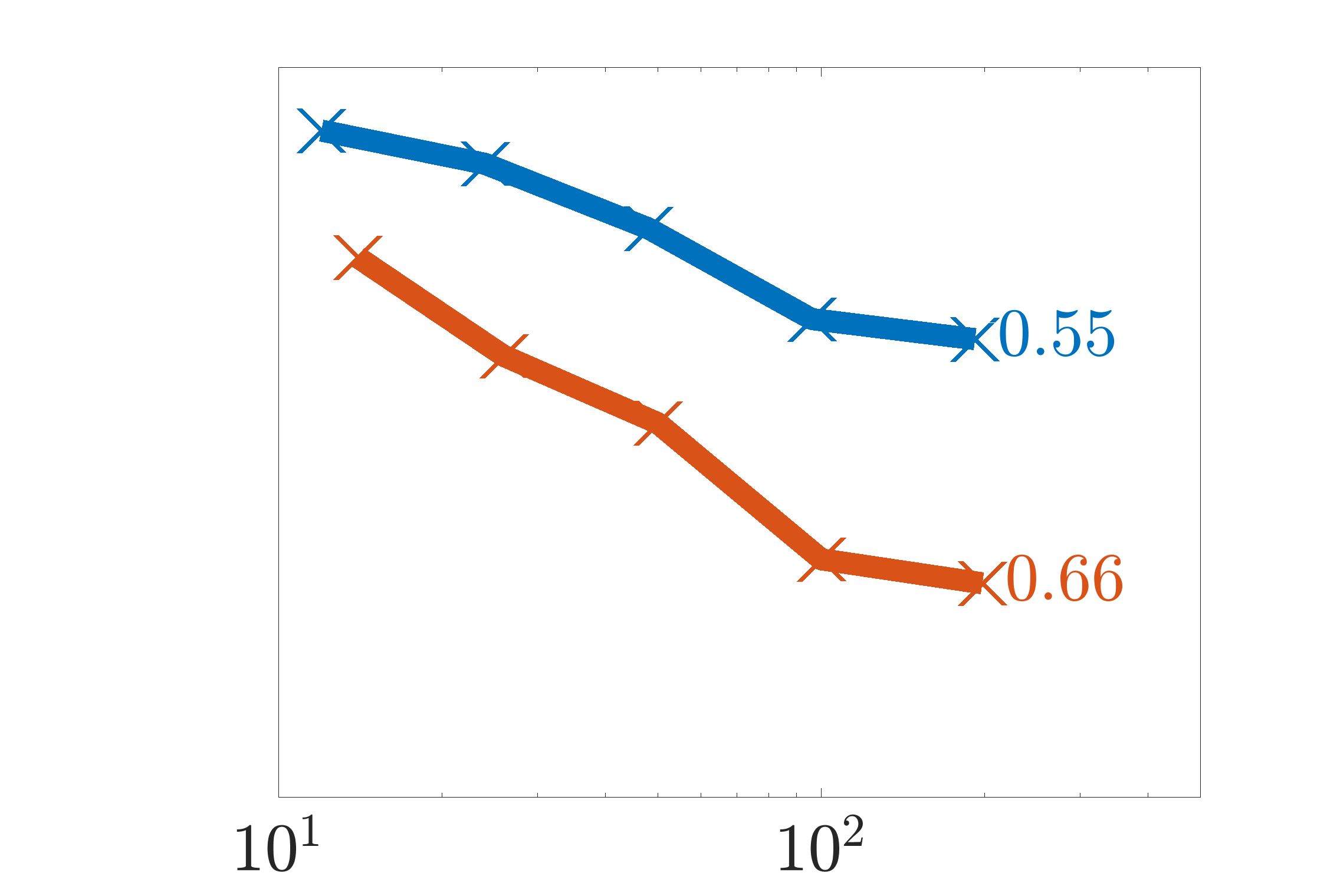}}
\subfigure[$t = 0.75$]{\includegraphics[width=2.1in]{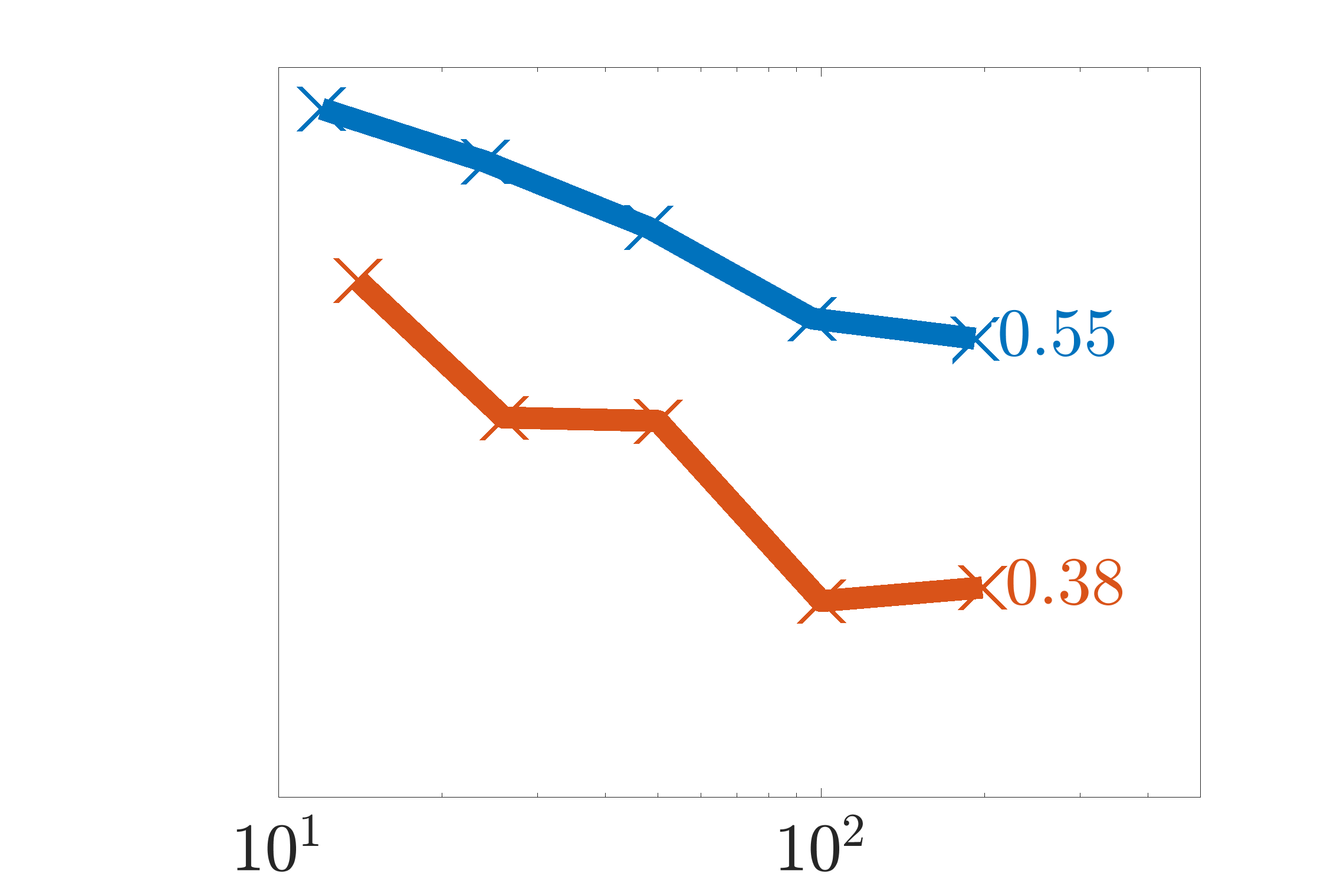}}
\end{subfigmatrix}
\caption{Convergence the relative $L_2$ integral norm for the shock problem with $\nu = \frac{1}{500}$}
\label{fig:Example2_L2_vsdofs_nu1over500}
\end{center}
\end{figure}

\begin{figure}[ht!]
\begin{center}
\begin{subfigmatrix}{6}
\subfigure[$t = 0$]{\includegraphics[width=2.1in]{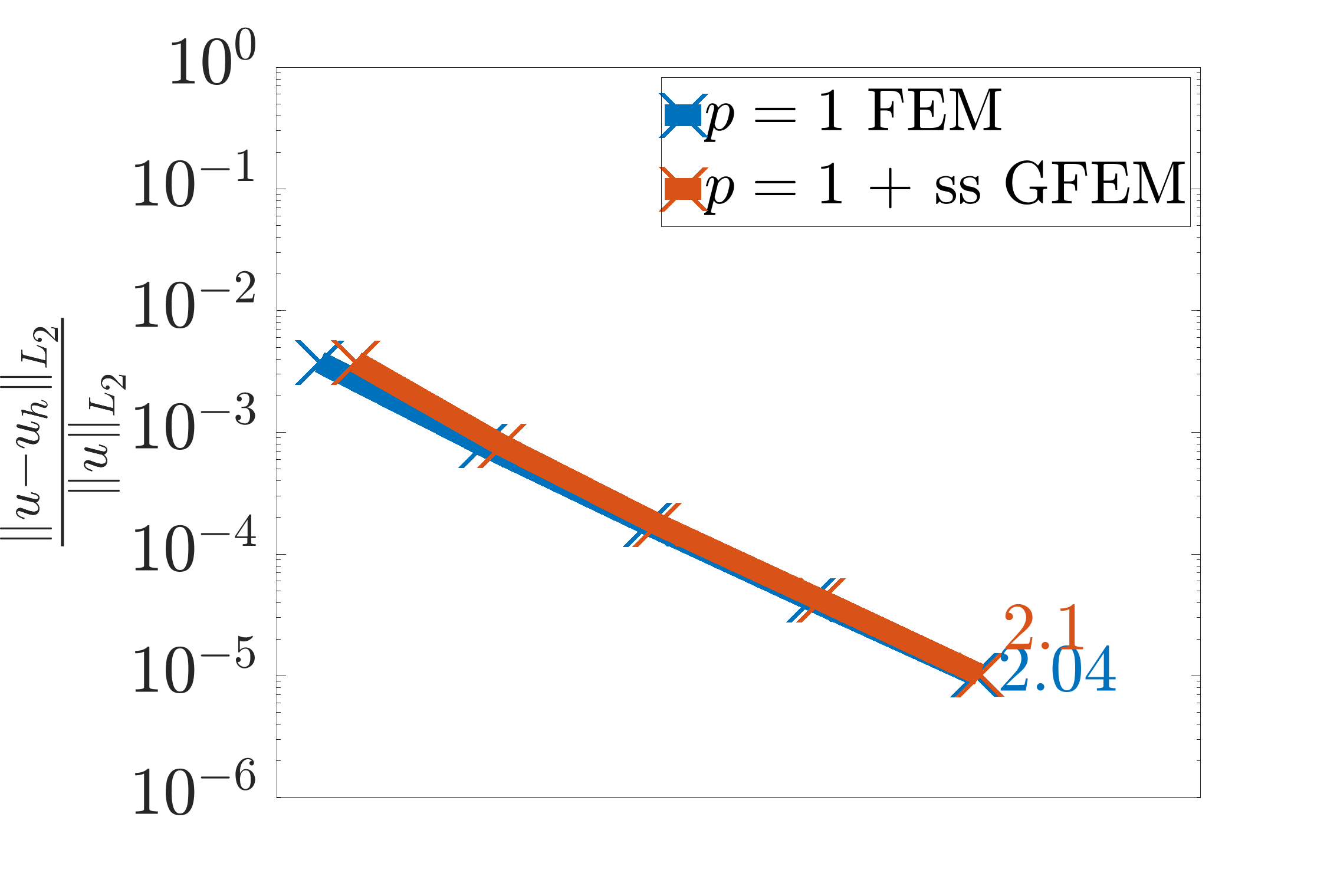}}
\subfigure[$t = 0.25$]{\includegraphics[width=2.1in]{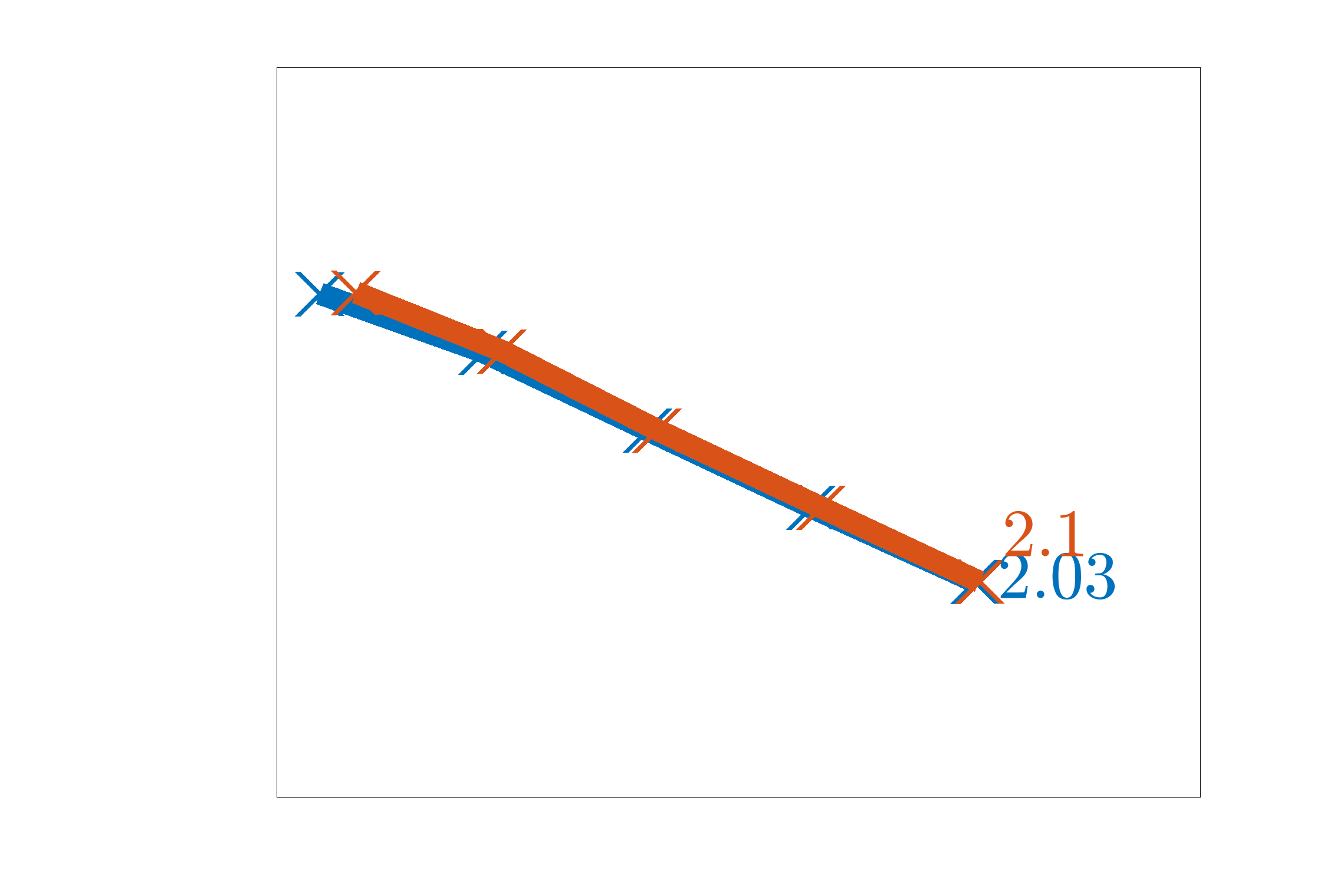}}
\subfigure[$t = 0.3$]{\includegraphics[width=2.1in]{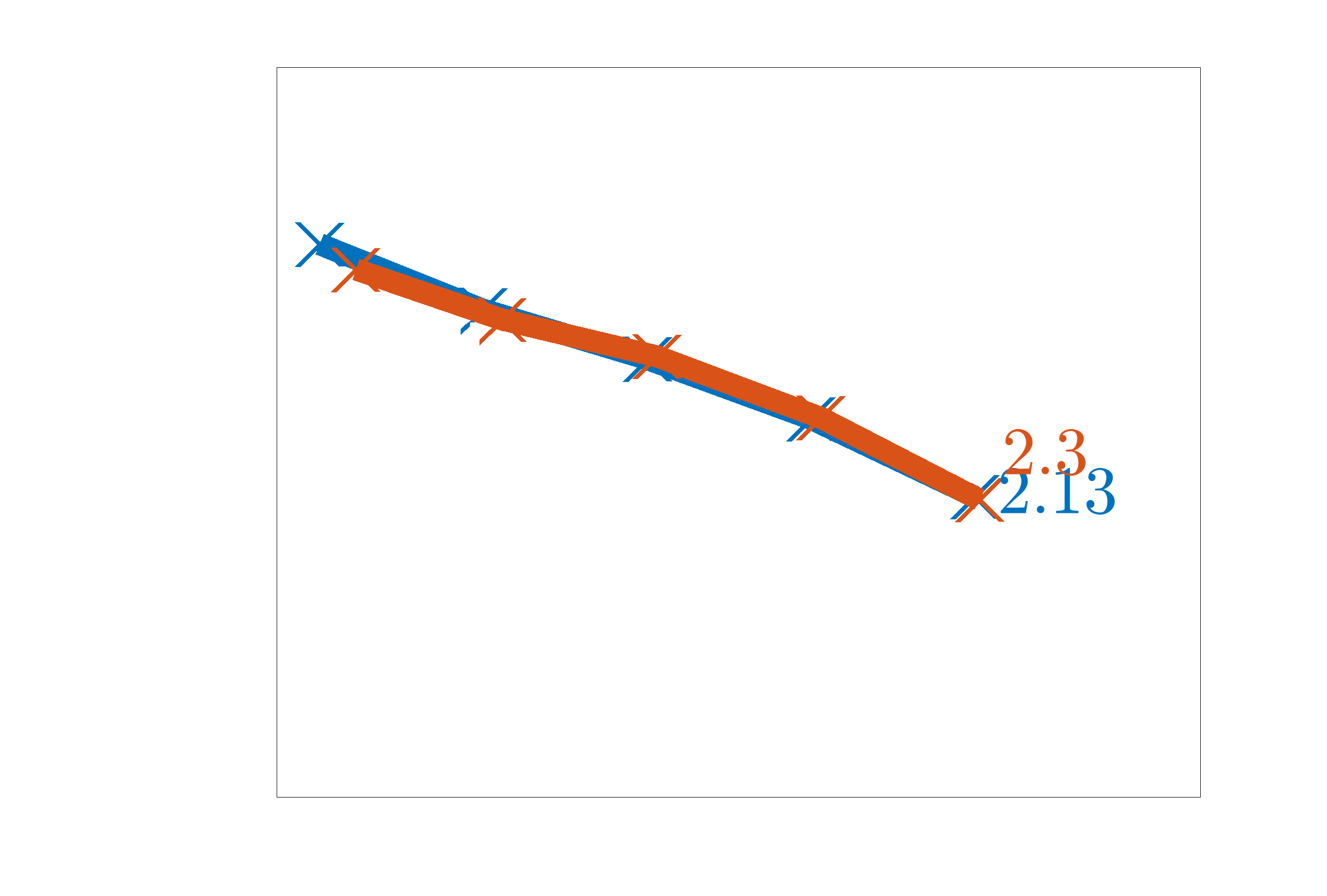}}
\subfigure[$t = 0.35$]{\includegraphics[width=2.1in]{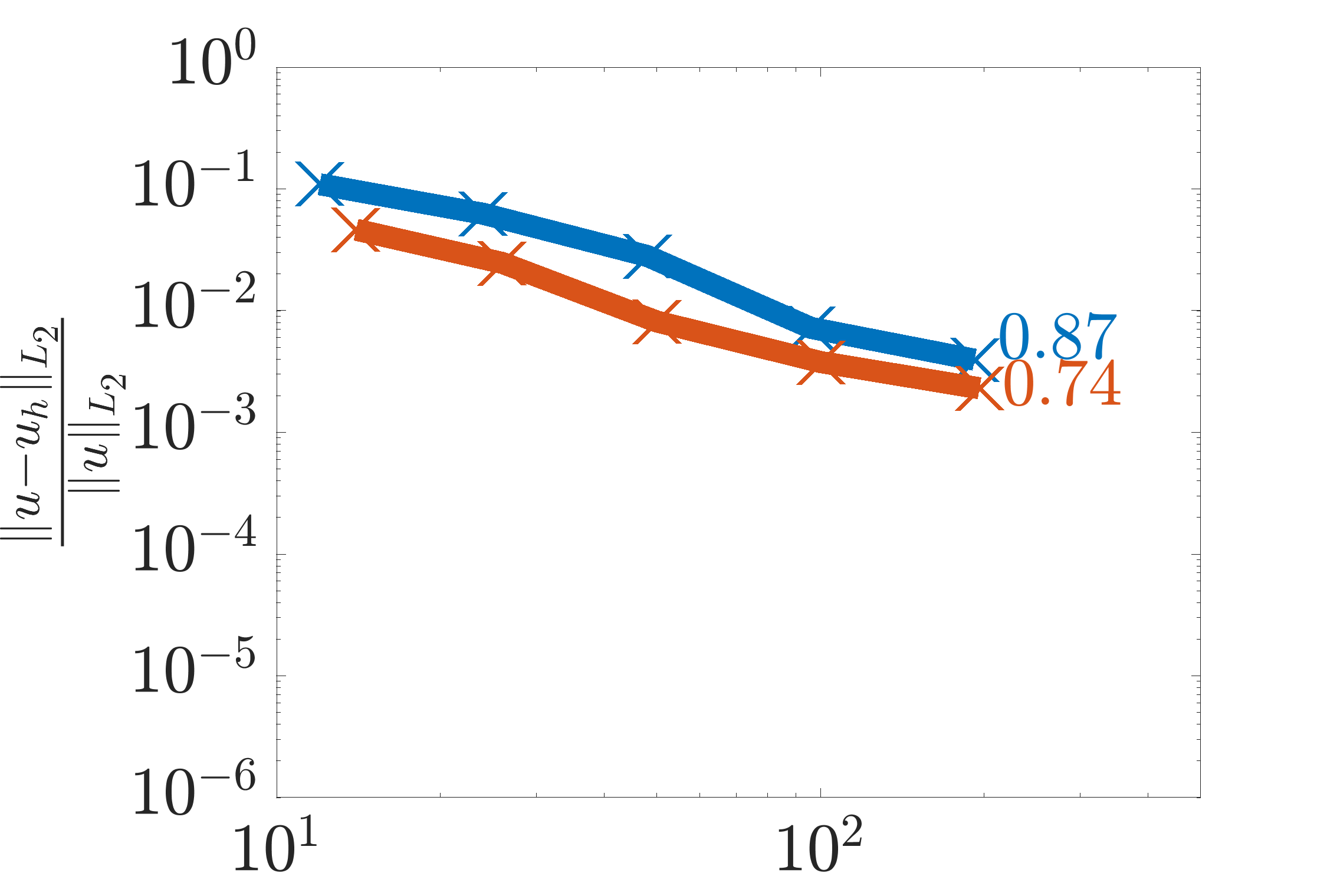}}
\subfigure[$t = 0.5$]{\includegraphics[width=2.1in]{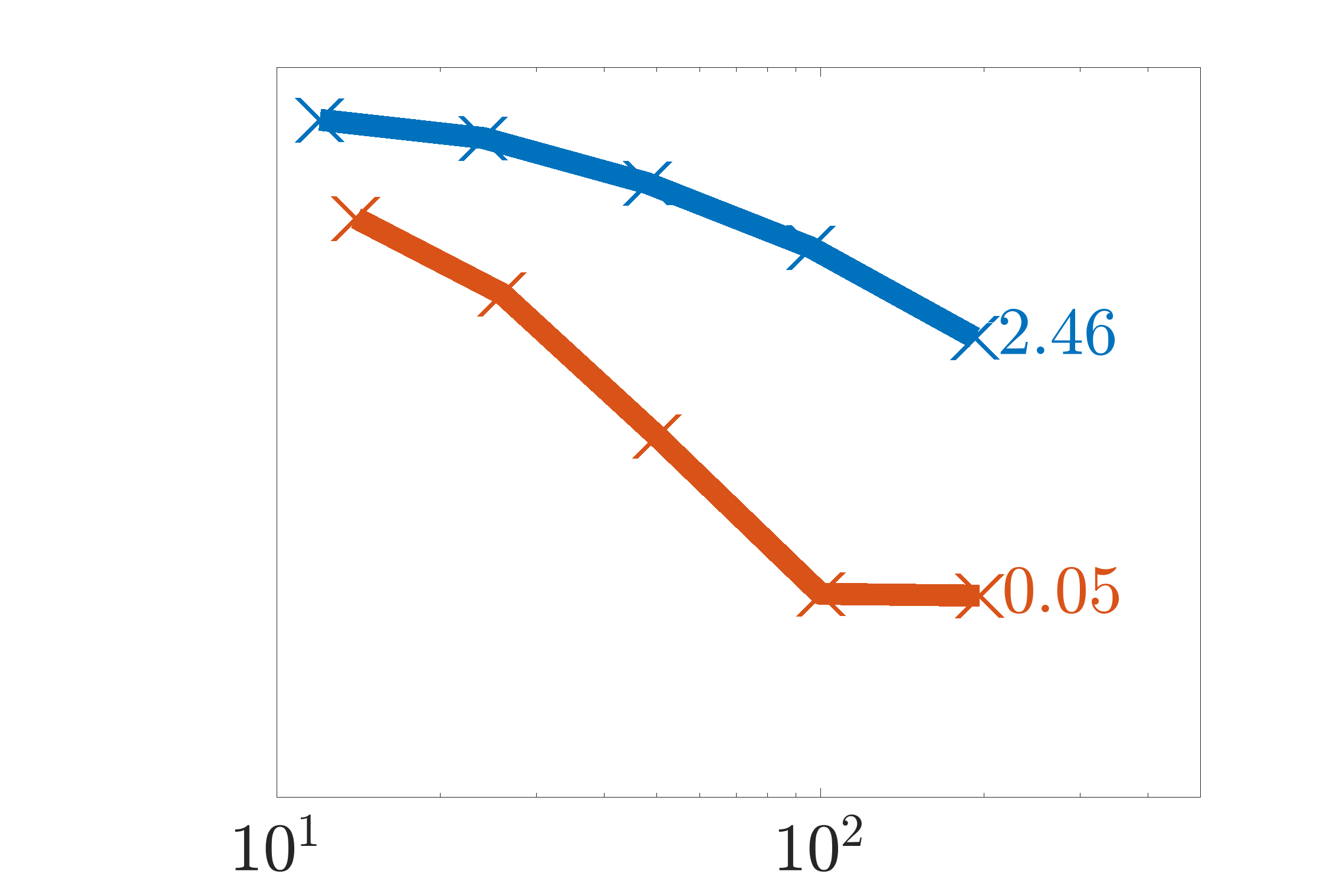}}
\subfigure[$t = 0.75$]{\includegraphics[width=2.1in]{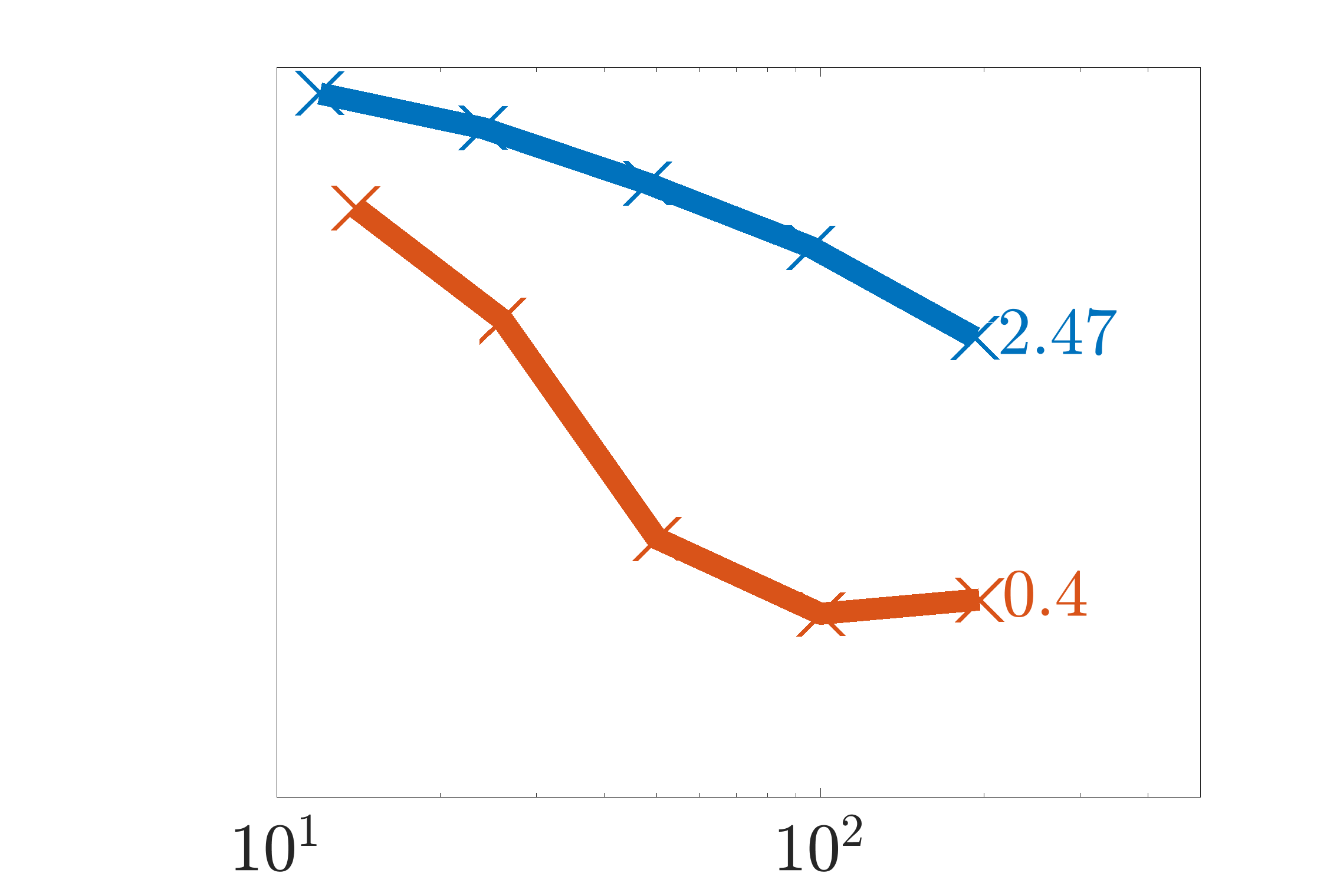}}
\end{subfigmatrix}
\caption{Convergence the relative $L_2$ integral norm for the shock problem with $\nu = \frac{1}{1000}$}
\label{fig:Example2_L2_vsdofs_nu1over1000}
\end{center}
\end{figure}

\begin{figure}[ht!]
\begin{center}
\begin{subfigmatrix}{6}
\subfigure[$t = 0$]{\includegraphics[width=2.1in]{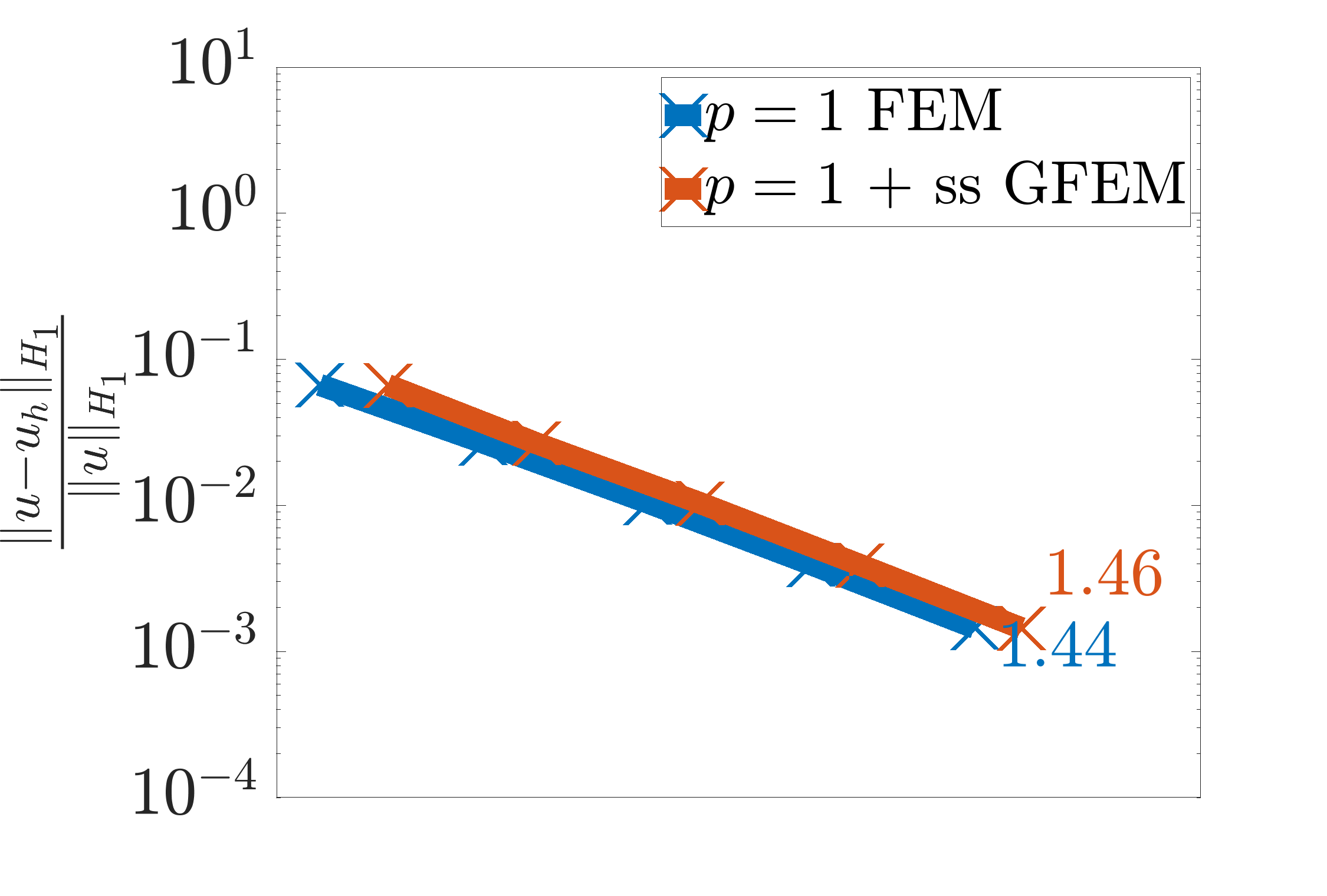}}
\subfigure[$t = 0.25$]{\includegraphics[width=2.1in]{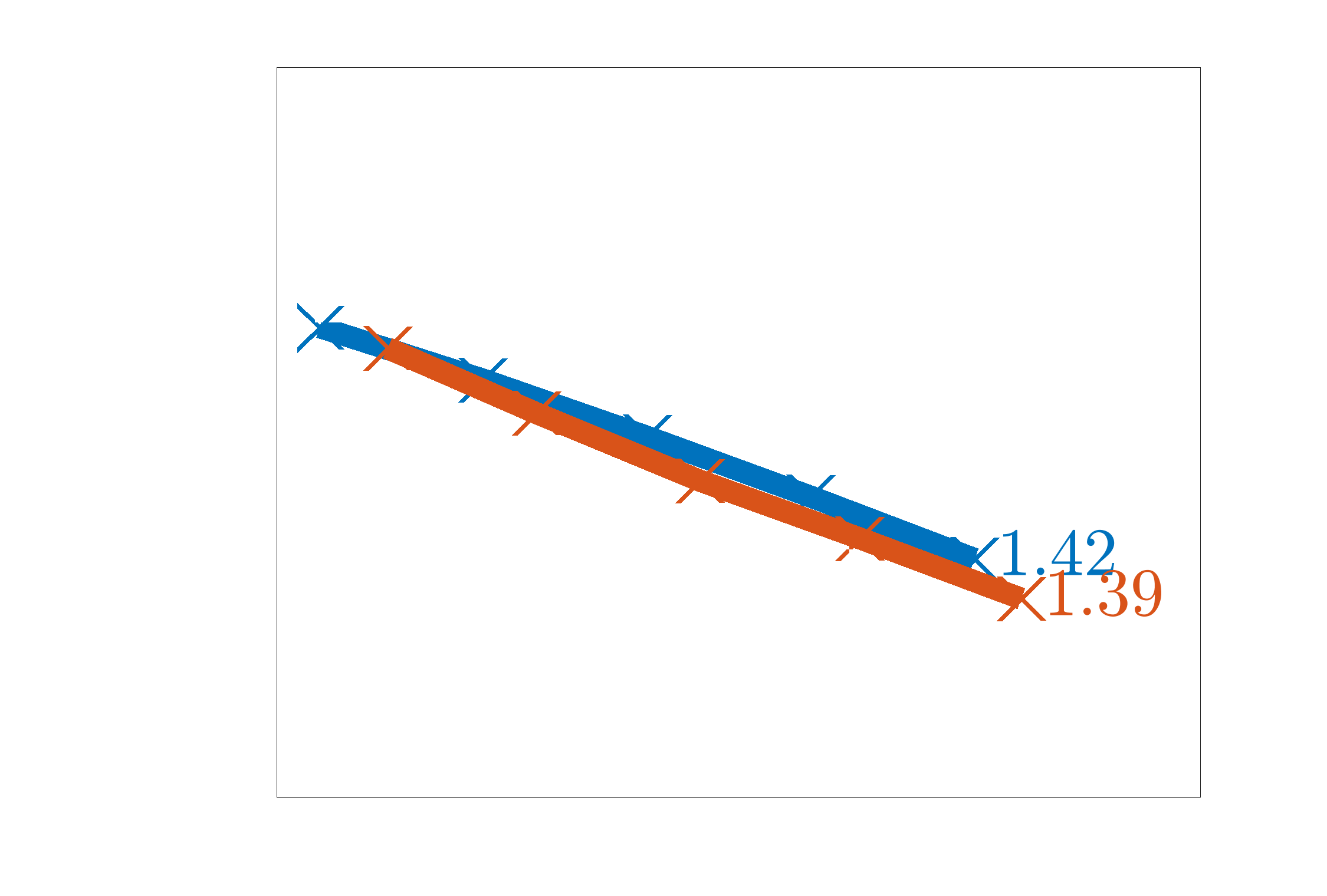}}
\subfigure[$t = 0.3$]{\includegraphics[width=2.1in]{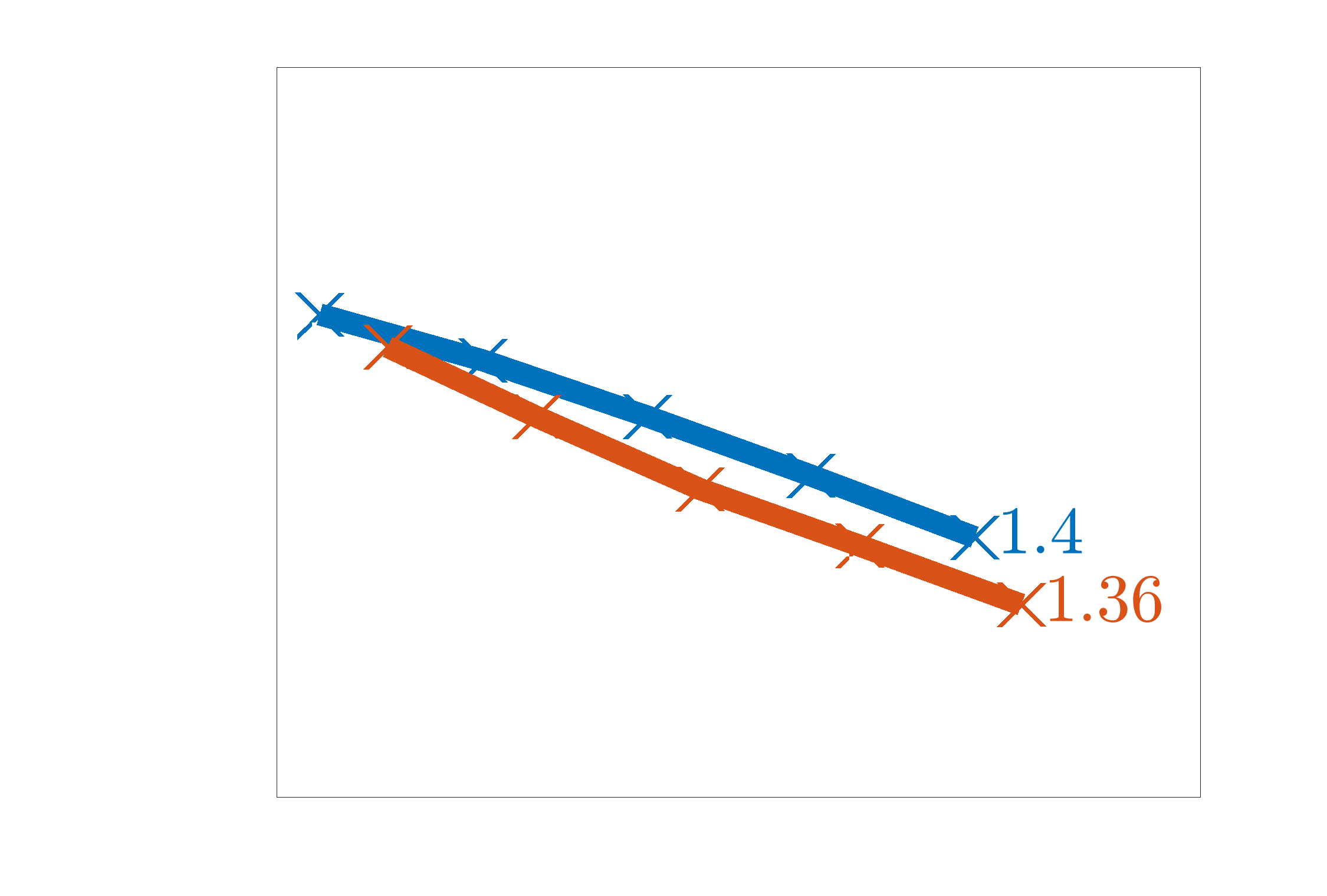}}
\subfigure[$t = 0.35$]{\includegraphics[width=2.1in]{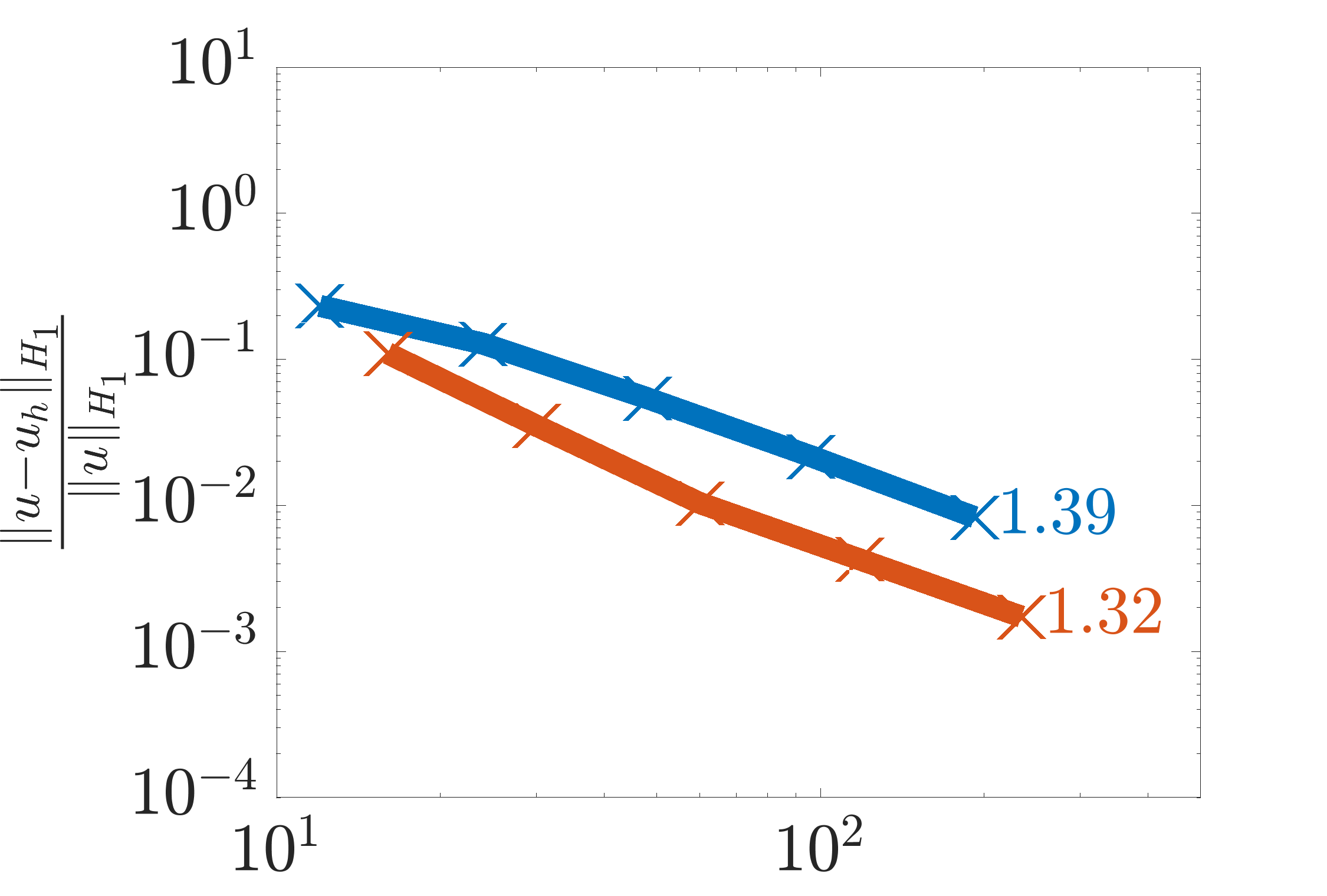}}
\subfigure[$t = 0.5$]{\includegraphics[width=2.1in]{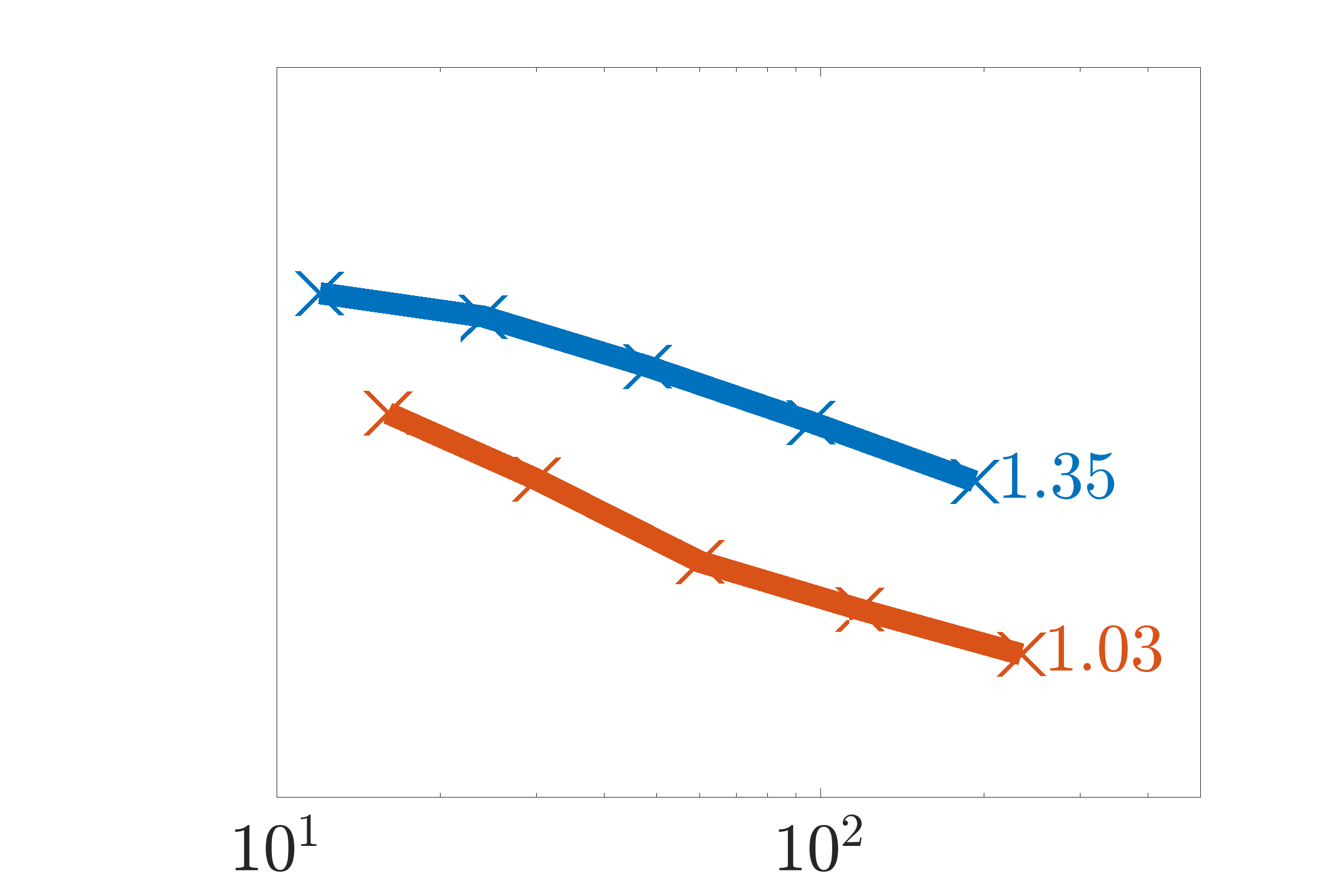}}
\subfigure[$t = 0.75$]{\includegraphics[width=2.1in]{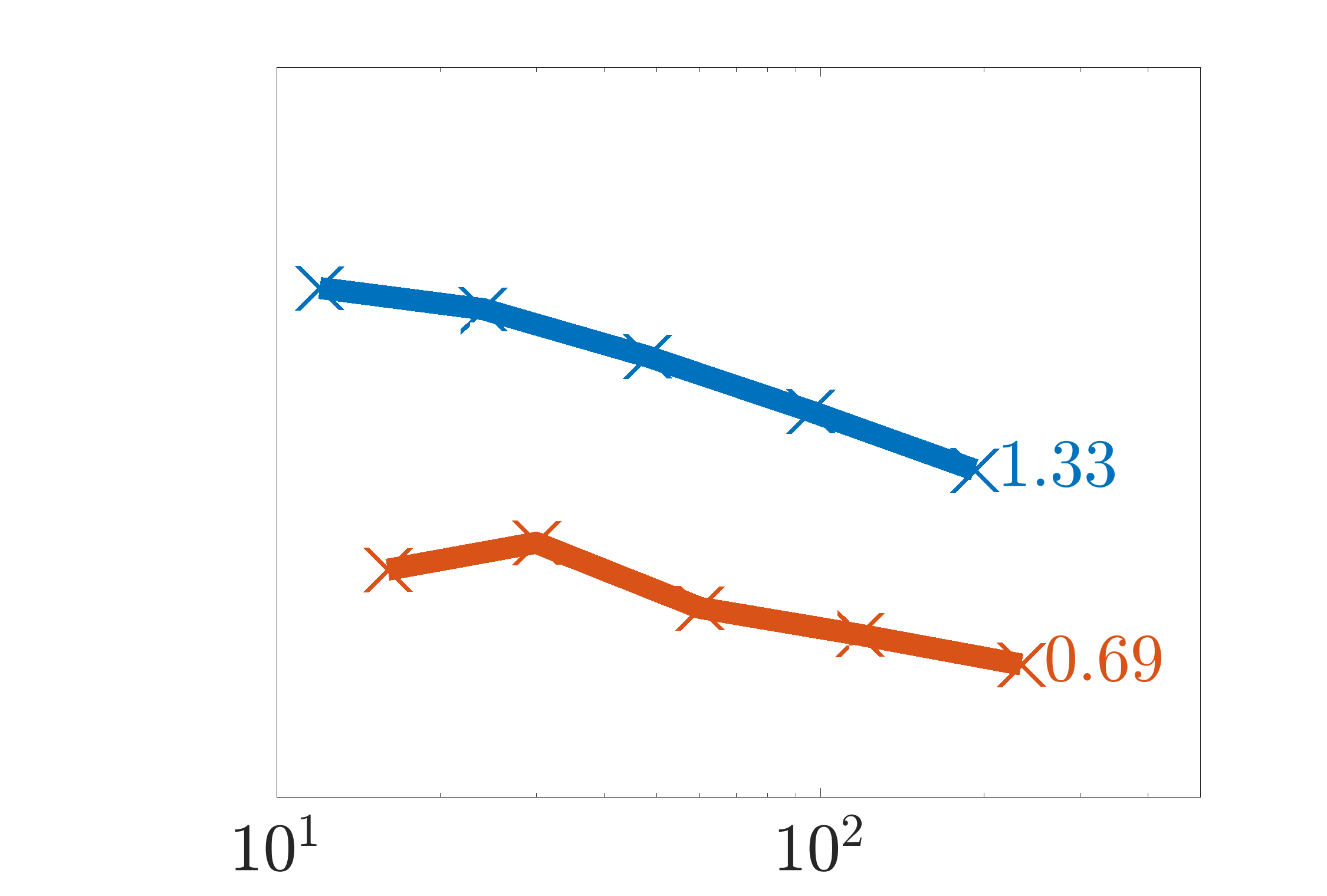}}
\end{subfigmatrix}
\caption{Convergence the relative $H_1$ integral norm for the shock problem with $\nu = \frac{1}{50}$}
\label{fig:Example2_H1_vsdofs_nu1over50}
\end{center}
\end{figure}

\begin{figure}[ht!]
\begin{center}
\begin{subfigmatrix}{6}
\subfigure[$t = 0$]{\includegraphics[width=2.1in]{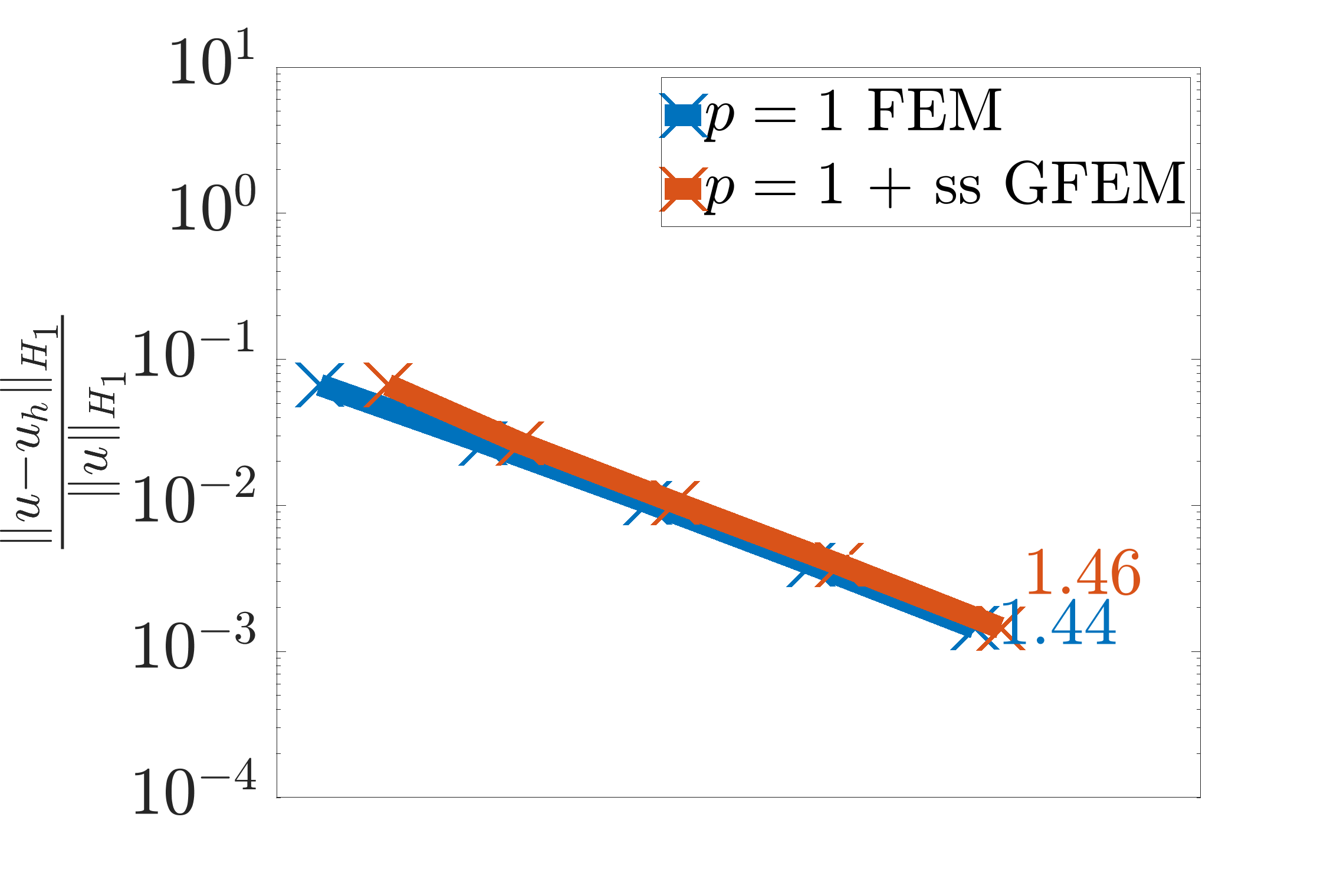}}
\subfigure[$t = 0.25$]{\includegraphics[width=2.1in]{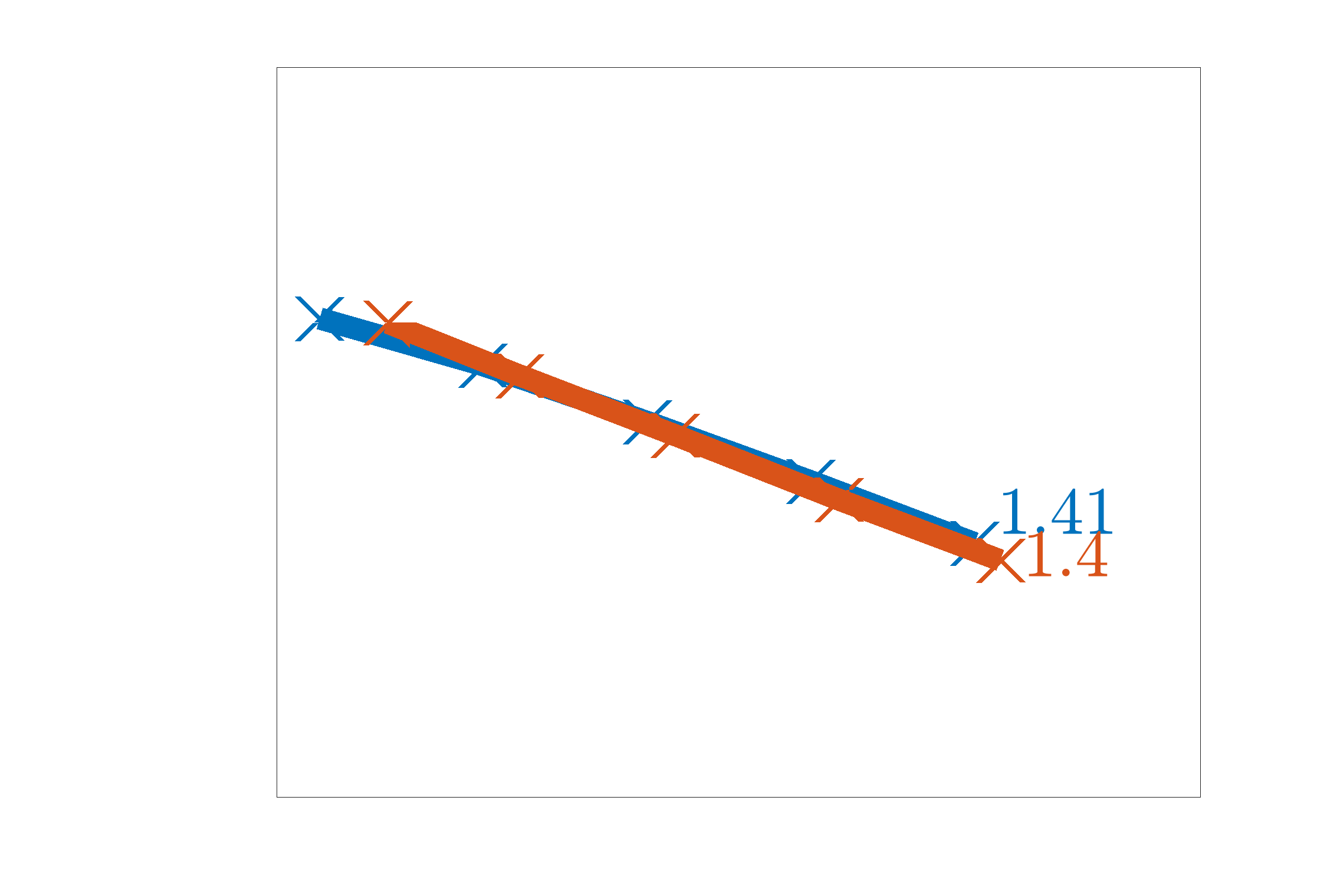}}
\subfigure[$t = 0.3$]{\includegraphics[width=2.1in]{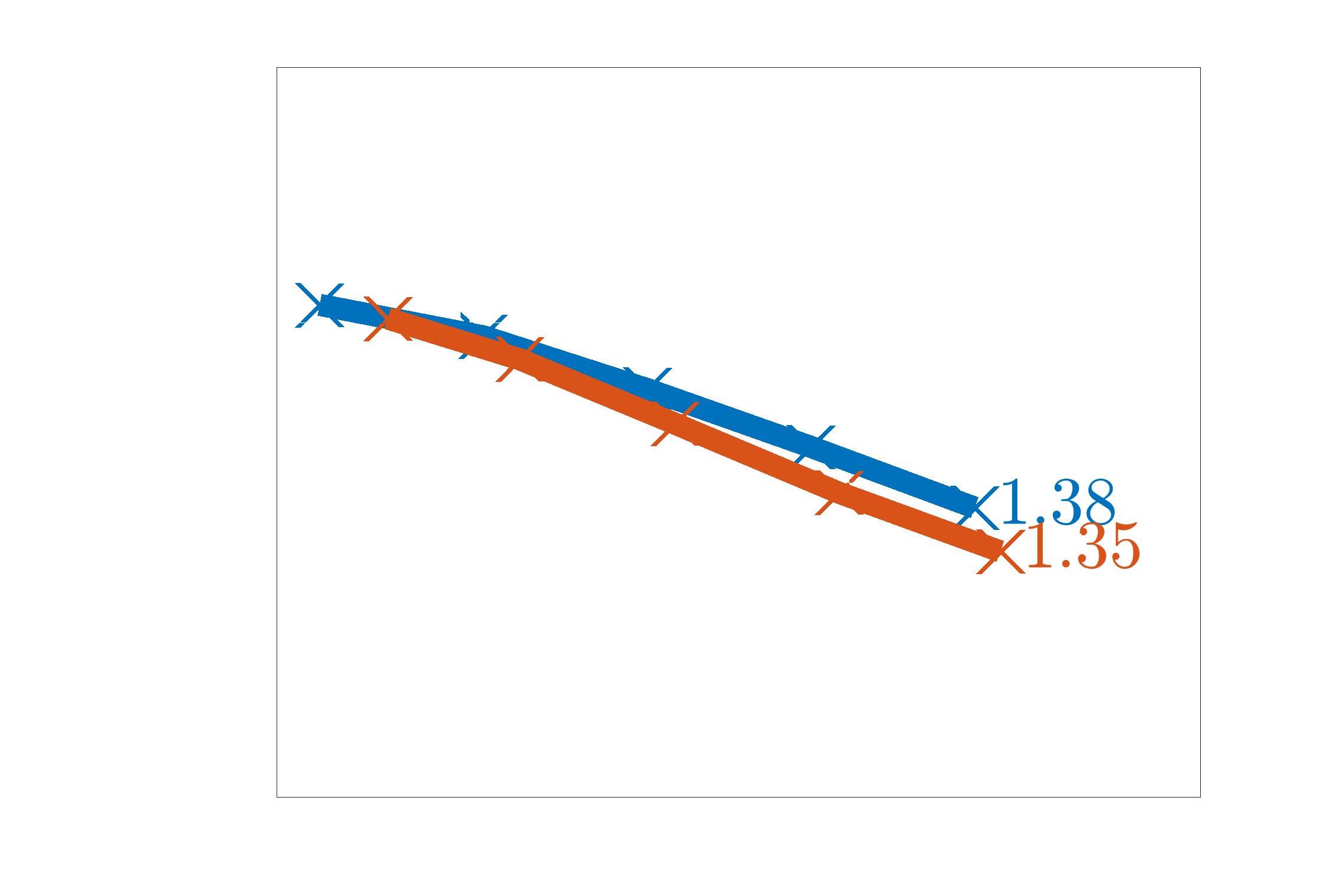}}
\subfigure[$t = 0.35$]{\includegraphics[width=2.1in]{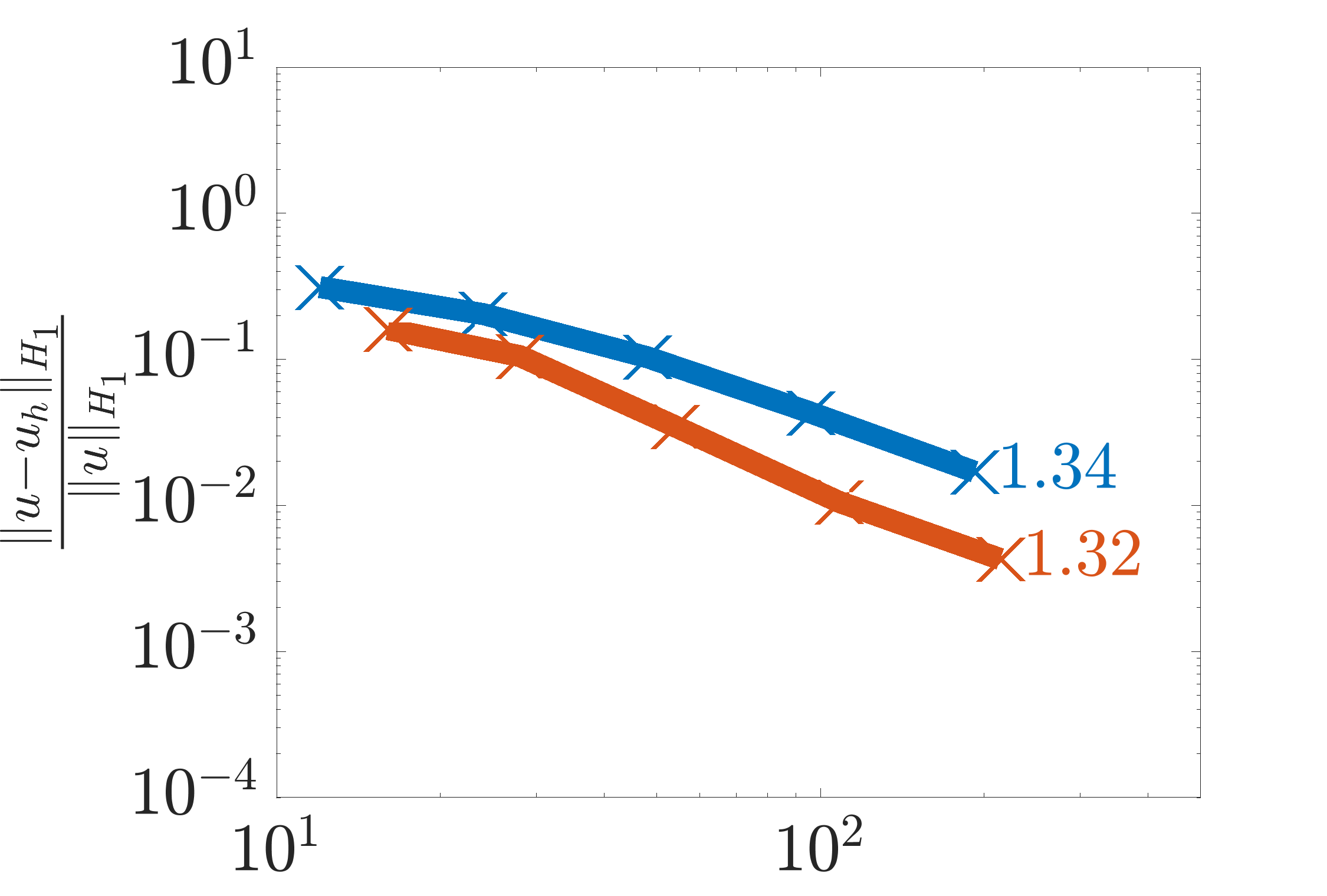}}
\subfigure[$t = 0.5$]{\includegraphics[width=2.1in]{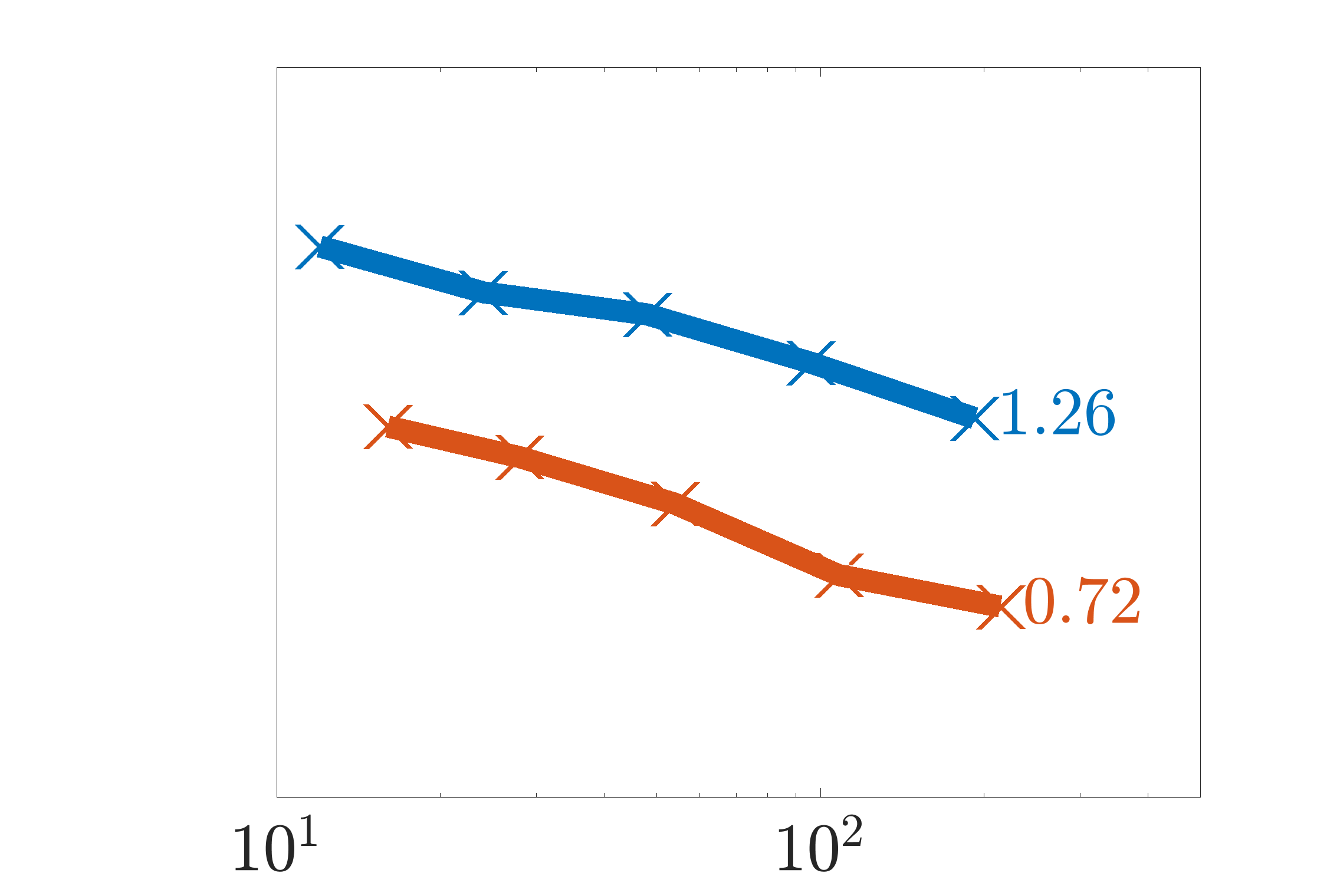}}
\subfigure[$t = 0.75$]{\includegraphics[width=2.1in]{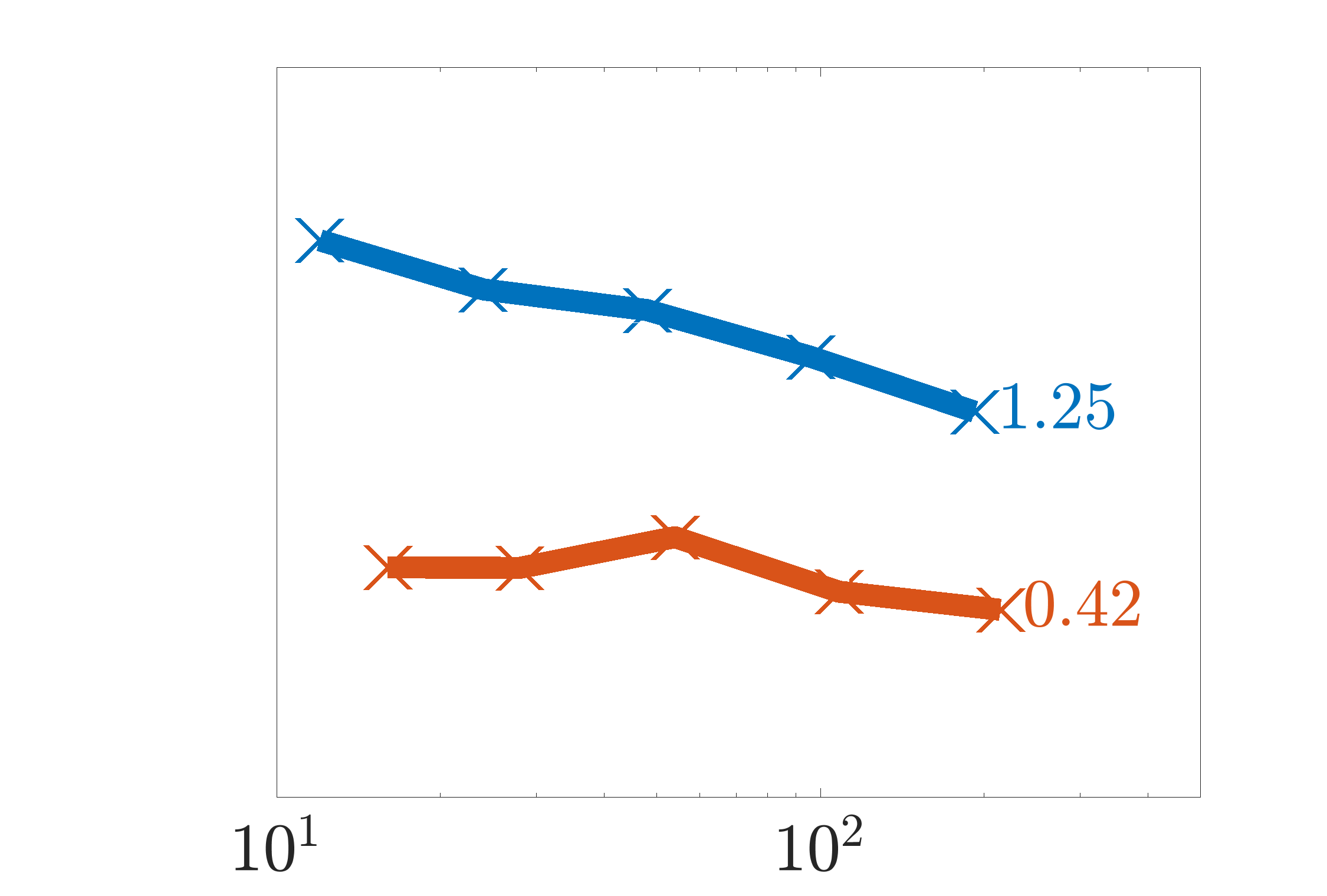}}
\end{subfigmatrix}
\caption{Convergence the relative $H_1$ integral norm for the shock problem with $\nu = \frac{1}{100}$}
\label{fig:Example2_H1_vsdofs_nu1over100}
\end{center}
\end{figure}

\begin{figure}[ht!]
\begin{center}
\begin{subfigmatrix}{6}
\subfigure[$t = 0$]{\includegraphics[width=2.1in]{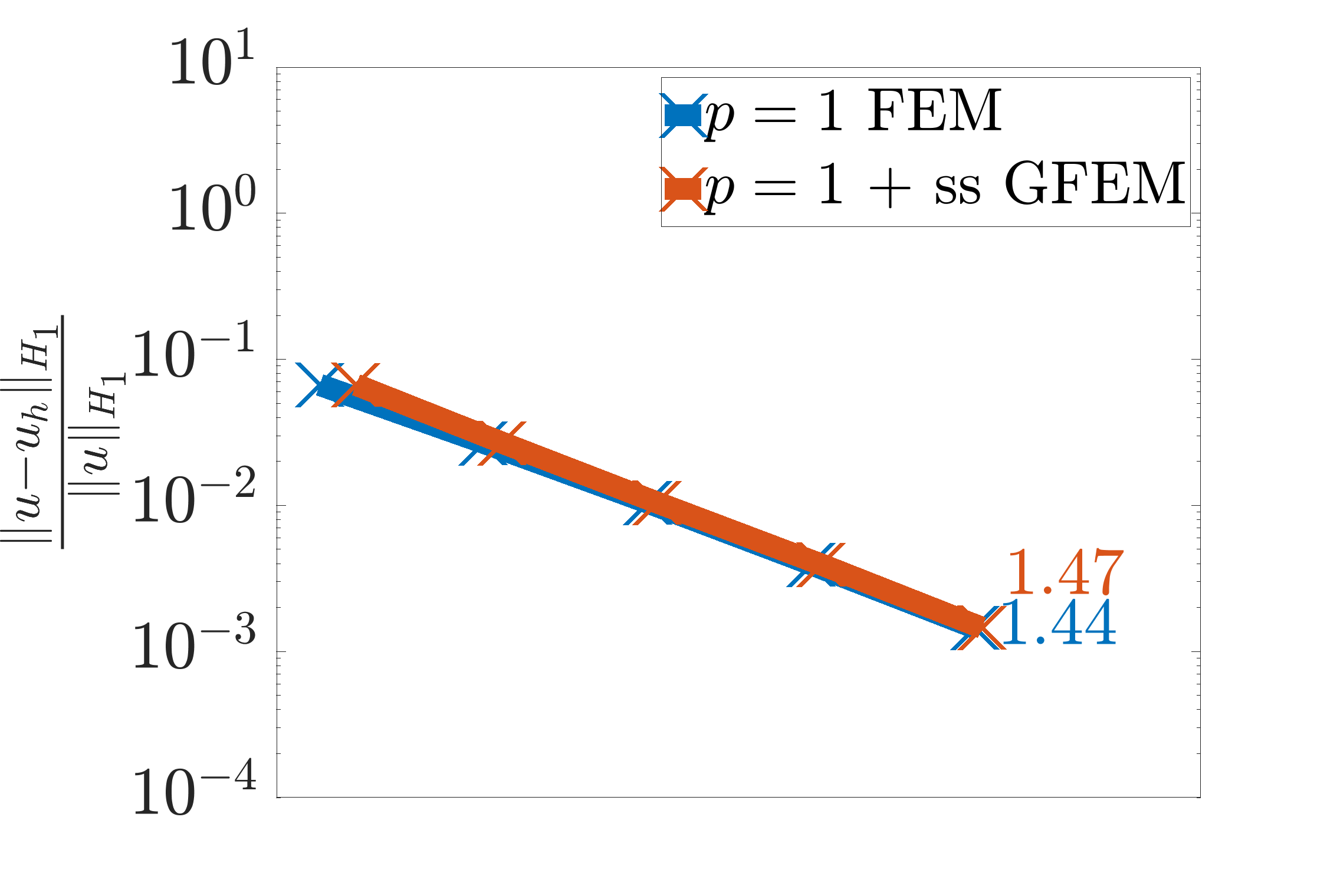}}
\subfigure[$t = 0.25$]{\includegraphics[width=2.1in]{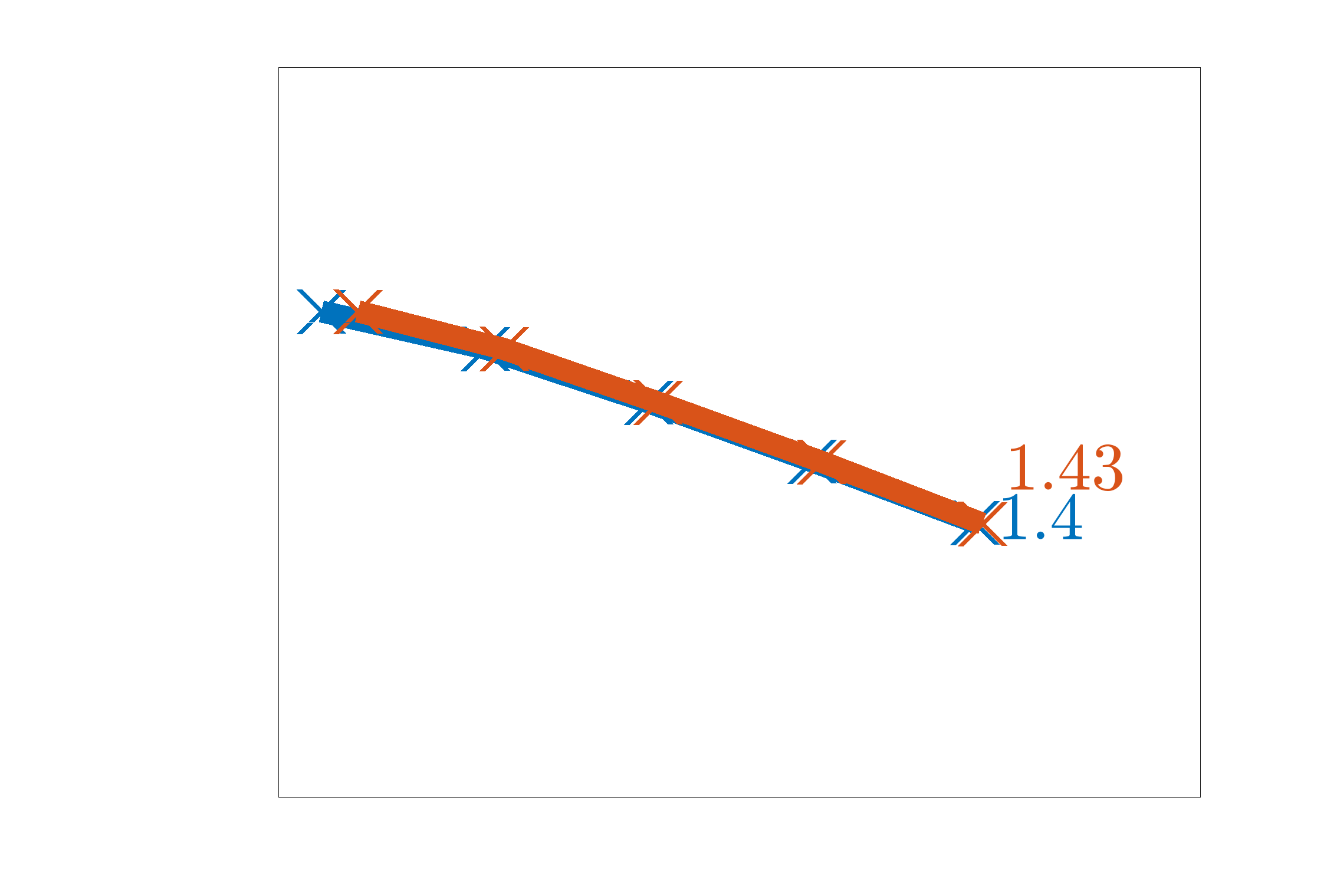}}
\subfigure[$t = 0.3$]{\includegraphics[width=2.1in]{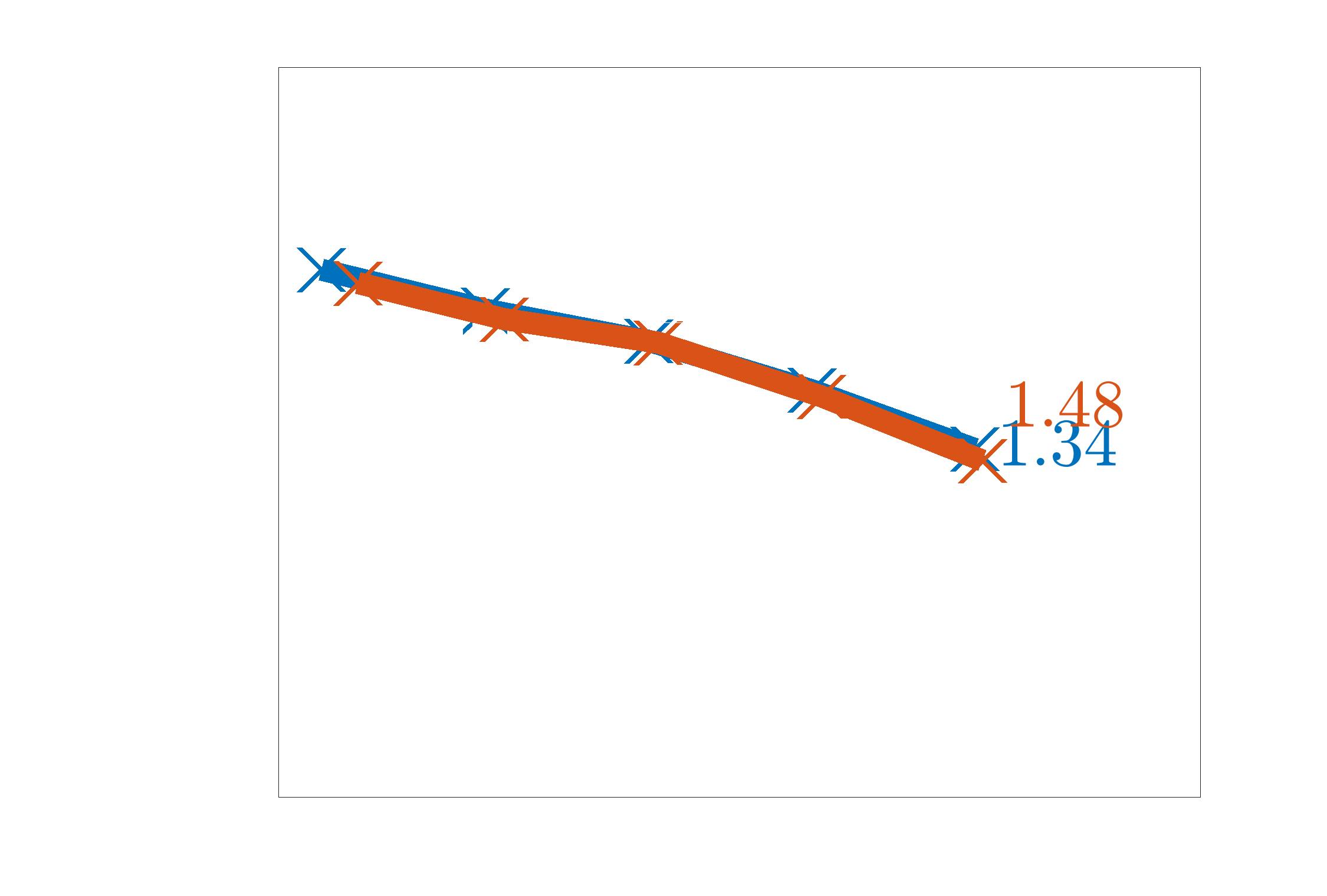}}
\subfigure[$t = 0.35$]{\includegraphics[width=2.1in]{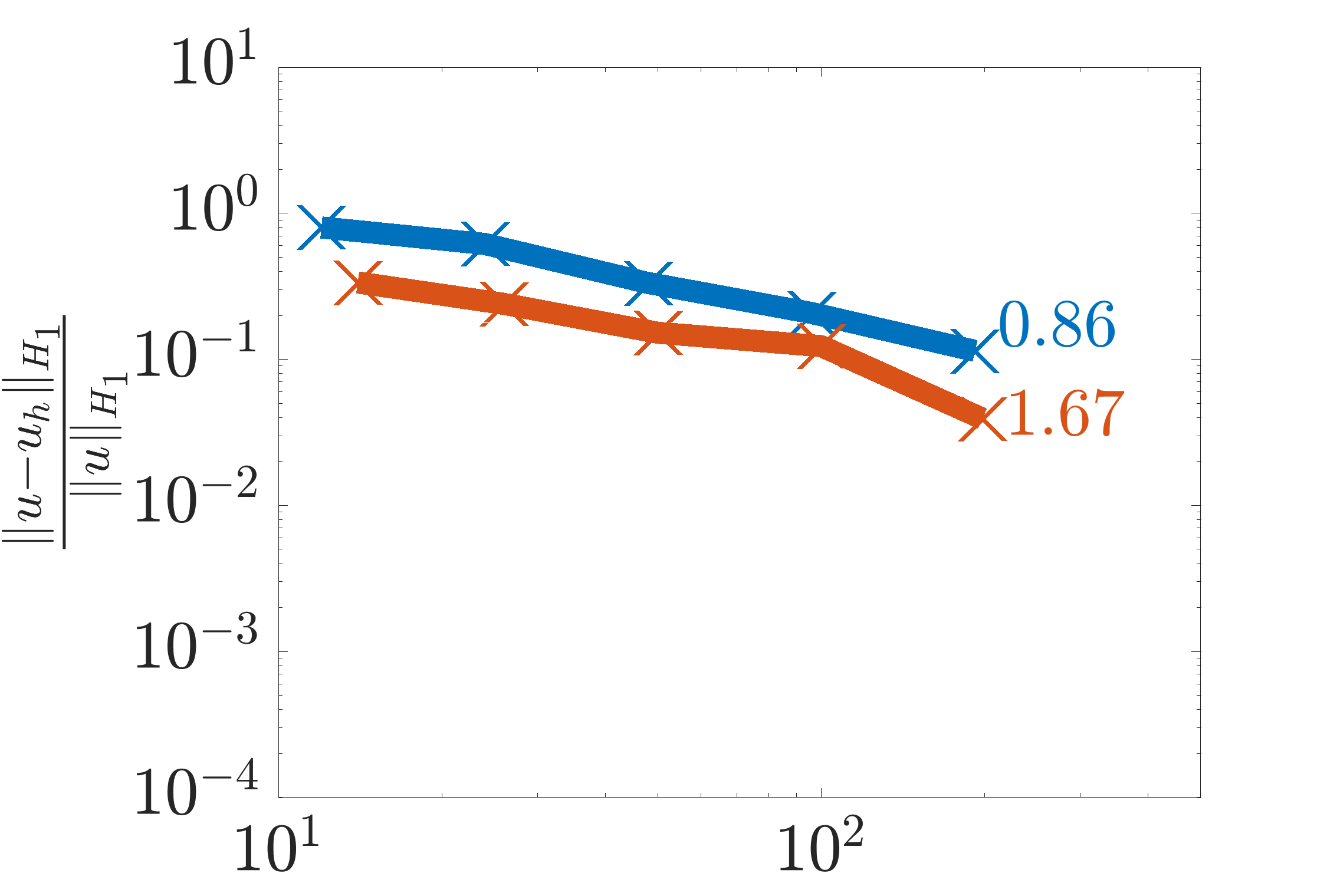}}
\subfigure[$t = 0.5$]{\includegraphics[width=2.1in]{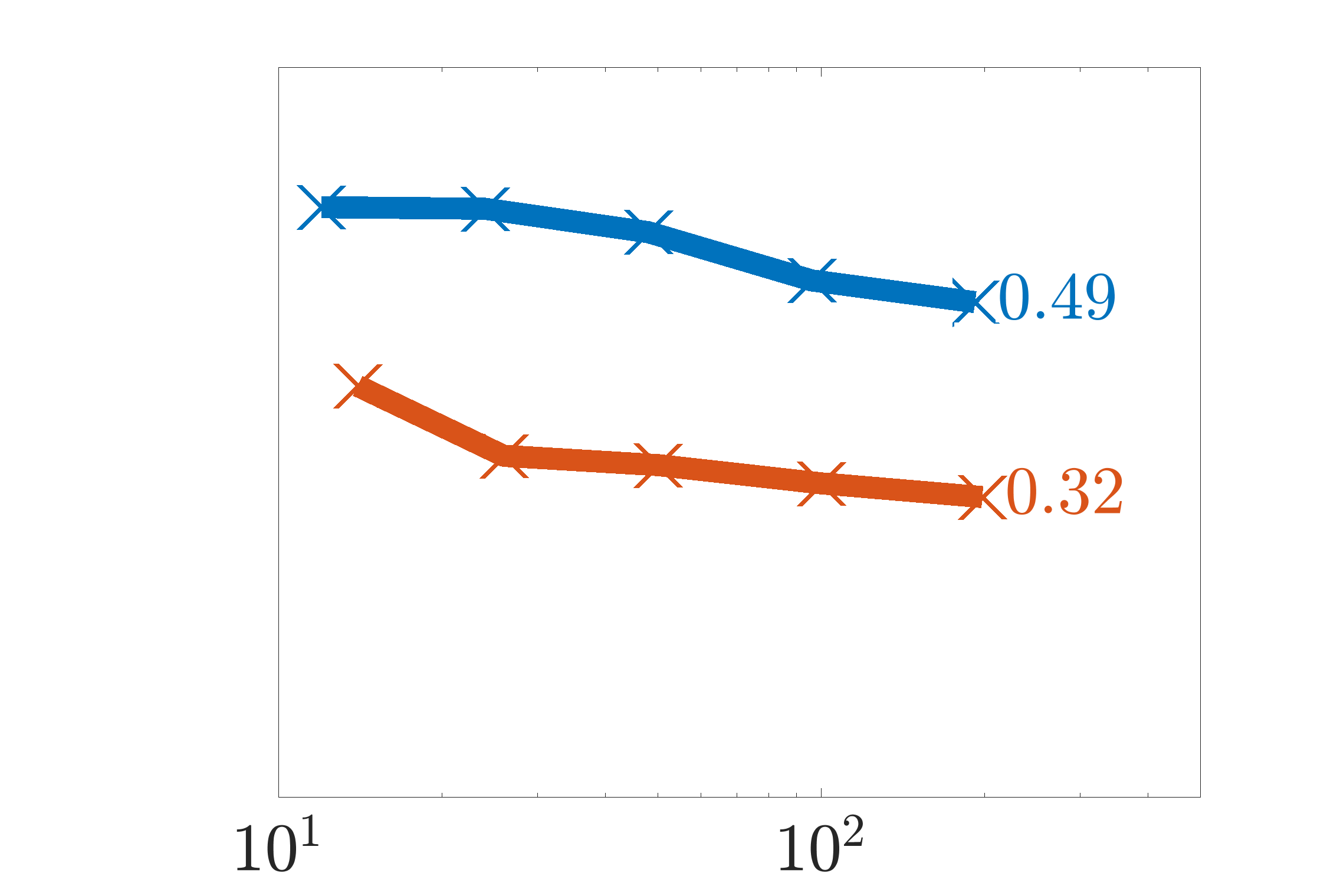}}
\subfigure[$t = 0.75$]{\includegraphics[width=2.1in]{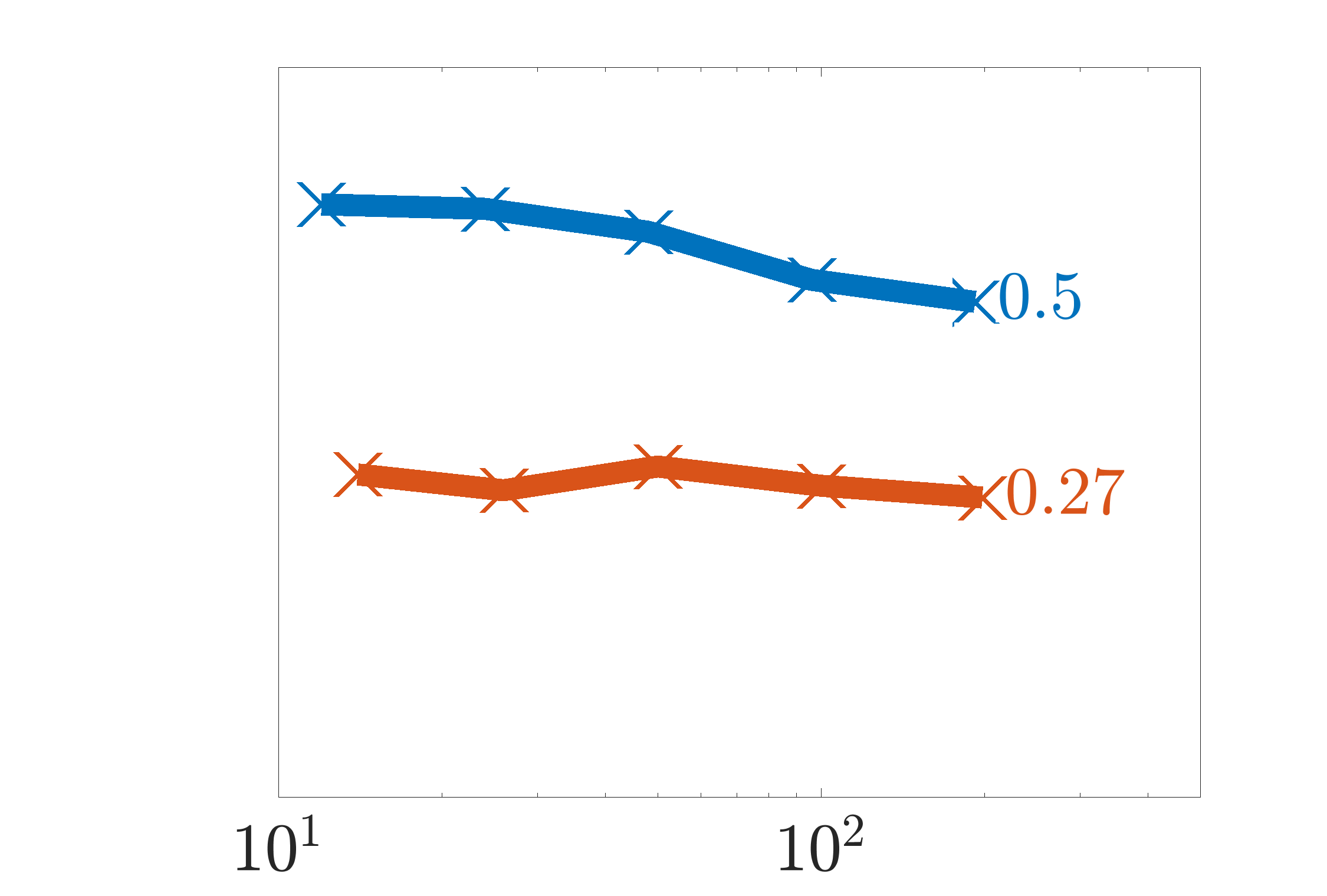}}
\end{subfigmatrix}
\caption{Convergence the relative $H_1$ integral norm for the shock problem with $\nu = \frac{1}{500}$}
\label{fig:Example2_H1_vsdofs_nu1over500}
\end{center}
\end{figure}

\begin{figure}[ht!]
\begin{center}
\begin{subfigmatrix}{6}
\subfigure[$t = 0$]{\includegraphics[width=2.1in]{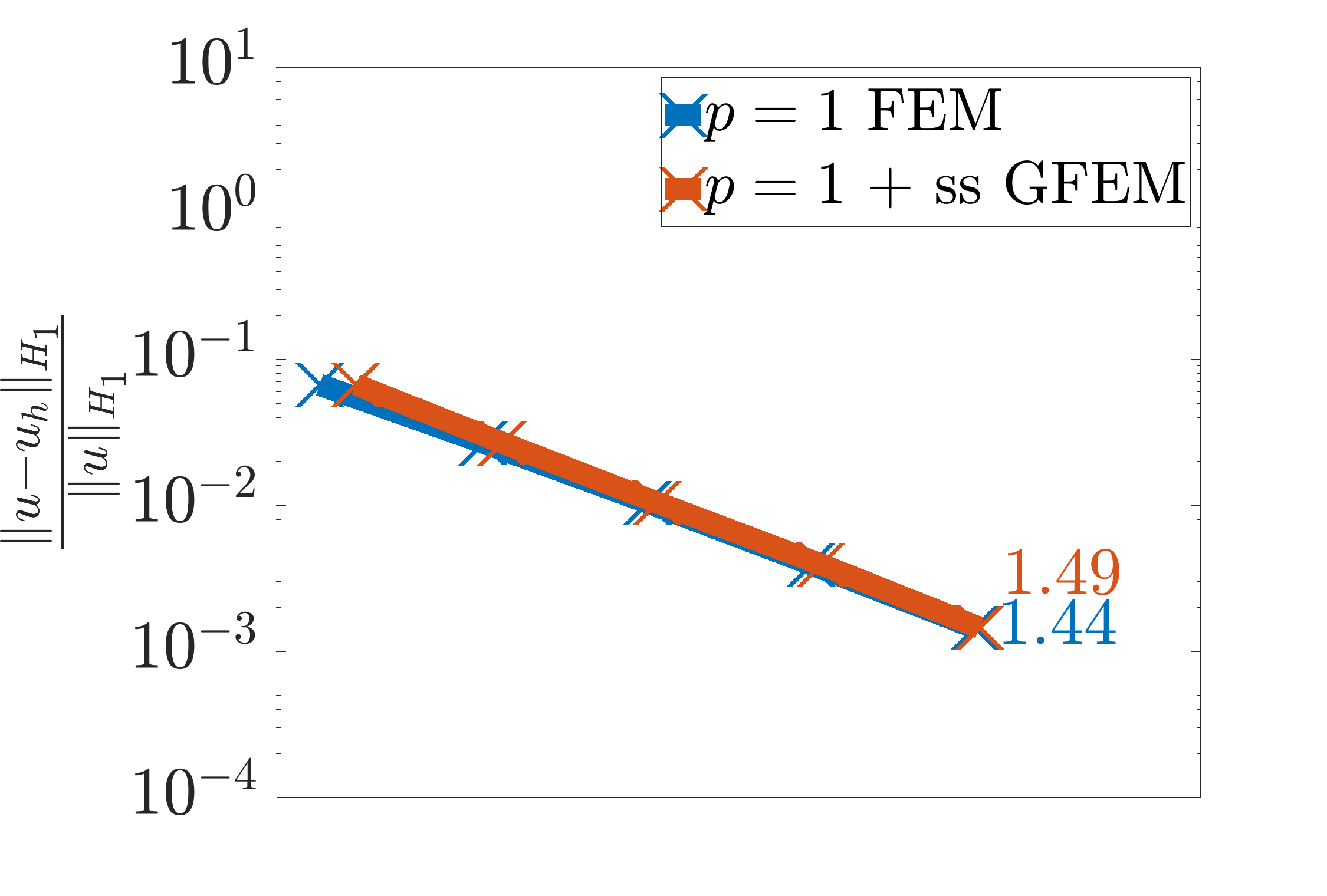}}
\subfigure[$t = 0.25$]{\includegraphics[width=2.1in]{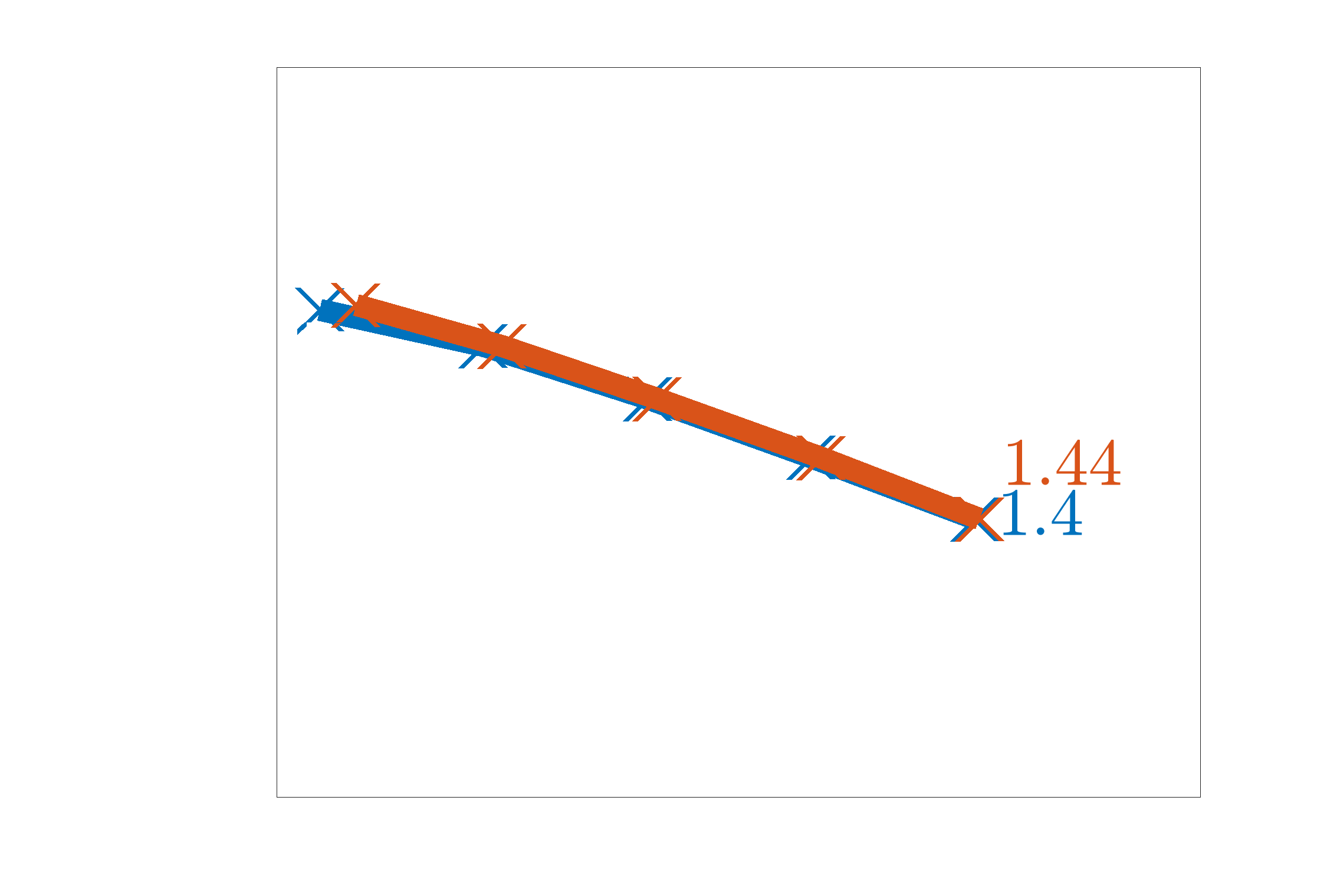}}
\subfigure[$t = 0.3$]{\includegraphics[width=2.1in]{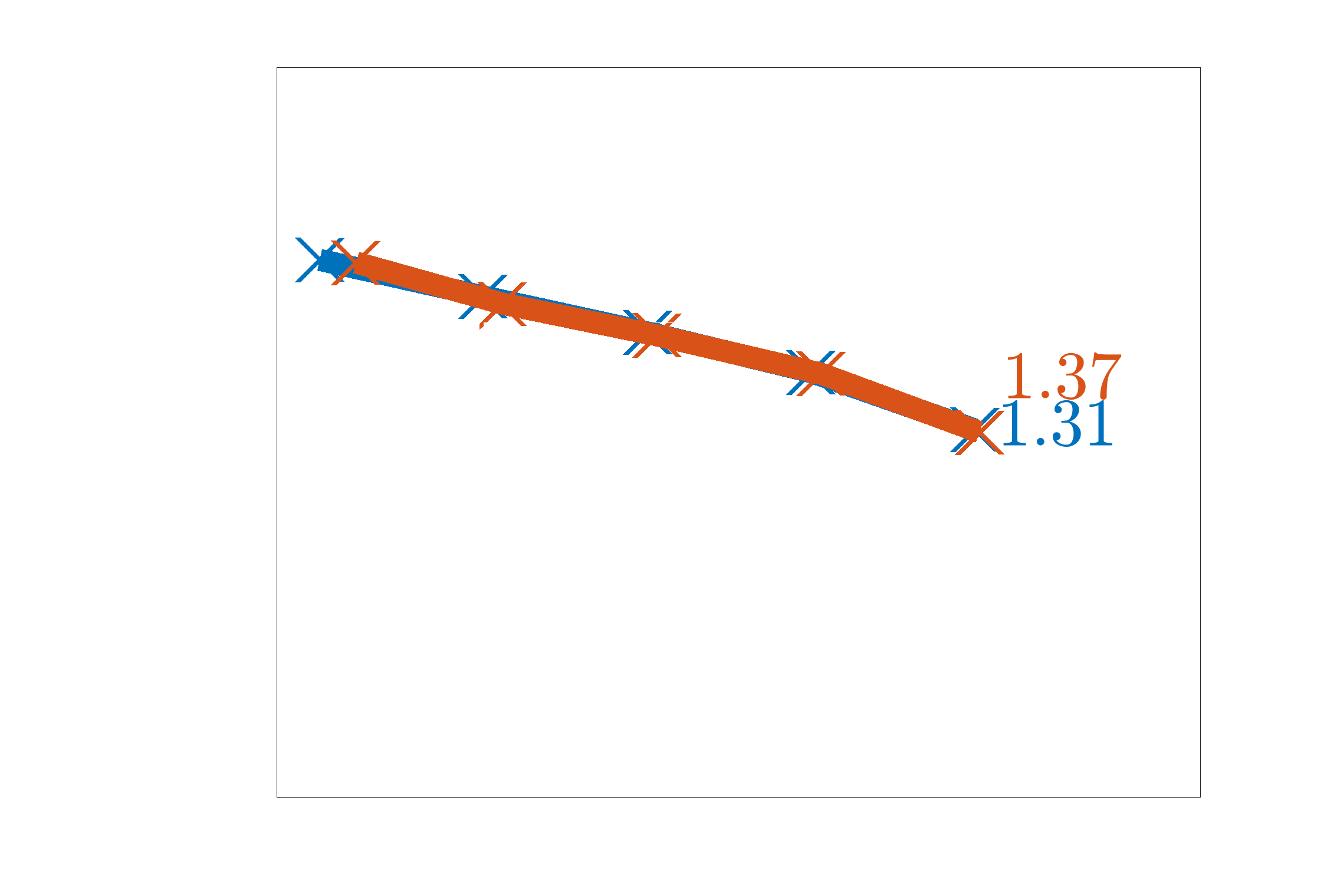}}
\subfigure[$t = 0.35$]{\includegraphics[width=2.1in]{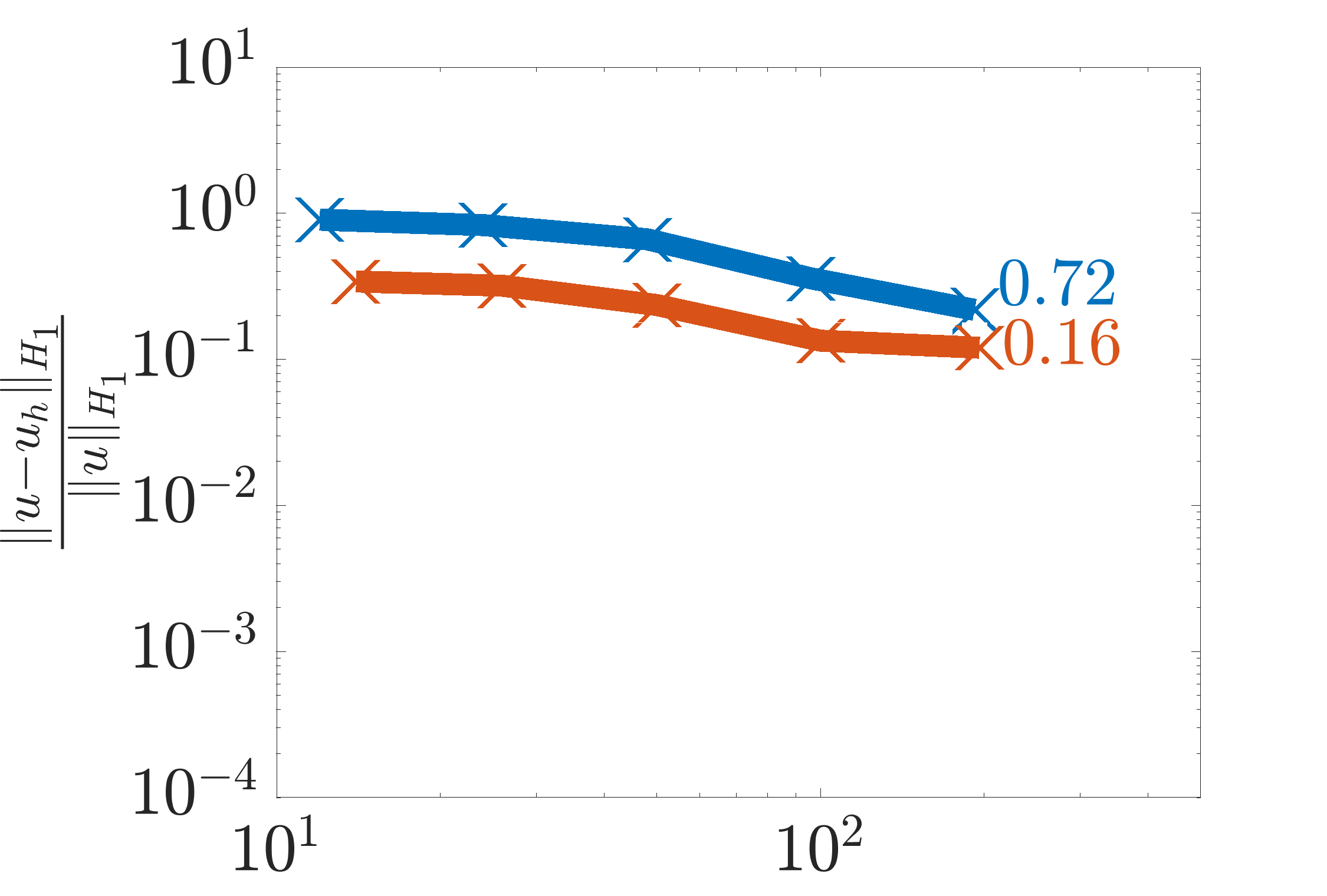}}
\subfigure[$t = 0.5$]{\includegraphics[width=2.1in]{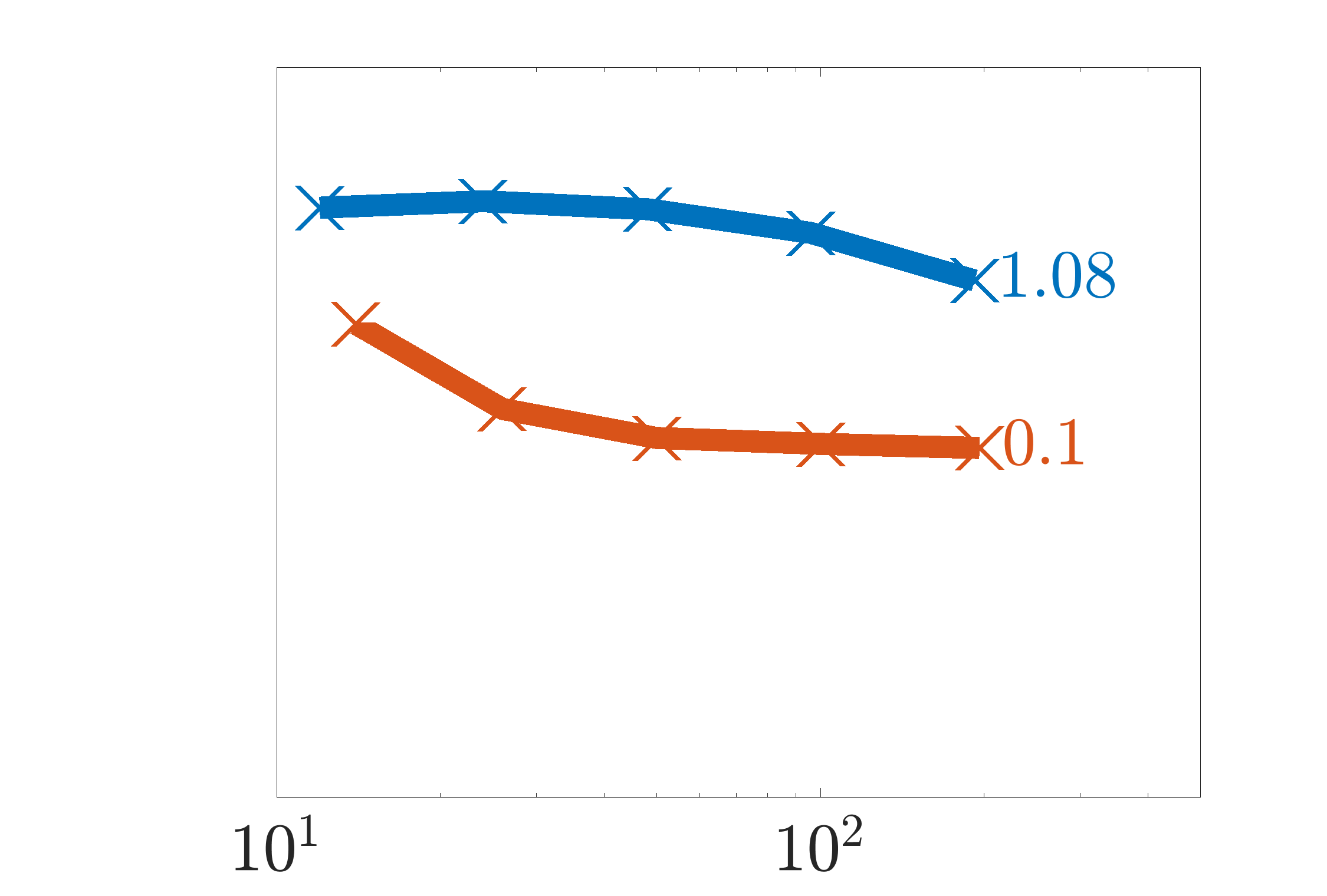}}
\subfigure[$t = 0.75$]{\includegraphics[width=2.1in]{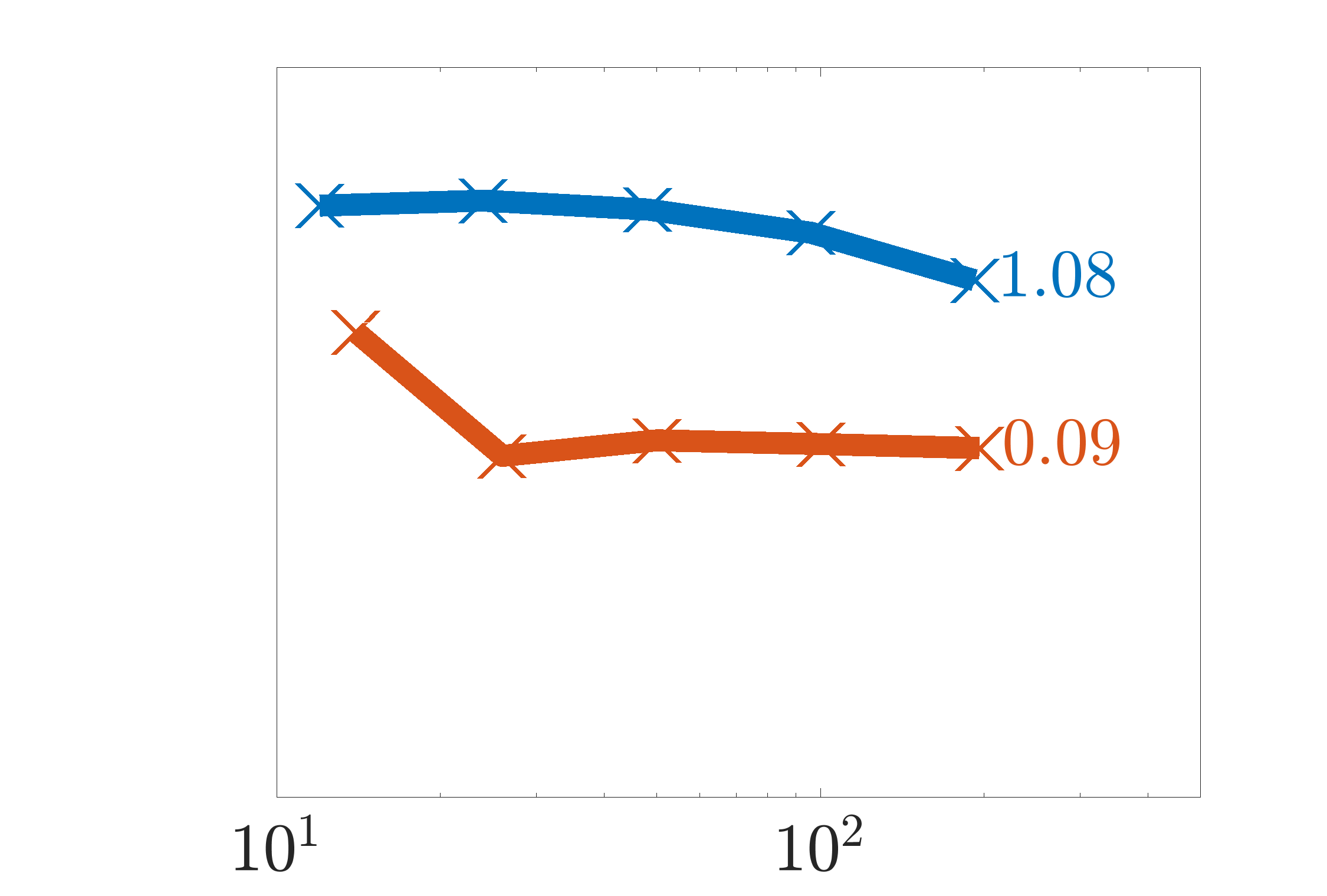}}
\end{subfigmatrix}
\caption{Convergence the relative $H_1$ integral norm for the shock problem with $\nu = \frac{1}{1000}$}
\label{fig:Example2_H1_vsdofs_nu1over1000}
\end{center}
\end{figure}

Convergence plots in the relative $L_2$ norm versus total degrees of freedom at times $t = [0, 0.25, 0.3, 0.35, 0.5, 0.75]$ for $\nu = \Big[\frac{1}{50}, \frac{1}{100}, \frac{1}{500}, \frac{1}{1000} \Big]$ are shown in Figs. \ref{fig:Example2_L2_vsdofs_nu1over50} - \ref{fig:Example2_L2_vsdofs_nu1over1000}, respectively. Similarly, convergence plots for the relative $H_1$ norm are shown in Figs. \ref{fig:Example2_H1_vsdofs_nu1over50} -  \ref{fig:Example2_H1_vsdofs_nu1over1000}. Note since the steady state solution is closely correlated with the transient solution post-shock, the $p = 1$ + ss GFEM solutions are on the same order of numerical precision as the 5000-element, $p = 1$ FEM reference solution. This is observed in the convergence plots as the $p = 1$ + ss GFEM solution convergence begins to flatten out after around $t = \frac{1}{\pi}$. Before shock formation, both FEM and GFEM converge similarly in all norms studied. Additionally, $p = 1$ FEM performs slightly better since the steady state as an enrichment is not correlated with the initial transient solution. However, around $t = 0.25$ and persisting through $t = 0.35$, the formation of the shock results in a shift in the error, as well as sub-optimal convergence in both the FEM and GFEM solutions. After the shock has mostly formed around $t = 0.35$, error levels in the FEM solutions for the 95-element solution are larger than 30\% in the $H_1$ norm. Error levels in the GFEM solutions at the same degrees of freedom are less than 2\%. However, the effect of the oscillations which arise in the GFEM solution at earlier time steps greatly affects the convergence rate in the GFEM solution in the $H_1$ norm, and sub-optimal convergence is observed over all grids studied.  

\subsubsection{Capturing intermediate solution features}
\begin{figure}[ht!]
\centering{\includegraphics[width=3in]{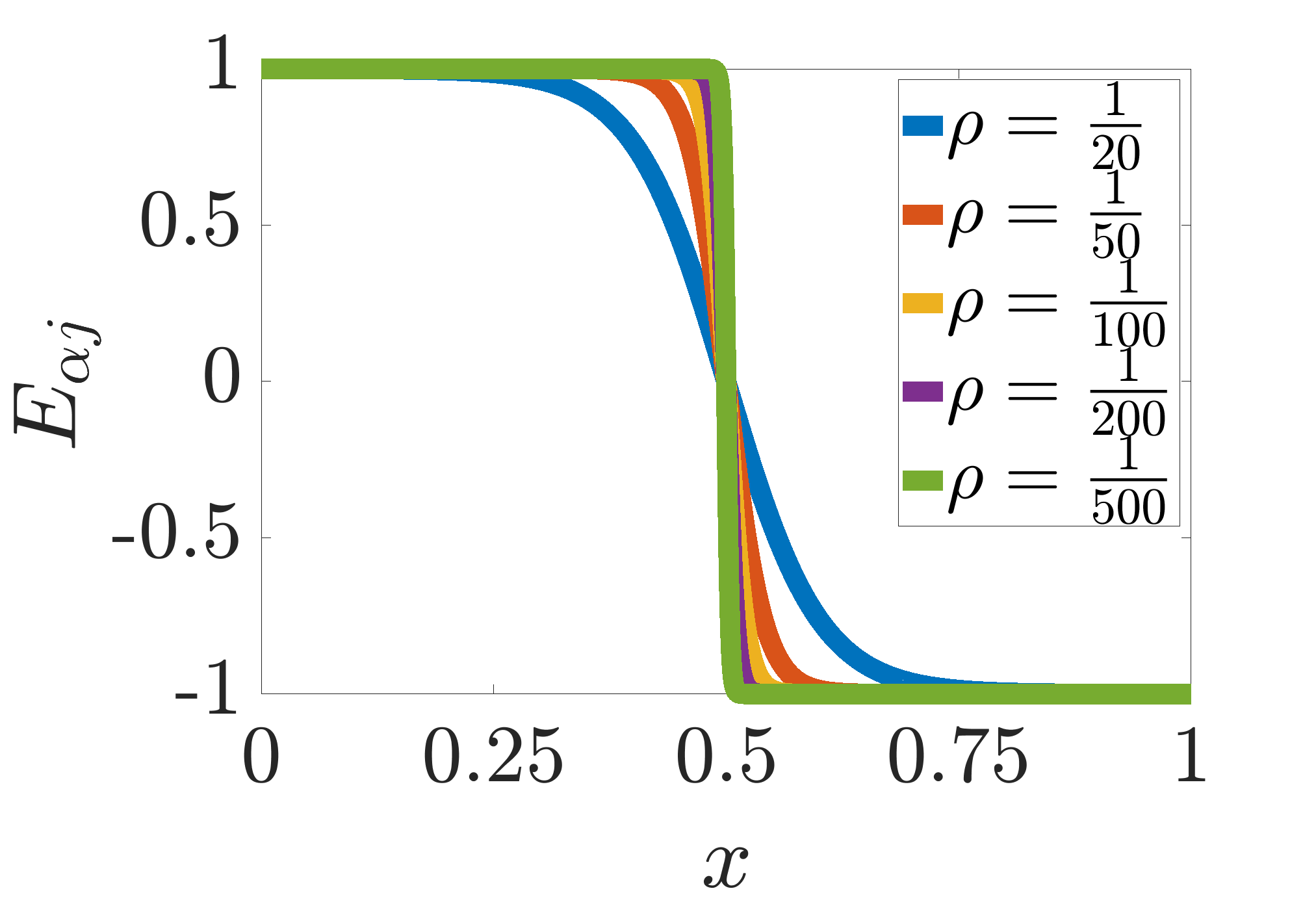}}
\caption{Set of shock enrichments $\Big( E_{\alpha j} = \tanh{ \Big[\frac{1}{2 \rho} \Big( \frac{1}{2} - x\Big) \Big]} \Big)$ for various $\rho$}
\label{fig:Example2_shock_enrichments}
\end{figure}

\begin{figure}[ht!]
\begin{center}
\begin{subfigmatrix}{2}
\subfigure[Relative $L_2$ integral norm]{\includegraphics[width=3in]{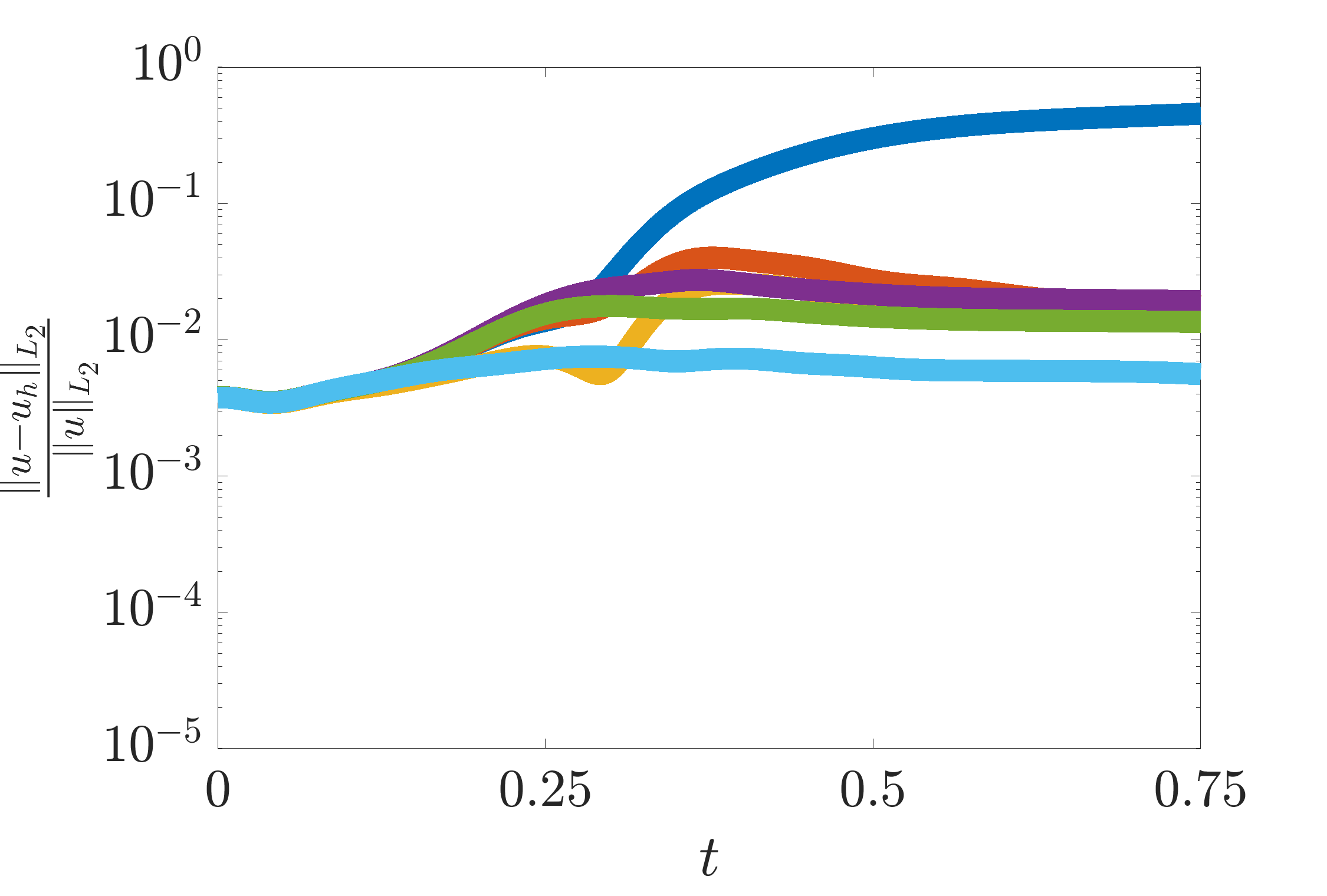}}
\subfigure[Relative $H_1$ integral norm]{\includegraphics[width=3in]{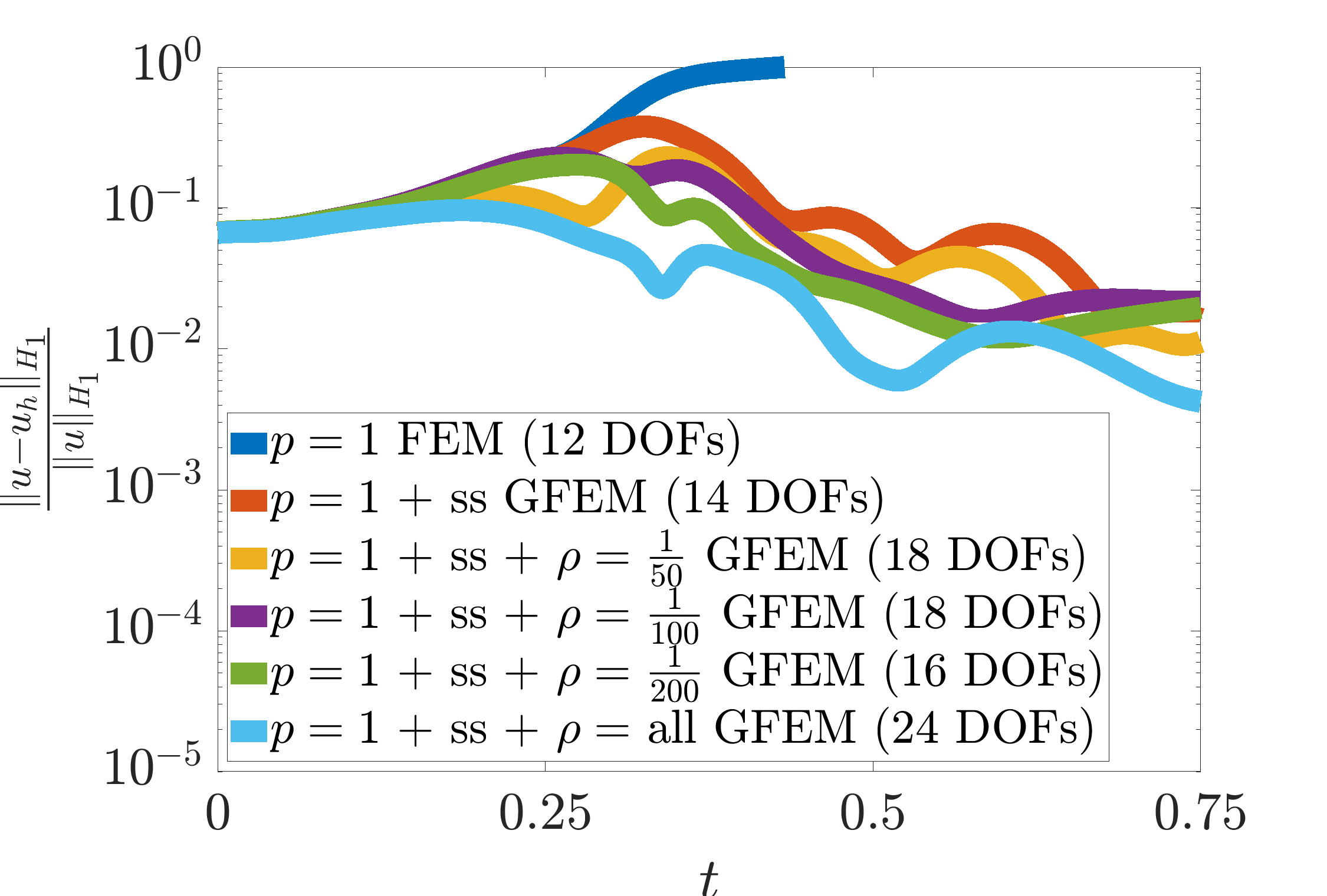}}
\end{subfigmatrix}
\caption{Relative $L_2$ and $H_1$ integral norms versus time for 11-element FEM and GFEM solutions to the shock formation problem when $\nu = \frac{1}{500}$}
\label{fig:Example2_11element_L2H1vstime_nu1over500_intermediatescales}
\end{center}
\end{figure}

\begin{figure}[ht!]
\begin{center}
\begin{subfigmatrix}{6}
\subfigure[$t = 0$]{\includegraphics[width=2.1in]{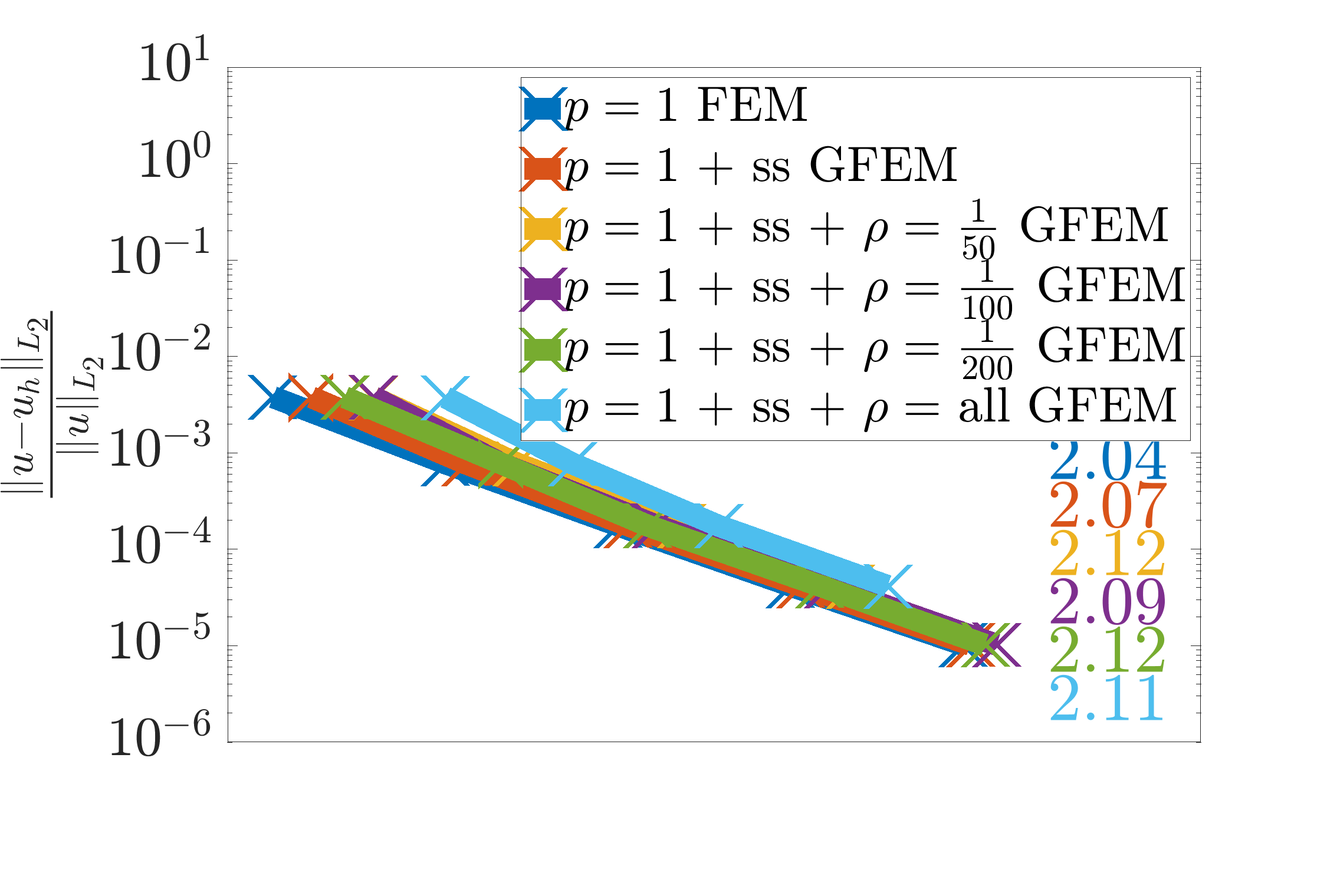}}
\subfigure[$t = 0.25$]{\includegraphics[width=2.1in]{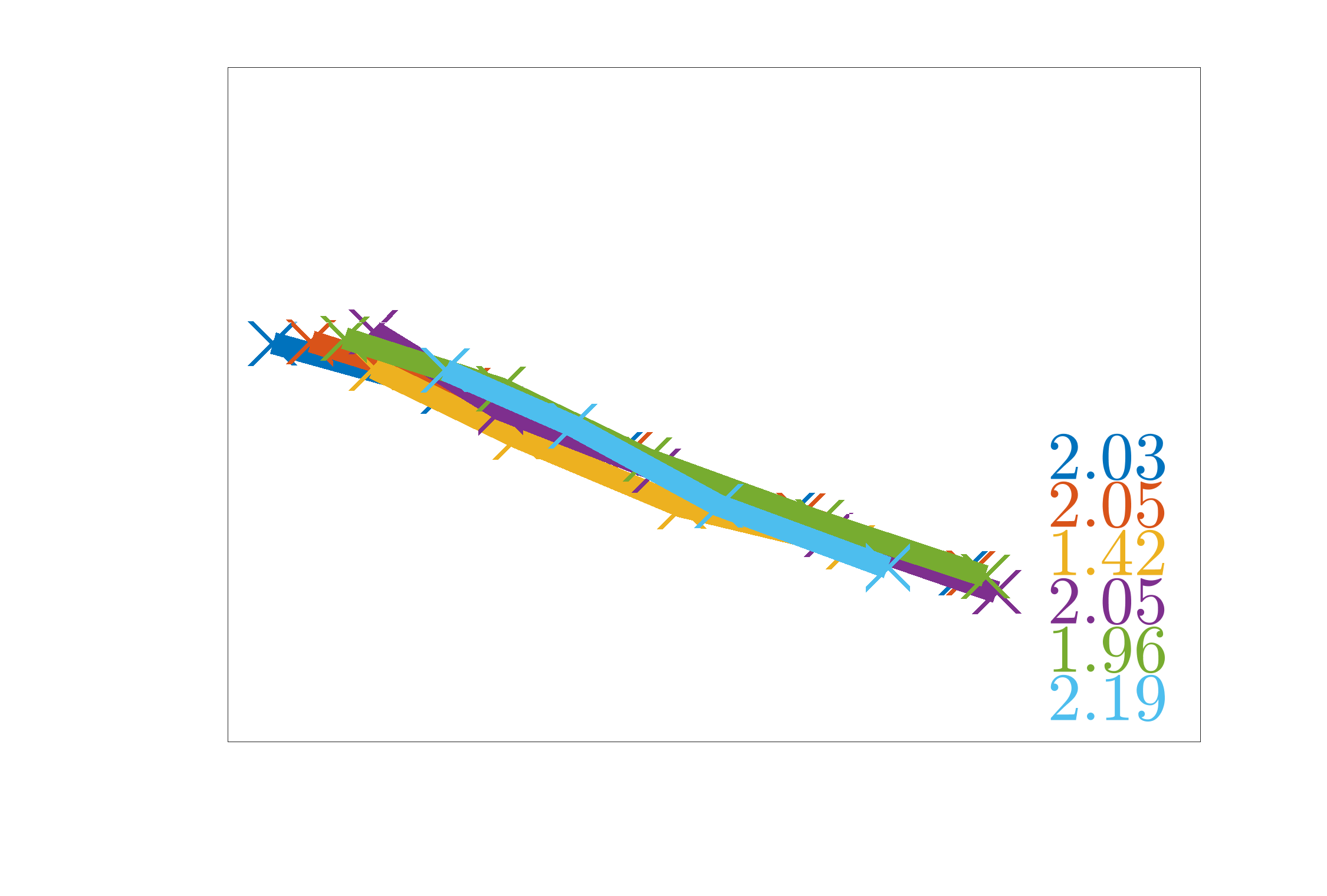}}
\subfigure[$t = 0.3$]{\includegraphics[width=2.1in]{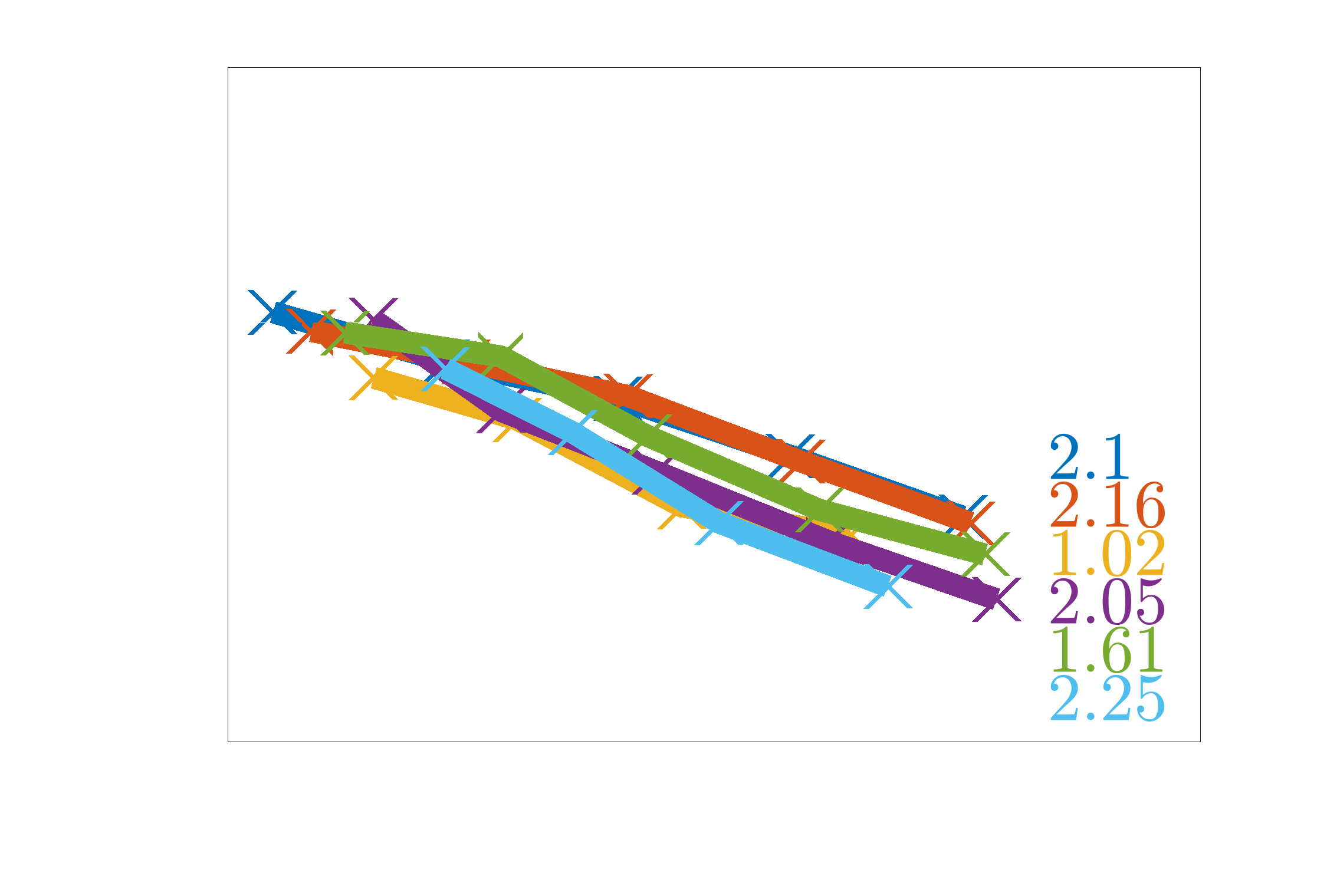}}
\subfigure[$t = 0.35$]{\includegraphics[width=2.1in]{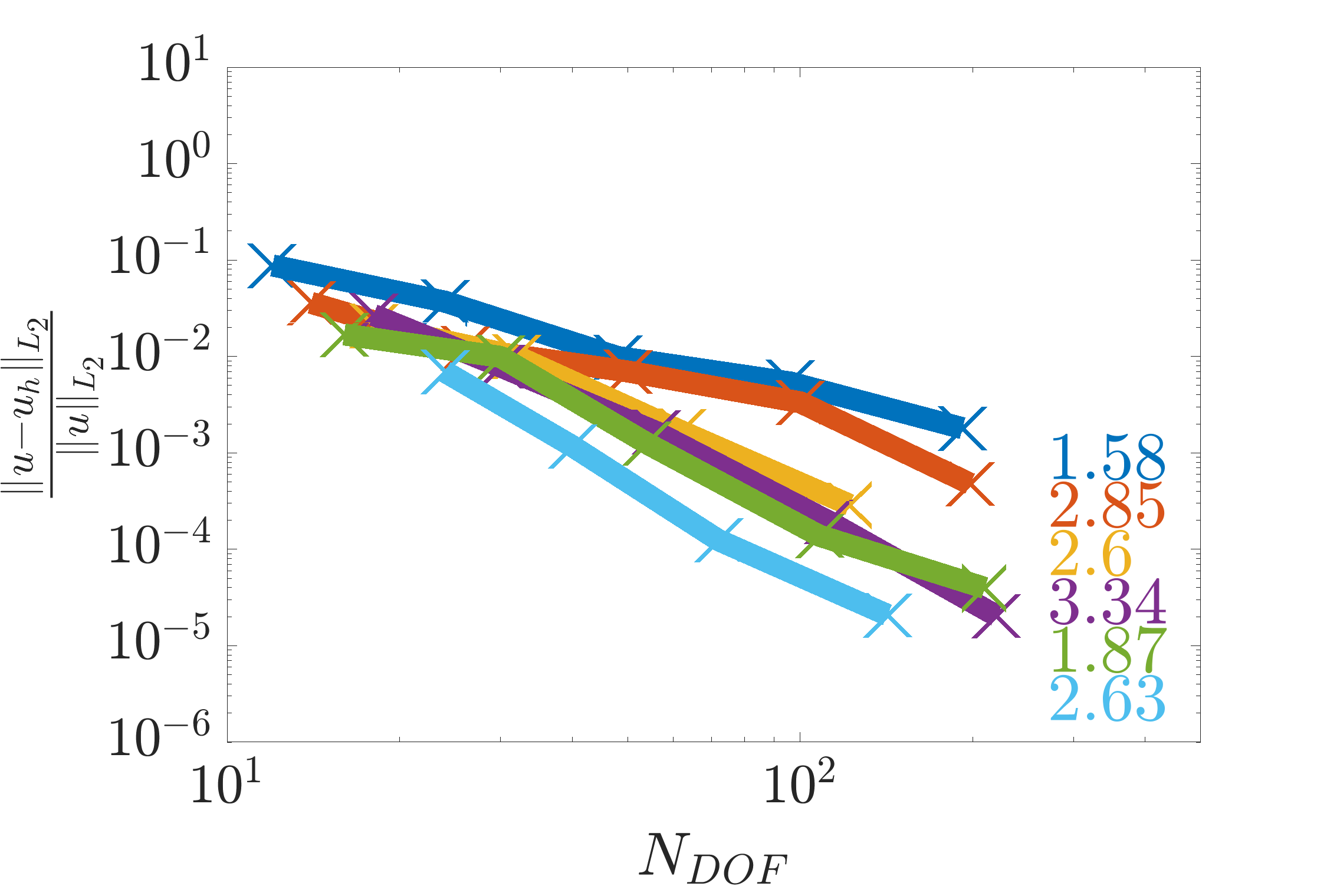}}
\subfigure[$t = 0.5$]{\includegraphics[width=2.1in]{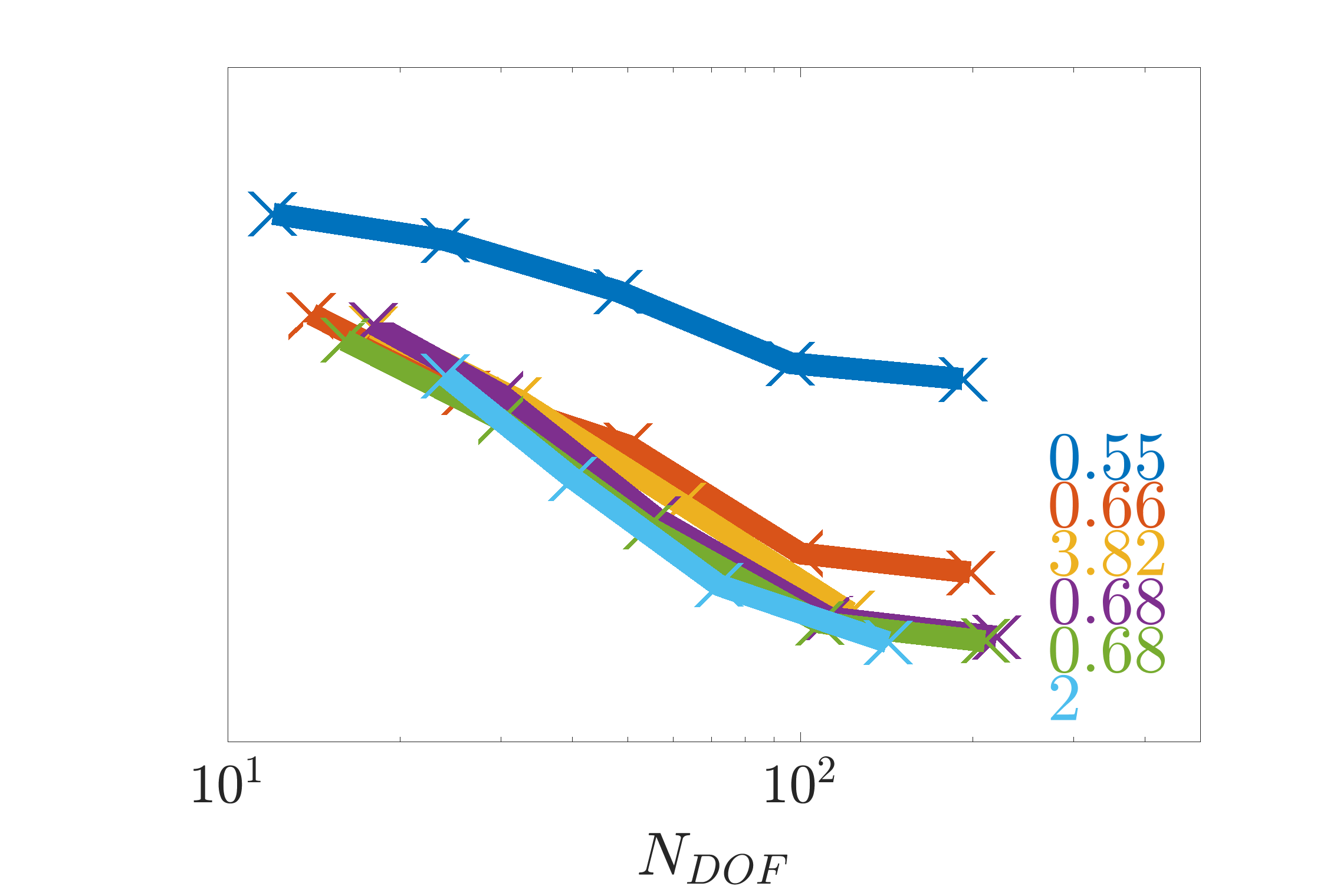}}
\subfigure[$t = 0.75$]{\includegraphics[width=2.1in]{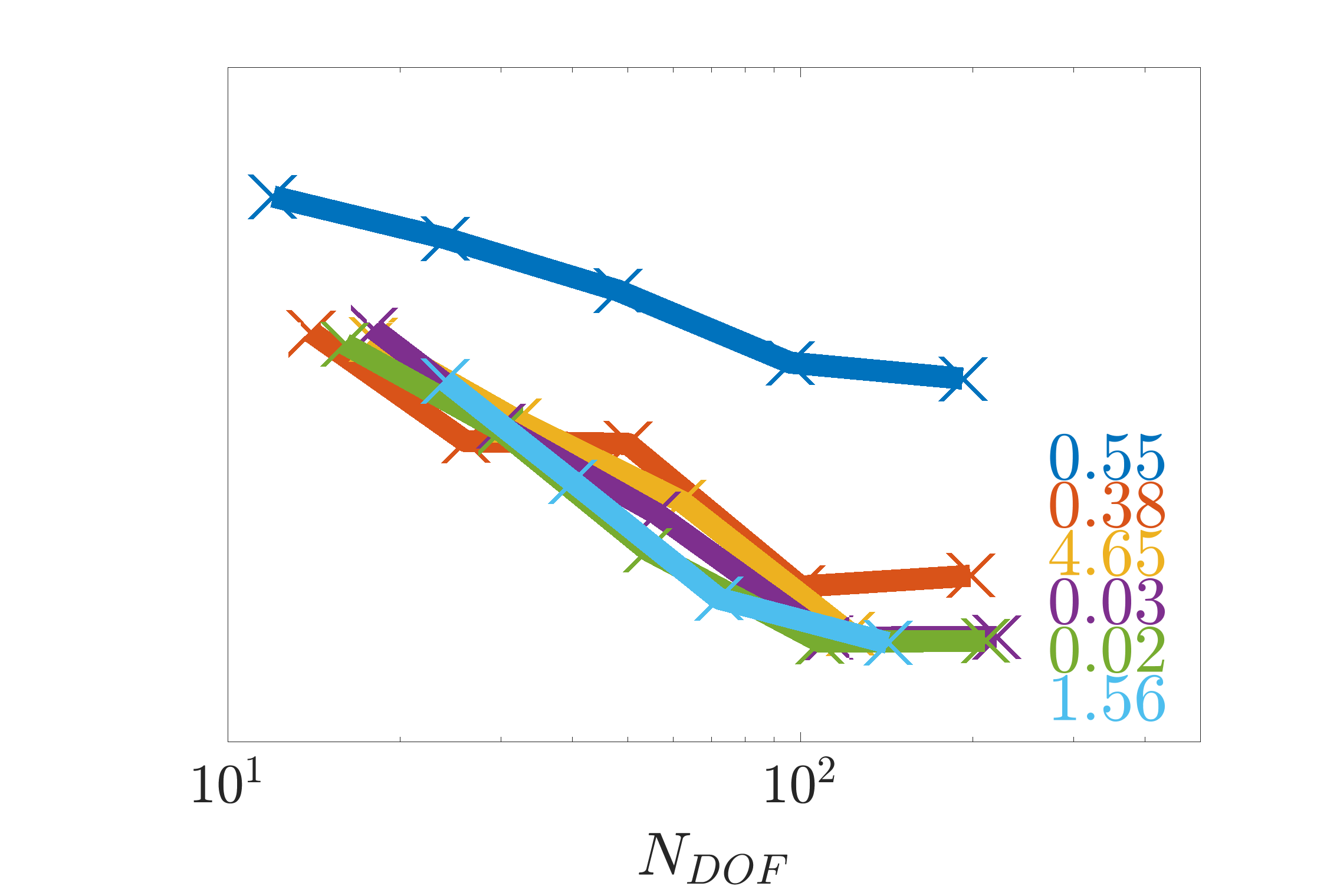}}
\end{subfigmatrix}
\caption{Convergence the relative $L_2$ integral norm for the shock problem with $\nu = \frac{1}{500}$}
\label{fig:Example2_L2_vsdofs_nu1over500_withrho}
\end{center}
\end{figure}

\begin{figure}[ht!]
\begin{center}
\begin{subfigmatrix}{6}
\subfigure[$t = 0$]{\includegraphics[width=2.1in]{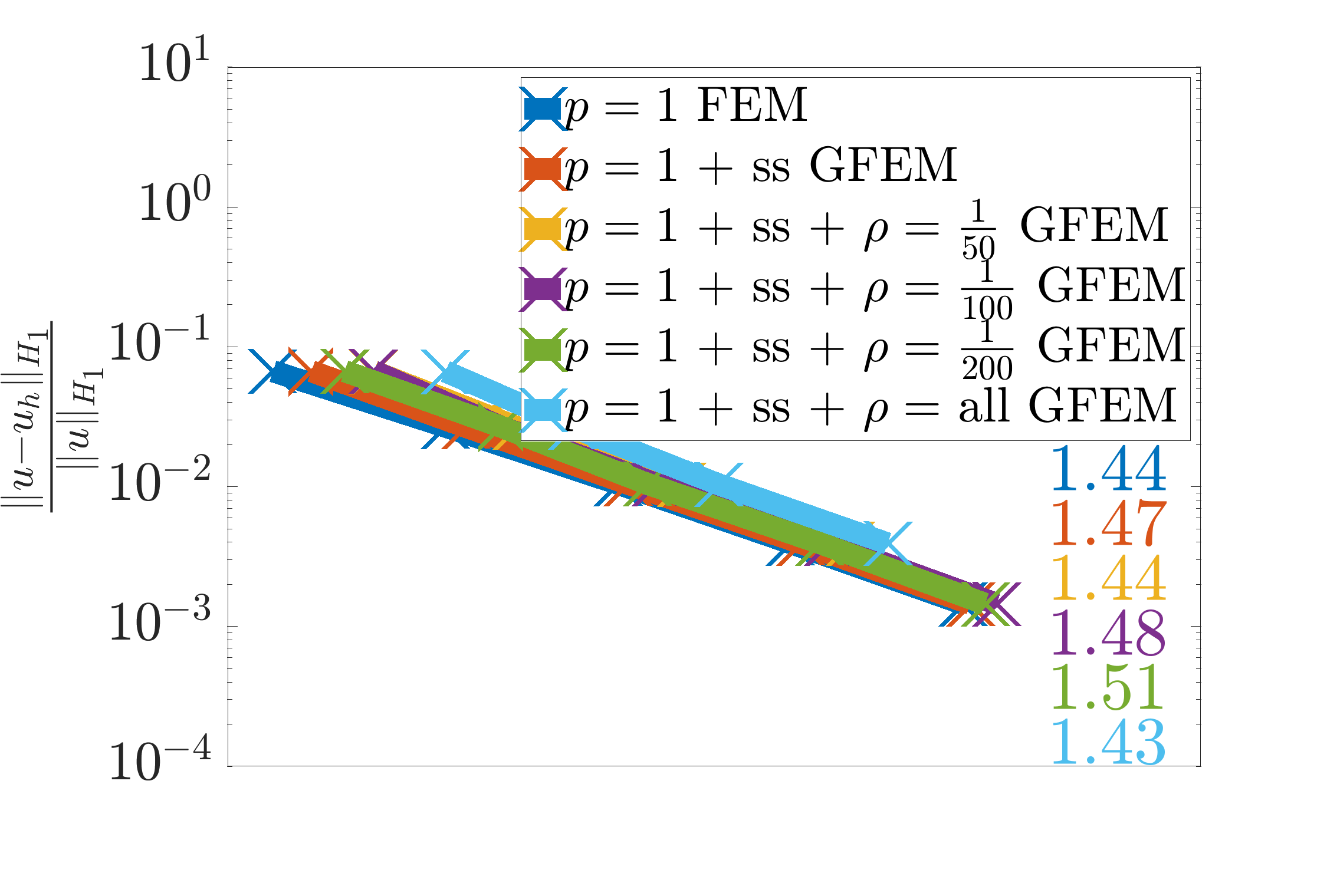}}
\subfigure[$t = 0.25$]{\includegraphics[width=2.1in]{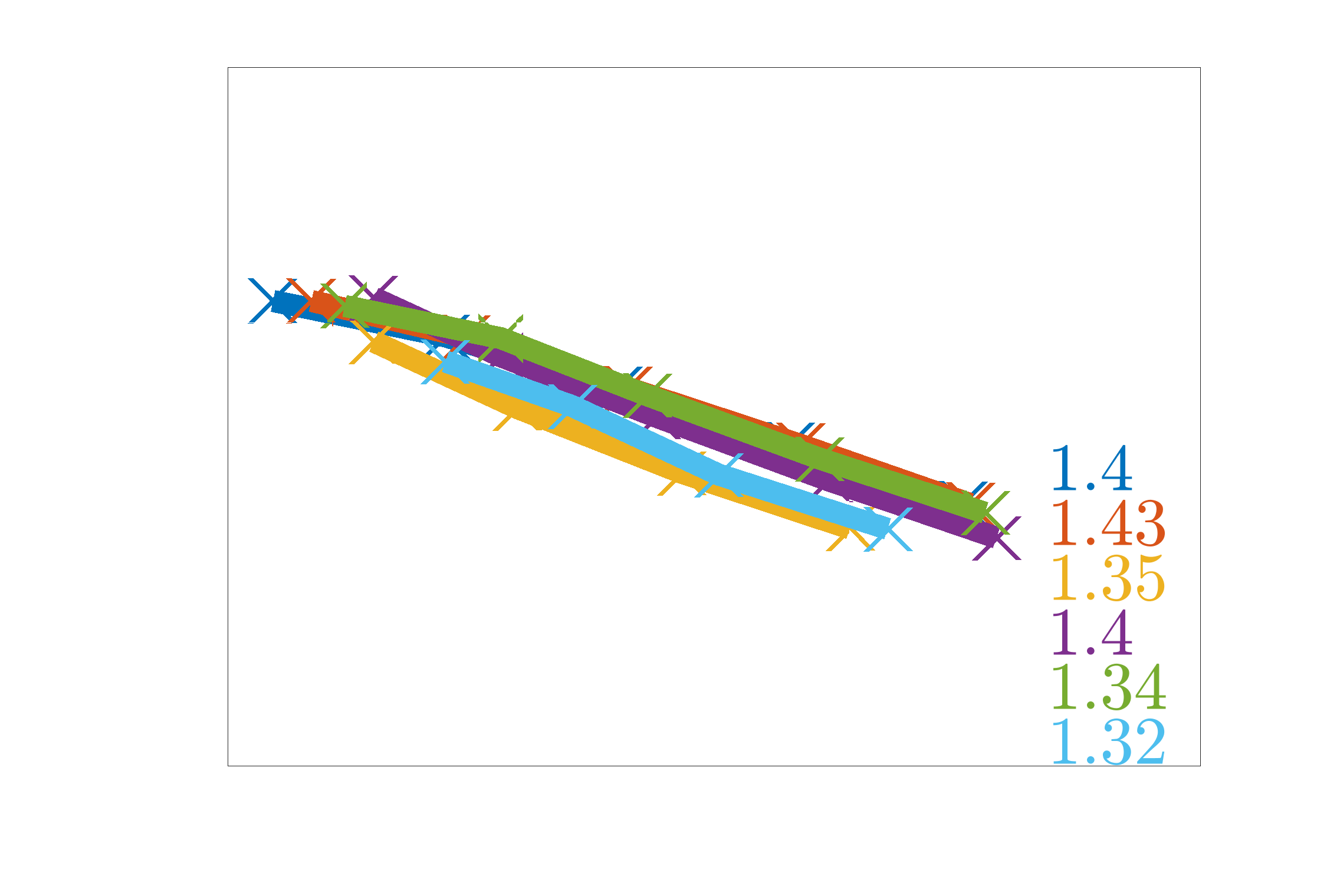}}
\subfigure[$t = 0.3$]{\includegraphics[width=2.1in]{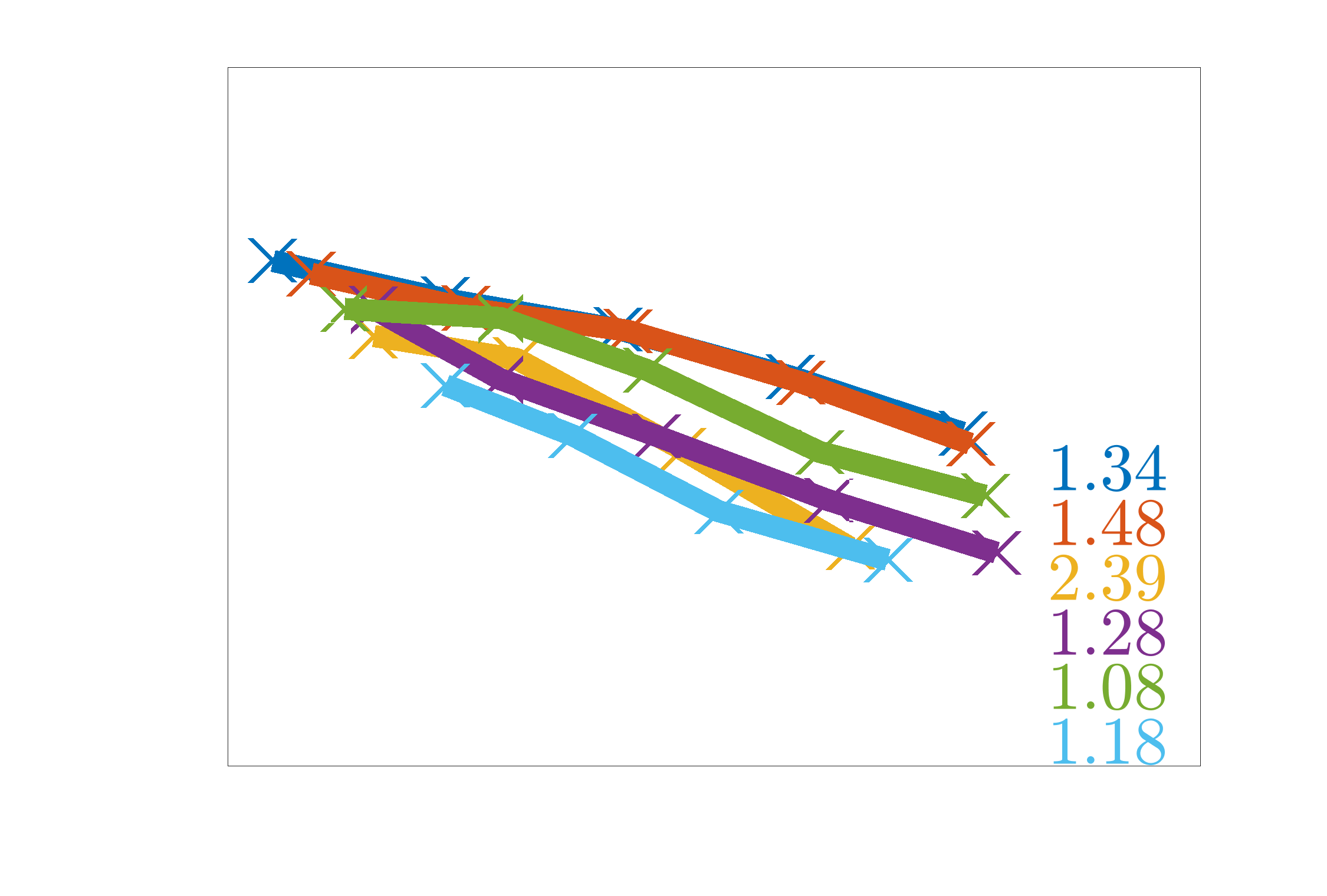}}
\subfigure[$t = 0.35$]{\includegraphics[width=2.1in]{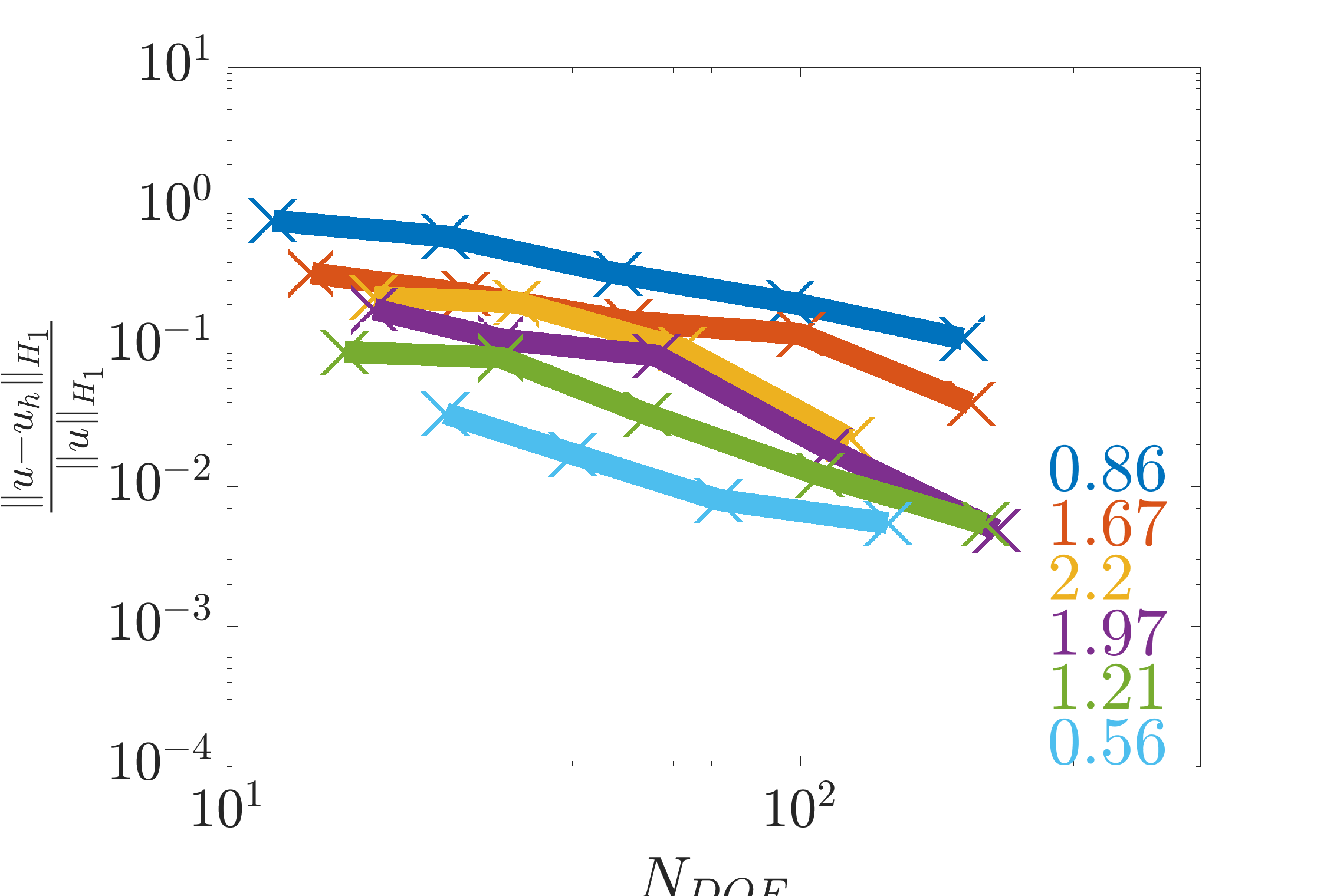}}
\subfigure[$t = 0.5$]{\includegraphics[width=2.1in]{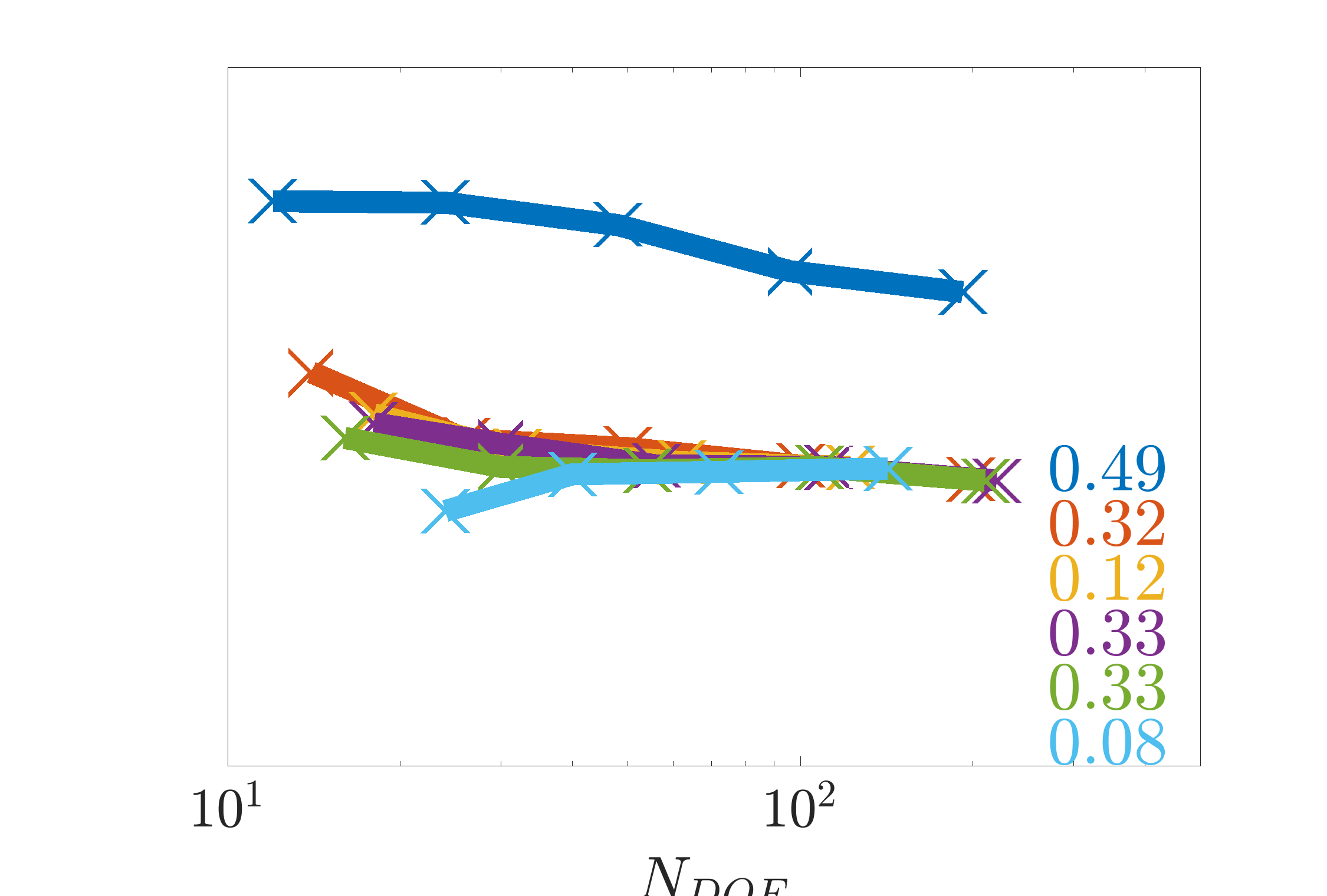}}
\subfigure[$t = 0.75$]{\includegraphics[width=2.1in]{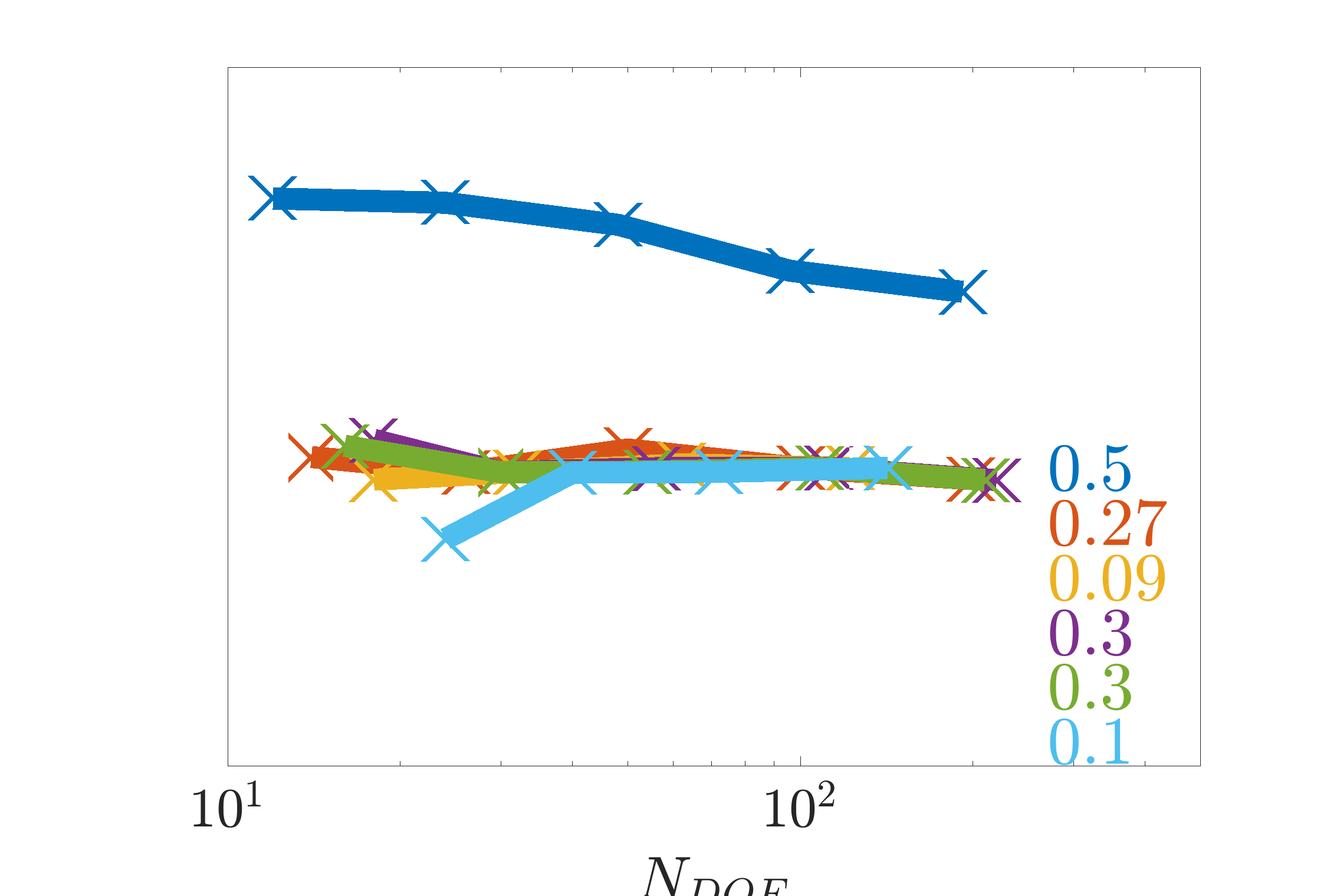}}
\end{subfigmatrix}
\caption{Convergence the relative $H_1$ integral norm for the shock problem with $\nu = \frac{1}{500}$}
\label{fig:Example2_H1_vsdofs_nu1over500_withrho}
\end{center}
\end{figure}

For $\nu = \frac{1}{500}$, GFEM solutions are further improved during shock formation by enriching the domain with additional shock enrichments, $E_{\alpha j} = \tanh{ \Big[\frac{1}{2 \rho} \Big( \frac{1}{2} - x\Big) \Big]}$, where $\rho$ controls the shock thickness. The local domain(s) these enrichments are applied are given by $\Omega_{local} = \Big[\frac{1}{2} - 2\rho \tanh^{-1}{0.99}-h_{e}, \frac{1}{2} + 2\rho \tanh^{-1}{0.99} + h_{e}\Big]$, which is the region where $|E_{\alpha j}| \leq 0.99$. Plots of the shock enrichments for various $\rho$ are shown in Fig. \ref{fig:Example2_shock_enrichments}.
The $p = 1$ + ss GFEM solution is further enriched using $\rho = \frac{1}{50}$, $\rho = \frac{1}{100}$, $\rho = \frac{1}{200}$, $\rho = \Bigg[\frac{1}{50}, \frac{1}{100}, \frac{1}{200}\Bigg]$. These solutions are presented as $p = 1$ + ss + $\rho = \frac{1}{50}$, $p = 1$ + ss + $\rho = \frac{1}{100}$, $p = 1$ + ss + $\rho = \frac{1}{200}$, and $p = 1$ + ss + $\rho =$ all GFEM solutions, respectively. Relative $L_2$ and $H_1$ norm versus time plots are shown in Fig. \ref{fig:Example2_11element_L2H1vstime_nu1over500_intermediatescales}.
 for 11-element solutions. Here the addition of various shock enrichments reduces the maximum error in the GFEM solutions, with the $p = 1$ + ss + $\rho =$ all GFEM providing the largest reduction of error. The maximum error in the $L_2$ and $H_1$ norm for 11-element $p = 1$ + ss GFEM is 4\% and 37.9\%, respectively. For 11-element $p = 1$ + ss + $\rho =$ all GFEM, the maximum error in the $L_2$ and $H_1$ norm is 0.75\% and 9.6\%, respectively, providing a 4-5 times reduction of error. Error levels improve further with grid refinement, as shown in Fig. \ref{fig:Example2_11element_L2H1vstime_nu1over500_intermediatescales} for 47-element solutions. Here, the maximum error in the $L_2$ and $H_1$ norm for 47-element $p = 1$ + ss GFEM is 0.75\% and 17.1\%, respectively. For 47-element $p = 1$ + ss + $\rho =$ all GFEM, the maximum error in the $L_2$ and $H_1$ norm is 0.031\% and 1.55\%, respectively, providing over a 10 times reduction of error. Lastly, convergence in the $L_2$ and $H_1$ norms versus total degrees of freedom are shown in Figs. \ref{fig:Example2_L2_vsdofs_nu1over500_withrho} and \ref{fig:Example2_H1_vsdofs_nu1over500_withrho}, respectively. Here, the addition of multiple shock enrichments improves overall convergence, specifically during shock formation between $t = 0.25$ and $t = 0.35$.  

\begin{figure}[ht!]
\begin{center}
\begin{subfigmatrix}{2}
\subfigure[Relative $L_2$ integral norm]{\includegraphics[width=3in]{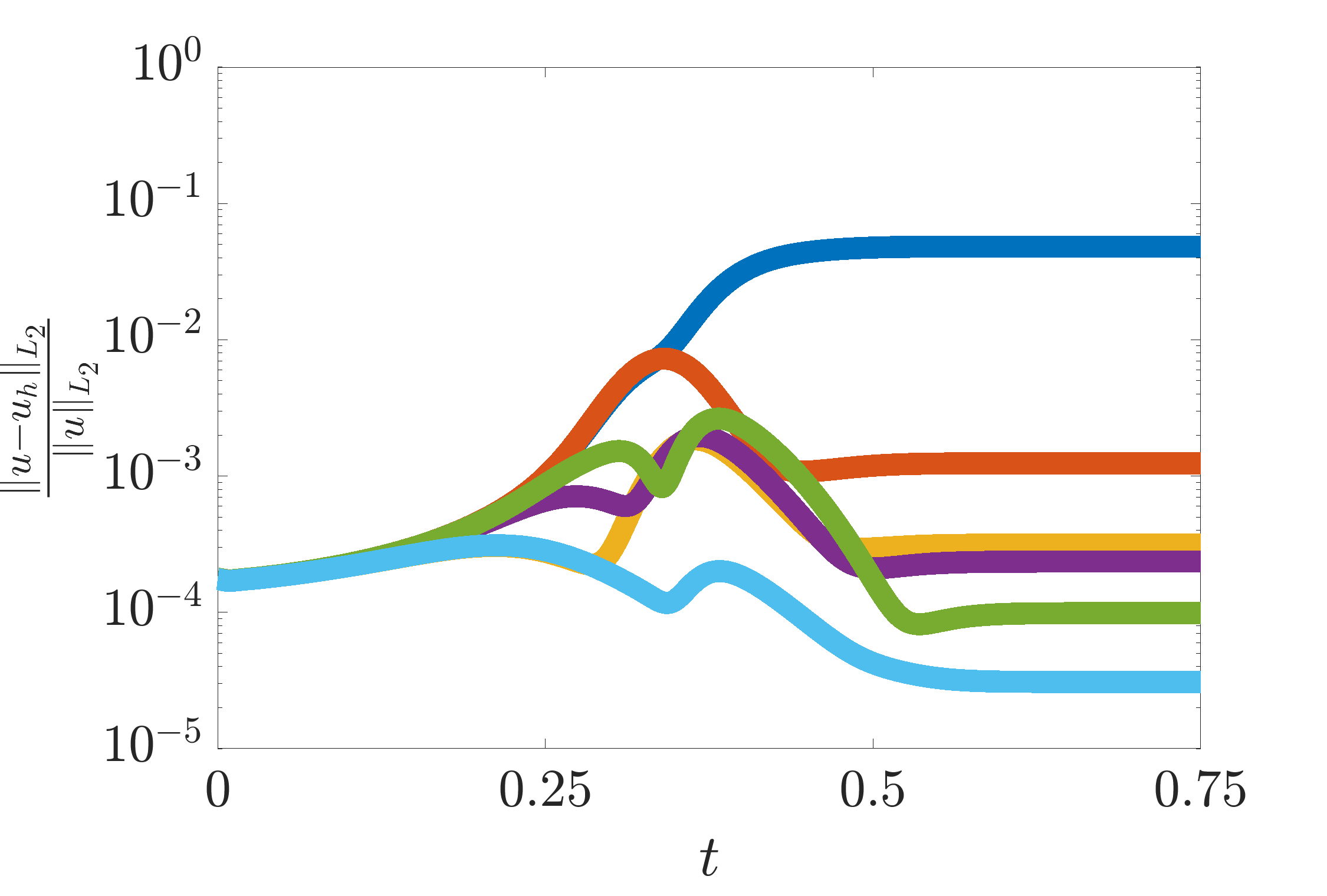}}
\subfigure[Relative $H_1$ integral norm]{\includegraphics[width=3in]{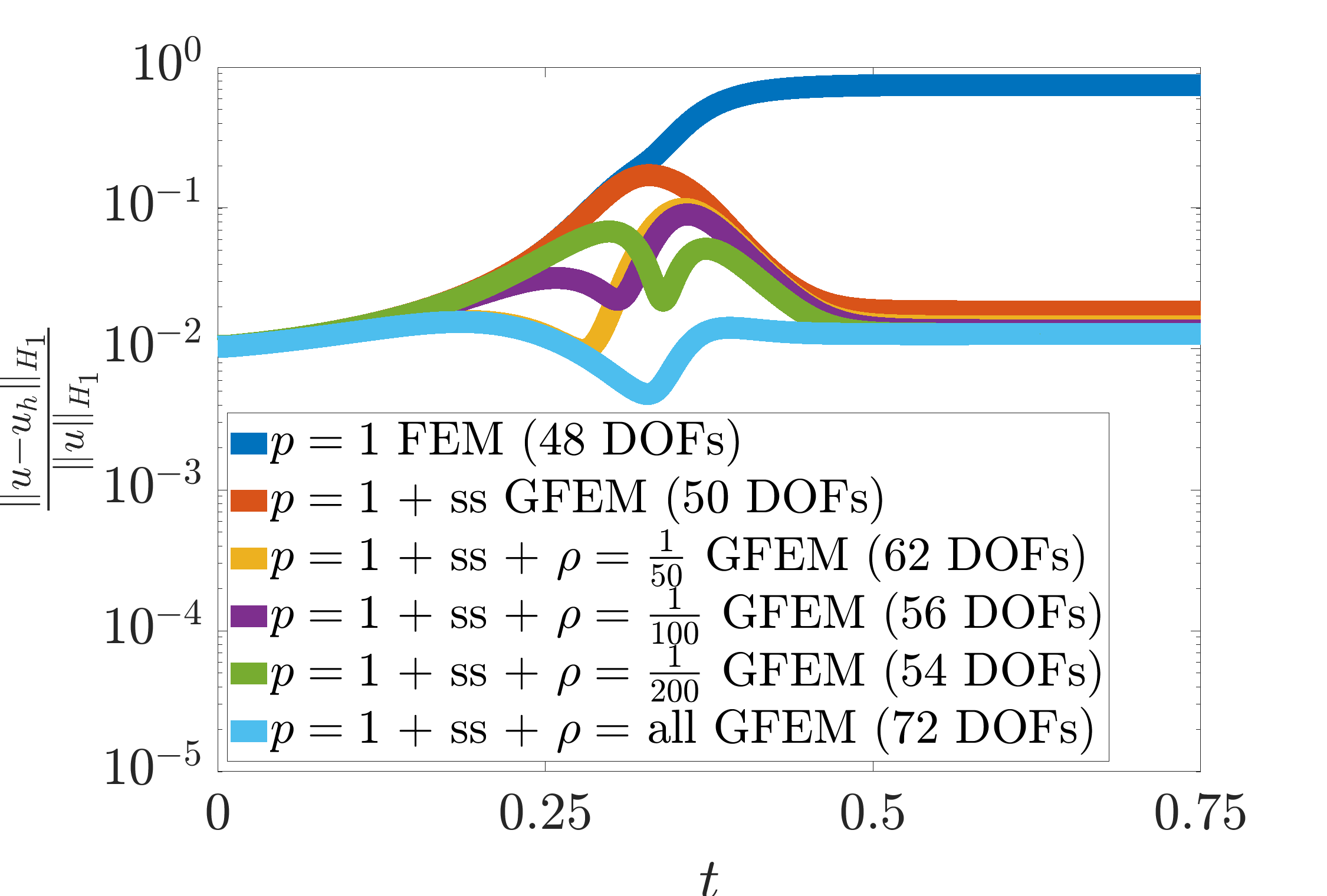}}
\end{subfigmatrix}
\caption{Relative $L_2$ and $H_1$ integral norms versus time for 47-element FEM and GFEM solutions to the shock formation problem when $\nu = \frac{1}{500}$}
\label{fig:Example2_47element_L2H1vstime_nu1over500_intermediatescales}
\end{center}
\end{figure}

\section{Conclusions}
This work presents a stable, numerical solution of the one-dimensional, unsteady Burgers' equation for a boundary layer and shock formation problem over a range of small kinematic viscosities. Compared to linear FEM, the GFEM using solution-tailored enrichments yields a significant error reduction at the same number of degrees of freedom. For the boundary layer problem, the exponential enrichments obtained in \citep{Shilt2021} are sufficient for capturing the formation of steep boundary layers. For the shock formation problem, hyperbolic tangent functions capture the thin shock forming in the domain. For both examples presented, the enrichments effectively capture the local phenomena up to relatively small kinematic viscosities and reduce errors significantly regarding linear FEM in both the relative $L_2$ and $H_1$ norms. However, as the kinematic viscosity approaches extremely small values, the intermediate solution features impact the GFEM solution stability. Specifically, the boundary layer and shock formation in both examples exhibit a range of steep gradients over intermediate time scales when $\nu << 1$, which neither the linear interpolation nor initially presented solution-tailored enrichments are sufficient at capturing. The result is spurious oscillations in the GFEM solution during the formation of the boundary layer/shock. These oscillations propagate through later time steps and affect $H_1$ and $L_2$ norm convergence. Although oscillations exist in the GFEM solutions over extremely small kinematic viscosities, the oscillations are small over coarse grids, and the errors of the GFEM solutions are still significantly reduced compared to the linear FEM. However, to further improve the shock formation problem results, a set of shock enrichments were introduced to capture the range of scales, resulting in a further reduction of error in both the $L_2$ and $H_1$ norms. Specifically, roughly greater than 100 times error reduction is observed compared to the linear FEM at the same number of degrees of freedom. Capturing the intermediate, transitional solution features will likely be an important challenge for solving more complex flow field problems using the presented GFEM framework. Such problems may demand a set of enrichments that capture various scales of the flow as presented in the shock example or time-dependent enrichments as presented in \citep{Ohara2011}. This analysis is beyond the scope of the current work and a subject for future studies. 

\section*{Acknowledgments}

This work is made possible through an Ohio State Presidential Fellowship to Troy Shilt. The authors thank Prof. C. Armando Duarte, University of Illinois Urbana-Champaign, and Dr. Rohit Deshmukh, The Ohio State University, for their insights on the GFEM. 

\bibliography{references}
\end{document}